\documentclass[a4paper,11pt]{article}
\usepackage{amssymb,amsthm,graphicx,inputenc,mathtools,verbatim}
\usepackage{hyperref,url,enumitem,bm,esint,color}
\usepackage[margin=30mm]{geometry}

\usepackage{tikz}
\usepackage{cite}
\usepackage{booktabs}

\allowdisplaybreaks

\definecolor{rosso}{rgb}{0.85,0,0}
\definecolor{azzurro}{rgb}{0.13, 0.67, 0.8}
\definecolor{mygray}{gray}{0.6}

\def\Accorpa #1#2 #3 {\gdef #1{\eqref{#2}--\eqref{#3}}%
  \wlog{}\wlog{\string #1 -> #2 - #3}\wlog{}}

\newcommand{\ov}[1]{\overline{#1}}

\renewcommand{\hat}[1]{\widehat{#1}}
\let\badeps\epsilon
\newcommand{\eps}{\varepsilon}

\newcommand{\E}{\mathcal{E}}

\newcommand{\pd}{\partial}

\newcommand{\dt}{\partial_t}
\newcommand{\dn}{\partial_{\bnn}}
\def\bnn{{\boldsymbol n}}
\def\bnu{{\boldsymbol \nu}}

\def\norma #1{\mathopen \| #1\mathclose \|}
\def\<#1>{\mathopen\langle #1\mathclose\rangle}
\def\iO {\int_\Omega}
\def\intinf{\int_{-\infty}^{+ \infty}}

\renewcommand{\tilde}{\widetilde}

\newcommand{\ph}{{{\varphi}}}

\newcommand{\phz}{\phi_0}
\def\non{\notag}

\def\erre{{\mathbb{R}}}
\def\enne{{\mathbb{N}}}

\def\Vp{{(\Hx1)^*}}

\newcommand{\bigchi}{\ensuremath{\mathrm{\mathcal{X}}}}
\newcommand{\charfcn}[1]{\bigchi_{#1}} 

\theoremstyle{plain}
\newtheorem{thm}{Theorem}[section]

\newtheorem{prop}{Proposition}[section]

\newtheorem{remark}{Remark}[section]

\numberwithin{equation}{section}

\def\genspazio #1#2#3#4#5{#1^{#2}(#5,#4;#3)}
\def\spazio #1#2#3{\genspazio {#1}{#2}{#3}T0}

\def\L {\spazio L}
\def\H {\spazio H}
\def\C #1#2{C^{#1}([0,T];#2)}

\def\Lx #1{L^{#1}(\Omega)}
\def\Hx #1{H^{#1}(\Omega)}

\def\lhs{left-hand side}
\def\rhs{right-hand side}

\def\jump#1{[#1]^+_-}
\def\Sp{{S_+}}
\def\Sm{{S_-}}
\def\Kp{{K_+}}
\def\Km{{K_-}}

\def\th{\vartheta}

\def\Il{\mathbb{I}_\ell}
\def\Kl{\mathbb{K}_\ell}
\def\Iz{\mathbb{I}_0}
\def\Kz{\mathbb{K}_0}

\def\rad{Y}

\def\multibold #1{\def\arg{#1}%
  \ifx\arg\pto \let\next\relax
  \else
  \def\next{\expandafter
    \def\csname #1#1\endcsname{{\bf #1}}%
    \multibold}%
  \fi \next}

\def\pto{.}

\def\multimathbb #1{\def\arg{#1}%
  \ifx\arg\pto \let\next\relax
  \else
  \def\next{\expandafter
    \def\csname #1#1#1\endcsname{{\mathbb #1}}%
    \multimathbb}%
  \fi \next}

\def\multical #1{\def\arg{#1}%
  \ifx\arg\pto \let\next\relax
  \else
  \def\next{\expandafter
    \def\csname cal#1\endcsname{{\cal #1}}%
    \multical}%
  \fi \next}

\def\multimathop #1 {\def\arg{#1}%
  \ifx\arg\pto \let\next\relax
  \else
  \def\next{\expandafter
    \def\csname #1\endcsname{\mathop{\rm #1}\nolimits}%
    \multimathop}%
  \fi \next}

\multibold
qwertyuiopasdfhjklzxcvbnmQWERTYUIOPASDFGHJKLZXCVBNM.

\multimathbb
QWERTYUIOPASDFGHJKLZXCVBNM.

\multical
QWERTYUIOPASDFGHJKLZXCVBNM.

\multimathop
diag dist div dom mean meas sign supp .

\def\nn{\boldsymbol{n}}

\def\qstar{{q^*}}

\def\mustpl{\mu^{*}}

\def\bfeta{\boldsymbol{\eta}}
\def\bfJ{\boldsymbol{J}}

\def\phd{\ph_\delta}
\def\mud{\mu_\delta}
\def\mdel{m_\delta}

\def\phz{\ph_0}
\def\muz{\mu_0}

\def\esp{\theta}
\def\Source{S}

\newcommand{\sigmaO}{o}
\newcommand{\Snh}{\mathcal S^n_h}
\newcommand{\SmD}{\mathcal S^n_{h,-}}
\newcommand{\Whm}{V(\Sigma^n)}
\newcommand{\Vhm}{\underline{V}(\Sigma^n)}
\newcommand{\bB}{\mathbb{B}}
\title{{One-phase models for chemically active droplets}
}
\date{}

\author{Harald Garcke \footnotemark[1] \and Kei Fong Lam \footnotemark[2] \and Robert N\"urnberg \footnotemark[3] \and Andrea Signori \footnotemark[4]}

\begin{document}
\maketitle

\begin{abstract} We study one-phase models for chemically active droplets arising from a Cahn--Hilliard equation with one-sided degenerate mobility. The degeneracy confines diffusion to only one of the two pure phases and in the sharp-interface limit leads to a reaction-modified one-phase Mullins--Sekerka problem. Compared with the corresponding two-phase formulation, the one-phase models involve fewer material and reaction parameters. This reduction simplifies analytical investigations and allows for more transparent and computationally efficient parameter studies.
	
We establish the existence of weak solutions to the diffuse-interface problem by regularizing the degenerate mobility and deriving uniform energy and entropy estimates. Formally matched asymptotic expansions are then used to identify the associated one-phase sharp-interface limits. Depending on the pure phase in which the  mobility vanishes, a \emph{plus-phase model} and a \emph{minus-phase model} are obtained. 
	 
For radially symmetric configurations in two and three spatial dimensions, we investigate the existence and radial stability of stationary droplets and analyse their stability with respect to non-radial perturbations. The resulting amplification factors identify parameter regimes in which stationary droplets become morphologically unstable. Numerical computations complement the analytical results and illustrate the influence of surface tension, chemical reactions, and confinement on the stability and evolution of active droplets. \end{abstract} 
\medskip \noindent
\textbf{Keywords:} Cahn--Hilliard equation, degenerate mobility, one-phase model, Mullins--Sekerka problem, free boundary problem, active droplets, sharp-interface limit,  stability analysis. 

\medskip \noindent\textbf{AMS Subject Classification:} 35B25, 35K55, 35K61, 35R35, 76T30, 92D25.


\renewcommand{\thefootnote}{\fnsymbol{footnote}}
\footnotetext[1]{Fakult{\"a}t f\"ur Mathematik, Universit{\"a}t Regensburg, 93040 Regensburg, Germany
({\texttt harald.garcke@ur.de}).}
\footnotetext[2]{Department of Mathematics, Hong Kong Baptist University, Kowloon Tong, Hong Kong ({\texttt akflam@hkbu.edu.hk}).}
\footnotetext[3]{Department of Mathematics, University of Trento, Trento, Italy
({\texttt robert.nurnberg@unitn.it}).}
\footnotetext[4]{Department of Mathematics, Politecnico di Milano, 20133 Milano, Italy ({\texttt andrea.signori@polimi.it}), Alexander von Humboldt Research Fellow.}
%


%
%

\renewcommand{\thefootnote}{\arabic{footnote}}
\newcommand{\typo}{\footnote{Typo}}

\section{Introduction}
Phase separation in the presence of chemical reactions can produce dynamics that differ substantially from the classical coarsening behaviour of passive systems. In particular, chemical reactions may suppress Ostwald ripening and permit the coexistence of droplets of controlled size \cite{zwickerostwald}. Depending on the balance between diffusion, surface tension, and reaction rates, growing droplets may also lose their stability and undergo division and subsequent growth of the daughter droplets. The resulting cycles of growth and splitting have been proposed as a possible mechanism for the formation and reproduction of protocells; see, for instance, \cite{Active_drops}. Synthetic realisations of chemically active  droplets have also been developed, providing experimentally accessible model systems for membraneless organelles and protocells \cite{DonauSSetal20}. Chemically active droplet systems may also exhibit more complex morphologies, including non-equilibrium liquid shells \cite{bauermann2023formation,BergmannBBetal23}.

More generally, phase-separated droplets provide spatially localised compartments and therefore play an important role in the organisation of chemical and biological processes. A diffuse-interface description of chemically active droplets can be given by a Cahn--Hilliard equation supplemented by suitable reaction terms. In \cite{OOL}, we studied such a model from an analytical and numerical point of view. In particular, we established well-posedness of the phase-field system and used formally matched asymptotic expansions to derive a Mullins--Sekerka-type free boundary problem with bulk and interfacial reactions. The resulting sharp-interface model differs from the classical Mullins--Sekerka problem through affine reaction terms in the quasi-static bulk equations and an additional source term in the law for the normal velocity. Numerical computations for the diffuse-interface formulation demonstrated, among other phenomena, the suppression of Ostwald ripening, droplet growth and division, and the formation of shell-like structures. 
 
The corresponding two-phase sharp-interface problem was subsequently investigated in detail in \cite{OOL2}. There, radially symmetric stationary solutions and their stability with respect to radial and non-radial perturbations were analysed. The analysis also included multilayered configurations and a numerical method was introduced which was capable of resolving topological changes such as droplet splitting and merging. These results show that the reaction-modified Mullins--Sekerka model captures a wide range of complex droplet dynamics. At the same time, the two-phase formulation involves separate diffusion, mobility, and reaction parameters in both phases. This substantial number of parameters renders a systematic exploration of the model, and in particular the identification of stability boundaries, highly non-trivial.
 
The purpose of the present work is to derive and analyse one-phase versions of the active Mullins--Sekerka model. These models arise naturally when the mobility in one of the pure phases vanishes. At the diffuse-interface level, we therefore consider a Cahn--Hilliard equation with a one-sided degenerate mobility. In the sharp-interface limit, diffusion is confined to the phase in which the mobility remains positive, while the chemical potential in the other phase no longer appears as an independent unknown. As a consequence, the limiting free boundary problem contains only one bulk equation coupled to the interfacial conditions. This one-phase formulation, arising from a mathematically natural limit of the diffuse-interface description, has several practical advantages. It contains substantially fewer material and reaction parameters than the corresponding two-phase model. This reduction removes several parameter ratios and makes it easier to distinguish the effects of surface tension, bulk reactions, interfacial reactions, and the size of the domain in which the evolution problem is formulated in. Hence, the exploration of the parameter space becomes both more transparent and computationally less expensive. 
The one-phase models preserve the principal mechanisms responsible for stationary droplets and morphological instabilities, while providing a simpler framework for their analysis. 
  
There are two possible one-phase limits, depending on the pure phase at which the mobility degenerates. If the mobility vanishes at $\ph=-1$, diffusion takes place only in the \emph{plus-phase}, leading to what we call the \emph{plus-phase model}. Conversely, degeneracy at $\ph=1$ gives rise to the \emph{minus-phase model}, in which the chemical potential is determined in the \emph{minus-phase}. The analytical arguments leading to the two formulations are closely related. We establish the existence theory for both choices and derive their sharp-interface limits within a common framework. 
  
The subsequent analysis, however, is centered on the \emph{minus-phase model}. The \emph{plus-phase} formulation is included mainly to display the symmetry of the diffuse-interface construction and to clarify the relation between the two possible degeneracies. The \emph{minus-phase model} is particularly suitable for the study of a droplet surrounded by a chemically active exterior phase. In this setting, the quasi-static diffusion equation is posed in the region outside the droplet and is supplemented by a homogeneous Neumann condition on the external boundary. The location of this boundary can influence both the existence and the stability of stationary droplets. This makes confinement an intrinsic part of the model and distinguishes the resulting problem from a purely local interfacial evolution law. 
  
Cahn--Hilliard equations with degenerate mobility have been studied extensively, with the analysis typically relying on regularization arguments and entropy estimates; see, e.g.,\cite{ElliottG96, Agosti2017, ebenbeck2021cahn, garcke2025anisotropic, Carrillo}. We first establish the existence of weak solutions to the Cahn--Hilliard system with one-sided degenerate mobility. The proof is based on a regularization of the mobility and on estimates that are uniform with respect to the regularization parameter. In addition to the usual energy estimate, an entropy-type estimate is employed to control the solution near the pure phase at which the mobility vanishes. This estimate allows us to pass to the degenerate limit and ensures that the order parameter remains on the physically relevant side of the corresponding pure-phase value. Our existence proof follows the regularization and entropy approach developed for Cahn--Hilliard equations with degenerate mobility; see \cite{ElliottG96} and, for a recent related treatment, \cite{ebenbeck2021cahn,garcke2025anisotropic}. We then use formally matched asymptotic expansions to derive the associated one-phase free boundary problems as the interfacial thickness tends to zero. For positive mobility, the convergence of the Cahn--Hilliard equation to the Mullins--Sekerka, or Hele--Shaw, free boundary problem was rigorously established in \cite{AlikakosBC94}. The sharp-interface dynamics of Cahn--Hilliard equations with degenerate mobility can differ substantially from the classical Mullins--Sekerka evolution. For mobilities vanishing in both pure phases, surface-diffusion-type limits have been derived and analysed in \cite{CahnEN-C96, LeeMS16}. Here, we instead consider a one-sided degeneracy: the mobility vanishes in only one pure phase and remains positive in the other. As a result, diffusion is suppressed in one phase but persists in the complementary phase, leading to a one-phase Mullins--Sekerka problem with bulk and interfacial reactions. We refer to \cite{Glasner03} for the case without chemical reactions.
  
In the \emph{minus-phase case}, the limiting chemical potential satisfies a {reaction-diffusion} equation only in the \emph{minus-phase}. Its trace on the interface is proportional to the mean curvature, while the normal velocity is determined by the one-sided diffusive flux together with the interfacial reaction rate. The reduction from a two-sided flux jump to a single flux is the characteristic structural feature of the one-phase limit. After non-dimensionalisation, the model is governed by a small set of parameters measuring the relative influence of capillarity, interfacial reactions, and confinement. Our analysis of the sharp-interface problem focuses on radially symmetric droplets in two and three spatial dimensions. We derive explicit representations of the chemical potential in terms of modified Bessel functions and obtain an ordinary differential equation for the droplet radius. This reduction allows us to investigate the existence and radial stability of stationary droplets. Depending on the parameters, the radial evolution may admit multiple stationary radii with different stability properties. In particular, the interaction between the exterior reaction-diffusion field and the outer boundary can stabilise droplets whose radius would otherwise decrease or grow. 
  
We next study the response of stationary spherical droplets to non-radial perturbations. Stability and long-time behaviour of equilibria for the classical Mullins--Sekerka problem have been investigated by several authors; see, for instance, \cite{MullinsS63,ESCHER1998267}. By decomposing perturbations into Fourier modes in two dimensions and spherical harmonics in three dimensions, we derive an amplification factor for every mode. The resulting expressions identify the competition between the destabilising effect of the chemically driven flux and the stabilising effect of surface tension. They also permit direct parameter studies of the onset of morphological instabilities. Owing to the reduced number of parameters in the one-phase formulation, these stability regimes can be explored more systematically than in the full two-phase model in \cite{OOL2}. 
  
Finally, we complement the analytical results with numerical computations for the diffuse and the sharp-interface problem using approaches developed in \cite{BarrettBG99,bgnreview}. The simulations are used to examine the radial dynamics, verify the predictions of the linear stability analysis, and investigate the nonlinear evolution beyond the onset of instability. In particular, we explore how changes in the capillary and reaction parameters, as well as in the size of the outer domain, affect the transition from stable droplets to deformed or splitting configurations. 
  
The remainder of the paper is organised as follows. In Section~\ref{sec:models}, we formulate the diffuse-interface problem with one-sided degenerate mobility and the corresponding one-sided Mullins--Sekerka sharp interface problems. In Section~\ref{sec:anal} we establish the existence of weak solutions to the Cahn--Hilliard equation with reactions and one-sided degenerate mobility. In Section~\ref{sec:Sharp}, we derive the one-phase sharp-interface limits and introduce their non-dimensional formulations. In Section~\ref{sec:radial}, we derive radially symmetric solutions and investigate the existence and radial stability of stationary droplets. In Section~\ref{sec:linstab}, we perform the linear stability analysis with respect to non-radial perturbations. Section~\ref{sec:num} introduces the numerical method and presents numerical computations illustrating the analytical findings and the nonlinear behaviour of the model.

\section{The mathematical models}\label{sec:models}
\subsection{The Cahn--Hilliard model}
Consider a bounded domain $\Omega \subset \mathbb{R}^d$, with $d\in \{1,2,3\}$, and smooth boundary $\partial \Omega$, modelling a mixture of two chemical species. The evolution of the system is described by the order parameter $\ph$, representing the normalised concentration difference between the two components. This variable satisfies the following Cahn--Hilliard system incorporating chemical reactions \cite{Active_drops} (cf.~\cite[(2.1)]{OOL}):
\begin{subequations}\label{eq:sys21}
\begin{alignat}{2}
	\label{sys:1}
	\dt \ph  &= \div(m(\ph)\nabla  \mu) + S_\eps(\ph)
	\qquad && \text{in $Q := (0,T) \times \Omega $,}
	\\
	\label{sys:2}
	 \mu &= \beta (-\eps \Delta \ph + \tfrac 1 \eps \psi'(\ph))
	\qquad && \text{in $Q$,}
	\\
	\label{sys:3}
	\dn \mu  &= \dn \ph =0
	\qquad && \text{on $\Gamma := (0,T) \times \partial \Omega $,}
	\\
	\label{sys:4}
	 \ph(0) &=\phz
	\qquad && \text{in $\Omega$.}
\end{alignat}
\end{subequations}
\Accorpa\Sys {sys:1} {sys:4}
In this formulation, $\mu$ denotes the chemical potential, while $m: \erre \to \erre_{\geq 0}$ represents a mobility function that depends on the local concentration. The positive constant $\beta$ is linked to the surface energy density, and $\eps>0$ is a small parameter characterizing the thickness of the diffuse interface. The function $\psi: \erre \to \erre_{\geq 0}$ stands for a double-well potential, $\phz$ provides the initial condition for the order parameter, and $\dn$ indicates the outward normal derivative on $\partial \Omega$ with respect to the unit normal vector $\bnn$. The system is analysed over the space-time cylinder $Q$ for an arbitrary but fixed final time $T>0$. This model has recently been investigated by the authors in \cite{OOL, OOL2}, where the well-posedness of the phase-field formulation was rigorously established and the corresponding sharp-interface limit was derived through an asymptotic analysis as $\eps \to 0$. Furthermore, those studies include a linear stability analysis of the resulting sharp-interface problem, complemented by numerical simulations. {In those studies, the mobility function $m$ is taken to be strictly positive over $\erre$ and is non-degenerate.}

The chemical reaction mechanism is incorporated through the source term $S_\eps:\erre \to \erre$, which is decomposed as
\begin{equation}\label{def:source}
S_\eps(r)= S_1(r)+\tfrac 1\eps S_2(r){,\quad \forall r \in \erre}.
\end{equation}
The $\frac 1\eps$ scaling of $S_2$ is specifically designed to model rapid reaction dynamics localised within the interfacial layer. For a fixed threshold parameter $r_c \in (0,1]$ and given constants $S_+$, $S_-$, $K_+$, $K_-$, and $L$, we define the bulk and interfacial contributions as follows:
\begin{align}
	\label{def:source1} 
	S_1(r)= 
	\begin{cases}
		\Sp  & \quad \text{if $r \geq {r_c}$,}
		\\
		\Sm + G_1(r)(\Sp-\Sm) 
		& \quad \text{if $r\in (-{r_c},{r_c})$,}
		\\
		\Sm & \quad \text{if $r \leq -{r_c}$,}
	\end{cases}
\end{align}
and
\begin{align}
	\label{def:source2} 
	S_2(r)= 
	\begin{cases}
		- \Kp  (r-1) & \quad \text{if $r \geq {r_c}$,}
		\\
		\hat S_2(r)
		& \quad \text{if $r\in (-{r_c},{r_c})$,}
		\\
		- \Km   (r+1)& \quad \text{if $r \leq -{r_c}$,}
	\end{cases}
\end{align}
where we define for $r\in (-{r_c},{r_c})$
$$
\hat S_2(r)=		- {\Km}  G_2(r) -  {\Kp}  G_3(r)+ L G_4(r)
{- K_+(r_c-1)G_1(r) - K_-(1 - r_c)(1-G_1(r))}.
$$
The functions $G_1,G_2,G_3, G_4 :[-r_c,r_c] \to \erre$ are smooth interpolants that will be specified in the subsequent subsection. They are required to fulfill the matching conditions
\begin{align}
	& \label{ass:G:1}
	G_1({r_c})=1,
	\quad 
	G_1(-{r_c}) =0,
	\quad 
	G_2(\pm{r_c}) =G_3(\pm{r_c}) =G_4(\pm{r_c}) =0, 
	\\ 
	& \label{ass:G:2}
	G_1'(\pm{r_c})=G_4'(\pm{r_c})=0,
	\quad 
	G_2'({r_c})=G_3'(-{r_c})=0,
	\quad 
	G_2'(-{r_c})=G_3'({r_c})=1,
\end{align}
which guarantee the global differentiability of the source term $S_\eps$ across the whole real line. In the asymptotic analysis that follows, the choice $r_c=1$ is typically adopted as it significantly streamlines both the algebraic form of $\hat S_2$ and the subsequent matched asymptotic expansion procedure. Nevertheless, alternative values remain admissible; for instance, $r_c = \frac12$ was employed in \cite{Active_drops}.

The evolution equation \eqref{eq:sys21} is associated with the following \textit{free energy} functional
\begin{align*}
	{\cal E}(\ph) = \beta \Big(
	\frac \eps 2 \iO |\nabla \ph|^2
	+ \frac 1\eps \iO \psi(\ph)	
	\Big).
\end{align*}
Throughout this work, we assume that the potential $\psi$ is symmetric, i.e., $\psi(r)=\psi(-r)$ for all $r \in \erre$, and possesses  minima at $\pm 1$, meaning $\psi(\pm 1) = 0$ and $\psi''(\pm 1) \neq 0$. A standard choice satisfying these properties is the \textit{quartic double-well potential}
\begin{align}
	\label{quartic}
	\psi(r) = \tfrac 14 (1-r^2)^2, \quad \forall r \in \erre.
\end{align}

While the forthcoming theoretical framework remains valid for any set of interpolants satisfying the conditions \eqref{ass:G:1}--\eqref{ass:G:2}, it is convenient to construct explicit expressions tied to the underlying potential $\psi$. Specifically, for $r \in [-r_c, r_c]$, we define
\begin{align}
	\label{def:G1}
	G_1(r) { = \hat G_1(\tfrac{r}{r_c}), \quad \hat G_1(r) }
	& =
	\tfrac 34 (r+1)^2- \tfrac 14(r+1)^3,
	\\
	\label{def:G23}
	G_2(r) { = r_c \hat G_2(\tfrac{r}{r_c}), \quad \hat G_2(r) }
	& = -\tfrac 12\tfrac {1}{ \sqrt{\psi''(-1)}} (r-1) \sqrt{2 \psi(r)},
	\quad
	G_3 (r) = - G_2(-r),
\end{align}
which immediately yield the following boundary properties:
\begin{align*}
	\hat G_1(1)=1,
	\quad 
\hat 	G_1(-1) = \hat G_1'(\pm1)=0,	
	\quad 
\hat 	G_2(\pm 1) =0,
	\quad 
\hat 	G_2'(1)=0,
	\quad 
\hat 	G_2'(-1) = 1.
\end{align*}
By symmetry, the relation $\hat G_3(r) = - \hat G_2(-r)$ immediately guarantees $\hat G_3'(1) = 1$. Specializing to the quartic potential \eqref{quartic}, where $\psi''(-1) = 2$, these definitions simplify to
\[
\hat G_2(r) = - \tfrac{1}{4} (1-r^2)(r-1), \quad \hat G_3(r) = - \hat G_2(-r) = - \tfrac{1}{4}(1-r^2)(r+1),
\]
which readily satisfy the required matching conditions \eqref{ass:G:1} and \eqref{ass:G:2}. Finally, we define the fourth interpolant as
\begin{equation}\label{eq:G4}
	G_4(r)  = \hat G_4(\tfrac{r}{r_c}), \quad \hat G_4(r) 
	= 2\psi(r) ,\quad 
	\forall r\in \mathbb{R},
\end{equation}
a choice that trivially fulfills the boundary requirements \eqref{ass:G:1} and \eqref{ass:G:2}.

The mobility function $m$, which plays a central role in this work, is selected to degenerate in one of the two pure phases. 
A possible choice is the function
\begin{equation} \label{eq:mob}
m(r) = \begin{cases}
m_+ & \text{ for } r \geq 1, \\
m_+ \frac{1+r}{2} + m_- \frac{1-r}{2} & \text{ for } |r| \leq 1, \\
m_- & \text{ for } r \leq -1,
\end{cases}
\end{equation}
with constants $m_+ > 0$ and $m_-=0$ for the \emph{plus-phase model}, and with constants $m_- > 0$ and $m_+=0$ for the \emph{minus-phase model}.

\subsection{Mullins--Sekerka system with reactions}
As established in \cite{OOL} {for the case of non-degenerate mobility}, the phase-field system \eqref{eq:sys21} converges to a reaction-modified Mullins--Sekerka model in the sharp-interface limit $\eps \to 0$. We now present the resulting free boundary problem, omitting the technical details of the asymptotic derivation.

The domain $\Omega$ is partitioned into two disjoint, time-dependent subdomains $\Omega^+(t)$ and $\Omega^-(t)$, separated by a smooth, evolving interface $\Sigma(t) := \partial \Omega^+(t) \cap \partial \Omega^-(t)$. In the canonical configuration, $\Omega^+(t)$ represents a droplet fully contained within the bulk phase $\Omega^-(t)$, so that $\Sigma(t) \cap \partial \Omega = \varnothing$. 

The sharp-interface dynamics of active droplets are governed by a Mullins--Sekerka system with bulk reactions. Given an initial interface $\Sigma(0) = \Sigma_0$, the evolution for $t \in (0,T)$ is described by the following free boundary problem: find the moving interface $\Sigma(t)$ and the associated chemical potential $\mu(t,\cdot): \Omega \to \mathbb{R}$ such that (cf.~\cite[(2.1)]{OOL2})
\begin{subequations}\label{SharpI:original}
\begin{alignat}{2}
	\label{sharpI:1:p}
	 - m_+ \Delta \mu &= \Sp - \rho_+ \mu  \quad &&\text{in $\Omega^+(t)$,}
		\\
        \label{sharpI:1:m} - m_- \Delta \mu & =\Sm - \rho_- \mu  \quad && \text{in $\Omega^-(t)$,}
	\\
	\label{sharpI:2}
	\mu &=  	{\alpha}	\kappa
	\qquad && \text{on $\Sigma(t)$,}
	\\
	\label{sharpI:4}
	 - 2 {\cal V} & = \jump{m\nabla \mu} \cdot {\bnu}
	+ S_I
	\qquad && \text{on $\Sigma(t)$,}
	\\
	\label{sharpI:5}
	 \dn \mu & =0
	\qquad && \text{on $\partial \Omega${.}}
\end{alignat}
\end{subequations}
Here, $S_+ < 0$ denotes a sink term, whereas $S_- > 0$, $\rho_\pm > 0$, $m_\pm > 0$, and $\alpha > 0$ are given physical constants, and $S_I \in \erre$ quantifies the interfacial reaction rate. 

The unit normal $\bnu$ points from $\Omega^-$ into $\Omega^+$, and $\kappa = -\div_\Sigma \bnu$ denotes the mean curvature, taken as positive for convex interfaces. The symbol $\mathcal{V}$ represents the normal velocity of $\Sigma(t)$ in the direction of $\bnu$. 

Assuming the classical quartic double-well potential \eqref{quartic}{,} 
the macroscopic parameters $\alpha$, $\rho_{\pm}$, and $S_I$ governing the free boundary problem \eqref{SharpI:original} are explicitly determined by the phase-field coefficients in \eqref{eq:sys21} through the following identities:
\[
 \alpha=\frac{ \sqrt{2} \beta}3,
 \quad
 \rho_\pm= \frac {K_\pm}{2\beta}, 
 \quad
 S_I = \frac{1}{\sqrt{2}} \Big (K_{+} - K_{-} + \frac{4}{3} L \Big ).
\]
For the reader's convenience, we collect here the one-phase sharp-interface models that will be formally derived in the subsequent sharp-interface limit of \eqref{eq:sys21} using one-sided degenerate mobilities. Both the plus-phase and minus-phase formulations are presented in their dimensional and nondimensional forms.

\vspace{1.2em}
\noindent
\begin{minipage}[t]{0.48\textwidth}
\centering
\textbf{Dimensional Plus-Phase Model}
\begin{align*}
    - m_+ \Delta \mu &= \Sp - \rho_+ \mu  &&\text{in $\Omega^+$,} \\
    \mu &= \alpha\kappa  &&\text{on $\Sigma$,} \\
    - 2 \mathcal{V} &= m_+\nabla \mu \cdot \bnu + S_I  &&\text{on $\Sigma$,} \\
    m_+ \dn \mu &= 0  &&\text{on $\partial \Omega^+ \cap \partial \Omega$.}
\end{align*}
\end{minipage}%
\hfill
\begin{minipage}[t]{0.48\textwidth}
\centering
\textbf{Dimensional Minus-Phase Model}
\begin{align*}
    - m_- \Delta \mu &= \Sm - \rho_- \mu  &&\text{in $\Omega^-$,} \\
    \mu &= \alpha\kappa  &&\text{on $\Sigma$,} \\
    - 2 \mathcal{V} &= -m_-\nabla \mu \cdot \bnu + S_I  &&\text{on $\Sigma$,} \\
    m_- \dn \mu &= 0  &&\text{on $\partial \Omega$.}
\end{align*}
\end{minipage}

\vspace{1.2em}

\noindent
\begin{minipage}[t]{0.48\textwidth}
\centering
\textbf{Nondimensional Plus-Phase Model}
\begin{align*}
    - \Delta \mu &= -1 - \mu  &&\text{in $\Omega^+$,} \\
    \mu &= \alpha\kappa  &&\text{on $\Sigma$,} \\
    - 2 \mathcal{V} &= \nabla \mu \cdot \bnu + S_I  &&\text{on $\Sigma$,} \\
    \dn \mu &= 0  &&\text{on $\partial \Omega^+ \cap \partial \Omega$.}
\end{align*}
\end{minipage}%
\hfill
\begin{minipage}[t]{0.48\textwidth}
\centering
\textbf{Nondimensional Minus-Phase Model}
\begin{align*}
    - \Delta \mu &= 1 - \mu  &&\text{in $\Omega^-$,} \\
    \mu &= \alpha\kappa  &&\text{on $\Sigma$,} \\
    - 2 \mathcal{V} &= -\nabla \mu \cdot \bnu + S_I  &&\text{on $\Sigma$,} \\
    \dn \mu &= 0  &&\text{on $\partial \Omega$.}
\end{align*}
\end{minipage}

\section{Existence of weak solutions to the Cahn--Hilliard model}\label{sec:anal}

\subsection{Notation}
For any Banach space $X$, its norm is represented as $\| \cdot \|_X$, its dual space as $X^*$, and the duality pairing between $X^*$ and $X$ is denoted by $\langle \cdot, \cdot \rangle_X$. In the case where $X$ is a Hilbert space, the inner product is denoted by $(\cdot,\cdot)_X$. For each $1 \leq p \leq \infty$, $k \geq 0$, and $s > 0$, the standard Lebesgue and Sobolev spaces defined on $\Omega$ are denoted as $L^p(\Omega)$ and $W^{k,p}(\Omega)$, with their respective norms $\| \cdot \|_{L^p(\Omega)}$ and $\| \cdot \|_{W^{k,p}(\Omega)}$. 
We adopt the convention $H^k(\Omega) := W^{k,2}(\Omega)$ for all $k \in \mathbb{N}$, and denote the generalised mean value of a functional $h \in \Vp$ as
\[
(h)_\Omega := \frac{1}{|\Omega|} \langle h,1 \rangle_{H^1(\Omega)}.
\]

\subsection{Mathematical setting and main result}

We state the following structural assumptions for the domain, potential, mobility, physical parameters, and the source term. 
\begin{enumerate}[label=$(\mathrm{A \arabic*})$, ref = $(\mathrm{A \arabic*})$,topsep=0ex]
\item \label{ASS:Omega}
Let $\Omega\subset \erre^d$ be a
bounded domain, with $d\in \{1,2,3\}$, and smooth boundary denoted by $\partial \Omega$, and suppose that $\beta, \eps > 0$ are fixed. 
\item  \label{ASS:pot}
The double-well potential $\psi:\erre \to [0,+\infty)$ is twice differentiable and can be decomposed as $\psi = \psi_1+ \psi_2,$ where $\psi_1 $ is convex, and $\psi_2 $ is twice differentiable and satisfies
\begin{align}
	\label{pot:1}
| \psi_2''(r) | \le c_0, \quad \forall \, r \in \erre ,
\end{align}
for some $c_0>0$. Without loss of generality we also suppose that $\psi_1(0)=\psi_1'(0)=0$. Moreover, it holds that
\begin{align}
	\label{pot:2}
	\forall \eta>0 \, \exists \,C_\eta:  \quad |\psi'(r)|  \leq \eta \psi(r) + C_\eta,
	\quad \forall \, r \in \erre.
	\end{align}
In analogy to \eqref{pot:2}, we assume the following growth condition:
			\begin{align}
			\label{pot:3}
			\forall \eta>0 \, \exists \,C_\eta:  \quad |\psi''(r)|  \leq \eta \psi(r) + C_\eta,
			\quad \forall \, r \in \erre.
		\end{align}
\end{enumerate}
Since $\psi_1'' = \psi'' - \psi_2'' \ge 0$ and $|\psi_2''| \le c_0$, condition \eqref{pot:3} is equivalent to the same condition for $\psi_1''$ in place of $|\psi''|$. It is satisfied by every polynomial double-well potential, and more generally whenever $|\psi''(r)|/\psi(r) \to 0$ as $|r| \to \infty$.	Some growth condition of this type is needed because, in contrast to the case of a two-sided degenerate mobility, the order parameter is not bounded from above under the mobility assumption \ref{ASS:mob:plus}, and respectively not bounded from below under the mobility assumption \ref{ASS:mob}, so that the term $\psi''(\ph)\nabla\ph$ in the flux identity of Theorem~\ref{THM:degenerate:plus} below has to be controlled by the a priori estimates alone.
\begin{enumerate}[label=$(\mathrm{A \arabic*})$, ref = $(\mathrm{A \arabic*})$,topsep=0ex,start=3]
\item   \label{ASS:source}
The source term $\Source_\eps : \erre \to \erre$ is Lipschitz continuous. Hence, there exists $C_S > 0$ such that
\begin{align*}
  |\Source_\eps(r)| \leq C_S (|r| + 1), \quad \forall \, r \in \erre.
\end{align*}
Moreover, it holds that
\begin{align*}
	\Source_\eps(r) \geq 0, \quad \forall r \in (-\infty, -1],
	\quad 
	\text{ and }
	\quad  
	\Source_\eps(r) \leq 0, \quad \forall r \in [1,+\infty).
\end{align*}
\end{enumerate}
We observe that assumptions \ref{ASS:pot}--\ref{ASS:source} encompasses the quartic potential defined in \eqref{quartic}, as well as the source term presented in the preceding section.

The assumption on the mobility function is stated separately, depending on whether the degeneracy occurs at the pure phase $-1$ or $1$.
\begin{enumerate}[label=$(\mathrm{A \arabic* ^+})$, ref = $(\mathrm{A \arabic* ^+})$,topsep=0ex, start=4]
\item  \label{ASS:mob:plus}
The mobility function $m: [-1,+\infty) \to [0,+\infty)$ is 
non-negative, bounded, piecewise continuously differentiable and globally Lipschitz continuous. In particular, we suppose that for positive constants $m_+$ and $m_*$, it holds that $0\leq m(r) \le m_+$ for all $r \in [-1,+\infty) $ and $m(r) \geq m_*$ for all $r \in [0,{+}\infty)$. Moreover, we assume that it degenerates at $-1$ in the following sense: there exists a function $m_0 \in C^1([-1,0);\mathbb{R})$, positive constants $m_-$ and $M$ and an exponent $\esp \geq 1$ such that 
\begin{align}
	\label{def:deg:mob:minus}
	m(r) = (1+r)^\esp m_0(r),
	\quad \forall \, r \in [-1,0),
\end{align}
with 
\begin{align*}
0<m_- \le m_0(r),  \quad \forall r \in [-1,0),
\qquad \|m_0\|_{C^1([-1,0))} \le M{.}
\end{align*}
The mobility is then extended continuously to the whole real line by setting $m(r) = 0$ for $r < -1$.
\end{enumerate}
\begin{enumerate}[label=$(\mathrm{A \arabic* ^-})$, ref = $(\mathrm{A \arabic* ^-})$,topsep=0ex, start=4]
\item \label{ASS:mob}
The mobility function  $m: (-\infty,1] \to [0,+\infty)$  is non-negative, bounded, piecewise continuously differentiable and globally Lipschitz continuous. In particular, we suppose that for positive constants $m_+$ and $m_*$, it holds that $0\leq m(r) \le m_+$ for all $r \in (-\infty,1]$ and $m(r) \geq m_*$ for all $r \in (-\infty,0]$. Moreover, we assume that it degenerates at $1$ in the following sense: there exists a function $m_0 \in C^1((0,1];\mathbb{R})$, positive constants $m_-$ and $M$ and an exponent $\esp \geq 1$ such that 
\begin{align}
	\label{def:deg:mob}
	m(r) = (1-r)^\esp m_0(r),
	\quad \forall \, r \in (0,1],
\end{align}
with
\begin{align*}
0<m_- \le m_0(r),  \quad \forall r \in (0,1],
\qquad \|m_0\|_{C^1((0,1])} \le M{.}
\end{align*}
The mobility is then extended continuously to the whole real line by setting $m(r) = 0$ for $r > 1$.
\end{enumerate}
We note that \ref{ASS:mob:plus} is satisfied by the example mobility 
\eqref{eq:mob} with $m_+>0$ and $m_-=0$, while \ref{ASS:mob} is 
satisfied by \eqref{eq:mob} with $m_->0$ and $m_+=0$; here $m_\pm$ denote the parameters in \eqref{eq:mob} and not the bounds in \ref{ASS:mob:plus} or \ref{ASS:mob}.

Finally, we define the entropy-like function $\Phi$, whose precise form depends on whether \ref{ASS:mob:plus} or \ref{ASS:mob} is adopted. Specifically, we set:
\begin{align*}
	r \mapsto \Phi(r) := \int_{0}^{r} \int_{0}^{s} \frac{1}{m(\xi)} \, d\xi \, ds,
\end{align*}
and have $\Phi : (-1,+\infty) \to [0,+\infty)$ in the \emph{plus-phase model}, whereas $\Phi : (-\infty, 1) \to [0,+\infty)$ in the \emph{minus-phase model}, respectively. In both cases, under either \ref{ASS:mob:plus} or \ref{ASS:mob}, the above properties of $m$ imply that $\Phi$ is non-negative, convex, twice continuously differentiable, and it satisfies for all $r$ in the domains of definition
\begin{align*}
\Phi''(r) = \frac{1}{m(r)},  \qquad
\text{and} \qquad   \Phi'(0) = 0, \quad \Phi(0) = 0.
\end{align*}
Moreover, owing to the uniform bounds on the mobility postulated in \ref{ASS:mob:plus} or \ref{ASS:mob}, it follows that $\Phi$ satisfies a quadratic growth condition and its derivative $\Phi'$ is at most linearly growing for $r\rightarrow +\infty$ under \ref{ASS:mob:plus} and, respectively, for $r\rightarrow -\infty$ under \ref{ASS:mob}. In what follows, we use the simplified notation $S$ for the source term, instead of $S_\varepsilon$, since $\varepsilon$ is fixed. We will revert to the precise notation in the asymptotic analysis carried out in Section~\ref{sec:Sharp}.

We now state the existence result that will be proved for the one-sided degenerate mobility case.

\begin{thm}[Existence for the plus-phase model] \label{THM:degenerate:plus}
Suppose that \ref{ASS:Omega}--\ref{ASS:source}, and \ref{ASS:mob:plus} are fulfilled,
and let the initial datum $\phz $ be such that
\begin{align}
	\label{ass:data}
	\phz \in \Hx1, \quad \text{$\phz \ge -1$\, $a.e.$ in $\Omega$,} \quad 
	\psi(\phz) \in \Lx1, \quad \Phi(\phz) \in \Lx1.
\end{align}
Then, there exists a pair $(\ph, \bfJ)$ such that
\begin{align*}
& \ph \in \H1{\Vp} \cap L^\infty(0,T; H^1(\Omega)) \cap L^2(0,T; H^2(\Omega)), \\
& \ph(0) = \phz \text{ in } \Lx2, \quad \dn \ph  = 0 \text{ a.e. on } \Gamma, \quad 
  \ph \ge -1  \text{ a.e.~in } Q,\\
&
\bfJ \in L^2(0,T; L^2(\Omega)^d) \text{,} \quad  \sqrt{\psi_1''( \ph)} \nabla\ph\in  L^2(Q)^d\text{,} \quad \Phi( \ph) \in  L^\infty(0,T;L^1(\Omega)),
	\\
	& 
\psi(\ph)
	\in  L^\infty(0,T;L^1(\Omega)) \text{,} \quad\psi''(\ph) \in L^1(Q).
\end{align*}
Moreover, these functions satisfy the weak formulation:
\begin{align}\label{weak1}
\int_0^T \langle \partial_t \ph, \xi \rangle_{H^1(\Omega)}
&= \int_Q \bfJ \cdot \nabla \xi
+\int_Q  \Source(\ph) \xi ,
\end{align}
for all  $\xi \in L^2(0,T;H^1(\Omega))$, and  
\begin{align}\label{weak2}
\int_Q \bfJ \cdot \bfeta 
= - \beta \eps \int_Q  \Delta \ph \div\big(m(\ph) \bfeta\big) 
- \frac \beta  \eps \int_Q  m(\ph) \psi''(\ph) \nabla \ph  \cdot \bfeta,
\end{align}
for all $\bfeta \in L^2(0,T; H^1(\Omega)^d) \cap L^\infty(0,T; L^\infty(\Omega)^d)$ with $ \bfeta \cdot \nn = 0$ a.e. on $\Gamma$.
\end{thm}
We remark that the variational equality for the flux function $\bfJ$ comes from the formal definition
\[
\bfJ = -m(\ph) \nabla \mu = -\beta  m(\ph) \nabla \Big( - \varepsilon \Delta \ph + \frac 1\eps \psi'(\ph) \Big).
\]
As the chemical potential $\mu$ cannot be fully controlled in the setting with a degenerate mobility, we eliminate its presence in the variational formulation by using \eqref{sys:2}.

Similarly, provided that assumption \ref{ASS:mob} holds for the mobility, the same result as above can be derived when the degeneracy occurs at the pure phase $1$.

\begin{thm}[Existence for the minus-phase model] \label{THM:degenerate:minus}
Suppose that \ref{ASS:Omega}--\ref{ASS:source}, and \ref{ASS:mob} are fulfilled, and let the initial datum $\phz$ 
be such that
\begin{align}
	\label{ass:data:minus}
	\phz \in \Hx1, \quad \text{$\phz \leq 1$\, $a.e.$ in $\Omega$,
} \quad 
	\psi(\phz) \in \Lx1, \quad \Phi(\phz) \in \Lx1.
\end{align}
Then, there exists a pair $(\ph, \bfJ)$ such that
\begin{align*}
& \ph \in \H1{\Vp} \cap 
L^\infty(0,T; H^1(\Omega)) \cap L^2(0,T; H^2(\Omega)), \\
& \ph(0) = \phz \text{ in } \Lx2, \quad \dn \ph  = 0 \text{ a.e. on } \Gamma, \quad 
 \ph \leq 1  \text{ a.e.~in } Q,\\
&
\bfJ \in L^2(0,T; L^2(\Omega)^d) \text{,} \quad  \sqrt{\psi_1''( \ph) }\nabla\ph\in  L^2(Q)^d\text{,} \quad \Phi( \ph) \in  L^\infty(0,T;L^1(\Omega)),
	\\
& 
	 \psi(\ph)
	  \in L^\infty(0,T;L^1(\Omega))\text{,} \quad\psi''(\ph)\in  L^1(Q) .
\end{align*}
Moreover, these functions satisfy the weak formulation:
\begin{align*}
\int_0^T \langle \partial_t \ph, \xi \rangle_{H^1(\Omega)}
&= \int_Q \bfJ \cdot \nabla \xi
+\int_Q  \Source(\ph) \xi ,
\end{align*}
for all $\xi \in L^2(0,T;H^1(\Omega))$, and 
\begin{align*}
\int_Q \bfJ \cdot \bfeta 
= - \beta \eps \int_Q  \Delta \ph \div\big(m(\ph) \bfeta\big) 
- \frac \beta  \eps \int_Q  m(\ph) \psi''(\ph) \nabla \ph  \cdot \bfeta,
\end{align*}
for all $\bfeta \in L^2(0,T; H^1(\Omega)^d) \cap L^\infty(0,T; L^\infty(\Omega)^d) $ with $ \bfeta \cdot \nn = 0$ a.e. on $\Gamma$.
\end{thm}

These results are obtained by adapting the ideas of \cite{ElliottG96} (see also \cite{garcke2025anisotropic,ebenbeck2021cahn}) to handle the one-sided degenerate mobility $m$. The main strategy is to approximate the degenerate mobility by non-degenerate mobilities $\mdel$ indexed by $\delta \in (0,1]$, so that the resulting problem becomes regular and classical theory can be applied. Within this regularised framework, one can then apply previously established well-posedness results, which are provided in \cite{OOL} in our specific case. Several regularity estimates are derived in this setting, including an entropy estimate that involves a suitable regularization $\Phi_\delta$ of the entropy function $\Phi$  defined above (see \cite{ElliottG96}), which allows passing to the limit as $\delta \to 0$ and recovering a solution for the original problem. Furthermore, we note that the regularity of the mobility required in \ref{ASS:mob:plus} and \ref{ASS:mob} is essential to ensure that the weak formulation is mathematically rigorous.

We note that Theorems \ref{THM:degenerate:plus} and \ref{THM:degenerate:minus} share the same structure and can be proved by essentially identical arguments. Consequently, we will henceforth provide an extensive proof of just Theorem \ref{THM:degenerate:plus}. 
We adopt the previously outlined strategy, starting with the formulation of the approximate problem.

\subsection{Approximation scheme}
Let $\delta \in (0,1]$. We introduce a positive regularised mobility $\mdel:\mathbb{R}\to[0,+\infty)$ defined by
\begin{align*}
\mdel(r) :=
\begin{cases}
m(-1 + \delta), & \text{if } r \le -1 + \delta,\\
m(r), & \text{if } r > -1 + \delta,
\end{cases}
\end{align*}
and we define the corresponding regularised entropy function 
$\Phi_\delta:\mathbb{R}\to[0,+\infty)$ by setting:
\begin{align*}
\Phi_\delta(r) := \int_0^r \int_0^s \frac{1}{\mdel(\xi)}\, d\xi\, ds, \quad \forall r \in \mathbb{R}.
\end{align*}
Owing to the above definition of $\mdel$, it follows that $\Phi_\delta \in C^2(\mathbb{R})$, is non-negative and convex, and satisfies
\begin{align}\label{Phi:delta:prop}
0<\Phi_\delta''(r) = \frac{1}{\mdel(r)}, \quad \forall  r \in \mathbb{R}, 
\qquad
\Phi_\delta(r) = \Phi(r), \quad \forall  r \geq -1+\delta,
\end{align}
together with the normalization properties $\Phi_\delta(0)=0, $ and $ \Phi_\delta'(0)=0$. Moreover, there exists a constant $c_\delta>0$ depending on $\delta$ such that
\begin{align*}
|\Phi_\delta''(r)| \le c_\delta, 
\qquad
|\Phi_\delta'(r)| \le c_\delta(|r|+1),
\qquad \forall r \in \mathbb{R}.
\end{align*}
Thus, the approximated system to \eqref{eq:sys21} reads as: 
\begin{subequations}\label{eq:approx:sys}
\begin{alignat}{2}
	\label{sys:approx:1}
	\dt \phd  &= \div(\mdel(\phd)\nabla  \mud) + \Source(\phd)
	\qquad && \text{in $Q  $,}
	\\
	\label{sys:approx:2}
	 \mud &= \beta (-\eps \Delta \phd + \tfrac 1 \eps \psi'(\phd))
	\qquad && \text{in $Q$,}
	\\
	\label{sys:approx:3}
	\dn \mud  &= \dn \phd =0
	\qquad && \text{on $\Gamma  $,}
	\\
	\label{sys:approx:4}
	 \phd(0) &=\phz
	\qquad && \text{in $\Omega$.}
\end{alignat}
\end{subequations}
\Accorpa\Sysapprox {sys:approx:1} {sys:approx:4}
With this  choice of $\mdel$, for every $\delta \in (0,1]$, there exists a weak solution $(\phd, \mud)$
solving \eqref{eq:approx:sys} in the corresponding natural weak setting with
\begin{equation}\label{weaksoln:reg}
	\begin{aligned}
	\phd & \in \H1 \Vp \cap \L\infty {\Hx1} \cap \L2 {\Hx2},
	\\
	\mud & \in \L2 {\Hx1},
\end{aligned}
\end{equation}
as guaranteed by \cite[Thm. 3.1]{OOL}.

\begin{remark} We point out some inaccuracies in the proof of \cite[Thm. 3.1]{OOL} that do not affect the validity of the result. The proof of \cite[Thm. 3.1]{OOL} derives entropy estimates for Galerkin solutions $(\ph_n,\mu_n)$ to a Galerkin approximation of \eqref{eq:approx:sys} by testing the analogue of \eqref{sys:approx:1} with $\Phi_\delta'(\ph_n)$, which however is not an admissible test function at the Galerkin level. The rigorous procedure is to consider a two-step approximation where we first apply an affine linear truncation $\psi_{1,M}$ of $\psi_1$, defined for $M \geq 1$ as follows:
	\[
	\psi_{1,M}(r) = \begin{cases}
		\psi_1(M) + \psi_1'(M) (r-M) & \text{ for } r > M, \\
		\psi_1(r) & \text{ for } |r| \leq M, \\
		\psi_1(-M) + \psi_1'(-M)(r+M) & \text{ for } r < -M,
	\end{cases}
	\]
	and set $\psi_M := \psi_{1,M} + \psi_2$ and $\mathcal{E}_M(\ph) :=  \iO \big( \frac{\beta\eps}{2} |\nabla \ph|^2 + \frac\beta\eps \psi_M(\ph) \big)$. Due to the Lipschitz nonlinearity, a standard Galerkin approach  without the need of entropy estimates yields the existence of a weak solution $(\ph_{M,\delta}, \mu_{M,\delta})$ with regularities stated in \eqref{weaksoln:reg} to the following approximation problem:
	\begin{subequations}
		\begin{alignat}{2}
		\label{Galerkin:2:ph} \dt \ph_{M,\delta}  &= \div(\mdel(\ph_{M,\delta})\nabla  \mu_{M,\delta}) + \Source(\ph_{M,\delta})
		\qquad && \text{in $Q  $,}
		\\
		\label{Galerkin:2:mu} \mu_{M,\delta} &= \beta (-\eps \Delta \ph_{M,\delta} + \tfrac 1 \eps \psi_M'(\ph_{M,\delta}))
		\qquad && \text{in $Q$,}
		\\
		\dn \mu_{M,\delta}  &= \dn \ph_{M,\delta} =0
		\qquad && \text{on $\Gamma  $,}
		\\
		\ph_{M,\delta}(0) &=\phz
		\qquad && \text{in $\Omega$.}
		\end{alignat}
	\end{subequations}
A priori estimates for the Galerkin solutions can be obtained by differentiating $Y_n := \E_M(\ph_n) + \frac K2 \|\ph_n\|^2_{L^2(\Omega)}$, with $K:=\frac{2\beta c_0}{\eps} + 1$ and $c_0$ being defined in \eqref{pot:1}, with respect to time and using the equations on the Galerkin level. For fixed $\delta> 0$, $\Phi_\delta'(\ph_{M,\delta})$ is now an admissible test function for \eqref{Galerkin:2:ph}, and we can follow the proof of \cite[Thm. 3.1]{OOL} to derive a priori estimates uniform in $M$. The first is the entropy estimate obtained from testing \eqref{Galerkin:2:ph} with $\Phi_\delta'(\ph_{M,\delta})$ and \eqref{Galerkin:2:mu} with $- \Delta \ph_{M,\delta}$, leading to
	\[
	\norma{\Phi_\delta(\ph_{M,\delta})}_{L^\infty(0,T;L^1(\Omega))} + \norma{\ph_{M,\delta}}_{L^\infty(0,T;L^2(\Omega)) \cap L^2(0,T;H^2(\Omega))} \leq C_\delta
	\]
for a positive constant $C_\delta$ independent of $M$. Then, we perform a standard energy estimate by testing \eqref{Galerkin:2:ph} with $\mu_{M,\delta}$ and \eqref{Galerkin:2:mu} with $- \dt \ph_{M,\delta}$, and summing the resulting identities to obtain
		\begin{align*}
		& \frac d{dt} {{\cal E}_M}
		(\ph_{M,\delta})
		+ \iO m_\delta (\ph_{M,\delta})
		 |\nabla \mu_{M,\delta}|^2
		= \iO \Source(\ph_{M,\delta}) \mu_{M,\delta}
	\\
	& \quad = \iO \Source(\ph_{M,\delta}) (\mu_{M,\delta}-(\mu_{M,\delta})_\Omega)
		+ \iO \Source(\ph_{M,\delta})(\mu_{M,\delta})_\Omega = I_1 + I_2.
	\end{align*}
Now, employing the Young inequality and \ref{ASS:source}, we infer that
	\begin{align*}
		|I_1| &\leq C ( \norma{\ph_{M,\delta}}_{L^2(\Omega)}+1)\norma{\nabla \mu_{M,\delta}}_{L^2(\Omega)}
		\leq 
		\frac{c_\delta}2 \norma{\nabla \mu_{M,\delta}}^2_{L^2(\Omega)}
		+ C_\delta(\norma{\ph_{M,\delta}}^2_{L^2(\Omega)}+1) \\
			& \leq  \frac{c_\delta}2
 \norma{\nabla \mu_{M,\delta}}^2_{L^2(\Omega)}
			+ C_\delta,
		\\
		|I_2| & \leq C
		( \norma{\ph_{M,\delta}}_{L^2(\Omega)}+1)|(\mu_{M,\delta})_\Omega|
		\leq
		C |(\mu_{M,\delta})_\Omega|,
	\end{align*}	
	with $c_\delta := \min \{
	\delta^\theta m_-, m_*
	\} \leq \mdel (r) $ for all $r \in \erre.$ To control the mean value, we observe that \eqref{pot:2} is not inherited by $\psi_M$ with constants independent of $M$. However, for every $\eta > 0$ there exists a constant $C_\eta > 0$, independent of $M$, such that
	\begin{align}
		\label{pot:2:M}
		|\psi_M'(r)| \le \eta \, \psi_{1,M}(r) + C_\eta \big( 1 + r^2 \big), \qquad \forall \, r \in \erre .
	\end{align}
Indeed, \eqref{pot:1} gives $|\psi_2'(r)| \le C(1+|r|)$ and $|\psi_2(r)| \le C(1+r^2)$. For $|r| \le M$ we have $\psi_M = \psi$ and $\psi_{1,M} = \psi_1$, so that \eqref{pot:2:M} follows from \eqref{pot:2} and $\psi(r) \le \psi_1(r) + C(1+r^2)$. For $r > M$ we have $\psi_{1,M}'(r) = \psi_1'(M) \ge 0$ and, again by \eqref{pot:2},
	\begin{align*}
		\psi_1'(M) = \psi'(M) - \psi_2'(M)
		\le \eta \, \psi_1(M) + C_\eta \big( 1 + M^2 \big)
		\le \eta \, \psi_{1,M}(r) + C_\eta \big( 1 + r^2 \big),
	\end{align*}
	since $\psi_{1,M}$ is non-decreasing on $[0,+\infty)$; the case $r < -M$ is analogous. Subsequently, integrating \eqref{Galerkin:2:mu} over $\Omega$, and using $\dn \ph_{M,\delta} = 0$ and \eqref{pot:2:M} with $\eta = 1$, we obtain
	\begin{align*}
		|(\mu_{M,\delta})_\Omega|
		= \frac{\beta}{\eps |\Omega|} \, \Big| \iO \psi_M'(\ph_{M,\delta}) \Big|
		& \le C \big( \|\psi_{1,M}(\ph_{M,\delta})\|_{\Lx1} + \|\ph_{M,\delta}\|_{\Lx2}^2 + 1 \big)
		\\
		& \le C_\delta \big( \|\psi_{1,M}(\ph_{M,\delta})\|_{\Lx1} + 1 \big).
	\end{align*}
	Moreover, $0 \le \psi_{1,M}(r) \le \psi_1(r)$ and $\psi_2(r) \ge - C(1+r^2)$ imply $\E_M(\phz) \le \E(\phz)$ and
	\begin{align*}
		\E_M(\ph_{M,\delta}) \ge \frac{\beta\eps}{2} \|\nabla \ph_{M,\delta}\|_{\Lx2}^2
		+ \frac{\beta}{\eps} \|\psi_{1,M}(\ph_{M,\delta})\|_{\Lx1} - C_\delta.
	\end{align*}
Thus, combining the previous estimate with Gr\"onwall’s lemma, we conclude that
	\begin{align*}
		\norma{\ph_{M,\delta}}_{\L\infty {\Hx 1}}
		+\norma{\psi_{1,M}(\ph_{M,\delta})}_{\L\infty {\Lx 1}}
		+ \norma{\nabla \mu_{M,\delta}}_{\L2 {{\Lx2}^d}}
		\leq C_\delta.
	\end{align*}
	This also implies $(\mu_{M,\delta})_\Omega \in L^\infty(0,T)$ and by the Poincar\'e--Wirtinger inequality we have
		\[
		\norma{\mu_{M,\delta}}_{\L2{\Hx 1}} \leq C_\delta.
		\]
	Lastly, by comparison of terms in \eqref{Galerkin:2:ph} we deduce also that 
	\[
	\norma{\dt \ph_{M,\delta}}_{\L2{\Vp}} \leq C_\delta.
	\]
	Passing to the limit $M \to \infty$ then yields limit functions $(\phd, \mud)$ satisfying the regularities stated in \eqref{weaksoln:reg} and is a weak solution to the approximate system \eqref{eq:approx:sys}. 
\end{remark}

\subsection{Uniform estimates}
The next step is to derive estimates that are uniform with respect to $\delta$.

\begin{prop}[Uniform estimates] \label{LEM:apriori_delta}
Suppose that \ref{ASS:Omega}--\ref{ASS:source}, and \ref{ASS:mob:plus} are fulfilled. There exists $\delta_0 \in (0,1]$ such that, for all $\delta \in (0, \delta_0]$, the weak solution 
$(\phd, \mud)$ introduced above satisfies the following estimates, 
with a constant $C$ independent of $\delta$:  
\begin{align*}
& (i) \quad  \mathop{\textrm{ess sup}}_{0 \le t \le T} \Big( \|\phd(t)\|_{\Hx1}^2  
+ \|\psi_{1}(\phd(t))\|_{\Lx1} 
+ \|\Phi_\delta(\phd(t))\|_{\Lx1} 
\Big) 
\\
& \quad \qquad 
+ \int_0^T 
\Big( 
\|\sqrt{\mdel(\phd)} \nabla \mud\|_{\Lx2^d}^2 
 + \|\Delta \phd\|_{\Lx2}^2 
 +\norma{\sqrt{\psi_{1}''(\phd)} \nabla \phd}^2_{\Lx2^d}
 \Big)
  \le C.
\\
&(ii) \quad 
\norma{\dt \phd}_{\L2 \Vp} +
\norma{\bfJ_\delta}_{\L2 {\Lx2^d} } \le C, \quad \text{ with } \bfJ_\delta := -\mdel(\phd) \nabla \mud.
\\
&(iii)\quad  \mathop{\textrm{ess sup}}_{0 \le t \le T} \int_\Omega \big(-\phd(t) - 1 \big)^2_+ \le C \delta^\esp,
\end{align*}
where $\esp$ is the exponent in \eqref{def:deg:mob:minus}, and $(\cdot)_+$ denotes the positive part.
\end{prop}

\begin{proof}

We begin with the preliminary observation that, throughout the proof, the phrase ``$\delta$ sufficiently small'' is understood in the following precise sense: there exists $\delta_0 \in (0,1]$ such that the stated property holds for all $\delta \in (0,\delta_0]$, as specified in each occurrence.

We then proceed with an entropy-type estimate  along the lines of \cite{ElliottG96} (see also \cite{garcke2025anisotropic,ebenbeck2021cahn}). 
First, owing to the above construction, notice that $\Phi_\delta''(\phd) $ is uniformly bounded. We then test \eqref{sys:approx:1} with $\Phi'_\delta(\phd)$ and \eqref{sys:approx:2} with $- \Delta \phd$, integrate by parts, rearrange the terms, and sum the resulting equations to obtain:
\begin{align*}
	 & \frac d {dt} \norma{\Phi_\delta(\phd)}_{\Lx1}
	+ \beta \eps \norma{\Delta \phd}^2_{\Lx2}
	+ \frac \beta \eps \iO \psi_{1}''(\phd)|\nabla \phd|^2
	\\  &\quad 
	= \iO \Source(\phd) \Phi_\delta'(\phd)
	- \frac \beta \eps \iO \psi_2''(\phd)|\nabla \phd|^2,
\end{align*}
for almost every $t \in (0,T)$. Next, integrating over $(0,t)$, for an arbitrary $t \in (0,T)$, and rearranging terms, we infer that 
\begin{equation}\label{est:entropy:1}
\begin{aligned}
& \int_\Omega \Phi_\delta(\phd(t))
+ \beta\eps \int_0^t \int_\Omega |\Delta \phd|^2
+ \frac \beta\eps \int_0^t \int_\Omega \psi_{1}''(\phd)\,|\nabla \phd|^2 \\
&\quad 
= \int_\Omega \Phi_\delta(\ph_{0})
- \frac \beta\eps \int_0^t \int_\Omega \psi_2''(\phd)\,|\nabla \phd|^2
+ \int_0^t \int_\Omega \Source(\phd)\Phi_\delta'(\phd).
\end{aligned}
\end{equation}
We define the function $G$ such that 
\[
G(0) = 0, \quad G'(0) = 0, \quad G''(r) = (1+r)^{-\esp},
\]
so that for all $r \in (-1,0]$ the following relations hold between $G''$ and $\Phi''$:
\[
G''(r) = (1+r)^{-\esp} = \frac{m_0(r)}{m(r)} = m_0(r) \Phi''(r) \leq M \Phi''(r).
\]
Then, for $r \in (-1+\delta,0]$, we have
\begin{align*}
\Phi_\delta''(r) = \Phi''(r) = \frac{1}{m(r)} = \frac{(1+r)^{-\esp}}{m_0(r)} \leq \frac{G''(r)}{m_-}.
\end{align*}
On the other hand, for $r \in (-1,-1+\delta]$, we estimate as follows:
\begin{align*}
\Phi_\delta''(r) = \frac{1}{m(-1+\delta)} = \frac{(1+(-1+\delta))^{-\esp}}{m_0(-1+\delta)} \leq \frac{(1+r)^{-\esp}}{m_-} = \frac{G''(r)}{m_-}.
\end{align*}
Hence, we find that for all $r \in (-1,0]$ it holds that 
\[
\Phi_\delta''(r) \leq \frac{1}{m_-} G''(r) \leq \frac{M}{m_-} \Phi''(r).
\]
Integrating twice and using  $\Phi_\delta(0) = \Phi_\delta'(0) = 0$ and $\Phi(0) = \Phi'(0) = 0$ yields
\[
\Phi_\delta(r) \leq \frac{M}{m_-} \Phi(r), \quad \forall r \in [-1,0]{.}
\]
Then, combining with the second property in \eqref{Phi:delta:prop} and the assumption on the initial datum in \eqref{ass:data}, we have for the first term on the \rhs\ of \eqref{est:entropy:1}
\begin{align*}
\int_\Omega \Phi_\delta(\ph_0) = \int_{\{ \ph_0(x) <-1+\delta\}} \Phi_\delta(\ph_0) + \int_{\{ \ph_0(x) \geq -1+\delta\}} \Phi(\ph_0) \leq \Big (\frac{M}{m_-} + 1 \Big ) \int_\Omega \Phi(\ph_0) \leq C.
\end{align*}
The second integral on the \rhs\ of \eqref{est:entropy:1} can be estimated using interpolation as
\begin{align*}
& \left| \int_0^t \int_\Omega \psi_2''(\phd)\,|\nabla \phd|^2 \right|
\le \|\psi_2''\|_{L^\infty(\erre)} 
\int_0^t \int_\Omega |\nabla \phd|^2 
\\  & \quad 
\le \frac{\beta\eps}2 \int_0^t \int_\Omega |\Delta \phd|^2 + C \int_0^t (\norma{\phd}^2_{\Lx2}+1),
\end{align*}
for all $t \in (0,T)$.

Next, we take note that the convexity of $\Phi_\delta$ (hence monotonicity of $\Phi_\delta'$), together with $\Phi_\delta'(0) = 0$, implies $\Phi_\delta'(r) \ge 0$ for $r \ge 0$ and $\Phi_\delta'(r) \le 0$ for $r \le 0$. To handle the third integral on the \rhs\ of \eqref{est:entropy:1}, we first establish the following bound holding uniformly in $\delta \in (0,1]$:
\begin{align}
		\label{Phi:delta:lin}
		(1+r)\, |\Phi_\delta'(r)| \le C_\Phi \big( 1 + \Phi_\delta(r) \big), \qquad \forall \, r \in [-1,0],
\end{align}
with a positive constant $C_\Phi$ depending only on $\esp$, $m_-$ and $M$. We aim to study the behavior of $\Phi_\delta'$ near the degeneracy point $r=-1$. Since $m(r) \to 0$ as $r \to -1$, it holds that $\Phi_\delta'(r) $ is not uniformly bounded in $\delta$ as $r \to -1$, so that the factor $(1+r)$ in \eqref{Phi:delta:lin} is essential. To carry out this analysis, it is convenient to perform a change of variable that localizes the singularity at the origin. To show this, we define for $t \in (0,1]$ the function $g(t) := \Phi_\delta(-1+t)$, so that $g(1) = g'(1) = 0$ and $g''(t) = \frac{1}{\mdel(-1+t)}$. Under this transformation, \eqref{Phi:delta:lin} reduces to showing that
\begin{align}\label{g:bound:relation}
	t\, |g'(t)| \le C_\Phi (1 + g(t)), \qquad \forall \, t \in (0,1].
	\end{align}
To this end, we aim to control the \lhs\ in terms of $g$ and an antiderivative of $g''$. By \ref{ASS:mob:plus} we have $m_- (1+r)^\esp \le m(r) \le M (1+r)^\esp$ for $r \in [-1,0)$, and hence, by the definition of $\mdel$,
	\begin{align}
		\label{g:bounds}
		g''(u) \le \frac{u^{-\esp}}{m_-}, \quad \forall u \in (0,1],
		\qquad
		g''(u) \ge \frac{u^{-\esp}}{M}, \quad \forall u \in [\delta,1].
	\end{align}
Non-negativity of $\Phi_\delta$ and Taylor's formula with remainder term at $t = 1$ gives
	\begin{align}
		\label{g:taylor}
		|g'(t)| = \int_t^1 g''(u)\, du, \qquad g(t) = \int_t^1 (u-t)\, g''(u)\, du \qquad \text{for } t \in (0,1].
	\end{align}
Let $\tau := \max\{t, \delta\}$. The introduction of $\tau$ serves to unify the two distinct regimes: the truncation regime where $t \leq \delta$ and $m_\delta$ is constant, and the  regime where $t > \delta$ and $m_\delta$ coincides with $m$. Since $u - t \ge t$ for $u \ge 2t$, \eqref{g:taylor} yields
	\begin{equation}\label{g:bounds:2}
	\begin{aligned}
		t\, |g'(t)| &\le t \int_t^{\min\{2t,1\}} g''(u)\, du + \int_{\min\{2t,1\}}^1 (u-t)\, g''(u)\, du
		\\
		&\le t^2 \Big (\max_{u \in [t,\, \min\{2t,1\}]} g''(u)\Big ) + \int_{t}^1 (u-t) \, g'' (u) \, du \le \frac{\tau^{2-\esp}}{m_-} + g(t),
	\end{aligned}
	\end{equation}
where we used the property $g''(u) \le \frac{\tau^{-\esp}}{m_-}$ on $[t,\min\{2t,1\}]$ by \eqref{g:bounds}. Indeed, if $u < \delta$ (which forces $\tau = \delta$), as $g''(u) = g''(\delta) \leq \frac{\delta^{-\theta}}{m_-}= \frac{\tau^{-\theta}}{m_-}$, by \eqref{g:bounds}. On the other hand if $u \geq \delta$, then as $u \in [t,2t]$ we have $u \geq \tau$, so that $g''(u) \leq \frac{u^{-\theta}}{m_-} \le \frac{\tau^{-\esp}}{m_-}$ on $[t,2t]$ as well.
		
From \eqref{g:bounds:2}, if $\esp \le 2$, then $\tau^{2-\esp} \le 1$ and this leads to \eqref{g:bound:relation}. Otherwise, if $\esp > 2$ and $2 \tau > 1$, we simply have $\tau^{2-\theta} \leq 2^{\theta-2}$ which leads to \eqref{g:bound:relation}. On the other hand, if $\esp > 2$ and $2\tau \le 1$, then, as $g$ is non-increasing, $g(t) \ge g(\tau)$, together with \eqref{g:bounds} and \eqref{g:taylor} give us
	\begin{align*}
		g(\tau) \ge \int_{2\tau}^1 (u-\tau)\, g''(u)\, du \ge \frac{1}{2M} \int_{2\tau}^1 u^{1-\esp}\, du = \frac{(2\tau)^{2-\esp} - 1}{2M(\esp-2)}.
	\end{align*}
Rearranging leads to 
	\[
	\tau^{2-\esp} \le 2^{\esp-2} \big( 1 + 2M(\esp-2)\, g(t) \big),
	\]
so that we obtain \eqref{g:bound:relation} from \eqref{g:bounds:2}. In all cases the desired bound \eqref{Phi:delta:lin} follows, with a positive constant $C_\Phi$ depending only on $\esp$, $m_-$ and $M$.
		
We now turn to the last term on the right-hand side of \eqref{est:entropy:1}. Let $L_S$ denote the Lipschitz constant of $\Source$. We ignore the time integral for now and split
	\begin{align*}
		\int_\Omega \Source(\phd)\Phi_\delta'(\phd)
		= \bigg( \int_{\{\phd < -1\}} + \int_{\{-1 \le \phd \le 0\}} + \int_{\{\phd > 0\}} \bigg) \Source(\phd)\Phi_\delta'(\phd)
		=: \JJJ_1 + \JJJ_2 + \JJJ_3 .
	\end{align*}	On $\{\phd < -1\}$ we have $\Source(\phd) \ge 0$ by \ref{ASS:source} and $\Phi_\delta'(\phd) \le 0$, so that $\JJJ_1 \le 0$. On $\{-1 \le \phd \le 0\}$ we have $\Phi_\delta'(\phd) \le 0$. At points where $\Source(\phd) \ge 0$ the integrand is non-positive, and at points where $\Source(\phd) < 0$ the Lipschitz continuity of $\Source$ and $\Source(-1) \ge 0$ give $-\Source(\phd) \le \Source(-1) - \Source(\phd) \le L_S (1+\phd)$. Hence, by \eqref{Phi:delta:lin},
	\begin{align*}
		\JJJ_2 \le L_S \int_{\{-1 \le \phd \le 0\}} (1+\phd)\, |\Phi_\delta'(\phd)| \le L_S C_\Phi \int_\Omega \big( 1 + \Phi_\delta(\phd) \big).
	\end{align*}
On $\{\phd > 0\}$, finally, $\Phi_\delta'(\phd) = \Phi'(\phd) = \int_0^{\phd} \frac{1}{m(s)}ds \in [0, \frac{\phd}{m_*}]$, so that $|\JJJ_3| \le \frac{C_S}{m_*} \int_\Omega (1+\phd)\, \phd \le C \int_\Omega (1 + |\phd|^2)$. Let us highlight that the sign conditions in \ref{ASS:source} eliminate the need for a truncation of the source term, as employed in other works (see, e.g., \cite{ebenbeck2021cahn,EbenbeckL21}). Thanks to the uniform quadratic growth from below of $\Phi_\delta$ with respect to $\delta$: 
	\[
	\Phi_\delta(r) \geq \int_0^r \int_0^s \frac{1}{m_+} \, d\xi \, ds = \frac{1}{2 m_+} r^2, \quad \forall r \in \mathbb{R},
	\]
the right-hand side of \eqref{est:entropy:1} can be bounded as
\begin{equation}\label{est:entropy:2}
\begin{aligned}
	\non
 & \int_\Omega \Phi_\delta(\phd(t))
+ \frac{\beta\eps}{2} \int_0^t \int_\Omega |\Delta \phd|^2
+ \frac \beta\eps \int_0^t \int_\Omega \psi_{1}''(\phd)\,|\nabla \phd|^2 \notag
\\ & \quad 
 \leq 
 \int_\Omega \Phi_\delta(\ph_0) +  C \int_0^t (\norma{\phd}^2_{\Lx2}+1)+ C \int_0^t \int_\Omega (\Phi_\delta(\phd) + 1)
 \\ & \quad  \leq 
 C + C \int_0^t \int_\Omega (\Phi_\delta(\phd) + 1),
\end{aligned}
\end{equation}
with $C$ independent of $\delta$. We apply Gr\"onwall's lemma to conclude that
\begin{align}	\label{est:approx:0}
	\norma{\Phi_\delta(\phd)}_{\L\infty {\Lx1}} + \norma{\Delta \phd}_{\L2 {\Lx2}} 
	 +\norma{\sqrt{\psi_{1}''(\phd)} \nabla \phd}_{L^2(Q)^d} \le C.
\end{align}
Next, elliptic regularity theory yields
\begin{align}
	\norma{\phd}_{\L\infty {\Lx2} \cap \L2 {\Hx2}} \le C.
	\label{est:approx:1}
\end{align}
Finally, we can operate the standard energy estimate by testing \eqref{sys:approx:1} with $\mud$ and \eqref{sys:approx:2} with $-\dt \phd$. Summing and using \eqref{sys:approx:2} leads us to
\begin{align*}
	\frac d{dt} {\cal E}(\phd)
	+ \norma{\sqrt{\mdel(\phd)}\nabla \mud}^2_{\Lx2^d}
	= \iO \Source(\phd) (\beta (-\eps \Delta \phd + \tfrac 1 \eps \psi' (\phd)).
\end{align*}
To control the term on the \rhs, we notice that 
\begin{align*}
	& \iO \Source(\phd) \big(\beta (-\eps \Delta \phd + \tfrac 1 \eps \psi' (\phd)\big)
	\leq C (\norma{ \phd}^2_{\Hx2}+1)
	+ C\iO \Source(\phd)\psi' (\phd).
\end{align*}%
Next, to control the last integral on the \rhs,
we proceed arguing similarly as above, using the sign properties of the source term $S$ and the convex part $\psi_1$ of the potential. Namely, recalling \ref{ASS:pot} and \ref{ASS:source}, we observe that $S(r) \psi_1'(r) \leq 0$ for all $r \notin [-1,1]$. 
 Then,
\begin{align*}
	& \iO \Source(\phd)\psi' (\phd) = \iO \Source(\phd) \psi_1'(\phd) + \iO \Source(\phd) \psi_2'(\phd)
    \\ & \quad \leq \Big |\int_{\{\phd(x) \in [-1,1]\}} S(\phd) \psi_1'(\phd) \Big | + \Big |\iO S(\phd) \psi_2'(\phd) \Big | \leq C (\norma{\phd}^2_{\Lx2}+1) \leq C.
\end{align*}
Applying the bound \eqref{est:approx:1} once more and gathering the above estimates, we readily obtain
\begin{align}
	\label{est:approx:3}
	\norma{\phd}_{\L\infty {\Hx1}}
	+\norma{\psi_{1}(\phd)}_{\L\infty {\Lx1}}
	+ \norma{\sqrt{\mdel(\phd)}\nabla \mud}_{\L2 {\Lx2^d}}
	\le C.
\end{align}
Combining \eqref{est:approx:0}--\eqref{est:approx:3} concludes the proof of $(i)$.

Once estimate $(i)$ is established, estimate $(ii)$ follows directly by observing that
\begin{align*}
	\norma{\bfJ_\delta}_{\L2 {\Lx2^d} }^2 \le  m_+ \int_Q \mdel(\phd) |\nabla \mud|^2 \leq C,
\end{align*}
and, from the weak formulation of \eqref{sys:approx:1}, by comparison of terms we also gain
\begin{align}
	\label{est:approx:4}
	\norma{\dt \phd}_{\L2 \Vp} \le C.
\end{align}

Next, combining estimates $(i)$ and $(ii)$, we obtain a quantitative estimate on the deviation of $\phd$ from the pure phase $-1$, with an explicit dependence on $\delta$, as stated in estimate $(iii)$. To make this claim precise, we perform a Taylor expansion of $\Phi_\delta$ around $-1+\delta$. For any $z<-1$, there exists $\xi \in [z, -1+\delta]$ such that
\begin{align}
	\label{def:Phi}
\Phi_\delta(z)
= \Phi_\delta(-1+\delta)
+ \Phi_\delta'(-1+\delta)(z+1-\delta)
+ \frac12 \Phi_\delta''(\xi)(z+1-\delta)^2.
\end{align}
Furthermore, using the above definition of $\mdel$ and $\Phi_\delta$, we realise that
\begin{align*}
\Phi_\delta(-1+\delta) & = \Phi(-1+\delta) \ge 0, \quad
\Phi_\delta'(-1+\delta) = \Phi'(-1+\delta) \le 0,
\end{align*}
and hence the first two summands on the right-hand side of \eqref{def:Phi} are non-negative.
Moreover, since $\xi \in [z, -1+\delta]$, it holds that
\begin{align*}
\Phi_\delta''(\xi) & = \frac{1}{\mdel(\xi)} = \frac{1}{m(-1+\delta)}.
\end{align*}
Collecting the above properties yields 
\begin{align*}
	\Phi_\delta(z)
& \ge \frac{(z+1-\delta)^2}{2m(-1+\delta)}
\ge \frac{(z+1)^2}{2(1-1+\delta)^\esp
M}
 = \frac1{2 \delta^{\esp} 
 	M}\big(-z-1\big)^2, \quad 
 \forall z<-1,
\end{align*}
whence we deduce that
\begin{align*}
	(-z - 1)^2 \le 2\delta^\esp M \Phi_\delta(z) \leq C \delta^\esp \, \Phi_\delta(z), \quad \forall z \le -1.
\end{align*}
Finally, upon integration and a small manipulation, we conclude using estimate $(i)$ that
\begin{align*}
\mathop{\textrm{ess sup}}_{0 \le t \le T} \int_\Omega (-\phd(t) - 1)^2_+ 
\le C \delta^\esp \, \mathop{\textrm{ess sup}}_{0 \le t \le T} \int_\Omega \Phi_\delta(\phd(t)) 
\le C \delta^\esp,
\end{align*}
which proves estimate $(iii)$.
\end{proof}

The next step consists of passing to the limit as $\delta \to 0$, and showing that, up to subsequences, the resulting limit yields a solution to the limit problem.

\subsection{Passing to the limit $\delta \to 0$}

With the estimates above at our disposal and using classical compactness results, we can now pass to the limit $\delta \to 0$.  First, for $\delta \in (0,\delta_0]$,  we have from Proposition \ref{LEM:apriori_delta} that up to non-relabelled subsequences as $\delta \to 0$, there exists a limiting pair $(\ph,\bfJ)$ such that
\begin{subequations}
\begin{alignat}{3}
\notag \phd & \to \ph \quad \text{weakly-star in } \H1 {\Vp} \cap \L\infty {\Hx1} \cap \L2 {\Hx2}, 
	\\
\label{conv:strong:phd}	&  \hspace{1.75cm} \text{and strongly in } \C0 {\Lx \sigma}  \cap \L2 {W^{1, \sigma}(\Omega)},
		\\
\bfJ_\delta & \to \bfJ \quad \text{weakly  in }   \L2 {L^2(\Omega;\erre^d)}, 
\end{alignat}
\end{subequations}
for every $\sigma \in [1,6)$ in three dimensions and for every $\sigma \in [1,\infty)$ in one and two dimensions, where we recall that $\bfJ_\delta = -\mdel(\phd) \nabla \mud$.  
Here, we have employed the Banach–Alaoglu theorem together with classical compactness results, such as the Aubin–Lions lemma (see, e.g., \cite[Sect.~8, Cor.~4]{Simon87}). Moreover, from estimate $(iii)$ of Proposition~\ref{LEM:apriori_delta}, it readily follows that 
\begin{align*}
	\ph \geq -1 \quad \text{a.e.~in } Q.
\end{align*}
At this stage, we formulate the weak version of system \eqref{eq:approx:sys}, incorporating the flux term $\bfJ_\delta$.
As $S$ grows linearly passing to the limit in the approximate version of \eqref{weak1} is standard. As in \cite{ElliottG96} we introduce a weak formulation involving  the flux term. By \eqref{sys:approx:2}, \eqref{sys:approx:3} and $\phd \in \L2{\Hx2}$ we have $\mud = \beta(-\eps\Delta\phd + \tfrac1\eps\psi'(\phd))$ a.e.~in $Q$; since $\mdel$ is Lipschitz continuous and bounded, $\mdel(\phd)\bfeta \in \L2{\Hx1^d}$ with vanishing normal trace on $\partial \Omega$; and for a.e.~$t$ we have $\phd(t) \in \Hx2 \subset L^\infty(\Omega)$, so that $\psi'(\phd(t)) \in \Hx1$ with $\nabla \psi'(\phd(t)) = \psi''(\phd(t))\nabla\phd(t)$. Two integrations by parts in space, with no boundary contributions, thus yield the main terms to be handled, which appear in
\begin{align}
	\label{weak:2:delta}
	\int_Q \bfJ_\delta \cdot \bfeta
	=  -\beta \eps \int_Q  \Delta \phd \div\big(\mdel(\phd) \bfeta\big)
	-\frac \beta  \eps \int_Q  \mdel(\phd) \psi''(\phd) \nabla \phd  \cdot \bfeta .
\end{align}
As $\mdel$ is piecewise $C^1$ and globally Lipschitz, 
passing to the limit in $\int_Q  \Delta \phd \div\big(\mdel(\phd) \bfeta\big)$ can be carried out as in \cite{ElliottG96} (see also \cite{garcke2025anisotropic,ebenbeck2021cahn}), 
by using the chain rule for $f\circ u$ with a Sobolev function $u$ and a
piecewise $C^1$ function $f$ that is globally Lipschitz, see \cite[Theorem~7.8]{GilbargT83}. However, in contrast to \cite{ElliottG96}, passing to the limit in 
\[
\int_Q  \mdel(\phd) \psi''(\phd) \nabla \phd  \cdot \bfeta
\]
is more involved. This is due to the fact that the limit function $\ph$ is not bounded. Hence, we give details on how to pass to the limit in this term. We show that the integrand $\mdel(\phd) \psi''(\phd) \nabla \phd  \cdot \bfeta$ is integrable over $Q$. By estimate $(i)$ of Proposition~\ref{LEM:apriori_delta} and the growth condition $|\psi_2(r)| \le C(1+r^2)$ following from \eqref{pot:1}, we have $\psi(\phd) = \psi_1(\phd) + \psi_2(\phd) \in \L\infty{\Lx1}$ with a bound independent of $\delta$. Hence, by \eqref{pot:3},
\begin{align}
	\label{est:psi:Ltheta}
	\int_E |\psi''(\phd)| \le \eta \int_E \psi(\phd) + C_\eta |E| \le \eta \, C \, T + C_\eta |E|
\end{align}
for every measurable $E \subset Q$ and every $\eta > 0$, with $C$ independent of $\delta$. In particular, $\{\psi''(\phd)\}_\delta$ is bounded in $L^1(Q)$ and uniformly integrable on $Q$. Furthermore, with the help of Young's inequality it holds
\begin{align}
	\non
	|\psi''(\phd)| \, |\nabla \phd|
	& \le \psi_1''(\phd) \, |\nabla \phd| + c_0 |\nabla \phd|
	\\ & 
	\label{est:psi:L1}
	\le \tfrac12 \psi_1''(\phd) + \tfrac12 \psi_1''(\phd) |\nabla \phd|^2 + c_0 |\nabla \phd| \in L^1(Q),
\end{align}
where we used $0 \le \psi_1''(\phd) \le |\psi''(\phd)| + c_0$ and the entropy estimate in Proposition~\ref{LEM:apriori_delta}. 
Since $\psi_1'' = \psi'' - \psi_2''$ is continuous and non-negative, we have $\psi_1''(\phd) \to \psi_1''(\ph)$ a.e.~in $Q$, and by \eqref{est:psi:Ltheta} and \eqref{pot:1}, the family $\{\psi_1''(\phd)\}_\delta$ is uniformly integrable on $Q$. Vitali's convergence theorem thus gives $\psi_1''(\phd) \to \psi_1''(\ph)$ strongly in $L^1(Q)$ and, as $|\sqrt a - \sqrt b|^2 \le |a-b|$ for $a,b \ge 0$, as $\delta \to 0$,
\begin{align}
	\label{conv:psi}
	\sqrt{\psi_1''(\phd)} \to \sqrt{\psi_1''(\ph)} \quad \text{strongly in } L^2(Q),
\end{align}
which implies $\psi_1''(\ph) \in L^1(Q)$ and $\psi''(\ph) = \psi_1''(\ph) + \psi_2''(\ph) \in L^1(Q)$. Moreover, the entropy estimate in Proposition~\ref{LEM:apriori_delta} shows that $\sqrt{\psi_1''(\phd)}\, \nabla \phd$ is bounded in $L^2(Q)^d$, so that, along a further subsequence, as $\delta \to 0$,
\begin{align}
	\label{conv:entropy}
	\sqrt{\psi_1''(\phd)}\, \nabla \phd \to \sqrt{\psi_1''(\ph)}\, \nabla \ph \quad \text{weakly in } L^2(Q)^d .
\end{align}
Here the weak limit is identified as follows: by \eqref{conv:psi} and $\nabla \phd \to \nabla \ph$ in $L^2(Q)^d$, see \eqref{conv:strong:phd}, the product $\sqrt{\psi_1''(\phd)}\, \nabla \phd$ converges to $\sqrt{\psi_1''(\ph)}\, \nabla \ph$ strongly in $L^1(Q)^d$, and the $L^1$-limit coincides with the weak $L^2$-limit. In particular $\sqrt{\psi_1''(\ph)}\, \nabla \ph \in L^2(Q)^d$, which is the additional regularity stated in Theorem~\ref{THM:degenerate:plus}.
We now split $\psi'' = \psi_1'' + \psi_2''$ and write
\begin{align*}
& \mdel(\phd) \psi''(\phd) \nabla \phd \cdot \bfeta
	\\
	& \quad = \Big( \mdel(\phd) \sqrt{\psi_1''(\phd)}\, \bfeta \Big) \cdot \Big( \sqrt{\psi_1''(\phd)}\, \nabla \phd \Big)
	+ \mdel(\phd) \psi_2''(\phd) \nabla \phd \cdot \bfeta .
\end{align*}
Concerning the first term, the sequence $\{\mdel(\phd) \sqrt{\psi_1''(\phd)}\, \bfeta\}_\delta$ converges to $m(\ph) \sqrt{\psi_1''(\ph)}\, \bfeta$ strongly in $L^2(Q)^d$: indeed,
\begin{align*}
	& \Big\| \mdel(\phd) \sqrt{\psi_1''(\phd)}\, \bfeta - m(\ph) \sqrt{\psi_1''(\ph)}\, \bfeta \Big\|_{L^2(Q)^d}
	\\
	&  \quad \le m_+ \|\bfeta\|_{L^\infty(Q)^d} \Big\| \sqrt{\psi_1''(\phd)} - \sqrt{\psi_1''(\ph)} \Big\|_{L^2(Q)}
	\\
	& \qquad + \|\bfeta\|_{L^\infty(Q)^d} \Big\|\big(\mdel(\phd) - m(\ph)\big) \sqrt{\psi_1''(\ph)} \Big\|_{L^2(Q)}
\end{align*}
which converges to zero by \eqref{conv:psi} and by dominated convergence with the majorant $2 m_+ \sqrt{\psi_1''(\ph)} \in L^2(Q)$, as $(\mdel(\phd) - m(\ph)) \sqrt{\psi_1''(\ph)} \to 0$ a.e.~in $Q$. Together with the weak convergence \eqref{conv:entropy} this gives, as $\delta \to 0$,
\begin{align*}
	\int_Q \Big( \mdel(\phd) \sqrt{\psi_1''(\phd)}\, \bfeta \Big) \cdot \Big( \sqrt{\psi_1''(\phd)}\, \nabla \phd \Big)
	&\to \int_Q \Big( m(\ph) \sqrt{\psi_1''(\ph)}\, \bfeta \Big) \cdot \Big( \sqrt{\psi_1''(\ph)}\, \nabla \ph \Big)\\
	&= \int_Q m(\ph) \psi_1''(\ph) \nabla \ph \cdot \bfeta .
\end{align*}
Concerning the second term, $\psi_2''$ is continuous and bounded by $c_0$, so that $\mdel(\phd) \psi_2''(\phd) \bfeta \to m(\ph) \psi_2''(\ph) \bfeta$ a.e.~in $Q$ with the constant majorant $m_+ c_0 \|\bfeta\|_{L^\infty(Q)^d}$, hence strongly in $L^2(Q)^d$ by dominated convergence, while $\nabla \phd \to \nabla \ph$ strongly in $L^2(Q)^d$. Therefore, as $\delta \to 0$,
\begin{align*}
	\int_Q \mdel(\phd) \psi_2''(\phd) \nabla \phd \cdot \bfeta \to \int_Q m(\ph) \psi_2''(\ph) \nabla \ph \cdot \bfeta .
\end{align*}
Using the arguments developed in \cite{ElliottG96} together with this convergence then completes the proof of Theorem~\ref{THM:degenerate:plus}. 

\section{Sharp interface limits}\label{sec:Sharp}
In this section, we perform a formal asymptotic analysis of the system \eqref{eq:sys21} in the limit as $\eps \to 0$, considering potentials $\psi \in C^2(\erre)$ satisfying
\begin{align}
\label{ass:pot:sharp}
\psi(r)= \psi(-r), \quad \forall r \in \erre,
\quad
\psi(\pm1)=\psi'(\pm1)= \psi'(0)=0, \quad \psi''(\pm1)\neq 0,
\end{align}
such that $\psi >0$ in $(-1,1)$ and $\psi'(r) =0$ if and only if $r \in \{-1,0,1\}$. For the degenerate mobility, we work with the prototypical choice \eqref{eq:mob}. That is, in the \textit{plus-phase model} (degeneracy at $\varphi=-1$), it holds that
\begin{align}\label{degmob}
m(\varphi) = \tfrac{m_+}{2} (1+ \varphi),
\quad 
\text{if $\ph \in [-1,1]$,}
\end{align}
and so $m(-1)=0$ and $m(1)=m_+$.
More generally, any $m \in C^0(\erre)$ fulfilling the conditions $m(r) > 0$ for $r > -1$ and $m(r)=0$ for every $r\leq -1$ is admissible for our sharp interface analysis. Similarly, for the \textit{minus-phase model} (degeneracy at $\varphi=1$), we consider
\begin{align}\label{degmo:minus}
m(\varphi) = \tfrac{m_-}{2} (1- \varphi),
\quad 
\text{if $\ph \in [-1,1],$}
\end{align}
so that $m(-1)=m_-$ and $m(1)=0$. Again, any continuous mobility with the same properties as above can be employed.

Finally, for the source term, we adopt the explicit $\varepsilon$-dependent notation for the source $S=S_\varepsilon$ introduced in \eqref{def:source}. Since the asymptotic analysis proceeds analogously for both models, we restrict our detailed exposition to the \emph{plus-phase model}, merely stating the corresponding result for the \emph{minus-phase model}.

The method of formally matched asymptotic analysis integrates outer and inner expansions into the model equations, solving them stepwise, and defines a region for their matching. Further elaboration on the methodology can be found in, e.g.,~\cite{AbelsGG12,OOL,GLSS,BDGP}. The following assumptions and conventions are in order:
\begin{itemize}
	\item
		It is assumed that for small values of $\eps$, the domain $\Omega$ can be partitioned into two open subdomains, $\Omega_\pm=\Omega_\pm(\eps,t):=\{x \in \Omega : \ph_{\eps} (t,x) \gtrless 0\}$, separated by the zero level set denoted as $\Sigma_\eps=\Sigma_\eps(t) =\{x \in \Omega : \ph_{\eps}(t,x)=0\}$. We assume that the introduced spatial sets evolve over time, though we opt not to specify the temporal dependence for convenience. 
	\item
		It is assumed that there exists a family of solutions to \eqref{eq:sys21} that are sufficiently smooth and exhibit an asymptotic expansion in $\eps\in (0,1)$ within the bulk regions away from $\Sigma_\eps$ (referred to as the \textit{outer expansion}), and another expansion in the interfacial region adjacent to $\Sigma_\eps$ (referred to as the \textit{inner expansion}). Those will be denoted by $\{(\ph_\eps,\mu_\eps)\}_\eps$ and $\{(\Phi_\eps,\Lambda_\eps)\}_\eps$, respectively.
	\item
	We postulate that the level sets $\Sigma_\eps$, as $\eps \to 0$, converge to a limiting hypersurface $\Sigma_0$ with normal $\bnu$ that evolves with a normal velocity $\cal V$.
\end{itemize}
Although the setting is similar to that of \cite{OOL}, the one-sided degenerate mobility \eqref{degmob} introduces some technical considerations that require some detailed explanations.

\subsection{Outer expansion for the plus-phase model}
We suppose that the solution variables $\ph_\eps$ and $\mu_\eps$, far away from the interface, can be expressed as
\begin{align*}
	\ph_\eps = \sum_{i=0}^{\infty} \eps^i \ph_i,
	\quad
	\mu_\eps = \sum_{i=0}^{\infty} \eps^i \mu_i.
\end{align*}
Then, we consider equation \eqref{sys:2} to leading order $\eps^{-1}$ to infer that
\begin{align*}
	- \beta \psi'(\phz)=0.
\end{align*}
Here, we accounted for the conditions $\psi(\pm1)=\psi'(\pm1)=0$ in \eqref{ass:pot:sharp}. This shows that $\phz$ is a root of $\psi'$. 
Since $0$ is an unstable root, this leads to $\phz=\pm 1$. We then set
$
	\Omega^\pm:=\{x \in \Omega: \phz(x) = \pm 1\}.
$
To the next order $\eps^0$, equation \eqref{sys:2} yields
\begin{align*}
	\muz = \beta \psi''(\phz)\ph_1.
\end{align*}
This provides a relation between the leading order term $\muz$ and the first order term $\ph_1$. Since $\phz$ equals $1$ in $\Omega^+$ we infer from \eqref{sys:1} to order $\eps^0$ that 
\begin{align*}
	- \div (m_+\nabla \muz) = S_+ - K_+ \ph_1 \quad \text{ in } \Omega^+.
\end{align*}
Combining these two equations and setting $\rho_{+} = \frac{K_{+}}{\beta \psi''(1)}$ yields
\begin{align*}
- m_+ \Delta \mu_{0} = S_{+} - \rho_+ \muz  \quad \text{ in } \Omega^+.
\end{align*}
On the other hand, in $\Omega^-$ we infer that
\begin{align*}
0 = S_- - K_- \ph_1,
\end{align*}
and as $K_- \neq 0$, this yields the relation
\begin{align*}
  \muz = \frac{S_{-}}{K_{-}}\beta \psi''(-1)  \quad \text{ in } \Omega^-.
\end{align*}
It will turn out that due to the one-sided degeneracy of \eqref{degmob} the above relation for $\muz$ in $\Omega^-$ and the parameter $S_-$ will not play a role in the limiting free boundary problem. We remark that if $\pd \Omega^+ \cap \pd \Omega \neq \varnothing$, we also have the boundary condition
\[
m_+ \dn \mu_0 = 0 \quad \text{ on } \pd \Omega^+ \cap \pd \Omega.
\]

\subsection{Inner expansion and matching conditions for the plus-phase model}

To study the interface $\Sigma_\eps = \{\varphi_\eps = 0\}$ as $\eps \to 0$, we introduce the rescaled normal coordinate $z = d/\eps$, where $d$ is the signed distance to $\Sigma_0$ (positive in $\Omega^+$, negative in $\Omega^-$), so that the unit normal vector of the interface $\nabla d = \boldsymbol{\nu}$ points from $\Omega^-$ to $\Omega^+$. Parameterizing $\Sigma_0$ by arc-length as $g(t,s)$, we obtain the following reparametrization rule for any sufficiently smooth function $f$ defined in a tubular neighborhood
\[
f(x) = f(g(t,s) + \eps z \boldsymbol{\nu}(g(t,s))) =: F(t,s,z),
\]
leading to the asymptotic derivatives
\[
\partial_t f \approx -\frac{1}{\eps} \mathcal{V} \partial_z F, \quad 
\nabla_x f \approx \frac{1}{\eps} \partial_z F \boldsymbol{\nu} + \nabla_{\Sigma_0} F, \quad 
\Delta_x f \approx \frac{1}{\eps^2} \partial_{zz} F - \frac{1}{\eps} \kappa \partial_z F,
\]
where $\nabla_{\Sigma_0}$ stands for the surface gradient on $\Sigma_0$, $\kappa=-\div_{\Sigma_0} \bnu$ for the corresponding mean curvature, and the symbol $\approx$ indicates that we omitted higher order terms with respect to $\eps$.

The solution in the inner region is supposed to possess the following expansion
\begin{align*}
	\Phi_\eps
	 = \sum_{i=0}^{\infty} \eps^i \Phi_i,
	\quad
	\Lambda_\eps = \sum_{i=0}^{\infty} \eps^i \Lambda_i.
\end{align*}
The postulated convergence of the level sets $\Sigma_\eps$ to the hypersurface $\Sigma_0$ translates into the condition
\begin{align}
	\label{ass:phz}
	\Phi_0(t,s,z=0)=0.
\end{align}
Furthermore, we assume
\begin{align*}
	\Phi_\eps (t,s,z=\infty) =1,
	\quad
	\Phi_\eps (t,s,z=-\infty) =-1.
\end{align*}
With reference to \cite{GLSS}, we employ the following matching conditions:
\begin{alignat*}{2}
&\lim_{z\to \pm\infty} \Phi_0 = \varphi_0^\pm(t,x), \quad
&&\lim_{z\to +\infty} \Lambda_0 = \muz^+(t,x),\\
&\lim_{z\to \pm\infty} \partial_z \Phi_0 = 0, \quad
&&\lim_{z\to +\infty} \partial_z \Lambda_0 = 0,\\
&\lim_{z\to \pm\infty} \partial_z \Phi_1 = \nabla \varphi_0^\pm \cdot \boldsymbol{\nu}, \quad
&&\lim_{z\to +\infty} \partial_z \Lambda_1 = \nabla \muz^+ \cdot \boldsymbol{\nu},
\end{alignat*}
with $\varphi_0^\pm(t,x) = \lim_{\delta\to 0} \varphi_0(t,x \pm \delta \boldsymbol{\nu})$, and similarly for $\muz$. 

\begin{remark}
Due to the one-sided degeneracy of the mobility \eqref{degmob} we do not possess full control of the gradient $\nabla \mu_\eps$ in the whole domain $\Omega$, particularly in regions where $\ph_\eps$ is close to $-1$. This is reflected in the one-sided matching conditions for $\Lambda_0$ and $\Lambda_1$, as opposed to the two-sided matching conditions for $\Phi_0$ and $\Phi_1$.
\end{remark}

\subsection{Leading-order inner expansions for the plus-phase model}
To leading order $\eps^{-1}$ in \eqref{sys:2} we find
\begin{align*}
	\partial_{zz} \Phi_0 - \psi'(\Phi_0)=0{.}
\end{align*}
Using \eqref{ass:phz}, we observe that
this entails $\Phi_0$ is just a function of $z$, reducing to the ODE relation
\begin{align}
	\label{ODE:phiz}
	\Phi_0''(z)-\psi'(\Phi_0(z))=0,
	\quad
	\Phi_0(0)=0,
	\quad
	\lim_{z \to \pm \infty} \Phi_0(z) = \pm1.
\end{align}
The latter condition implies that $|\Phi_0(z)| < 1$ for all $z \in \erre$. Upon multiplying \eqref{ODE:phiz} by $\Phi_0'$, integrating and using the matching conditions for $\Phi_0$ and $\pd_z \Phi_0$, we obtain the equipartition of energy:
\begin{align}
	\label{eq:energy}
	\frac 12 |\Phi_0'(z)|^2 = \psi
	(\Phi_0(z)),
	\quad \forall z \in \erre.
\end{align}
Hence, we find the identity
\begin{align}
	\label{def:gamma}
	\intinf |\Phi_0'(z)|^2\,dz =
	\intinf 2 \psi(\Phi_0(z))\,dz
	= \int_{-1}^1 \sqrt{2 \psi(s)} \, ds =: \gamma.
\end{align}
{We remark that in the special case of $\psi$ being the quartic potential \eqref{quartic}, the solution $\Phi_0(z)$ to the ODE problem \eqref{ODE:phiz} is $\Phi_0(z) = \tanh(z/\sqrt{2})$, which coincides with the optimal profile minimizing the one-dimensional transition energy. The interfacial parameter $\gamma$ is then computed from this optimal profile, yielding
}
%
\begin{align*}
	\gamma=
	\frac 12 \intinf
		 {\rm sech}^4 \Big (\frac{z}{\sqrt{2}} \Big)
	\,dz=\frac {2\sqrt{2}}3{,}
\end{align*}
{and represents the minimal interfacial energy that determines the surface tension in the sharp interface limit.}
Next, to order $\eps^{-2}$, \eqref{sys:1} gives
\[
\frac{m_+}{2} \partial_z \Big( \Big(1 + \Phi_0 \Big) \partial_z \Lambda_0 \Big) = 0.
\]
Testing with $\Lambda_0$, integrating in $z$ and employing the matching conditions $\lim_{z \to - \infty} \Phi_0(z) = -1$ and $\lim_{z \to + \infty} \partial_z \Lambda_0 = 0$ we see that
\[
\frac{m_+}{2}\int_{-\infty}^{+\infty}  (1+\Phi_0(z)) |\partial_z \Lambda_0(t,s,z)|^2 \, dz = 0.
\]
Since $|\Phi_0 (z)|<1$, it follows that $\partial_z \Lambda_0(t,s,z) = 0$, whence $\Lambda_0$ is independent of $z$.

\subsection{Expansion to first order for the plus-phase model}
Moving to terms of zeroth order in \eqref{sys:2},  we derive that
\begin{align*}
	\Lambda_0 = \beta \psi''(\Phi_0)\Phi_1
	- \beta \partial_{zz}\Phi_1
	+ \beta \kappa\Phi_0'.
\end{align*}
Testing the above by $\Phi_0'$ and integrating from $-\infty$ to $+\infty$ produces
\begin{align*}
	\intinf
	\Lambda_0(t,s) \Phi_0'(z)\,dz
	= \beta \intinf \big(  (\psi'(\Phi_0(z)))' \Phi_1(z)
	-  \partial_{zz}\Phi_1(z) \Phi_0'(z)
	+  \kappa|\Phi_0'(z)|^2 \big) \,dz.
\end{align*}
Using the matching conditions $\lim_{z \to \pm \infty} \Phi_0(z) = \pm 1$ and $\lim_{z \to \pm \infty} \partial_z \Phi_0 = 0$, upon integrating by parts and using $\psi'(\phz)=\psi'(\pm1)=0$, we infer
\begin{align*}
	& \intinf
	  (\psi'(\Phi_0(z)))' \Phi_1(z) -
	\partial_{zz}\Phi_1(z) \Phi_0'(z)
	 \,dz
	 \\ & \quad
	  =
	 \big[\psi'(\Phi_0)\Phi_1 -
	\partial_{z}\Phi_1 \Phi_0'\big]^{+\infty}_{-\infty}
	 -\intinf
		\partial_z \Phi_1(z) (\psi'(\Phi_0(z))-\Phi_0''(z))	
	 \, dz=0,
\end{align*}
whence, recalling \eqref{eq:energy}, \eqref{def:gamma} and the matching condition $\lim_{z \to +\infty} \Lambda_0 = \muz^+$, this entails that 
\[
2 \muz^+ = 2\Lambda_0 = \gamma\beta \kappa.
\]
Next, to order $\eps^{-1}$ we find from \eqref{sys:1} that
\begin{align*}
	-{\cal V} \Phi_0' = \partial_{z} \Big ( \frac{m_+}{2} (1 + \Phi_0) \partial_z \Lambda_1 \Big ) -[ \Km G_2(\Phi_0)-\Kp G_2(-\Phi_0)]+L G_4(\Phi_0) .
\end{align*}
Integrating from $-\infty$ to ${+\infty}$ and using the matching conditions $\lim_{z \to +\infty} \partial_z \Lambda_1 = \nabla \muz^+ \cdot \bnu$ and $\lim_{z \to-\infty} \Phi_0(z) = -1$ leads to
\begin{align*}
		-2 {\cal V} & =
		m_+ \nabla \muz^+ \cdot \bnu +    \int_{- \infty}^{{+\infty}} \Kp G_2(-\Phi_0(z)) - \Km G_2(\Phi_0(z)) + L  G_4(\Phi_0(z))\, dz.
\end{align*}
Following \cite{OOL}, the integral on the \rhs\ can be simplified using \eqref{def:G23} and the property $-\Phi_0(z)=\Phi_0(-z)$. It holds that
\begin{align*}
	\non
	S_I & :=   \int_{- \infty}^{{+\infty}} \Kp G_2(-\Phi_0(z)) - \Km G_2(\Phi_0(z)) \,dz +
	 \int_{- \infty}^{{+\infty}} L G_4(\Phi_0(z))\, dz\\
	&=
	{(\Kp-  \Km)} \int_{- \infty}^{{+\infty}} G_2(\Phi_0(z)) \, dz  +
	\int_{- \infty}^{{+\infty}} L G_4(\Phi_0(z))\, dz,
\end{align*}
where, upon setting $\tilde c := \tfrac {1}{2 \sqrt{\psi''(-1)}}$, we realise with the help of the equipartition of energy \eqref{eq:energy} that
\begin{align*}
	& \int_{- \infty}^{{+\infty}} G_2(\Phi_0(z)) \, dz	
	=
	 -\tilde c  \int_{-\infty}^{{+\infty}}
	(\Phi_0(z)-1)\sqrt{2 \psi(\Phi_0(z))}\, dz
	\\ 
	& \quad
  =-\tilde c \int_{-\infty}^{{+\infty}}  (\Phi_0(z)-1)\Phi_0'(z)\, dz
	= -\tilde c \int_{-1}^1 (s-1)\,ds =  2 \tilde c.
\end{align*}
In addition, we notice that
\[
G_4(\Phi_0(z))=2 \psi(\Phi_0(z))=\sqrt{2 \psi(\Phi_0(z))}
\sqrt{2 \psi(\Phi_0(z))} =\sqrt{2 \psi(\Phi_0(z))} \Phi'_0(z),
\]
and compute
\[
\int_{- \infty}^{{+\infty}}  G_4(\Phi_0(z))\, dz= \int_{- \infty}^{{+\infty}} \sqrt{2 \psi(\Phi_0(z))} \Phi'_0(z) =
 \int_{-1}^1 \sqrt{2 \psi(s)} \, ds=\gamma.
 \]
Thus, recalling the definition of $\tilde c$, we find that
\begin{align}\label{defn:SI}
S_I=  \frac {{\Kp-  \Km}}{\sqrt{\psi''(-1)}}+\gamma L.
\end{align}
For the particular choice \eqref{quartic} of $\psi$, it holds that $\tilde c = \frac {\sqrt 2}4$ so that 
\begin{align}\label{defn:SI:eg}
S_I= \frac {\sqrt 2}2 ({\Kp-  \Km})+\frac{2\sqrt{2}}3  L.
\end{align}
Dropping the subscripts, the sharp interface limit takes the form
\begin{subequations}\label{SharpI}
\begin{alignat}{2}
	\label{sharpI:1:p:one}
	 - m_+ \Delta \mu &= \Sp - \rho_+ \mu  \quad &&\text{in $\Omega^+$,}
	\\
	\label{sharpI:2:one}
	\mu &=  \frac{\gamma \beta \kappa}2
	\qquad && \text{on $\Sigma$,}
	\\
	\label{sharpI:4:one}
	 - 2 {\cal V} & = m_+ \nabla \mu \cdot {\bnu}
	+ S_I
	\qquad && \text{on $\Sigma$,}\\
    \label{sharpI:bc}
    m_+ \dn \mu  &= 0 \quad && \text{on $\pd \Omega^+ \cap \pd \Omega$.}
\end{alignat}
\end{subequations}
\Accorpa\SharpIone {sharpI:1:p:one} {sharpI:5:one}

\subsection{Nondimensionalisation of the plus-phase model}
We now perform a nondimensionalisation argument in order to identify important dimensionless parameters.
Choosing units $\tilde x$, $\tilde t$ and $\tilde \mu$ for length, time and chemical potential we introduce
the nondimensional variables
\begin{equation*}
	\hat x =\frac x{\tilde x},\quad \hat t =\frac t{\tilde t}, \quad \hat \mu =\frac \mu{\tilde \mu},
\end{equation*}
and now consider the rescaled variant of system \eqref{SharpI} on the rescaled domain $\hat \Omega^+$ and the corresponding interface $\hat{\Sigma}$. Denoting by $\hat{\nabla}$ and $\hat{\Delta}$ the gradient and Laplacian with respect to $\hat{x}$, for the new nondimensional variables, we obtain from \eqref{SharpI} the following system
\begin{subequations}
\begin{alignat*}{3}
			- \hat\Delta \hat\mu &=
		\frac{(\tilde x)^2}{m_+ \tilde \mu }\Sp -  \frac{(\tilde x)^2\rho_+}{m_+} \hat \mu  \quad && \text{in $\hat \Omega^+$},
		\\
	 \hat \mu &=  \frac{\gamma \beta {\hat \kappa}}{2\tilde \mu \tilde x} 
	\qquad && \text{on $\hat \Sigma$,}
	\\
	 - 2 { \cal \hat V} & = \tfrac{\tilde t\tilde \mu } {(\tilde x)^2} m_+ \hat \nabla \hat\mu \cdot {\hat \bnu}
	+\tfrac{\tilde t}{\tilde x} S_I
	\qquad && \text{on ${\hat{\Sigma}}$,}
	\end{alignat*}
\end{subequations}
with $\hat {\calV}$ being the normal velocity of the interface $\hat \Sigma$ in the direction of the normal $\hat \bnu$. To obtain simple nondimensional equations, and on assuming that $\rho_+>0$ and $S_+ <0$,
we set
\begin{equation*}
	\tilde x =\sqrt{\frac{m_+}{\rho_+}}, \quad \tilde \mu = \frac{(\tilde x)^2 (-S_+) }{m_+}= \frac{-S_+}{\rho_+}, \quad \tilde t =
	 \frac{(\tilde x)^2}{\tilde \mu m_+} =\frac{-1}{S_+} .
\end{equation*}
We now define the nondimensional parameter
\begin{equation*}
	\beta^* =\frac{\beta}{\tilde x} \frac 1{\tilde \mu}=
	 \frac{\beta \rho_+^{\tfrac 32}}{-m_+^{\frac12} S_+}= \frac {c_{l}}{\tilde x},
\end{equation*}
where we call
\begin{equation}
	\label{cap:lenght}
	c_{l}=  \beta \frac{\rho_+}{-S_+},
\end{equation}
the modified \textit{capillary length}, and introduce
\[
\alpha^* = \frac{\gamma \beta^*}{2}.
\]
In addition, we introduce the nondimensional interface reaction term $S_I^*$ as follows
\begin{equation*}
	S_I^* =\frac{\tilde t}{\tilde x}S_I= \frac{-1}{S_+} \sqrt{\frac{\rho_+}{m_+}} S_I.
\end{equation*}
We then obtain, dropping the hat and starred notation for convenience, the nondimensional sharp interface model
\begin{subequations}\label{SharpI:ND}
\begin{alignat}{3}
		\label{sharpI:1:ND}
			- \Delta  \mu &=  -1 - \mu  \quad && \text{in $  \Omega^+$,}
	\\
	 \label{sharpI:2:ND} \mu & = \alpha \kappa 
	\qquad & &\text{on $ \Sigma$,}
	\\
	\label{sharpI:4:ND}
	 - 2 { \cal  V} &= \nabla \mu \cdot { \bnu} 
	+ S_I
	\qquad && \text{on $\Sigma$,}\\
    \label{sharpI:bc:ND}
    \dn \mu & = 0 \quad && \text{on $\pd \Omega^+ \cap \pd \Omega$.}
\end{alignat}
\end{subequations}

\subsection{Sharp interface analysis for the minus-phase model}
The sharp-interface one-phase model describing the dynamics exclusively in $\Omega^+$ derived above can be formally obtained by setting $m_- = S_- = \rho_- = 0$ in the original system \eqref{SharpI:original}. Assuming that the free boundary $\Sigma$ does not intersect the external boundary $\partial\Omega$, this yields a closed system posed in $\Omega^+$ and on $\Sigma$, denoted by \eqref{SharpI}. Its nondimensionalised counterpart is given by \eqref{SharpI:ND}. 

By entirely analogous arguments, now involving a degenerate mobility at the opposite extremum \eqref{degmo:minus}, a second sharp-interface model capturing the dynamics solely in $\Omega^-$ can be derived.
For brevity, we state the resulting system here without further details. Specifically,  setting $m_+ = S_+ = \rho_+ = 0$ in the original system formally yields a corresponding set of equations defined in $\Omega^-$ and on the interface $\Sigma$.  The system reads as:
\begin{subequations}\label{SharpI:single:m}
\begin{alignat}{2}
        \label{oneph:-:sharpI:1:m} - m_- \Delta \mu & =\Sm - \rho_- \mu  \quad && \text{in $\Omega^-$,}
	\\
	\label{oneph:-:sharpI:2}
	\mu &=  \alpha\kappa
	\qquad && \text{on $\Sigma$,}
	\\
	\label{oneph:-:sharpI:4}
	 - 2 {\cal V} & = -m_-\nabla \mu_- \cdot {\bnu}
	+ S_I
	\qquad && \text{on $\Sigma$,}
	\\
	\label{oneph:-:sharpI:5}
	m_- \dn \mu & =0
	\qquad && \text{on $\partial \Omega$,}
\end{alignat}
\end{subequations}
\Accorpa\SharpIp {oneph:-:sharpI:1:m} {oneph:-:sharpI:5}
along with its nondimensional formulation
\begin{subequations}\label{SharpI:single:ad:m}
\begin{alignat}{3}
	\label{oneph:-:sharpIIn:1:p}
	-  \Delta \mu &=
		 1 -   \mu \quad && \text{in $ \Omega^-$},
		\\
	\label{oneph:-:sharpIIn:2}
	 \mu & = \alpha\kappa
	\qquad & &\text{on $ \Sigma$,}
	\\
	\label{oneph:-:sharpIIn:4}
	 - 2 { \cal  V} &= 
-	\nabla \mu_- \cdot { \bnu}
	+ S_I
	\qquad && \text{on $\Sigma$,}
	\\
	\label{oneph:-:sharpIIn:5}
	 \dn \mu &=0
	\qquad && \text{on $\partial \Omega$.}
\end{alignat}
\end{subequations}

\section{Radial solutions}\label{sec:radial}

We now consider the free boundary problems \eqref{SharpI:ND} and \eqref{SharpI:single:ad:m} in the $d$-dimensional ball of radius $R$ centered at the origin:
$
\Omega ={\mathbb B_R^d(0)},
$
for some external radius $R>0$ and seek  a radially symmetric solution under the geometry
\begin{align*}
	\Omega^+(t)& = \mathbb{B}^d_{q(t)}(0) = \{ x \in \erre^d \, : \, | x| < q(t) \},
	\\
	\Omega^-(t) & = \Omega \setminus \ov{\mathbb{B}^d_{q(t)}(0)} = \{ x \in \erre^d \, : \, q(t) < |x| < R \}.
\end{align*}
Here, the function $q(t)$ encodes the location of the time-dependent interface $\Sigma = \partial \mathbb{B}^d_{q(t)}(0)$. At $p \in \Sigma$, the unit normal, normal velocity and mean curvature are given as
\begin{align*}
	\bnu (p) =- \frac p{|p|} =- \frac p{q(t)},
	\qquad
	{\cal V} = -\dot{q}(t)	,
	\qquad
	\kappa = \frac {d-1}{q(t)},
\end{align*}
with the dot denoting the time derivative. In this radially symmetric setting, we make the ansatz
\begin{align}\label{radialansatz}
\mu(t, x) = \hat{\mu}(t, |x|) = \hat{\mu}(t, r),
\quad \forall r >0.
\end{align}
In what follows, we drop the hat-notation for convenience, and use the prime symbol to denote differentiation with respect to {the radial variable} $r$. 

Before proceeding to the analysis of the two single-phase models, we introduce some useful preliminaries.
To address the analysis in both two and three-dimensional settings simultaneously, for every $\ell \in \mathbb{N} \cup \{0\}$, we introduce the notation
\begin{align}\label{order:ell:Bess}
	\Il(r):=
	\begin{cases}
	I_\ell(r) \quad  & \text{if $d=2$},
	\\
	i_\ell(r) \quad  & \text{if $d=3$},
	\end{cases}
	\quad
	\Kl(r):=
	\begin{cases}
	K_\ell(r) \quad  & \text{if $d=2$},
	\\
	k_\ell(r) \quad  & \text{if $d=3$},
	\end{cases}
	\quad \forall r \in \erre^+,
\end{align}
where $\{I_\ell,K_\ell\}$ and $\{i_\ell,k_\ell\}$ represent the modified Bessel functions of first and second kind to the order $\ell$ and the spherical modified Bessel functions of first and second kind to the order $\ell$, respectively. It is well known that $\{I_\ell, K_\ell\}$ forms a set of two linearly independent solutions to the modified Bessel differential equation
\begin{align*}
	r^2 y''(r) + r y'(r) - r^2 y(r) = \ell^2 y(r),
	\quad \forall r \in \erre^+,
\end{align*}
while $\{i_\ell, k_\ell\}$ forms a set of two linearly independent solutions to the modified spherical Bessel differential equation
\begin{align*}
	r^2 y''(r) + 2r y'(r) - r^2 y(r) = \ell (\ell+1) y(r),
	\quad \forall r \in \erre^+.
\end{align*}
For a precise definition and further properties of these Bessel functions, we refer to \cite{WatsonBessel}.  It is worth pointing out that $\III_\ell$, $\KKK_\ell$ and $\III'_\ell$ are always non{-}negative, whereas $\KKK_\ell'$ is non{-}positive. We list some properties and identities involving the zeroth and first order modified Bessel functions:
\begin{align}
	I_0(0)& =1,
	\quad
	 I_0'(r)=I_1(r), \quad \lim_{r \to 0^+} K_{0}(r)  = + \infty, \quad K_{0}'(r)  = - K_1(r),  \label{Bessel:2d} \\
	\label{Bessel:3d:1}
	i_0(r) &= \frac {\sinh(r)}{r},	\quad  i_1(r) =i_0'(r) = \frac {r \cosh(r) - \sinh(r)}{r^2},\\
	\label{Bessel:3d:2}
	k_0(r) & = \frac{e^{-r}}{r}, \quad k_1(r) =-k_0'(r) = \frac {e^{-r}(r+1)}{r^2},
\end{align}
to be considered for every $r \in \erre^+$, whenever it occurs.

\subsection{Plus-phase model}
The \emph{plus-phase model} \eqref{SharpI:ND} in the radial setting reads as
\begin{subequations}\label{SharpI:rad}
\begin{alignat}{3}
	\label{SharpI:rad:1}
	  -\mu'' -\frac {d-1}r \mu' &=
	-1 - \mu
	 \quad &&\text{in $\{0<r<q(t)\}$,}
	\\[1ex]
    \label{SharpI:rad:2}
	 \mu & = \alpha \frac{d-1}{q(t)}
	\quad && \text{on $\{r=q(t)\}$,}
	\\[1ex]
    \label{SharpI:rad:3}
	 2\dot{q} & = -\mu'	+  S_I
	\quad && \text{on $\{r=q(t)\}$.}
\end{alignat}
\end{subequations}
We further complement the above system with the boundary condition at the origin
\begin{align}
	\label{SharpI:rad:4}
	\mu(t,r=0)<\infty.
\end{align}
We consider the ansatz 
\begin{alignat}{3}
	\label{muplus:rad}
	\mu(t,r) & = b_{d,+}(q(t)) \Iz(r) + e_+, \quad && \text{ for }r \in (0,q(t)),
\end{alignat}
for suitable constant $e_{+}$ and function $b_{d,+}(q(t))$ of the interfacial location $q(t)$ to be determined. Note that the expression for the coefficient function $b_{d,+}$ depends on the spatial dimension $d$. From the condition \eqref{SharpI:rad:4} we note that the ansatz for $\mu$ in \eqref{muplus:rad} does not contain a term involving $\Kz(r)$ due to its singular {behaviour} at $r = 0$. 

Suppressing the explicit dependence on $t$ and upon substituting the ansatz \eqref{muplus:rad} into \eqref{SharpI:rad:1}, we infer that
\begin{align*}
1 + \mu(r) & = \mu''(r) + \frac{d-1}{r} \mu'(r) = b_{d,+}(q) \Iz(r) = \mu(r) - e_+,
\end{align*}
which leads to the identification of the constant $e_{+} = -1$. Next, from the interface condition \eqref{SharpI:rad:2}, we infer that
\begin{align*}
b_{d,+} (q) \Iz(q) - 1 &  = \alpha \frac {d-1}q,
\end{align*}
which yields the formula for $b_{d,+}(q)$ as
\begin{align*}
	b_{d,+} (q) = \Big ( \alpha \frac {d-1}q + 1 \Big ) \frac{1}{\Iz(q)}, \quad \text{ for } q > 0,
\end{align*}
and thus the analytical formula for the solution $\mu$ reads as
\begin{align}\label{radial:mu:sol}
\mu(t,r) = \Big ( \alpha \frac{d-1}{q(t)} + 1 \Big ) \frac{\Iz(r)}{\Iz(q(t))} - 1, \quad \text{ for } r \in (0,q(t)).
\end{align}
From this we can express the ordinary differential equation \eqref{SharpI:rad:3} satisfied by the radius function $q(t)$ as
\begin{align}
2 \dot{q} = -\Big(\alpha \frac {d-1}q  +1 \Big) \frac{\III_1 (q)}{\Iz(q)} + S_I	=: {\cal H}^+ (q).\label{radial:ode}
\end{align}
By utilising the following properties of the modified Bessel functions of the first kind $I_0, I_1$, and the modified spherical Bessel functions of the first kind $i_0, i_1$:
\begin{align*}
\lim_{r\to0^+} \frac{I_1(r)}{I_0(r)} = 0, \quad \lim_{r \to 0^+} \frac{I_1(r)}{r I_0(r)} = \frac{1}{2}, \quad \lim_{r \to \infty} \frac{I_1(r)}{I_0(r)} = 1, \quad \lim_{r \to \infty} \frac{I_1(r)}{rI_0(r)} = 0, \\
\lim_{r\to0^+} \frac{i_1(r)}{i_0(r)} = 0, \quad \lim_{r \to 0^+} \frac{i_1(r)}{r i_0(r)} = \frac{1}{3}, \quad \lim_{r \to \infty} \frac{i_1(r)}{i_0(r)} = 1, \quad \lim_{r \to \infty} \frac{i_1(r)}{ri_0(r)} = 0, 
\end{align*}
we see that the right-hand side $\mathcal{H}^+$ defined in \eqref{radial:ode} satisfies
\[
\lim_{q \to \infty} \mathcal{H}^+(q) = S_I - 1, \quad \lim_{q \to 0^+} \mathcal{H}^+(q) = \begin{cases}
S_I - \frac{\alpha}{2} & \text{ if } d = 2, \\
S_I - \frac{2 \alpha}{3} & \text{ if } d = 3.
\end{cases}
\]
Namely, due to the continuity of $\mathcal{H}^+$, the existence of a root of $\mathcal{H}^+$ inferred via the intermediate value theorem can be guaranteed under the following settings:
\begin{equation}\label{radial:para:root}
\begin{aligned}
\text{ Case (a): } &\begin{cases}
    \alpha \in (0,2), \quad S_I \in (\frac{\alpha}{2}, 1) & \text{ if } d = 2, \\
    \alpha \in (0, \frac{3}{2}), \quad S_I \in (\frac{2\alpha}{3}, 1) & \text{ if } d = 3,
    \end{cases} \\
\text{ Case (b): } & \begin{cases}
    \alpha > 2, \quad S_I \in (1, \frac{\alpha}{2}) & \text{ if } d = 2, \\
    \alpha > \frac{3}{2}, \quad S_I \in (1, \frac{2\alpha}{3}) & \text{ if } d = 3.
    \end{cases} 
\end{aligned}
\end{equation}
Furthermore, invoking the identities $\III_1(r) = \III_0'(r)$ and $\III_1'(r) = \III_0(r) - \frac{d-1}{r} \III_1(r)$, we infer the identity for $d\in \{2,3\}$
\[
(\mathcal{H}^+)'(q) = - \Big ( 1 - \frac{d-1}{q} \frac{\III_1(q)}{\III_0(q)} - \Big (\frac{\III_1(q)}{\III_0(q)} \Big )^2 \Big ) - \frac{\alpha (d-1)}{q} \Big ( 1 - \frac{d}{q}  \frac{\III_1(q)}{\III_0(q)} - \Big (\frac{\III_1(q)}{\III_0(q)} \Big )^2 \Big )
\]
and the behaviour of $({\cal H}^+)'$ near zero is
\begin{align*}
\lim_{q \to 0^+} (\mathcal{H}^+)'(q) = - \frac{1}{d}.
\end{align*}
In particular, for the parameter setting of Case (a), as $\mathcal{H}^+(0) > 0$ and $(\mathcal{H}^+)'(0) < 0$, we can expect the existence of a positive root for $\mathcal{H}^+$ closer to $0$ (corresponding to a very small stationary radius that is stable with respect to radial perturbations). On the other hand, for the parameter setting of Case (b), as $\mathcal{H}^+(0) < 0$ and $(\mathcal{H}^+)'(0) < 0$, we can expect the existence of a positive root for $\mathcal{H}^+$ further away from $0$ (corresponding to a larger stationary radius that may be unstable with respect to radial perturbations). For instance, we report an equilibrium radius $\qstar = 0.2064$ for parameters $\alpha = 1$ and $S_I = 0.6$ (Case (a)) and an equilibrium radius $\qstar = 6.2837$ for parameters $\alpha = 4$ and $S_I = 1.5$ (Case (b)). 

For the linear stability analysis we focus only on the parameter regimes corresponding to Case (a), as the associated roots of $\mathcal{H}^+$ are stable with respect to radial perturbations.

\subsection{Minus-phase model}
The \emph{minus-phase model} \eqref{SharpI:single:ad:m} in the radial setting
reads as
\begin{subequations}\label{SharpI:rad:m}
\begin{alignat}{3}
	\label{SharpI:rad:1:m}
	  -\mu'' -\frac {d-1}r \mu' &=
    1 - \mu
	 \quad &&\text{in $\{q(t)<r<R\}$,}
	\\[1ex]
    \label{SharpI:rad:2:m}
	 \mu & = \alpha \frac{d-1}{q(t)}
	\quad && \text{on $\{r=q(t)\}$,}
	\\[1ex]
    \label{SharpI:rad:3:m}
	 2\dot{q} & = \mu'	+  S_I
	\quad && \text{on $\{r=q(t)\}$,}\\[1ex]
    \label{SharpI:rad:4:m}
	 \mu'(t,r) & = 0
	\quad && \text{on $\{r=R\}$.}
\end{alignat}
\end{subequations}
We consider the ansatz 
\begin{alignat}{3}
	\label{muminus:rad}
	\mu(t,r) & = b_{d,-}(q(t)) \Iz(r) + c_{d,-}(q(t)) \Kz(r) + e_-, \quad && \text{ for }r \in (q(t),R),
\end{alignat}
for suitable constant $e_{-}$ and functions $b_{d,-}(q(t))$ and $c_{d,-}(q(t))$ of the interfacial location $q(t)$ to be determined. Note that the expressions for the coefficient function $b_{d,-}$ and $c_{d,-}$ depend on the spatial dimension $d$. 

Suppressing the explicit dependence on $t$ and upon substituting the ansatz \eqref{muminus:rad} into \eqref{SharpI:rad:1:m}, we infer that
\begin{align*}
-1 + \mu(r) & = \mu''(r) + \frac{d-1}{r} \mu'(r)  = \mu(r) - e_-,
\end{align*}
which leads to the identification of the constant $e_{-} = 1$. Next, the external boundary condition \eqref{SharpI:rad:4:m} yields
\begin{align}\label{minus:rad:bdy}
0 = \mu'(R) = b_{d,-}(q) \III_0'(R) + c_{d,-}(q) \KKK_0'(R) = b_{d,-}(q)\III_1(R) - c_{d,-}(q) \KKK_1(R).
\end{align}
From the interface condition \eqref{SharpI:rad:2:m}, we infer that
\begin{align*}
b_{d,-} (q) \Iz(q) + c_{d,-}(q) \Kz(q) + 1 &  = \alpha \frac {d-1}q,
\end{align*}
which together with \eqref{minus:rad:bdy} allows us to solve for $b_{d,-}(q)$ and $c_{d,-}(q)$ with the following linear system:
\begin{align*}
    \begin{pmatrix}
        \III_1(R) & - \KKK_1(R) \\
        \III_0(q) & \KKK_0(q)
    \end{pmatrix} \begin{pmatrix}
        b_{d,-}(q) \\ c_{d,-}(q)
    \end{pmatrix} = \begin{pmatrix}
        0 \\ \alpha \frac{d-1}{q} - 1
    \end{pmatrix},
\end{align*}
yielding the formula
\begin{align*}
c_{d,-}(q) & = \Big ( \alpha \frac{d-1}{q} - 1 \Big ) \frac{1}{\Iz(q) \frac{\KKK_1(R)}{\III_1(R)} + \Kz(q)} , \\
b_{d,-}(q) & =\frac{\KKK_1(R)}{\III_1(R)} c_{d,-}(q) = \Big ( \alpha \frac{d-1}{q} - 1 \Big ) \frac{\frac{\KKK_1(R)}{\III_1(R)}}{\Iz(q) \frac{\KKK_1(R)}{\III_1(R)} + \Kz(q)},
\end{align*}
and thus the analytical formula for the solution $\mu$ reads as
\begin{align}\label{radial:mu:sol:minus}
\mu(t,r) = \Big ( \alpha \frac{d-1}{q(t)} - 1 \Big ) \frac{\Iz(r)\frac{\KKK_1(R)}{\III_1(R)}+ \Kz(r)}{\Iz(q(t)) \frac{\KKK_1(R)}{\III_1(R)} + \Kz(q(t))} + 1, \quad \text{ for } r \in {(q(t),R)}.
\end{align}
From this we can express the ordinary differential equation \eqref{SharpI:rad:3:m} satisfied by the radius function $q(t)$ as
\begin{align}
2 \dot{q} = \Big(\alpha \frac {d-1}q  - 1 \Big) \frac{\III_1 (q) \frac{\KKK_1(R)}{\III_1(R)} - \KKK_1(q)}{\Iz(q) \frac{\KKK_1(R)}{\III_1(R)} + \Kz(q)} + S_I	=: {\cal H}^- (q).\label{radial:ode:minus}
\end{align}
We see that the right-hand side ${\cal H}^-$ defined in \eqref{radial:ode:minus} satisfies for $\alpha>0$
\begin{align*}
\lim_{q \to R} {\cal H}^-(q) = S_I, \quad \lim_{q \to 0^+} {\cal H}^{-}(q) = - \infty,
\end{align*}
the latter due to the following properties{:}
\begin{align*}
\lim_{r \to 0^+} \frac{\KKK_1(r)}{\KKK_0(r)} = +\infty, \quad \lim_{r \to 0^+} \frac{\III_1(r)}{\KKK_1(r)} = 0, \quad  \lim_{r \to 0^+} \frac{\III_0(r)}{\KKK_0(r)} = 0, 
\end{align*}
so that 
\begin{align*}
\lim_{q \to 0^+} {\cal H}^-(q) & = \lim_{q \to 0^+} \Big ( \alpha \frac{d-1}{q} - 1 \Big ) \frac{\KKK_1(q)}{\KKK_0(q)}  \frac{\frac{\III_1(q)}{\KKK_1(q)} \frac{\KKK_1(R)}{\III_1(R)} - 1}{\frac{\III_0(q)}{\KKK_0(q)} \frac{\KKK_1(R)}{\III_1(R)} + 1} 
\\ 
& = \lim_{q \to 0^+} \Big ( \alpha \frac{d-1}{q} - 1 \Big ) \frac{\KKK_1(q)}{\KKK_0(q)}  \frac{\frac{\III_1(q)}{\KKK_1(q)} - \frac{\III_1(R)}{\KKK_1(R)} }{\frac{\III_0(q)}{\KKK_0(q)}  + \frac{\III_1(R)}{\KKK_1(R)}} 
= -\infty.
\end{align*}
Moreover, a short calculation with the identities $\III_0' = \III_1$ and $\KKK_0' = - \KKK_1$ shows
\begin{align*}
({\cal H}^-)'(q) &  = \frac{\alpha (d-1) q - q^2}{q^2 (\III_0(q) \frac{\KKK_1(R)}{\III_1(R)}+ \KKK_0(q))^2} \\
&  \qquad \times  \left [ \Big (\III_1'(q) \frac{\KKK_1(R)}{\III_1(R)} - \KKK_1'(q) \Big ) \Big (\III_0(q) \frac{\KKK_1(R)}{\III_1(R)} + \KKK_0(q) \Big ) - \Big ( \III_1(q) \frac{\KKK_1(R)}{\III_1(R)} - \KKK_1(q) \Big )^2  \right ] \\
& \quad - \frac{\alpha (d-1) (\III_1(q) \frac{\KKK_1(R)}{\III_1(R)} - \KKK_1(q))}{q^2 (\III_0(q) \frac{\KKK_1(R)}{\III_1(R)} + \KKK_0(q))} ,
\end{align*}
where we observe
\begin{align*}
    ({\cal H}^-)'(R) = \frac{\alpha(d-1)R - R^2}{R^2(\III_0(R) \frac{\KKK_1(R)}{\III_1(R)} + \KKK_0(R))} \Big (\III_1'(R) \frac{\KKK_1(R)}{\III_1(R)} - \KKK_1'(R) \Big ){.}
\end{align*}
Using the well-known Wronskian identity for the modified Bessel functions $(d=2)$ and the modified spherical Bessel functions $(d=3)$:
\[
\III_n'(r) \KKK_n(r) - \III_n(r) \KKK_n'(r) = \frac{1}{r^{d-1}},
\]
we infer for the range $0 < \alpha < \frac{R}{d-1}$ the derivative of $\mathcal{H}^{-}$ at $R$ satisfies $(\mathcal{H}^-)'(R) < 0$. In this context, $n \in \mathbb{N} \cup \{ 0 \}$ denotes the angular mode index, where $n=0$ corresponds to the radially symmetric mode. Thus, as long as $S_I < 0$ and $\max_{q \in (0,R]} {\cal H}^{-}(q) > 0$, we infer the existence of at least two roots; one root close to $0$ that is unstable with respect to radial perturbations and another root close to $R$ that is stable with respect to radial perturbations.

\section{Linear stability of radial solutions}
\label{sec:linstab}

We now investigate the linear stability of the radially symmetric equilibrium configuration characterised by the stationary radius $\qstar$, corresponding to a stable root of the functional relation $\mathcal{H}^+$ established in \eqref{radial:ode} and for ${\cal H}^{-}$ in \eqref{radial:ode:minus}. In what follows, we detail the analysis for the \emph{plus-phase model}. The corresponding result for the \emph{minus-phase} case is stated directly thereafter with fewer details, as the underlying arguments proceed analogously.

\subsection{Plus-phase model}
We introduce a small perturbation to the interface position, defined by $w = \qstar + \badeps \rad $, with $\badeps \in (0,1)$ denoting a perturbation parameter. The perturbation profile $\rad=\rad(t,\th,\phi)$ depends on time and the angular coordinates, where $\th$ and $\phi$ represent the polar and azimuthal angles, respectively. In two spatial dimensions, the dependence reduces to $\rad=\rad(t,\th)$. Let $X(\th,\phi)$ denote the canonical parametrization of the unit sphere $\mathcal{S}^{d-1}$. The perturbed subdomain is then characterised as
\begin{align*}
	\Omega_{\rad,\badeps}^+ & := \{
	x \in \Omega\,\, : \,\,  x = |x| X(\th,\phi), \,
	|x|  < \qstar + \badeps \rad(t,\th,\phi)
	\}.
\end{align*}
The chemical potential is expanded as $\mu(r,t,\th,\phi) = \mu^*(t,r) + \badeps u(r,t,\th,\phi)$, where $\mu^*$ denotes the stationary solution to \eqref{SharpI:ND} associated with $\qstar$. Substituting this ansatz into the governing free boundary problem \eqref{SharpI:ND} yields the following system on the perturbed domain:
\begin{subequations}
\begin{alignat}{2}
	\label{lin:1}
		-  \Delta (\mustpl + \badeps u ) &= -1 -   (\mustpl + \badeps u ) \qquad &&\text{in $\Omega_{\rad,\badeps}^+$,}
	\\
	\label{lin:3}
	 \mustpl + \badeps u &=  \alpha \kappa
	\qquad && \text{on $\{r=w\}$,}
	\\
	\label{lin:5}
	 - 2 {\cal V} & = \nabla (\mu^* + \badeps u) \cdot {\bnu}
	+ S_I
	\qquad && \text{on $\{r=w\}$,}
\end{alignat}
\end{subequations}
where $\kappa$, $\mathcal{V}$, and $\bm{\nu}$ have the same meaning as above, associated with the deformed interface $\{ r = w \}$.
We perform a first-order Taylor expansion of the boundary conditions about the unperturbed interface $\{r = \qstar\}$ via
\begin{align*}
	\mu^*(w) =
	\mu^*(\qstar)
	+ (\mu^*)'|_{r=\qstar} (w-\qstar) + \text{ h.o.t.}=
	\mu^*(\qstar)
	+ \badeps (\mu^*)'|_{r=\qstar} Y +\text{ h.o.t.}{.}
\end{align*}
  Upon neglecting  higher order terms, we derive the following linearised system governing the evolution of $\rad$ and $u$:
\begin{subequations} \label{pertrad:general}
\begin{alignat}{2}
	\label{pertrad:general:1}
	&
	-  \Delta  u =  -   u  \qquad &&\text{in $\{0 <r <\qstar\}$,}
	\\
	\label{pertrad:general:3}
	& (\mu^*)' Y + u= -  \alpha 	\frac {d-1}{(\qstar)^2} \Big( \rad + \frac 1{d-1} \Delta_{{\cal S}^{d-1}} \rad \Big)
	\qquad && \text{on $\{r=\qstar\}$,}
	\\
	&  2 \dot{Y} =
	- (\mu^*)'' \rad
	\label{pertrad:general:5}
	-   u'
	\qquad && \text{on $\{r=\qstar\}$,}
\end{alignat}
\end{subequations}
where $\Delta_{{\cal S}^{d-1}}$ denotes the Laplace--Beltrami operator on the $(d-1)$-dimensional unit sphere ${\cal S}^{d-1}$. The linearised curvature term ${\tfrac {d-1}{(\qstar)^2}} \big( \rad + {\tfrac 1{d-1}} \Delta_{{\cal S}^{d-1}} \rad \big)$ follows from standard differential geometric expansions around a spherical interface, see, e.g., \cite[Lemma 3.1]{ESCHER1998267} (noting the differing normalization convention for the mean curvature) or \cite[Lemma 39]{bgnreview}. To ensure regularity at the center of the domain, we impose the boundedness condition at the origin
\begin{align}
	\label{pertrad:general:initial}
	u (r=0)<\infty.
\end{align}
The linear evolution law \eqref{pertrad:general:5} takes the abstract form $2 \dot{Y} = \mathcal{G}(\rad)$, where the operator $\mathcal{G}$ is determined by solving the elliptic boundary value problem \eqref{pertrad:general} subject to \eqref{pertrad:general:initial}. To diagonalise this system, we employ a separation of variables ansatz for the interfacial perturbation
\begin{align}\label{rad:ansatz}
	\rad = \delta (t) Z(\th,\phi),
\end{align}
where for $\ell \in \enne \cup \{ 0 \}$ and $ k \in \{-\ell , ..., \ell \}$, we set
\begin{alignat}{2}	
Z(\th,\phi)& =
	\cos(\ell \th),
	\; \text{ or }\;
	Z(\th,\phi) =
	\sin(\ell \th),
	 \quad &&\text{if $d=2$,} \label{perturb2d}
	\\
	Z(\th,\phi)& =
	Y_{\ell,k}(\th,\phi), \quad &&\text{if $d=3$.}	
	\label{perturb3d}
\end{alignat}
In three dimensions, $Y_{\ell,k}(\th,\phi)$ represents the standard spherical harmonics, with $\ell$ and $k$ corresponding to the degree and order, respectively, while $\th$ and $\phi$ retain their geometric interpretation as polar and azimuthal angles. Consistent with this spectral decomposition, we expand the bulk perturbation as
\begin{align}\label{perturb:u:ansatz}
	 u (r,t, \th,\phi)= \delta (t) Z(\th,\phi)U (r),
\end{align}
where the radial profile $U(r)$ satisfies a decoupled ordinary differential system derived below. Recalling the explicit representation of the Laplace--Beltrami operator on $\mathcal{S}^{d-1}$,
\begin{align*}
	\Delta_{{\cal S}^{d-1}} = \begin{cases}
	\pd^2_\th, \quad &\text{if $d=2$},
	\\
	\pd^2_\th + \cot(\th) \pd_\th + \frac 1{\sin^2(\th)} \pd^2_\phi, \quad & \text{if $d=3$},
	\end{cases}
\end{align*}
the full Laplacian in spherical coordinates decomposes as
\begin{align}
	\label{Delta:rad}	
	\Delta f = f'' + \frac {d-1}r f' + \frac 1{r^2} \Delta_{{\cal S}^{d-1}} f.
\end{align}
Substituting the spectral ansatz into the governing equations and exploiting the eigenfunction property
\begin{align}	\label{lap}
	\Delta_{{\cal S}^{d-1}}Z(\th,\phi)
	=\zeta_{\ell,d}Z(\th,\phi) ,
	\quad \text{with} \quad
	\zeta_{\ell,d} = \begin{cases}
	-\ell^2, \quad &\text{in $d=2$},
	\\
	-\ell(\ell+1), \quad & \text{in $d=3$},
	\end{cases}
\end{align}
we find that the linearised system {\eqref{pertrad:general} together with
\eqref{pertrad:general:initial}} reduces to the following radial boundary value problem:
\begin{subequations}\label{linstab}
\begin{alignat}{2}
	\label{linstab:1}
	&
	-  \Big( U '' + \frac {d-1}r U ' + \frac {\zeta_{\ell,d}}{r^2}  U  \Big)=  -    U \quad &&\text{in $\{{0<r<\qstar}\}$,}
	\\
	\label{linstab:3}
	& (\mu^*)'
	+ U  = -\alpha \frac {d-1}{r^2} \Big( 1 + \frac {\zeta_{\ell,d}}{d-1}\Big)
	\quad && \text{on $\{r=\qstar\}$,}
	\\
	\label{linstab:4}
	& 2 \dot{\delta}
	=
	-  \big( (\mu^* )'' +   U' \big)
	  \delta
	\quad && \text{on $\{r=\qstar\}$,}
\end{alignat}
\end{subequations}
and
\begin{align}
 \label{linstab:6} & U (r = 0) < \infty.
\end{align}
Within the inner region $r \in (0,\qstar)$, the radial profile is explicitly given by:
\begin{align*}
	U (r) & =
	\left (- \alpha \frac {d-1}{(\qstar)^2} \Big( 1+ \frac {\zeta_{\ell,d}}{d-1}\Big)
	-   \Big( \alpha \frac {d-1}\qstar +1 \Big)  \frac {\III_1(\qstar)}{\Iz (\qstar)}  \right ) \frac{\III_\ell(  r)}{\Il (\qstar)}.
\end{align*}
Further details concerning its derivation can be found in \cite[Appendix A.2]{OOL2}. Substituting this solution into the kinematic boundary condition yields a closed ordinary differential equation governing the temporal evolution of the perturbation amplitude $\delta(t)$:
\begin{align}
\frac{2\dot{\delta}}{\delta}   & =	-   \Big( \alpha \frac {d-1}{\qstar} +1 \Big)  \frac{1}{\Iz(\qstar)} \Big ( \III_1' (\qstar) - \III_1 (\qstar) \frac{\Il' (\qstar)}{\Il (\qstar)} \Big ) 
 +  \alpha \frac {d-1}{(\qstar)^2} \Big( 1+ \frac {\zeta_{\ell,d}}{d-1}\Big)
\frac{\Il' (\qstar)}{\Il (\qstar)}. 
\label{radial:ode:perturb}
\end{align}
Using ${\cal H}^+ (\qstar)=0$, see \eqref{radial:ode}, and the ODE $\mathbb{I}_1'(q) = \mathbb{I}_0(q) - \frac{d-1}{q} \, \mathbb{I}_1(q)$, we can then express \eqref{radial:ode:perturb} 
as
\begin{align}\label{decay}
	\frac{2\dot\delta}{\delta}
	= -\Big( 1 + \alpha\,\frac{d-1}{q^*} \Big) + \frac{(d-1)\,S_I}{q^*}
	+ \frac{\mathbb{I}_\ell'(q^*)}{\mathbb{I}_\ell(q^*)}
	\Big( S_I - \alpha\,\frac{(\ell-1)(\ell+d-1)}{(q^*)^2} \Big).
\end{align}	
To understand in more detail how the (in-)stability behavior depends on the mode number $\ell$ and the physical parameters, we now analyse the behavior of the third term on the right-hand side of \eqref{decay}.
	Setting $\nu=\ell+\frac{d-2}{2}$, the power-series representation of 
	\[\mathbb{I}_\nu(q)=
	\Big( \frac q2 \Big)^{\nu} \sum_{k=0}^{\infty}
		\frac{\big( \tfrac14 q^2 \big)^{k}}{k! \, \Gamma(\nu+k+1)} ,
	\]
	see \cite[Eq.~(10.25.2)]{NIST:DLMF}, and the identity $\Gamma(\nu+k+2)=(\nu+k+1)\Gamma(\nu+k+1) $ yield, by term-wise comparison, \[ 0<\mathbb{I}_{\nu+1}(q) <\frac{q}{2(\nu+1)}\mathbb{I}_\nu(q), \qquad \forall q > 0. \] Combining this with the recurrence relation $ \mathbb{I}_\nu'(q)=\mathbb{I}_{\nu+1}(q)+\frac{\nu}{q}\mathbb{I}_\nu(q)$, see \cite[Eq.~(10.29.2)]{NIST:DLMF}, gives 
		\begin{align*}
		\frac{\ell}{q} < \frac{\mathbb{I}_\ell'(q)}{\mathbb{I}_\ell(q)}
		< \frac{\ell}{q} + \frac{q}{2\ell + d},
		\qquad \forall q > 0, \ \ell \ge 0 .
	\end{align*}
This in particular implies that the dominating term on the right-hand side of \eqref{decay} is cubic in $\ell$ and has a negative sign.

\begin{remark}\label{rem:radial:perturb}
For sufficiently large $\ell$ the right hand side of \eqref{decay} 
 assumes increasingly negative values as the mode number $\ell$ grows. In particular, for large $\ell$ the right-hand side of \eqref{decay} is strictly monotonically decreasing in $\ell$. Consequently, for sufficiently large $\ell$, 
 we obtain exponential decay of the high-frequency perturbations, which implies linear stability at small spatial scales. The onset of instability is therefore contingent upon parameter regimes for which the amplification rate becomes positive for specific low-order modes $\ell$. In such cases, the interfacial perturbation $Y$ defined in \eqref{rad:ansatz} exhibits exponential growth, driving a morphological deformation of the spherical interface characterised by $\ell$-dependent oscillatory patterns. Heuristically, reducing the capillary regularization parameter $\alpha$ diminishes the stabilising curvature contribution, thereby expanding the range of unstable modes.
\end{remark}

\begin{figure}[h!]
	\centering
	\includegraphics[width=0.48\textwidth]{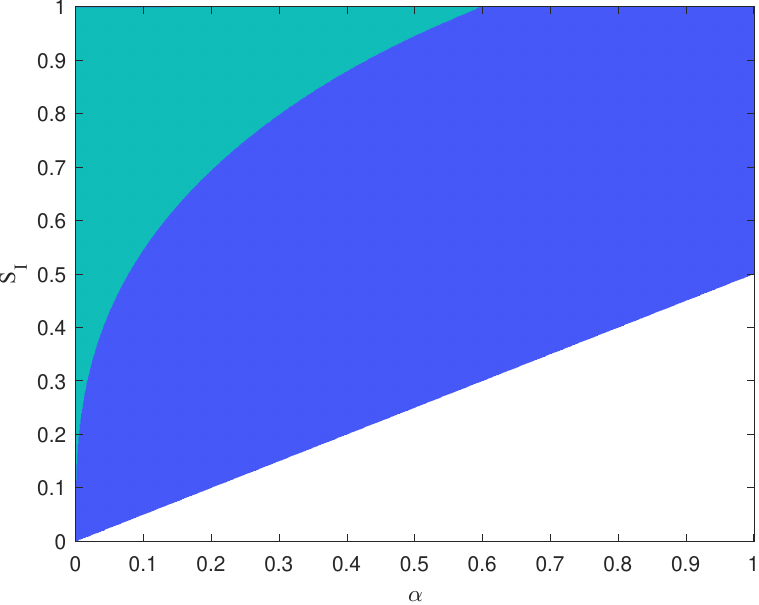}
	\includegraphics[width=0.48\textwidth]{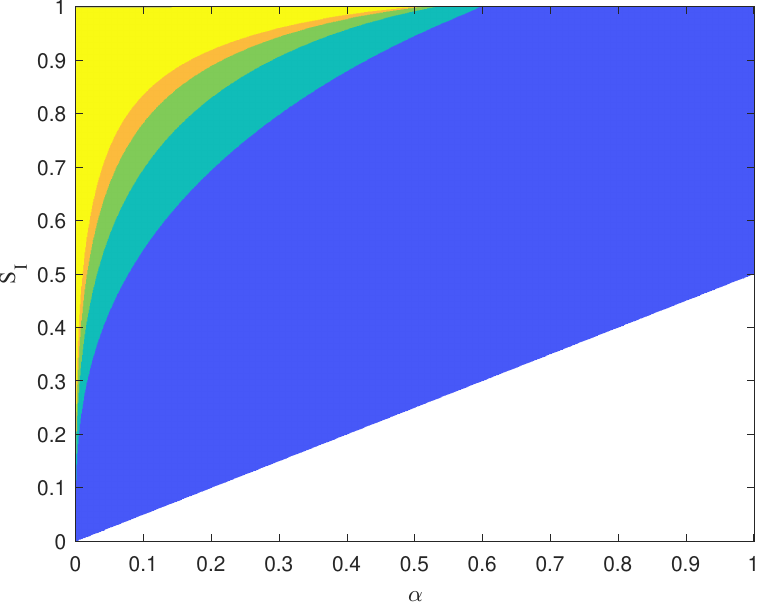}
	\caption{(Phase diagram for the plus-phase model) Radial setting $d = 2$ for the \emph{plus-phase model}.\\[1ex]
	(Left) A subdivision of the coloured region in the phase diagram where $\mathcal{H}^+$ in \eqref{radial:ode} admits a positive root into a region of stability (blue) and a region of instability (teal), both with respect to mode number $\ell = 2$ perturbations. The white region corresponds to the range of parameters where $\mathcal{H}^+$ may not exhibit a positive root.\\[1ex]
	(Right) A further subdivision of the coloured region where $\mathcal{H}^+$ in \eqref{radial:ode} admits a positive root into regions where the mode number $\ell \in \{2, 3, 4, 5\}$ has the largest positive value of the right-hand side of \eqref{radial:ode:perturb} amongst modes $2 \leq \ell \leq 5$. The colour correspondance is $\ell = 2$ (teal), $\ell = 3$ (green), $\ell = 4$ (orange) and $\ell = 5$ (yellow). The blue region corresponds to stability with respect to all mode numbers $\ell \in \{2, 3, 4, 5\}$.
	}
	\label{fig:rad2D_PD}
\end{figure}

For the two-dimensional setting $d= 2$ in Figure~\ref{fig:rad2D_PD} we show the partition of the phase space $(\alpha, S_I) \in [0,1] \times [0,1]$ into a white region where $\mathcal{H}^+$ in \eqref{radial:ode} may not exhibit a positive root, cf.~the range of parameters $S_I \leq \frac{\alpha}{2}$, compare Case (a) of \eqref{radial:para:root}, and a coloured region where $\mathcal{H}^+$ in \eqref{radial:ode} admits a positive root. In the left figure, we further subdivide the coloured region into a region of stability coloured in blue, and a region of instability coloured in teal, both with respect to perturbations with mode number $\ell = 2$. In particular, larger values of $\alpha$ stabilise the perturbations and higher values of $S_I$ combined with smaller values of $\alpha$ encourage the growth of perturbations.

In the right figure of Figure~\ref{fig:rad2D_PD} we further subdivide the coloured region into regions where the mode number $\ell \in \{2, 3, 4, 5\}$ has the largest positive value of the right-hand side of \eqref{radial:ode:perturb} in the range $2 \leq \ell \leq 5$.
Figure~\ref{fig:rad2D_PD} also show a monotonic growth with respect to $\ell$ which is in accordance with Remark \ref{rem:radial:perturb} where we observed the monotonic behavior for large $\ell$.

Next, we notice for dimension $d = 2$ and $\ell = 1$, the right-hand side of \eqref{radial:ode:perturb} is zero, thus perturbations with respect to mode $\ell = 1$ neither grow or decay. From the figure we observe some agreement with Remark \ref{rem:radial:perturb}, namely for certain parameter regimes perturbations with respect to higher mode numbers are stabilised, and only when $\alpha$ is reduced, the range of mode numbers for which perturbations are unstable can be expanded.
%

%

\subsection{Minus-phase model}
For the \emph{minus-phase model} a similar linear stability analysis would yield the following linearised system:
\begin{subequations} \label{pertrad:general:minus}
\begin{alignat}{2}
	\label{pertrad:general:1:m}
	&
	-  \Delta  u =  -   u  \qquad &&\text{in $\{\qstar <r < R\}$,}\\
	\label{pertrad:general:3:m}
	& (\mu^*)' Y + u= -  \alpha 	\frac {d-1}{(\qstar)^2} \big( \rad + \frac 1{d-1} \Delta_{{\cal S}^{d-1}} \rad \big)
	\qquad && \text{on $\{r=\qstar\}$,}
	\\
	& 	\label{pertrad:general:5:m}
    2 \dot{Y} =
	 (\mu^*)'' \rad
	+   u'
	\qquad && \text{on $\{r=\qstar\}$,}
	\\
	\label{pertrad:general:6:m}
	& u' =0
	\qquad &&
	\text{on $\{r=R\}$,}
\end{alignat}
\end{subequations}
so that with the ansatz \eqref{rad:ansatz} and \eqref{perturb:u:ansatz} we arrive at the following radial boundary value problem:
\begin{subequations}\label{linstab"minus}
\begin{alignat}{2}
	\label{linstab:1:m}
	&
	-  \Big( U '' + \frac {d-1}r U ' + \frac {\zeta_{\ell,d}}{r^2}  U  \Big)=  -    U \quad &&\text{in $\{{\qstar<r<R}\}$,}
	\\
	\label{linstab:3:m}
	& (\mu^*)'
	+ U  = -\alpha \frac {d-1}{r^2} \Big( 1 + \frac {\zeta_{\ell,d}}{d-1}\Big)
	\quad && \text{on $\{r=\qstar\}$,}
	\\
	\label{linstab:4:m}
	& 2 \dot{\delta}
	= \big( (\mu^* )'' +  U' \big)
	  \delta
	\quad && \text{on $\{r=\qstar\}$,}
	\\
	\label{linstab:5:m}
	& U' =0
	\quad && \text{on $\{r=R\}$.} 
\end{alignat}
\end{subequations}
The radial profile is explicitly given as
\begin{align*}
U(r) & = \frac{\Big (\Kl(r) - \frac{\Kl'(R)}{\Il'(R)} \Il(r) \Big )}{\Big ( \Kl(\qstar) - \frac{\Kl'(R)}{\Il'(R)} \Il(\qstar)\Big )} \\
& \quad \times \left ( - \alpha \frac{d-1}{(\qstar)^2} \Big ( 1 + \frac{\zeta_{\ell,d}}{d-1} \Big ) + \Big ( \alpha \frac{d-1}{\qstar} - 1 \Big ) \frac{\KKK_1(\qstar) - \frac{\KKK_1(R)}{\III_1(R)} \III_1(\qstar)}{\Kz(\qstar) + \frac{\KKK_1(R)}{\III_1(R)} \Iz(\qstar)} \right ),
\end{align*}
and the ordinary differential equation \eqref{linstab:4:m} for the temporal evolution of the perturbation amplitude $\delta(t)$ reads as
\begin{align}	
\non
\frac{2\dot{\delta}}{\delta}   & =	
- \Big (\alpha \frac {d-1}{\qstar}- 1 \Big)
	\frac { \KKK_1'(\qstar) - \frac{\KKK_1(R)}{\III_1(R)} \III_1'(\qstar)}{ \Kz(\qstar) + \frac {\KKK_1(R)}{\III_1(R)} \Iz (\qstar) }\\ 
\notag & \quad	+  \frac{\Kl'(\qstar) - \frac{\Kl'(R)}{\Il'(R)} \Il'(\qstar)}{\Kl(\qstar)-  \frac{\Kl'(R)}{\Il'(R)} \Il(\qstar)} \\
 & \qquad \times \left (- \alpha \frac {d-1}{(\qstar)^2} \Big( 1+ \frac {\zeta_{\ell,d}}{d-1}\Big)
	+ \Big(\alpha \frac {d-1}{\qstar}- 1 \Big)
	\frac {\KKK_1(\qstar)-  \frac {\KKK_1(R)}{\III_1(R)} \III_1(\qstar)}{ \Kz(\qstar) + \frac {\KKK_1(R)}{\III_1(R)} \Iz(\qstar)} \right ).
	\label{radial:ode:perturb:minus}
\end{align}
Let $A := \mathbb{K}_1(R)/\mathbb{I}_1(R)$. Since $q^*$ is a root of $\mathcal{H}^-$, \eqref{radial:ode:minus} gives
\begin{align}
	\label{root:minus}
	\Big( \alpha \, \frac{d-1}{q^*} - 1 \Big)
	\frac{\mathbb{K}_1(q^*) - A \, \mathbb{I}_1(q^*)}{\mathbb{K}_0(q^*) + A \, \mathbb{I}_0(q^*)} = S_I .
\end{align}
Hence the last bracket in \eqref{radial:ode:perturb:minus} equals
\begin{align*}
	- \alpha \, \frac{d-1}{(q^*)^2} \Big( 1 + \frac{\zeta_{\ell,d}}{d-1} \Big) + S_I
	= S_I + \alpha \, \frac{(\ell-1)(\ell+d-1)}{(q^*)^2} .
\end{align*}
For the first term in \eqref{radial:ode:perturb:minus} we use the identities
\[
	\mathbb{I}_1'(q) = \mathbb{I}_0(q) - \frac{d-1}{q} \, \mathbb{I}_1(q), \quad \mathbb{K}_1'(q) = - \mathbb{K}_0(q) - \frac{d-1}{q} \, \mathbb{K}_1(q),
\]
which give
\[
	\mathbb{K}_1'(q^*) - A \, \mathbb{I}_1'(q^*)
	= - \big( \mathbb{K}_0(q^*) + A \, \mathbb{I}_0(q^*) \big)
	- \frac{d-1}{q^*} \big( \mathbb{K}_1(q^*) - A \, \mathbb{I}_1(q^*) \big) .
\]
Together with \eqref{root:minus} we obtain
\begin{align*}
	- \Big( \alpha \, \frac{d-1}{q^*} - 1 \Big)
	\frac{\mathbb{K}_1'(q^*) - A \, \mathbb{I}_1'(q^*)}{\mathbb{K}_0(q^*) + A \, \mathbb{I}_0(q^*)}
	= \alpha \, \frac{d-1}{q^*} - 1 + \frac{(d-1) \, S_I}{q^*} .
\end{align*}
Inserting both expressions into \eqref{radial:ode:perturb:minus} yields
\begin{align}
	\label{decay:minus} \frac{2\dot\delta}{\delta}
	= \alpha \, \frac{d-1}{q^*} - 1 + \frac{(d-1) \, S_I}{q^*}
	+ G_\ell(q^*) \Big( S_I + \alpha \, \frac{(\ell-1)(\ell+d-1)}{(q^*)^2} \Big),
	\end{align}
	where 
	\[
	G_\ell(r) := 
	\frac{V'_\ell(r)}{V_\ell(r)}, \quad V_\ell := \mathbb{K}_\ell - \frac{\mathbb{K}_\ell'(R)}{\mathbb{I}_\ell'(R)} \, \mathbb{I}_\ell.
	\]
	One can now  use 
		$\mathbb{K}_\ell > 0$, $\mathbb{I}_\ell > 0$, $\mathbb{K}_\ell'(R) < 0$ and $\mathbb{I}_\ell'(R) > 0$,
		to deduce $V_\ell > 0$ in $(0,R]$. By \eqref{linstab:1:m}, $V_\ell$ satisfies
		\begin{align*}
			\big( r^{d-1} V_\ell'(r) \big)' = r^{d-1} \Big( 1 + \frac{\ell(\ell+d-2)}{r^2} \Big) V_\ell(r) > 0,
	\qquad \forall r \in  (0,R),
		\end{align*}
		where we used $-\zeta_{\ell,d} = \ell(\ell+d-2)$. Thus $r \mapsto r^{d-1} V_\ell'(r)$ is strictly increasing, and it vanishes at
		$r = R$ by \eqref{linstab:5:m}. Consequently $V_\ell' < 0$ in $(0,R)$, and therefore
		\begin{align}
			\label{G:sign}
			G_\ell(q) < 0, \qquad \forall q \in (0,R) .
		\end{align}
	A more detailed analysis gives that $G_\ell(q) = - \frac{\ell+d-2}{q} + O\big( \frac1\ell \big)$ as  $\ell \to \infty$ and 
that perturbations from large $\ell$ are stable.
This behavior contrasts with the \emph{plus-phase model} in two key aspects: 
first, due to $G_\ell<0$, the term $G_\ell S_I$ in \eqref{decay:minus} is destabilising if $S_I<0$, whereas in the \emph{plus-phase model} the term $\frac{\mathbb{I}_\ell'}{\mathbb{I}_\ell} S_I$ in \eqref{decay} is destabilising if $S_I>0$ as $\frac{\mathbb{I}_\ell'}{\mathbb{I}_\ell}>0$; second, $G_\ell$ explicitly depends on the outer boundary radius $R$, capturing exterior confinement effects, whereas the plus-phase counterpart is entirely independent of $R$.

\begin{figure}[h!]
\centering
\includegraphics[width=0.48\textwidth]{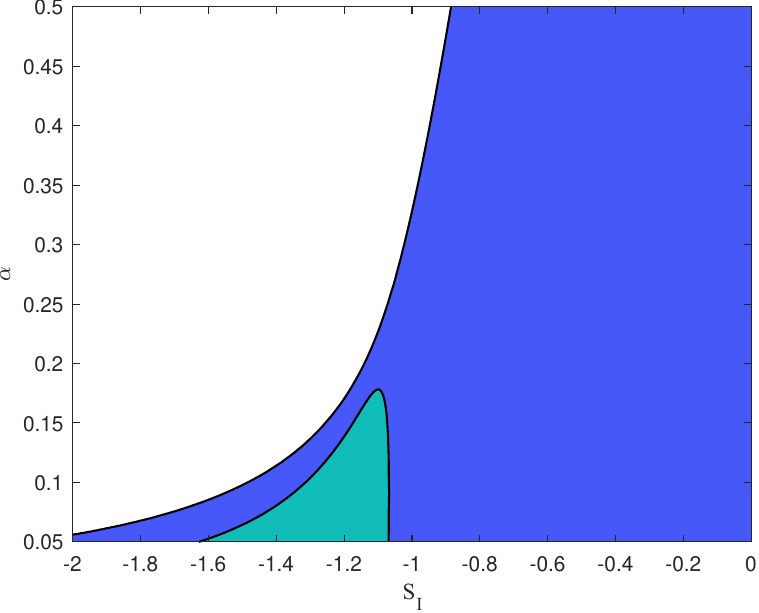}
\includegraphics[width=0.48\textwidth]{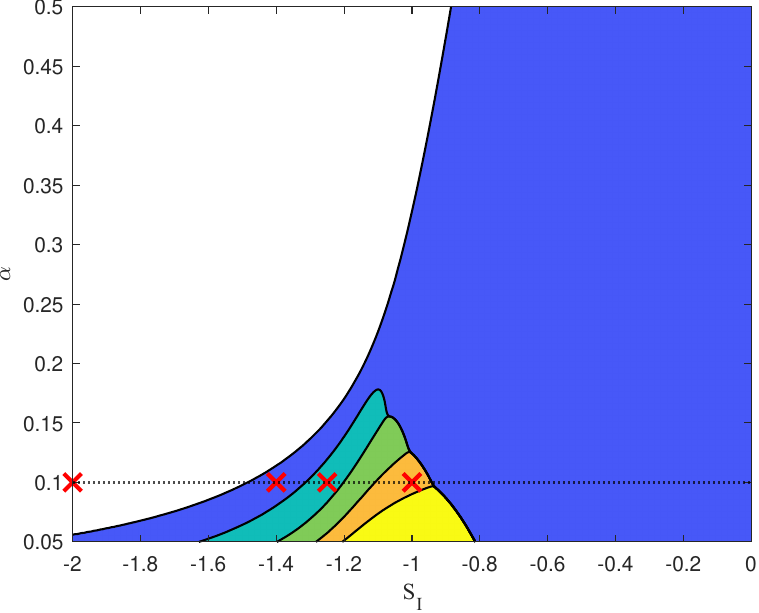}
\caption{(Phase diagram for the minus-phase model) Radial setting $d = 2$ for the \emph{minus-phase model} with $R = 4$.\\[1ex]
(Left) A subdivision of the coloured region in the phase diagram where $\mathcal{H}^{-}$ in \eqref{radial:ode:minus} admits a positive root into a region of stability (blue) and a region of instability (teal), both with respect to mode number $\ell = 2$ perturbations. The white region corresponds to the range of parameters where ${\cal H}^{-}$ may not exhibit a positive root. \\[1ex]
(Right) A further subdivision of the coloured region where $\mathcal{H}^-$ in \eqref{radial:ode:minus} admits a positive root into regions where the mode number $\ell \in \{2, 3, 4, 5\}$ has the largest positive value of the right-hand side of \eqref{radial:ode:perturb:minus} amongst modes $2 \leq \ell \leq 5$. The colour correspondance is $\ell = 2$ (teal), $\ell = 3$ (green), $\ell = 4$ (orange) and $\ell = 5$ (yellow). The blue region corresponds to stability with respect to all mode numbers $\ell \in \{2, 3, 4, 5\}$. A dotted horizontal line at $\alpha =0.1$ is drawn with four red crosses, corresponding to (from left to right) $(S_I, \alpha) = (-2, 0.1), (-1.4, 0.1), (-1.25, 0.1)$ and $(-1, 0.1)$, respectively. These correspond to later Figures~\ref{fig:pd_a01_SI-2_cp-4}, \ref{fig:pd_a01_SI-14_cp-4}, \ref{fig:pd_a01_SI-125_cp-4} and \ref{fig:pd_a01_SI-1_cp-4}.}
\label{fig:rad2D_minus}
\end{figure}

For the \emph{minus-phase model}, in Figure~\ref{fig:rad2D_minus} we plot the partitioning of the phase space $(S_I, \alpha) \in [-2, 0] \times [0.05, 0.5]$ into a white region where $\mathcal{H}^-$ in \eqref{radial:ode:minus} may not exhibit a positive root and a coloured region where $\mathcal{H}^-$ in \eqref{radial:ode:minus} admits a positive root. In the left figure, we further subdivide the coloured region into a region of stability coloured in blue and a region of instability coloured in teal, both with respect to perturbations with mode number $\ell = 2$. 
In contrast to the \emph{plus-phase model}, here we find that $S_I$ plays a more dominant role in matters of stability, where small values of $S_I$ stabilise the perturbations. For instability, a combination of small values of $\alpha$ and an intermediate range of negative values of $S_I$ is needed.

In right figure of Figure~\ref{fig:rad2D_minus} we further subdivide the coloured region into where the mode number $\ell \in \{2, 3, 4, 5\}$ has the largest positive value of the right-hand side of \eqref{radial:ode:perturb:minus} amongst modes $2 \leq \ell \leq 5$. We observe that instability for the higher mode numbers allows for a larger range of $S_I$, but the corresponding range of $\alpha$ shrinks.

\section{Numerical Computations}\label{sec:num}
\newcommand{\nabs}{\nabla_{\!s}}
\newcommand{\id}{{\rm id}}
\renewcommand{\vec}{\boldsymbol}

In this section, we present numerical computations for both the phase field model \eqref{eq:sys21}, as well as for the one-sided sharp interface models
\eqref{SharpI} and \eqref{SharpI:single:m}. Our numerical method for 
\eqref{eq:sys21} is the finite element approximation from \cite{OOL}, and so
we only state it briefly. The numerical method for the one-sided 
sharp interface models is an adaptation of the scheme from \cite{OOL2},
using the techniques devised in \cite{crystal} to deal with the one-sided nature of the models. Both methods are implemented with the help of the finite element toolbox ALBERTA, see \cite{Alberta}. To increase computational efficiency, we employ the same mesh refinement strategies as in \cite{OOL,OOL2}. In addition, for the sharp interface simulations we at times perform topological changes, once again using the same strategies as in \cite{OOL2}. 

In all our numerical simulations we will choose the plus phase to represent the interior of the modelled protocells. Moreover, we work exclusively with the nondimensional framework and therefore report only the values of $\alpha$ and $S_I$ for the sharp interface computations. Then, the corresponding parameters for the phase field models are: 
\[
\beta = \frac{3{\alpha}}{\sqrt{2}} ,\quad
{K_\pm = 2 \beta \rho_\pm,}
\quad
{L = \frac{3\sqrt{2}}4 \Big(S_I - \sqrt{2} \beta (\rho_+ - \rho_-) \Big),}
\] 
and 
\[
m_+ = 1, \quad  m_- = 0, \quad \rho_\pm = 1,
	\quad 
	S_+ = -1, \quad S_- =0, 
\]
for the \emph{plus-phase model} and
\[
m_+ = 0, \quad m_- = 1, \quad \rho_\pm = 1,
\quad 
	S_+ = 0, \quad S_- =1,  
\]
for the \emph{minus-phase model}, respectively.

\subsection{Finite element schemes}
\newcommand{\Sh}{\mathcal S_h}

\subsubsection{Phase field finite element scheme}

We assume that $\Omega$ is a polyhedral domain and let $\mathcal{T}_h$ be a regular triangulation of $\Omega$ into disjoint open simplices. Associated with $\mathcal{T}_h$ is the piecewise linear finite element space
\[
\Sh = \{ \zeta \in C^0(\overline{\Omega}) \, : \, \zeta \vert_{o} \in P_1(o) \quad \forall o \in \mathcal{T}_h\},
\]
where we denote by $P_1(o)$ the set of all affine linear functions on $o$,
see \cite{Ciarlet78}. 
Let $(\cdot,\cdot)_{\Omega}^h$ be the usual mass lumped $L^2$-inner product on $\Omega$ associated with $\mathcal{T}_h$, and let $\pi^h : C^0(\overline{\Omega}) \to \Sh$
be the standard interpolation operator. 
In addition, let $\Delta t$ denote a chosen uniform time step size. Then our finite element approximation of \eqref{eq:sys21} 
is given as follows. Let $\phi_h^0 \in \Sh$, e.g., $\phi_h^0 = \pi^h \varphi_0$
if $\varphi_0 \in C^0(\overline\Omega)$. 
Then, for $n \geq 0$, find $(\phi_h^{n+1}, \mu_h^{n+1}) \in \Sh \times \Sh$ such that
\begin{subequations} \label{eq:FEA}
\begin{alignat}{2}
\frac{1}{\Delta t}(\phi_h^{n+1}-\phi_h^n, \chi_h)_{\Omega}^h + (m(\phi_h^n) \nabla \mu_h^{n+1}, \nabla \chi_h)_{\Omega}^h = (S_\eps(\phi_h^n), \chi_h)_{\Omega}^h
\quad & \forall\chi_h\in \Sh, \\
\beta \varepsilon (\nabla \phi_h^{n+1}, \nabla \eta_h)_{\Omega}^h + \frac{\beta}{\eps} (\psi'(\phi_h^{n+1}), \eta_h)_{\Omega}^h - (\mu_h^{n+1}, \eta_h)_{\Omega}^h  = 0
\quad &\forall\eta_h\in \Sh.
\end{alignat}
\end{subequations}
The scheme \eqref{eq:FEA} is the same as \cite[(6.1)]{OOL}. 
The nonlinear system of equations arising at each time level of 
\eqref{eq:FEA} are solved with the help of Newton's method. 
The resulting linear systems at each iteration in two spatial dimensions
are solved by direct factorization using the package UMFPACK, 
see \cite{Davis04}, and in three spatial dimensions with a $V$-cycle multigrid solver using a block Gauss--Seidel smoother.

For the initial data $\varphi_0$ in general we choose a diffuse interface
representation of a desired sharp interface, with signed distance function
$d_0 : \overline\Omega \to \mathbb R$. In particular, unless otherwise stated, 
we let
\begin{equation} \label{eq:tanh0}
 \varphi_0(x) = \tanh \Big ( \frac{d_0(x)}{\eps \sqrt{2}} \Big ) .
\end{equation}

All our numerical simulations for the diffuse interface model \eqref{eq:sys21} employ the quartic potential \eqref{quartic} together with
the interpolation functions \eqref{def:G1}, \eqref{def:G23} and 
\eqref{eq:G4} for the fixed value $r_c=1$. 
For the mobility function $m$ we choose \eqref{eq:mob}.

\subsubsection{Sharp interface parametric finite element schemes}

Let us consider the \emph{minus-phase model} \eqref{SharpI:single:m}, the treatment
of the \emph{plus-phase model} \eqref{SharpI} is analogous.
Similarly to \cite{OOL2}, we first introduce the following 
natural variational formulation. 
Given $\Sigma(0)$, for $t > 0$ find $\Sigma(t)$
separating $\Omega$ into $\Omega^+(t)$ and $\Omega^-(t)$, as well as
$\mu \in H^1(\Omega^-(t))$ and $\kappa \in L^2(\Sigma(t))$ such that
\begin{subequations}
\label{eq:2}
\begin{align}
& m_- \left(\nabla\mu, \nabla{\phi}\right)_{\Omega^-(t)} 
+ \rho_-\left(\mu, {\phi}\right)_{\Omega^-(t)}
- 2\left\langle \mathcal{V},{\phi}\right\rangle_{\Sigma(t)} 
=S_-\left(1, {\phi}\right)_{\Omega^-(t)}+S_I \left\langle 1,{\phi}\right\rangle_{\Sigma(t)} , \nonumber \\ & \hspace{9cm}
\quad\forall {\phi}\in H^1(\Omega^-(t)),\label{eq:2aa} \\ &
 \left\langle \mu - \alpha\kappa , \chi \right\rangle_{\Sigma(t)} = 0,
\quad\forall \chi \in L^2(\Sigma(t))
,\label{eq:2bb} \\
& \left\langle \kappa\bnu,\vec\eta \right\rangle_{\Sigma(t)} 
+ \left\langle \nabs\vec\id,\nabs\vec\eta 
\right\rangle_{\Sigma(t)} = 0,
\quad\forall \vec\eta \in H^1(\Sigma(t))^d, \label{eq:2cc}
\end{align}
\end{subequations}
where we have also recalled a weak formulation for the curvature $\kappa$,
see, e.g., \cite[Remark~22]{bgnreview}. In addition, we let 
$(\cdot,\cdot)_{\Omega^-(t)}$ and
$\langle\cdot,\cdot\rangle_{\Sigma(t)}$ denote the {$L^2$-inner} products 
in $\Omega^-(t)$ and on $\Sigma(t)$, respectively. 
Moreover, $\vec\id$ denotes the identity function in $\mathbb R^d$.

Let $\Delta t > 0$ denote the time step size as before, 
and let $t_n = n\Delta t$, $n \geq 0$.
For $n\geq0$, let ${\cal T}^n_h$ 
be a regular partitioning of $\Omega$ into disjoint open simplices, so that 
$\overline{\Omega}=\cup_{\sigmaO\in{\cal T}^n_h}\overline{\sigmaO}$. 
As before, we associate with ${\cal T}^n_h$ the finite element space
$\Snh := \{\chi \in C^0(\overline{\Omega}) : \chi\!\mid_{\sigmaO} 
\text{ is affine } \ \forall \sigmaO \in {\cal T}^n_h\}$.
Let $K_\Omega^n$ be the number of nodes of ${\cal T}^n_h$ and  
let $\{\vec{p}^n_{j}\}_{j=1}^{K_\Omega^n}$ be the coordinates of these nodes.
Let $\{\phi_{j}^n\}_{j=1}^{K_\Omega^n}$ be the standard basis functions 
for $\Snh$.
We introduce $\pi^n_h:C^0(\overline{\Omega})\to \Snh$, the interpolation
operator, such that $(\pi^n_h \eta)(\vec{p}_k^n)= \eta(\vec{p}_k^n)$ 
for $k=1{,\dots,} K_\Omega^n$. 

Let $\Sigma^{n}$ be a polyhedral hypersurface, see
\cite[Definition~41]{bgnreview}, approximating the
surface $\Sigma(t_n)$, $n\geq0$.
In particular, let $\Sigma^n=\bigcup_{j=1}^{J^n_\Sigma} 
\overline{\sigma^n_j}$,
where $\{\sigma^n_j\}_{j=1}^{J^n_\Sigma}$ is a family of mutually disjoint open 
$(d-1)$-simplices, $J^n_\Sigma \in \enne$. 
Then let
\begin{equation*}
\Vhm := \{\vec\chi \in C^0(\Sigma^n)^d:\vec\chi\!\mid_{\sigma^n_j}
\mbox{ is affine}\ \forall\ j=1,\ldots, J^n_\Sigma\} 
=: \Whm^d,
\end{equation*}
where $\Whm \subset H^1(\Sigma^n)$ is the space of scalar continuous
piecewise linear functions on $\Sigma^n$.
For later purposes, we also introduce 
$\pi_{\Sigma^n}: C^0(\Sigma^n)\to \Whm$, the standard Lagrange
interpolation operator.

As in the continuous case, we let $\langle\cdot,\cdot\rangle_{\Sigma^n}$ denote
the {$L^2$-inner} product on $\Sigma^n$. In addition,
for scalar and vector valued functions $v,w$ that are piecewise continuous, 
with possible jumps across the edges of $\{\sigma_j^n\}_{j=1}^{J^n_\Sigma}$,
we introduce the mass lumped inner product
$\langle\cdot,\cdot\rangle_{\Sigma^n}^h$ as
\begin{equation*}
\langle v, w \rangle^h_{\Sigma^n} :=
\frac1d \sum_{j=1}^{J^n_\Sigma} |\sigma^n_j|\sum_{k=1}^{d} 
(v\cdot w)((\vec{q}^n_{j_k})^-),
\end{equation*}
where $\{\vec{q}^n_{j_k}\}_{k=1}^{d}$ are the vertices of $\sigma^n_j$,
and where we define $v((\vec{q}^n_{j_k})^-):=
\underset{\sigma^n_j\ni \vec{p}\to \vec{q}^n_{j_k}}{\lim}\, v(\vec{p})$.

Given $\Sigma^n$, we 
let $\Omega^n_+$ denote the interior of $\Sigma^n$ and let
$\Omega^n_-$ denote the exterior of $\Sigma^n$, so that
$\Sigma^n = \partial \Omega^n_+ = \overline\Omega^n_+ \cap 
\overline\Omega^n_-$. 
We then partition the elements of the bulk mesh 
$\mathcal{T}^n_h$ into interior, exterior and interfacial elements as follows.
Let
\begin{align}
\mathcal{T}^n_- & := \{ o \in \mathcal{T}^n_h : o \subset
\Omega^n_- \} , \qquad 
\mathcal{T}^n_+ := \{ o \in \mathcal{T}^n_h : o \subset
\Omega^n_+ \} , \nonumber \\
\mathcal{T}^n_{\Sigma^n} & := \{ o \in \mathcal{T}^n_h : o \cap
\Sigma^n \not = \emptyset \} . \label{eq:partT}
\end{align}
Clearly $\mathcal{T}^n_h = \mathcal{T}^n_- \cup \mathcal{T}^n_+ \cup
\mathcal{T}^n_{\Sigma^n}$ is a disjoint partition, which in practice
can easily be found{,} e.g.{,} with the Algorithm~4.1 in \cite{crystal}.
In order to define the one-sided bulk test and trial functions, we first let
\[
\Omega^{n,h}_{-}
= \{ \vec p \in \Omega : (\pi^n_h \charfcn{\overline{\Omega^n_-}})(\vec p) > 0\}
\ \text{ and } \
\Omega^{n,h}_{+} = \Omega \setminus \overline{\Omega^{n,h}_{-}},
\]
where $\charfcn{\Omega^n_-}$ denotes the characteristic function on 
$\Omega^n_-$, compare with \cite[(4.6)]{crystal}. We then let
\[
\SmD = \{ \chi \in \Snh : \chi(\vec p^n_j) = 0 \ \text{if} \
\supp(\phi^n_j) \subset \Omega^{n,h}_{+} \},
\]
compare with \cite[(3.5)]{crystal}.

In addition, we let $\vec{\nu}^n$ denote the piecewise constant unit normal 
to $\Sigma^n$, pointing into $\Omega^n_+$. We also
introduce the vertex normal vector $\vec\omega^{n}\in \Vhm$ as
the mass-lumped {$L^2$-}projection of $\vec\nu^n$ onto $\Vhm$, that is
\begin{equation} \label{eq:nuhomegah}
\langle \vec\omega^n, \vec\eta \rangle^h_{\Sigma^n} = 
\langle \vec\nu^n, \vec\eta \rangle_{\Sigma^n},
\quad\forall\vec\eta\in \Vhm.
\end{equation}

In addition to $\Sigma^n$, we will denote by $\mu_h^n$, and $\kappa_h^n$ 
the discrete approximations to $\mu(t_n)$ and $\kappa(t_n)$,
respectively.
We will parameterise the new discrete surface 
$\Sigma^{n+1}$ over $\Sigma^n$, with the help of a parameterization
$\vec{X}^{n+1} \in \Vhm$.

Our finite element scheme to approximate \eqref{eq:2} can then be formulated as
follows:
Given $\Sigma^0$, for $n\geq0$, find $(\mu_h^{n+1},\vec X^{n+1},\kappa_h^{n+1})
\in \SmD \times \Vhm \times \Whm$ such that 
$\vec X^{n+1} - \vec\id\!\mid_{\Sigma^n} \in \Vhm$ and
\begin{subequations} \label{eq:MSHG}
\begin{align}
& m_- \left(\nabla\mu_h^{n+1}, \nabla{\phi}\right)_{\Omega^{n,h}_{-}} 
+\rho_-\left(\mu_h^{n+1}, {\phi} \right)_{\Omega^{n,h}_{-}}
- 2 \left\langle \pi_{\Sigma^n}\left[
\frac{\vec X^{n+1}-\vec\id}{\Delta t} \cdot \vec\omega^n \right], {\phi}
\right\rangle_{\Sigma^n} 
\nonumber \\ & \hspace{6cm}
= S_- \left(1, {\phi}\right)_{\Omega^{n,h}_{-}}
+S_I\left\langle 1, {\phi} \right\rangle_{\Sigma^n},
\quad\forall {\phi} \in \SmD, \label{eq:MSHGa}\\
& \left\langle \mu_h^{n+1}, \chi \right\rangle_{\Sigma^n} 
- \alpha \left\langle \kappa_h^{n+1}, 
\chi\right\rangle_{\Sigma^n}^h =0, \quad\forall \chi \in \Whm
, \label{eq:MSHGb} \\
& \left\langle \kappa_h^{n+1}\vec\nu^n, \vec\eta \right\rangle_{\Sigma^n}^h + 
\left\langle \nabs\vec X^{n+1}, \nabs\vec\eta \right\rangle_{\Sigma^n}
= 0, \quad\forall \vec\eta \in \Vhm,
 \label{eq:MSHGc} 
\end{align}
\end{subequations}
and set $\Sigma^{n+1} = \vec X^{n+1}(\Sigma^n)$. 
The scheme \eqref{eq:MSHG} is the natural one-sided variant of 
\cite[(5.5)]{OOL2}, utilizing the techniques from \cite{crystal}. 

For the solution of the linear systems of equations arising from 
\eqref{eq:MSHG} we use the direct factorization package UMFPACK, 
see \cite{Davis04}, in two dimensions, and employ a Schur complement approach that is
solved with a preconditioned conjugate gradient method in three dimensions, see
\cite{dendritic} for details.

\subsection{Numerical computations for the plus-phase model}

\begin{figure}
\center
\includegraphics[angle=-0,width=0.45\textwidth]{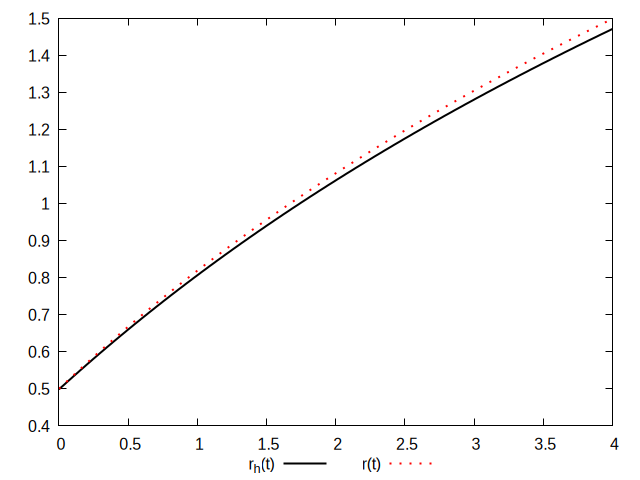}
\includegraphics[angle=-0,width=0.45\textwidth]{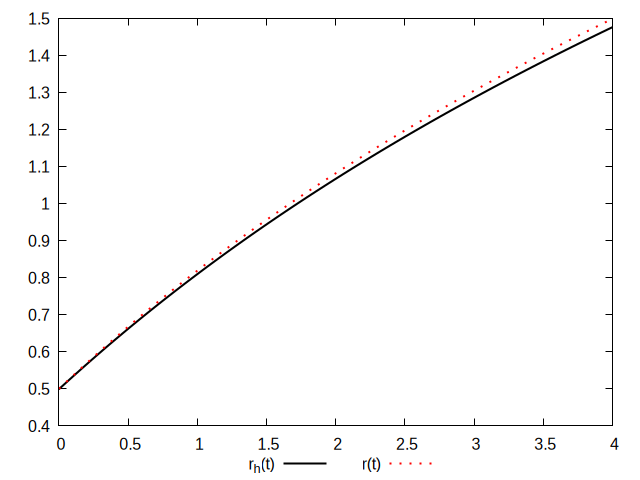}
\caption{(Plus-phase model: $\Omega = \bB_4^2(0)$, $\eps=\frac1{16\pi}, \frac1{32\pi}$) 
Comparison of $r(t)$ and $r_h(t)$ for $\alpha=0.1$ and $S_I=1$.
}
\label{fig:innercurved16pi}
\end{figure}

Some numerical results for radial evolutions in $\Omega = \bB_4^2(0)$
for the phase field model are shown in Figure~\ref{fig:innercurved16pi}.
Here we choose $\eps=\frac1{16\pi}$ or $\frac1{32\pi}$ and let the initial
radius be $r(0)=0.5$.
The physical parameters are $\alpha=0.1$ and $S_I=1$.
In particular, the plots in Figure~\ref{fig:innercurved16pi} compare
for a given $\phi^n_h$ the true radius $r(t_n)$, 
where $t_n = n \Delta t$, with the unique value $r_h(t_n) > 0$ such that 
$\phi^n_h(r_h(t_n), 0) = 0$. 

Notice that for the following sharp interface computations the size and shape 
of $\Omega$ play no role, as long as it is
sufficiently large to contain the moving interface. We nevertheless state the
computational domain used in all the presented simulations.

\begin{figure}
\center
\includegraphics[angle=-0,width=0.45\textwidth]{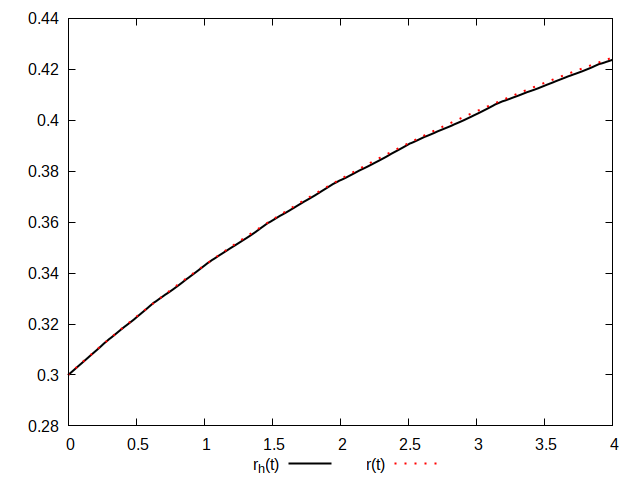}
\includegraphics[angle=-0,width=0.45\textwidth]{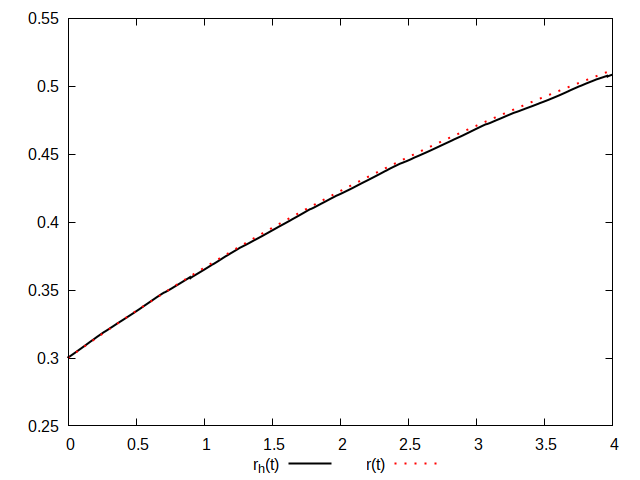}
\caption{(Plus-phase model: $\Omega = (0,8)^d$)
Comparison of $r(t)$ and $r_h(t)$ for $\alpha=0.01$ and $S_I=0.25$ 
with $d=2$ (left) and $d=3$ (right).
}
\label{fig:rad2}
\end{figure}

We begin with a single spherical interface inside the domain $\Omega = (0,8)^d$.
For the initial radius we choose $r(0) = 0.3$, while the physical
parameters are chosen as $\alpha=0.01$ and $S_I=0.25$.
A comparison between $r(t_n)$ and the radii $r_h(t_n)$ of the discrete
interfaces $\Sigma^n$, where now
\[
r_h(t_n) = \max\{|\vec z| : \vec z \in \Sigma^n \},
\]
is shown in Figure~\ref{fig:rad2}. We observe that in both two and three dimensions our
numerical algorithm matches the evolution of the true solution obtained from solving the ODE very well.

\begin{figure}
\center
\mbox{
\hspace*{-1cm}
\includegraphics[width=0.25\textwidth]{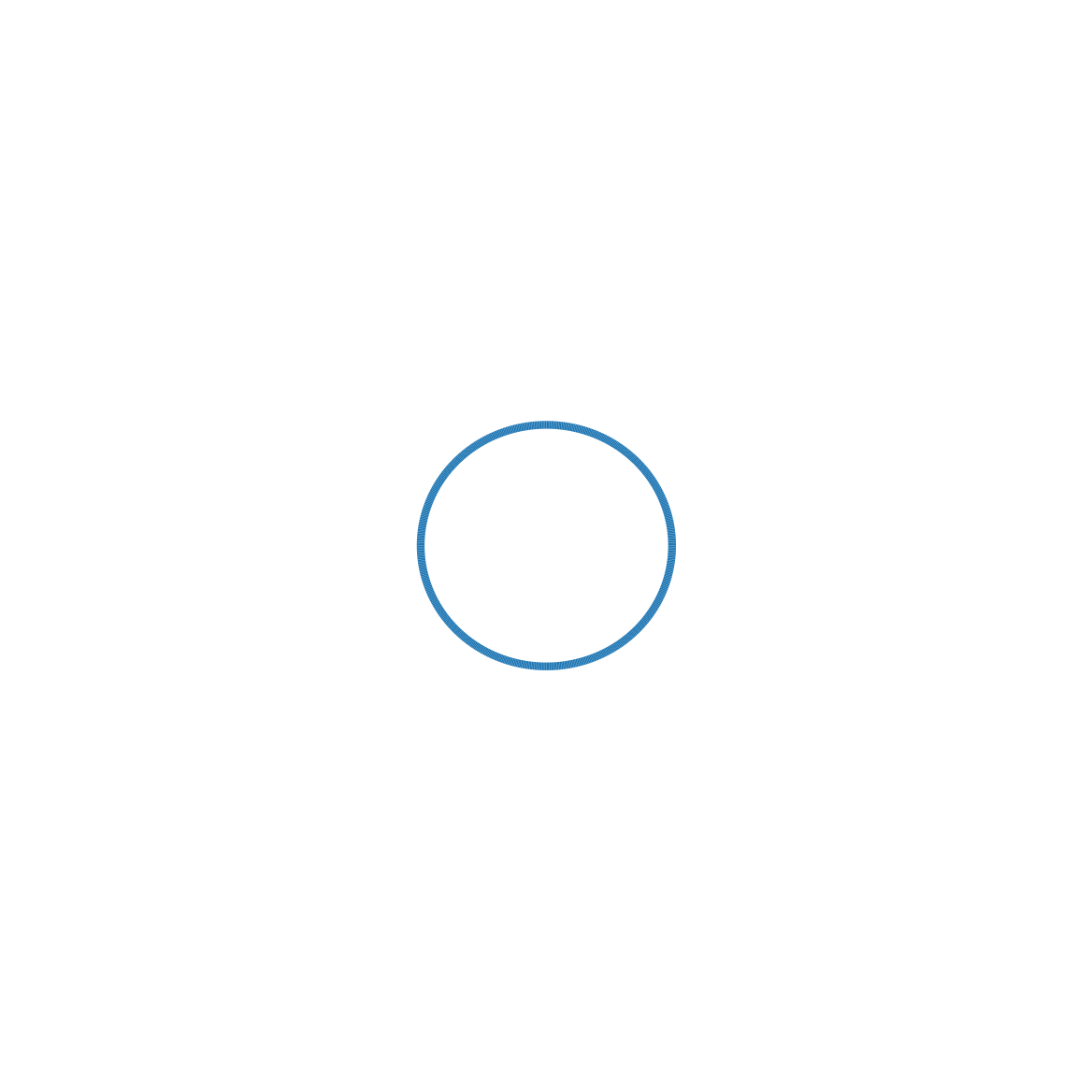}
\includegraphics[width=0.25\textwidth]{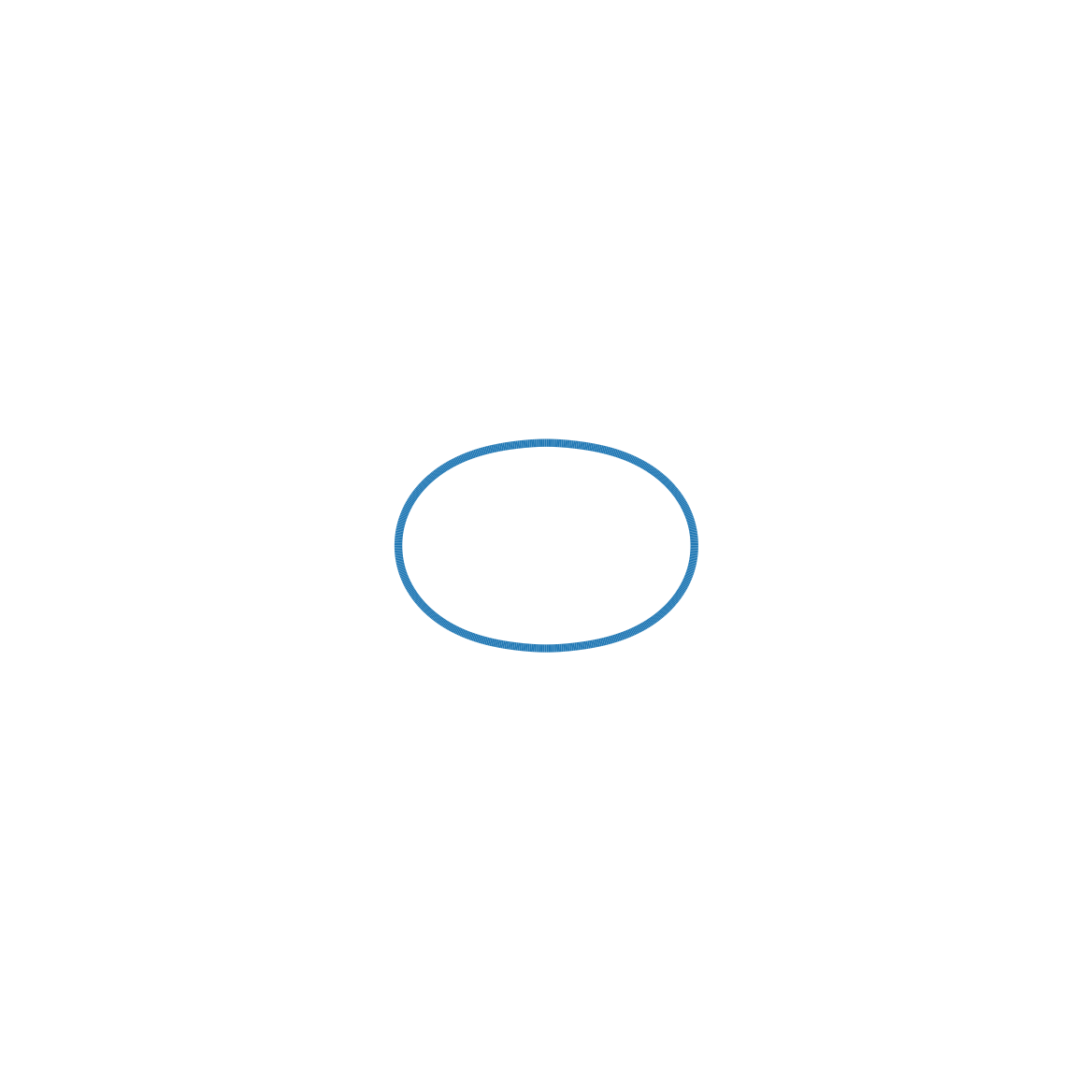}
\includegraphics[width=0.25\textwidth]{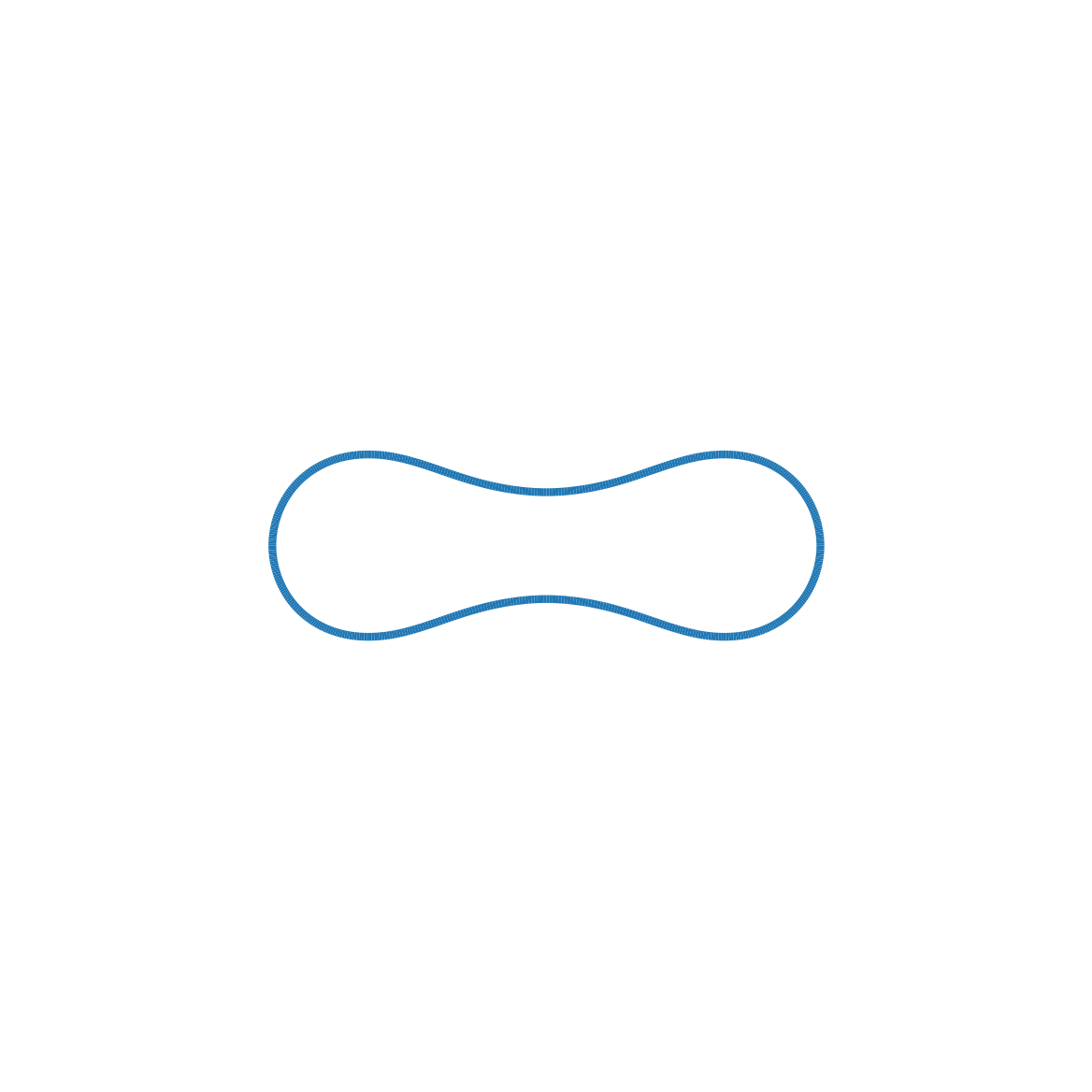}
\includegraphics[width=0.25\textwidth]{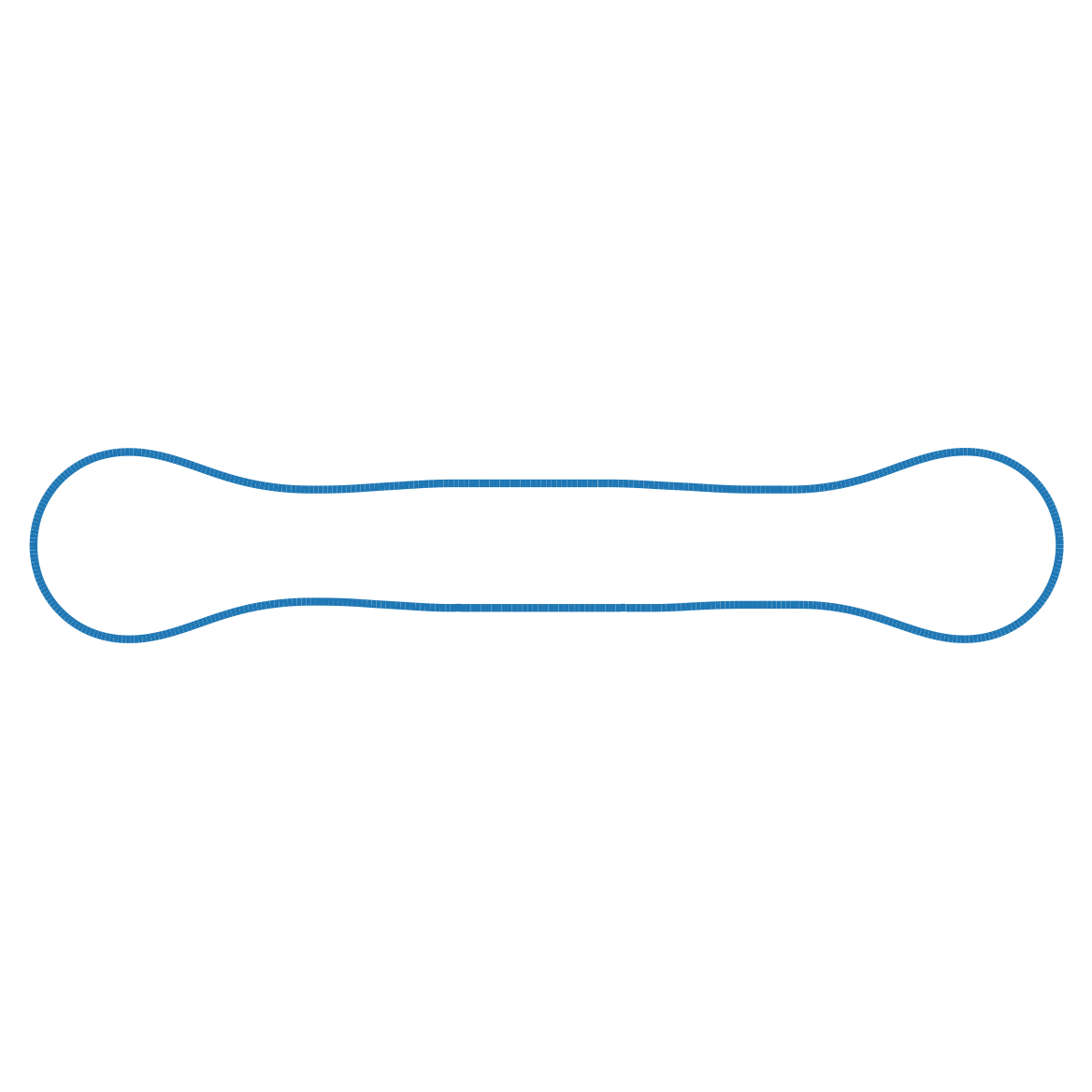}
}
\caption{(Plus-phase model: $\Omega = (0,8)^2$)
Evolution for a perturbed radial interface for $\alpha= 0.01$ and $S_I=0.3$.
We show the evolution at times $t=0,25,50,100$.
}
\label{fig:radial08pert2}
\end{figure}

For a radial interface with a mode-2 perturbation we show a simulation
in Figure~\ref{fig:radial08pert2}. In $\Omega = (0,8)^2$,
we let $\alpha=0.01$ and $S_I = 0.3$. The initial interface has total
dimensions $1.262\times1.212$, 
and so is a small perturbation of the stationary radius $0.6182$.
During the evolution we can observe the expected growth of
the $\ell=2$ perturbation. However, while the shape stretches significantly, 
no pinch-off takes place. 

\begin{figure}
\center
\mbox{
\hspace*{-1cm}
\includegraphics[angle=-0,width=0.2\textwidth]{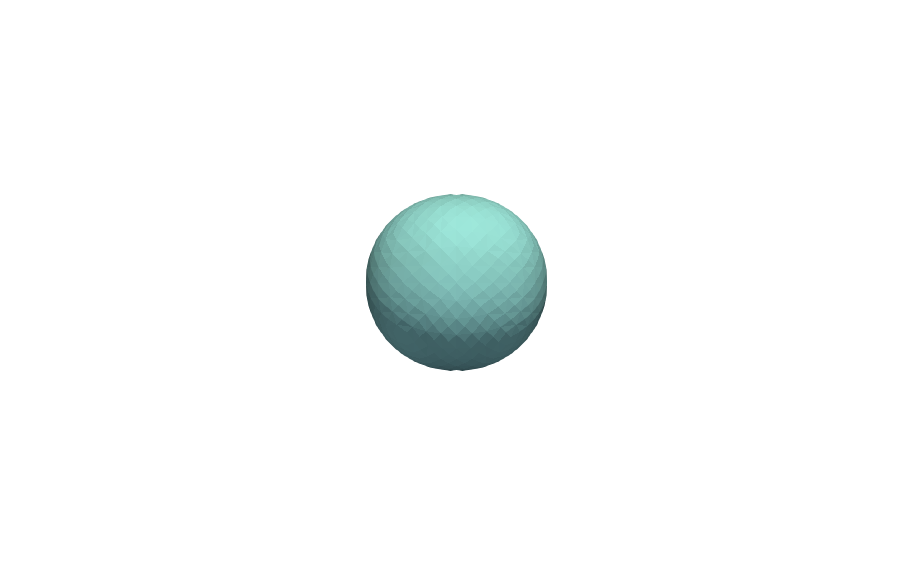}
\includegraphics[angle=-0,width=0.2\textwidth]{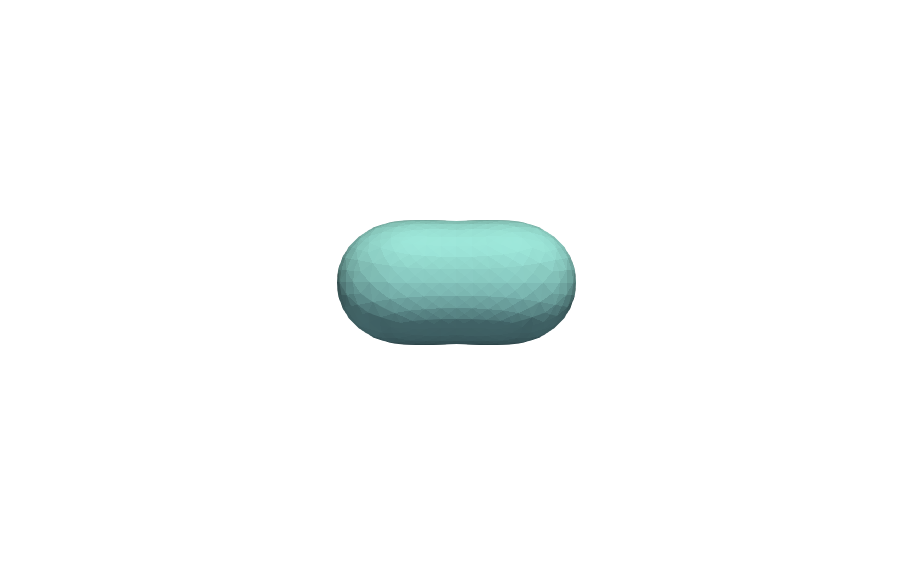}
\includegraphics[angle=-0,width=0.2\textwidth]{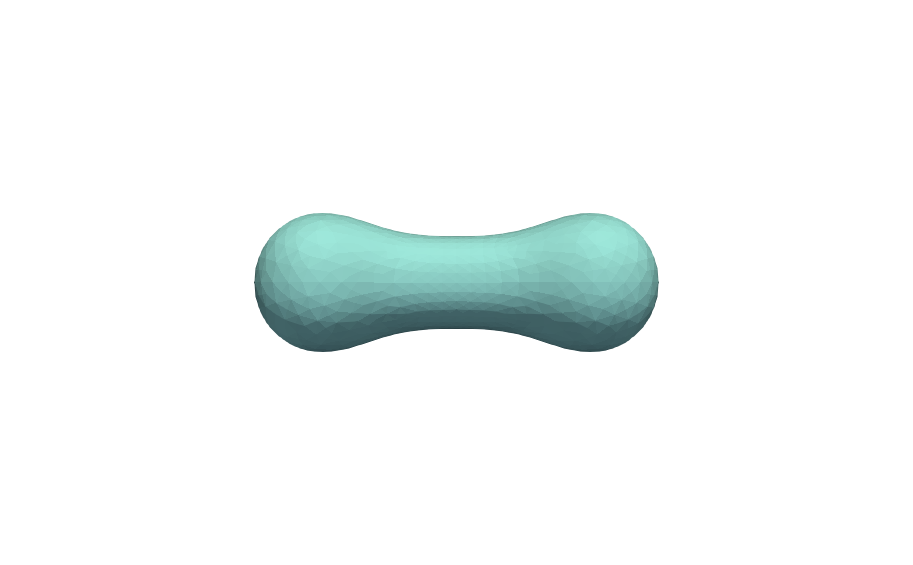}
\includegraphics[angle=-0,width=0.2\textwidth]{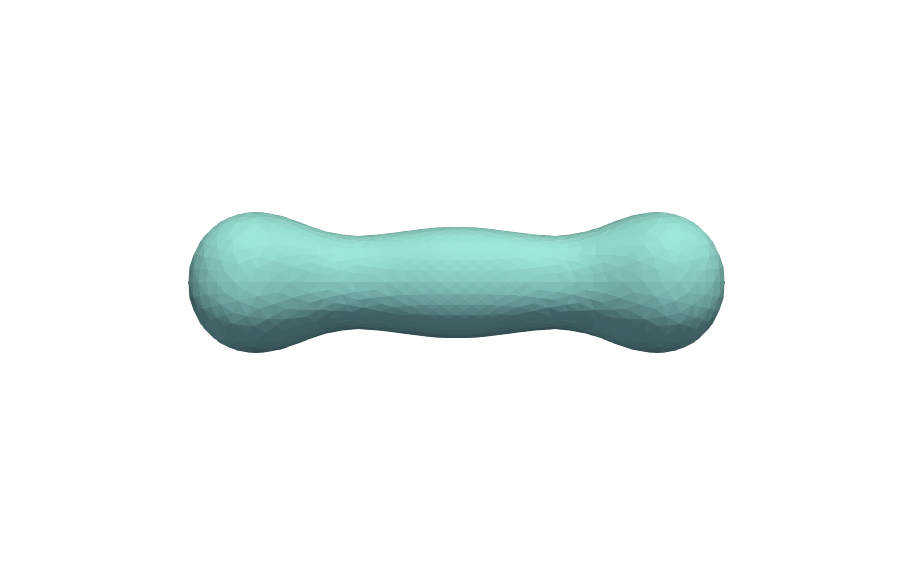}
\includegraphics[angle=-0,width=0.2\textwidth]{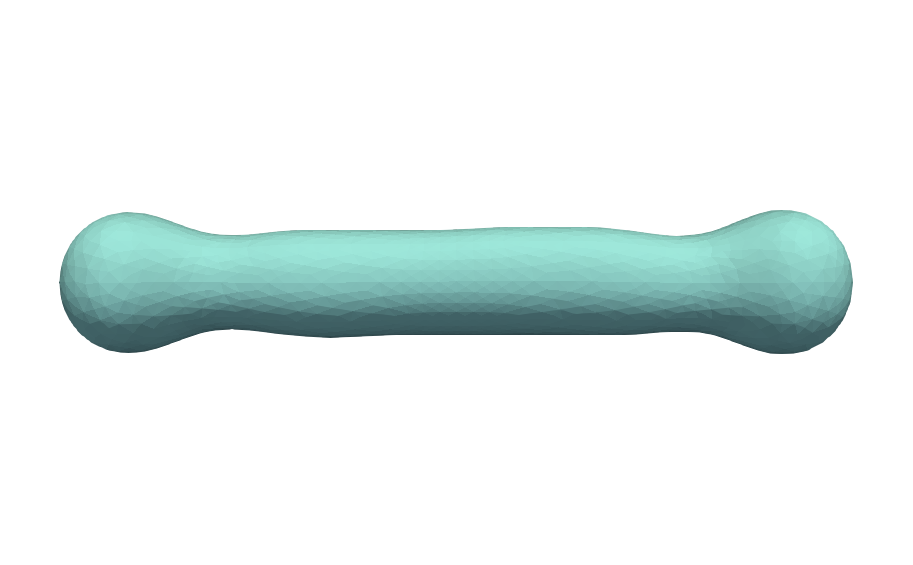}
}
\caption{(Plus-phase model: $\Omega = (0,8)^3$) 
Evolution for a perturbed radial interface for $\alpha=0.01$ and $S_I=0.25$.
We show the evolution at times $t=0, 50, 100, 150, 250$.
}
\label{fig:3dradpert2}
\end{figure}

\begin{figure}
\center
\mbox{
\hspace*{-1cm}
\includegraphics[angle=-0,width=0.2\textwidth]{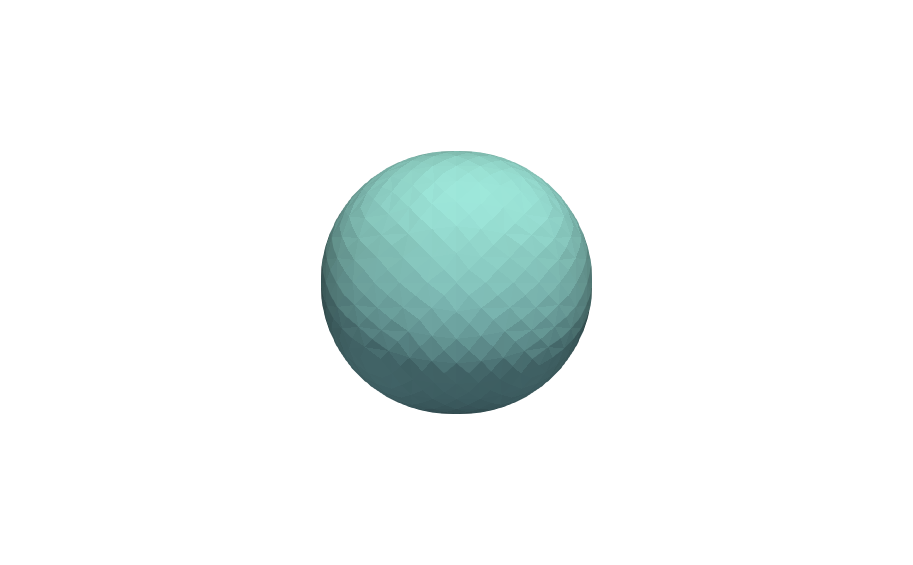}
\includegraphics[angle=-0,width=0.2\textwidth]{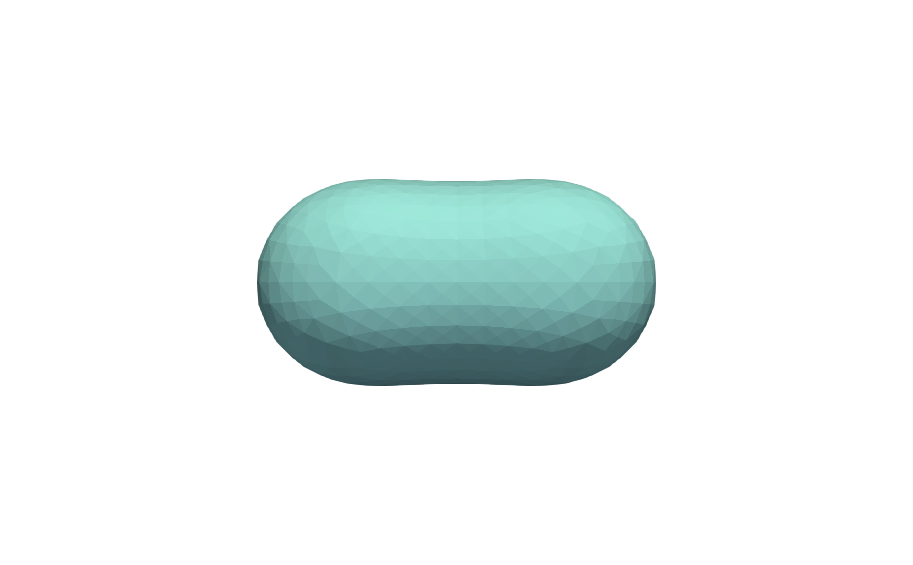}
\includegraphics[angle=-0,width=0.2\textwidth]{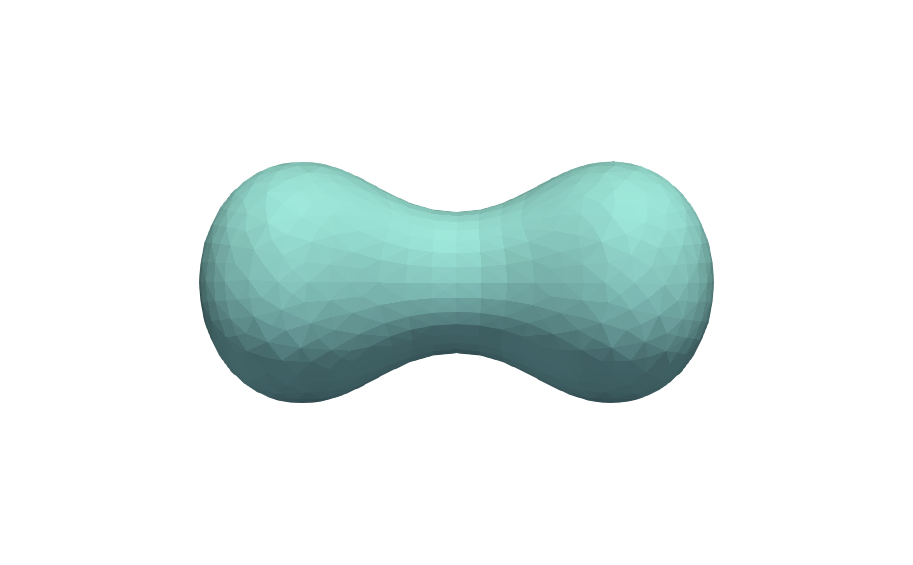}
\includegraphics[angle=-0,width=0.2\textwidth]{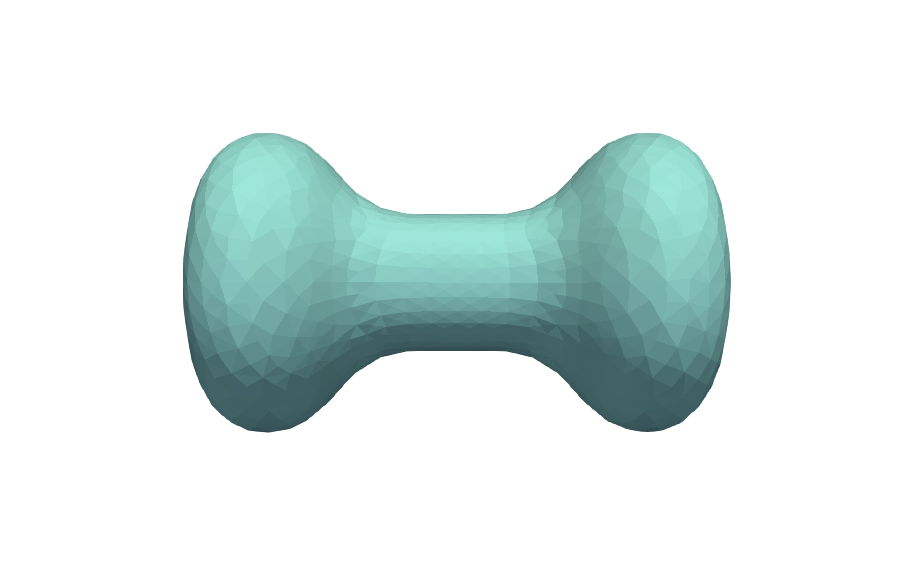}
\includegraphics[angle=-0,width=0.2\textwidth]{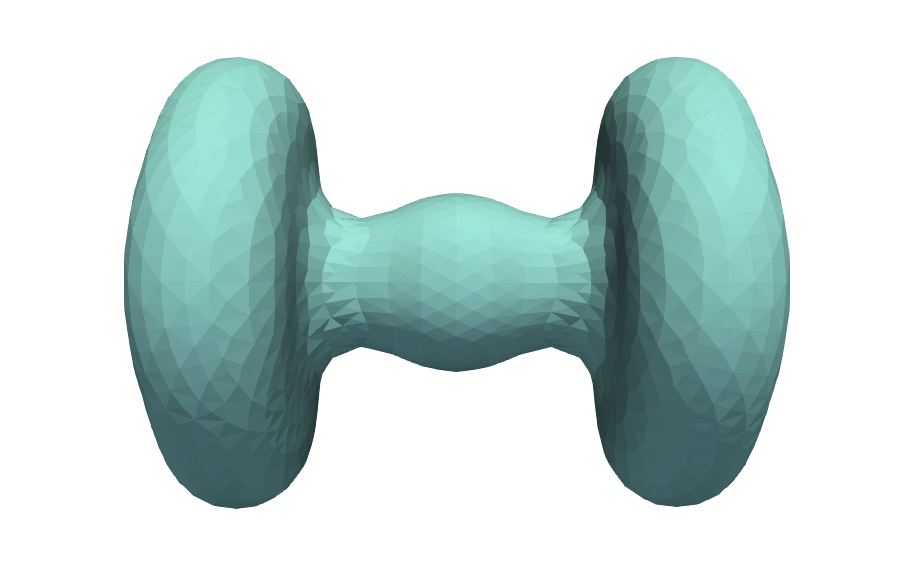}
}
\caption{(Plus-phase model: $\Omega = (0,8)^3$)
Evolution for a perturbed radial interface for $\alpha=0.009$ and $S_I=0.275$.
We show the evolution at times $t=0, 40, 60, 80, 100$.
}
\label{fig:3dradpert3}
\end{figure}

In $\Omega = (0,8)^3$,
we take $\alpha = 0.01$ and $S_I = 0.25$. Then the equilibrium radius is
$0.7583$. We expect small $\ell=2$ perturbations to grow. 
Using the initial data as a cigar-like shape of total dimension 
$1.75\times 1.691 \times 1.691$, shows the expected unstable growth, which
leads to an elongate shape, see Figure~\ref{fig:3dradpert2}. Once again, no
pinch-off occurs.
Repeating the same experiment, but now choosing $\alpha = 0.009$ and 
$S_I = 0.275$, once again leads to unstable growth without pinch-off in Figure~\ref{fig:3dradpert3}. 

\subsection{Numerical computations for the minus-phase model}

\begin{figure}
\center
\includegraphics[angle=-0,width=0.45\textwidth]{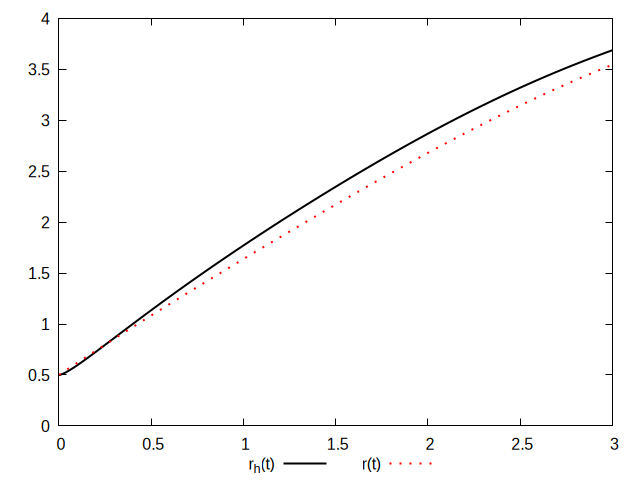}
\includegraphics[angle=-0,width=0.45\textwidth]{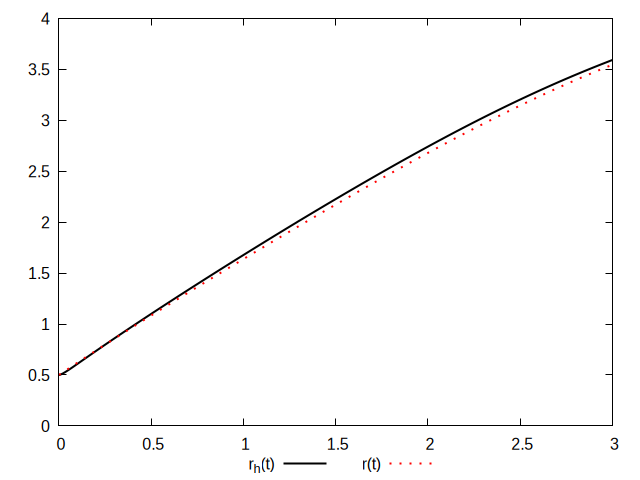}
\caption{(minus-phase model: $\Omega = \bB_4^2(0)$, $\eps=\frac1{16\pi}, \frac1{32\pi}$) 
Comparison of $r(t)$ and $r_h(t)$ for $\alpha=0.1$ and $S_I=1$.
}
\label{fig:outercurved16pi}
\end{figure}

A comparison similar to Figure~\ref{fig:innercurved16pi} for the phase field
model is shown
in Figure~\ref{fig:outercurved16pi}, once again choosing 
$\alpha=0.1$ and $S_I=1$.

\begin{figure}
\center
\includegraphics[angle=-0,width=0.24\textwidth]{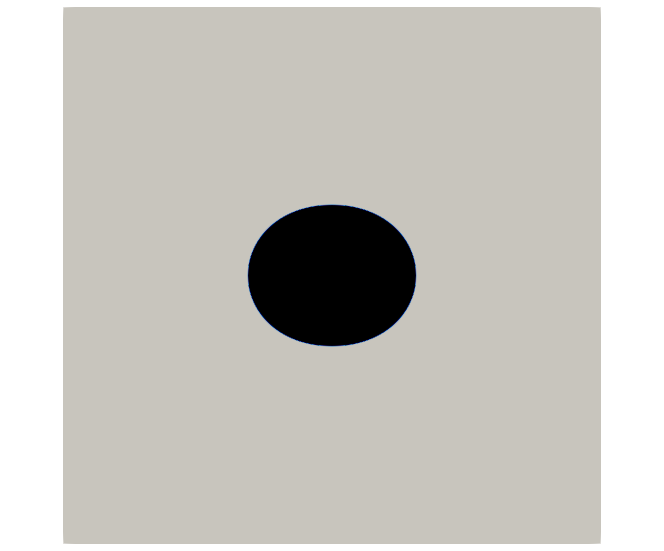}
\includegraphics[angle=-0,width=0.24\textwidth]{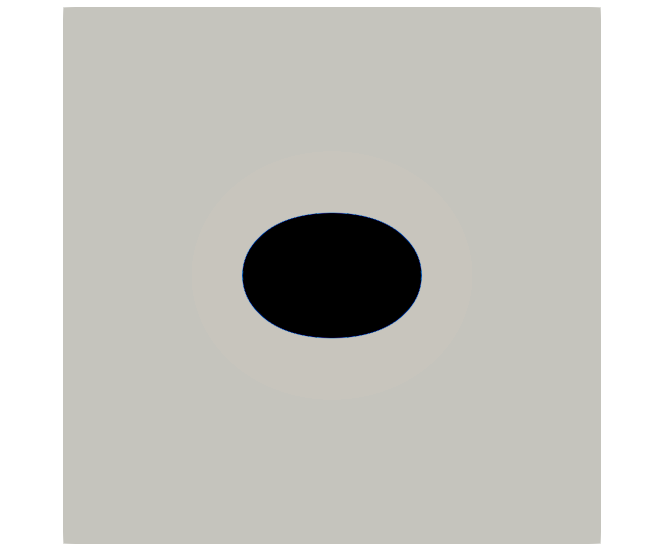}
\includegraphics[angle=-0,width=0.24\textwidth]{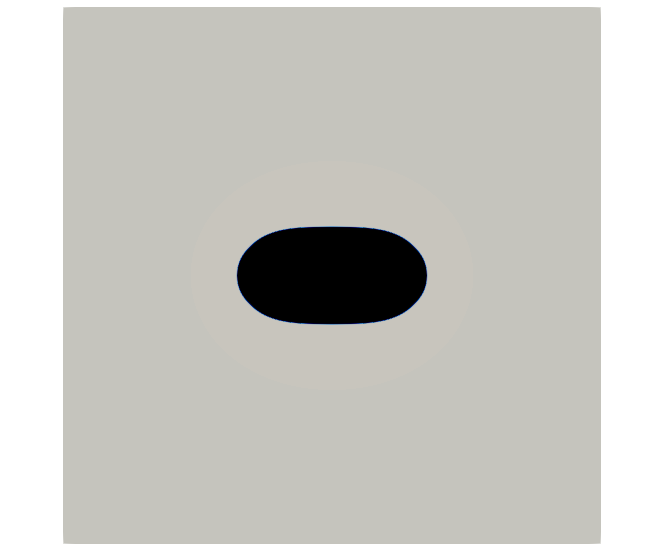}
\includegraphics[angle=-0,width=0.24\textwidth]{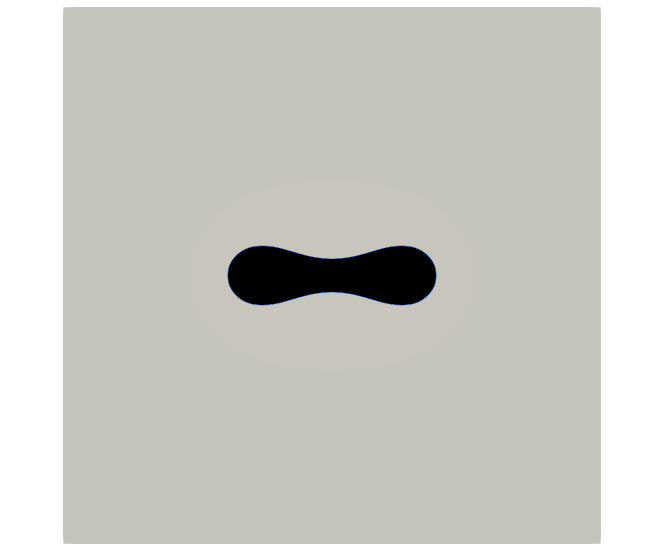}
\includegraphics[angle=-0,width=0.24\textwidth]{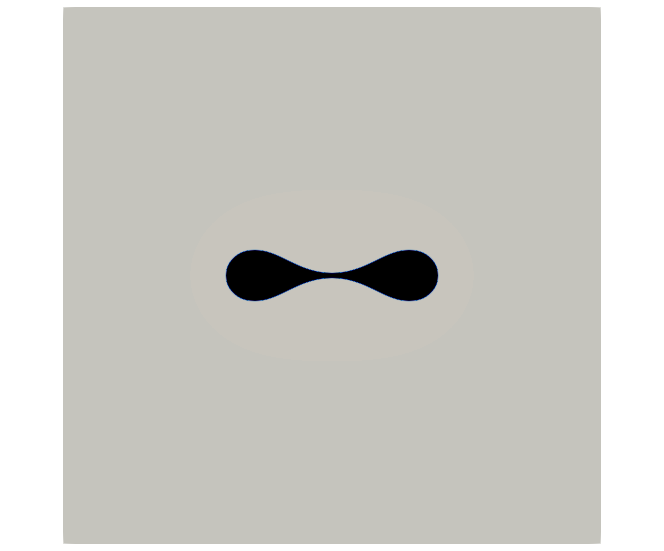}
\includegraphics[angle=-0,width=0.24\textwidth]{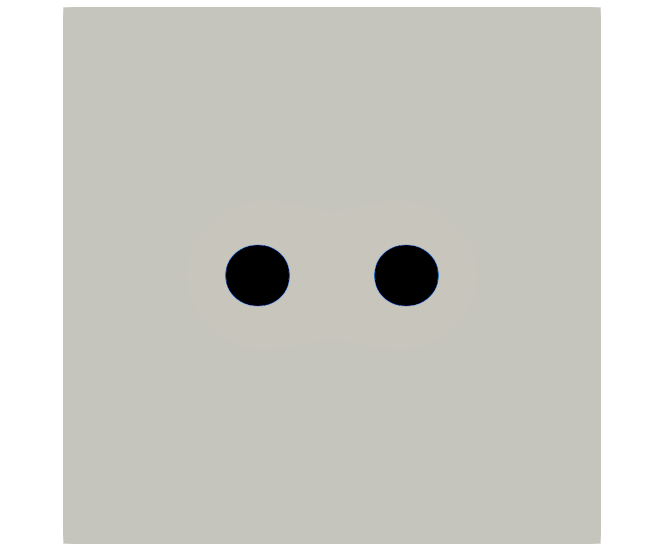}
\includegraphics[angle=-0,width=0.24\textwidth]{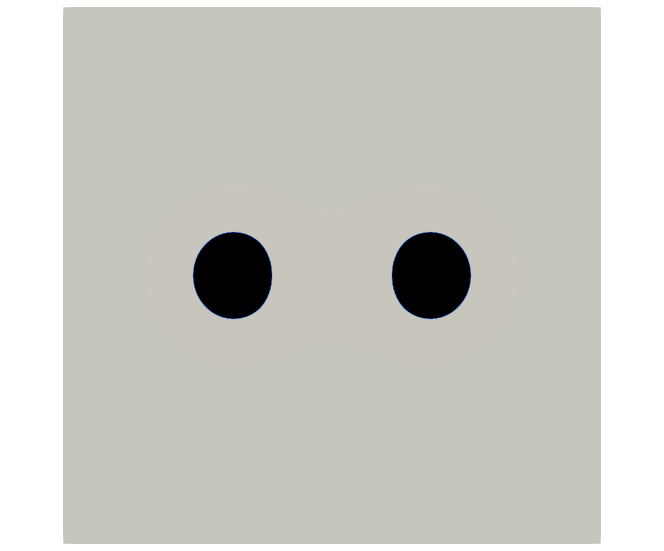}
\includegraphics[angle=-0,width=0.24\textwidth]{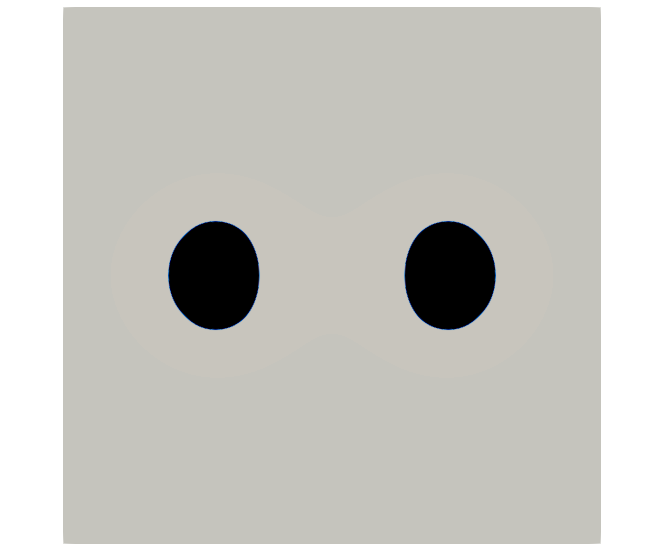}
\caption{(Minus-phase model: $\Omega = (0,8)^2$, $\eps=\frac1{64\pi}$) 
Pinch-off for $\alpha=0.1$ and $S_I=-1.25$.
We show the evolution at times $t=0$, 10, 20, 30, 32, 33, 40, 50.
}
\label{fig:2d64pi_a01_SI-125}
\end{figure}

In Figure~\ref{fig:2d64pi_a01_SI-125} we show a phase field simulation 
starting from an initial configuration with a small $\ell=2$ perturbation. In particular, 
the zero level set of the initial data has total dimension $2.5\times2.1$.
On choosing $\eps$ sufficiently small, here $\eps = \frac1{64\pi}$,
we do observe the expected pinch-off behaviour.

\begin{figure}
\center
\includegraphics[angle=-0,width=0.24\textwidth]{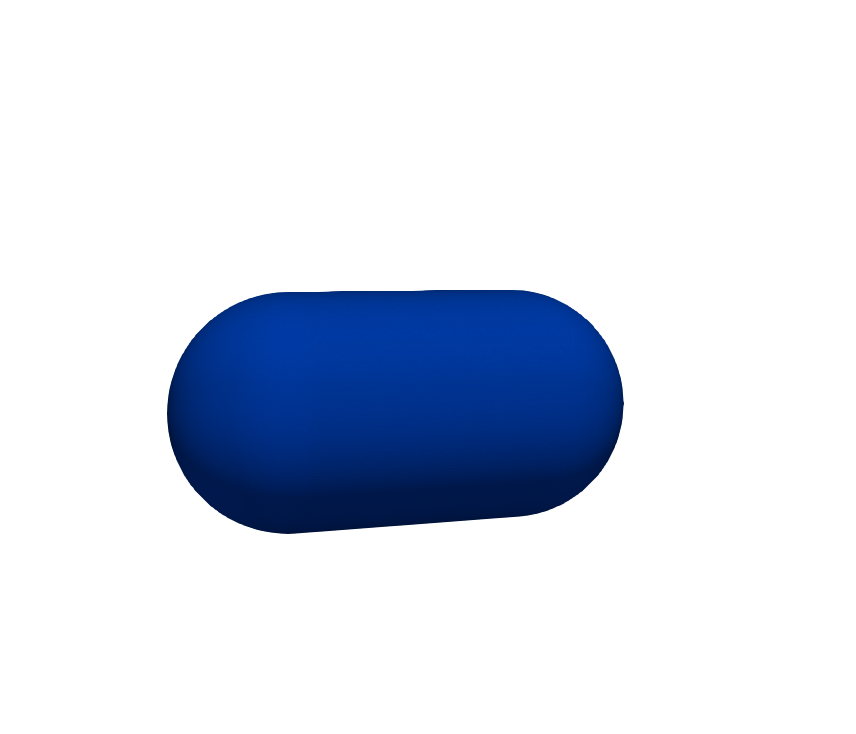}
\includegraphics[angle=-0,width=0.24\textwidth]{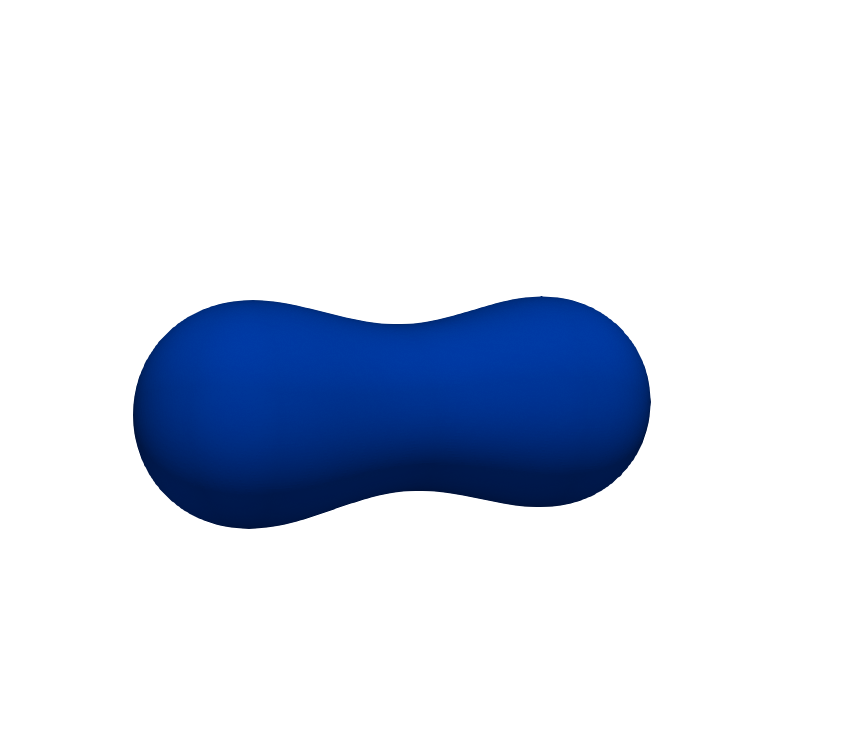}
\includegraphics[angle=-0,width=0.24\textwidth]{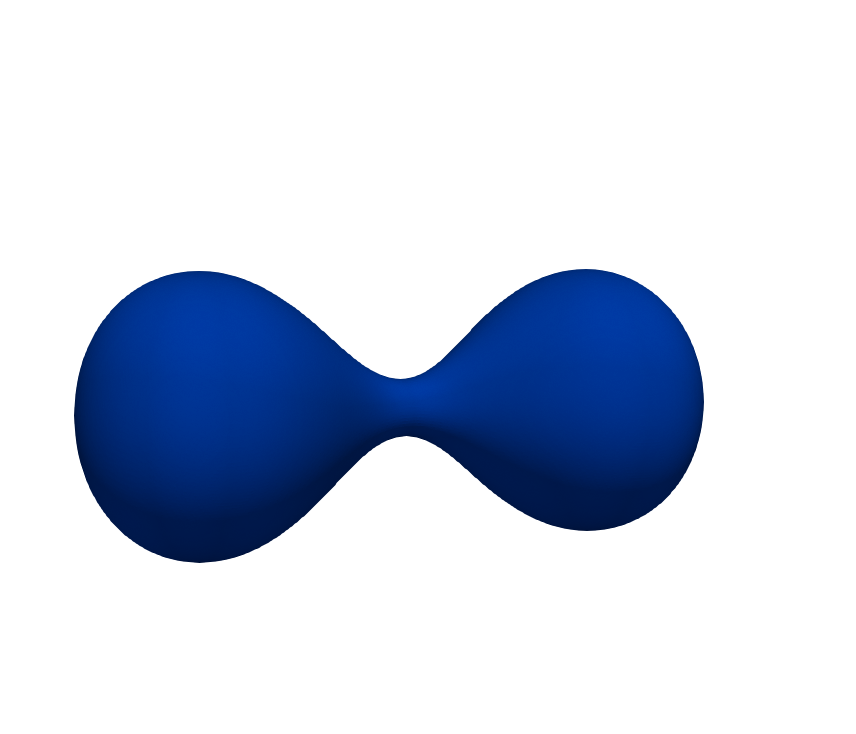}
\includegraphics[angle=-0,width=0.24\textwidth]{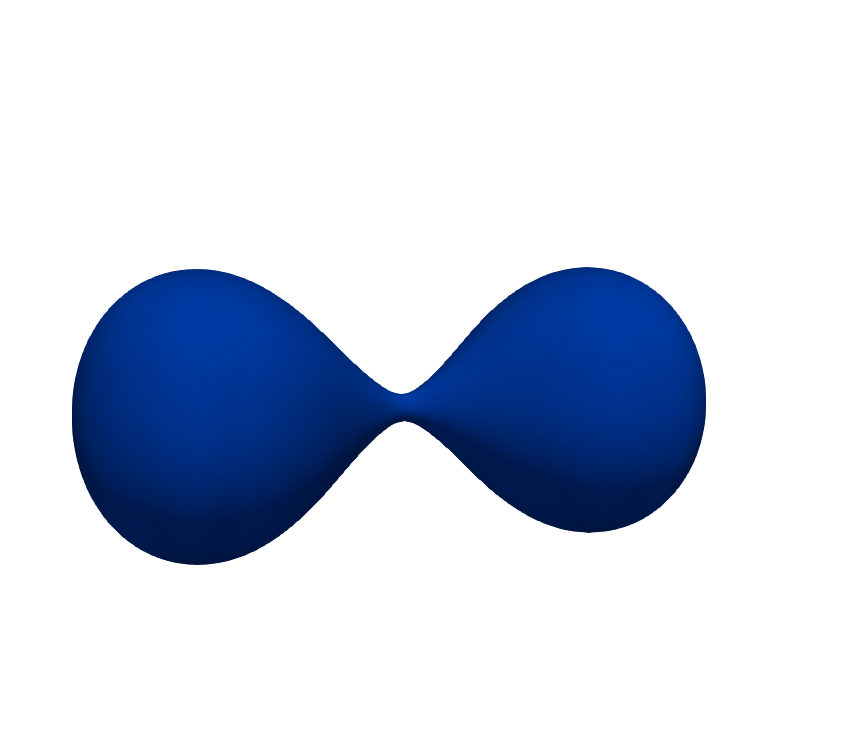}
\includegraphics[angle=-0,width=0.24\textwidth]{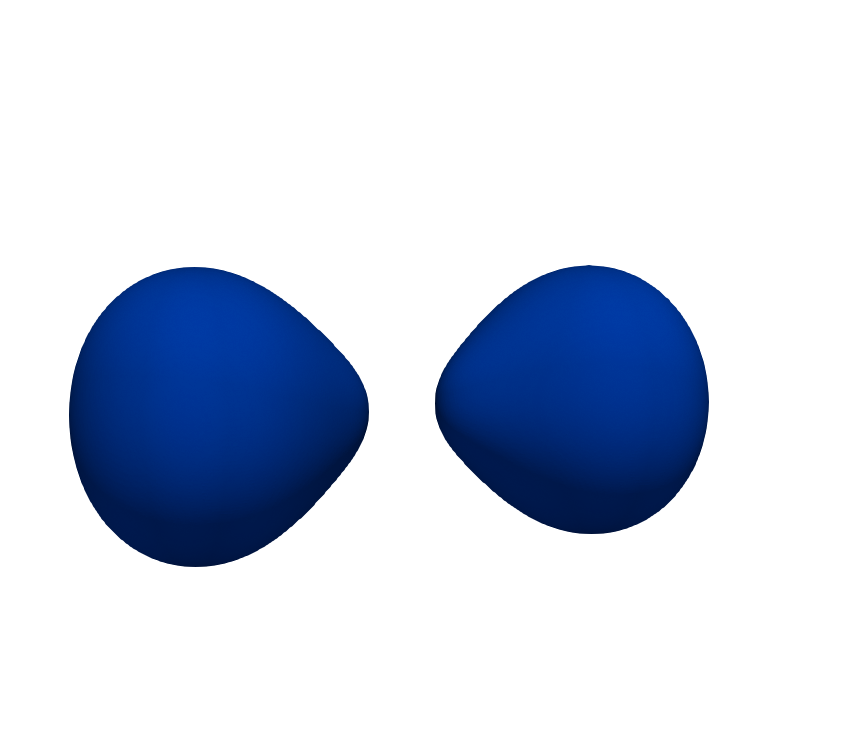}
\includegraphics[angle=-0,width=0.24\textwidth]{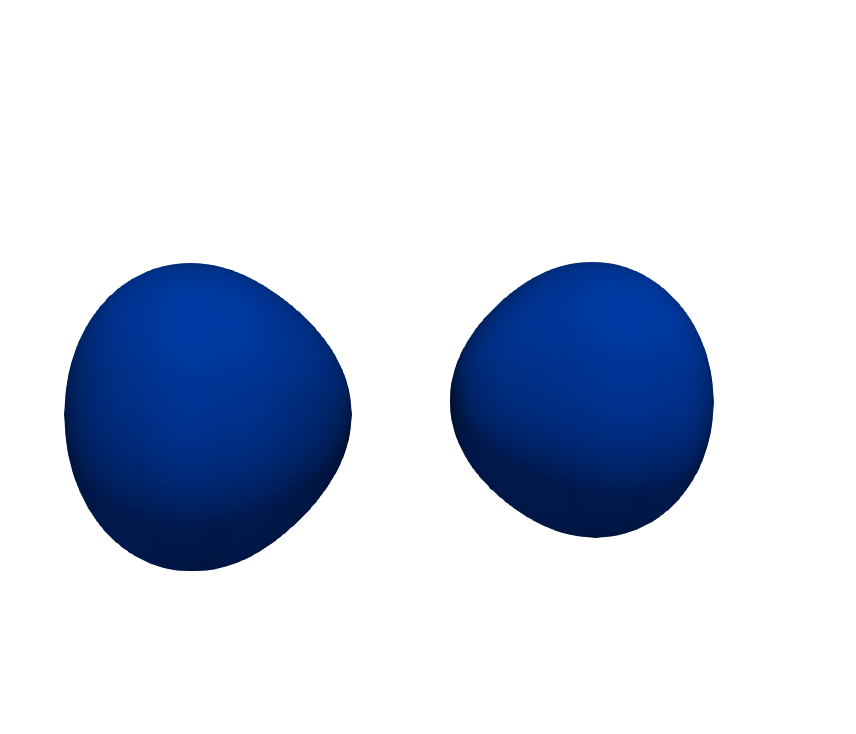}
\includegraphics[angle=-0,width=0.24\textwidth]{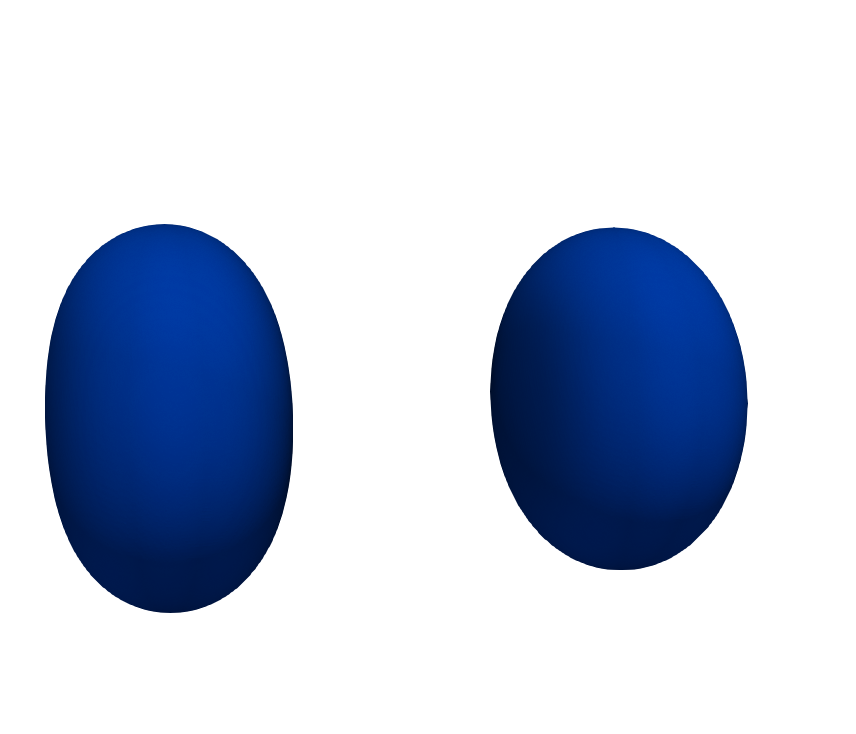}
\includegraphics[angle=-0,width=0.24\textwidth]{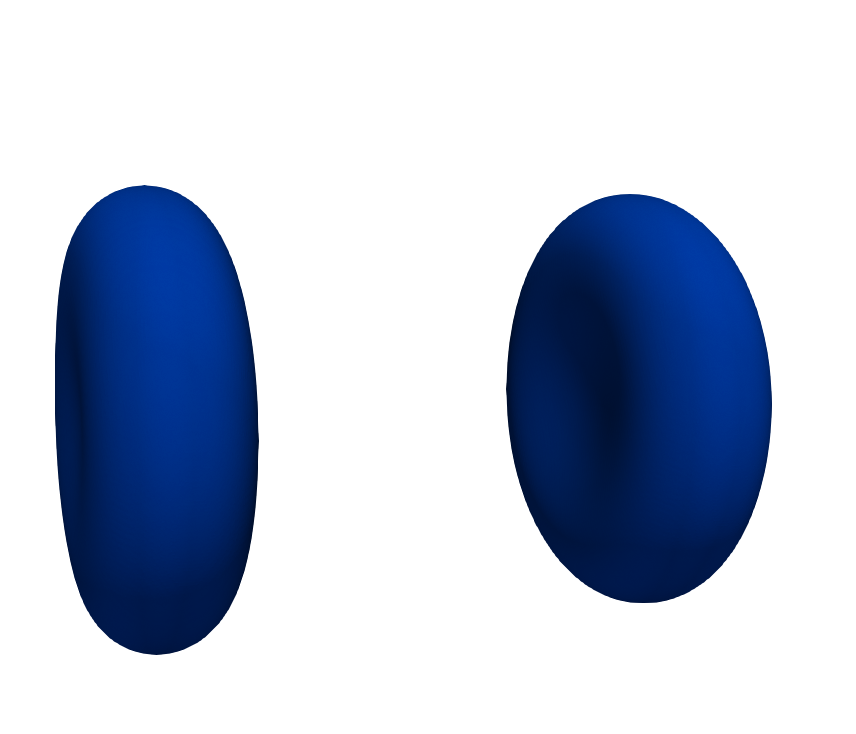}
\includegraphics[angle=-0,width=0.24\textwidth]{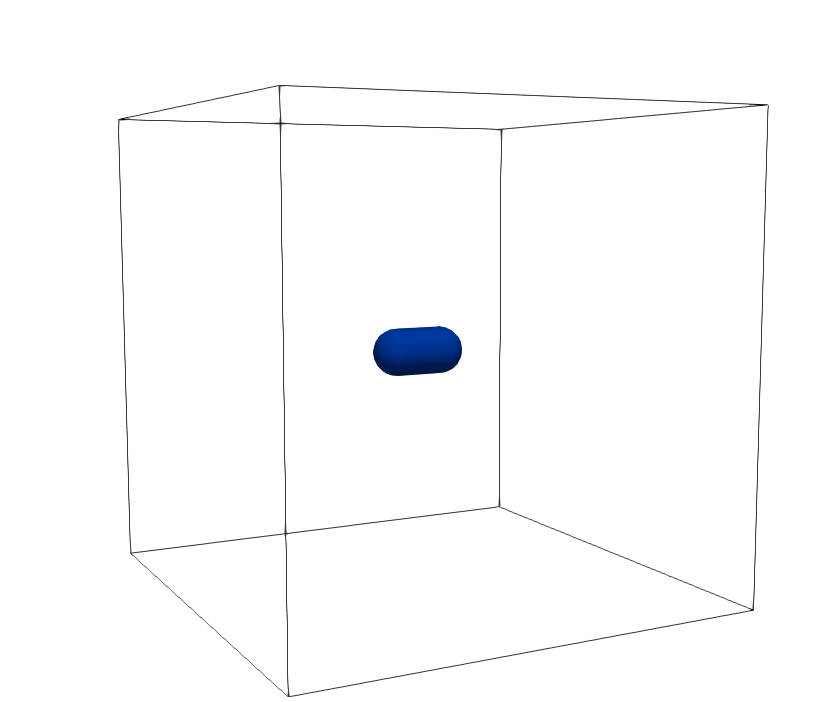}
\includegraphics[angle=-0,width=0.24\textwidth]{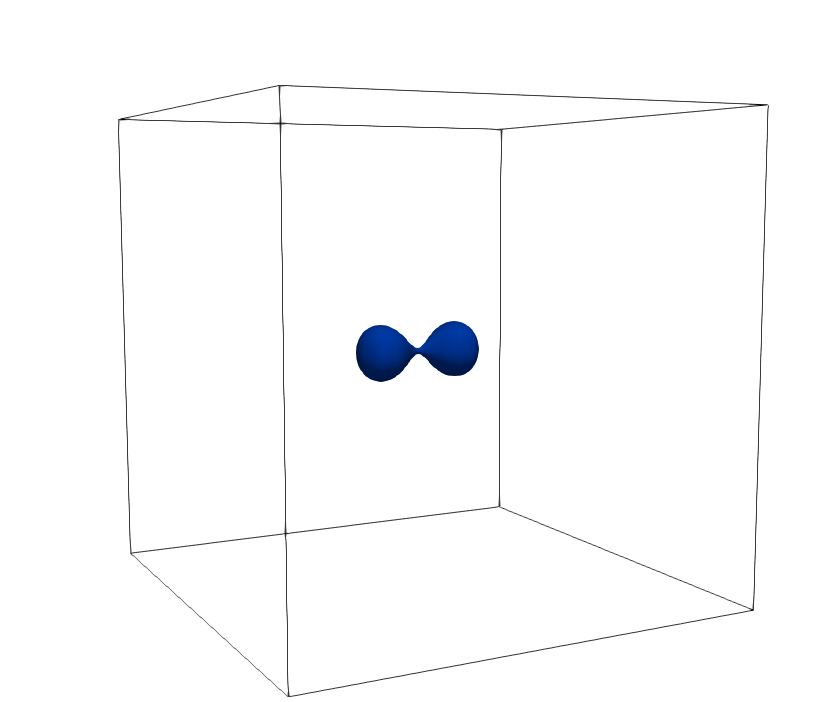}
\includegraphics[angle=-0,width=0.24\textwidth]{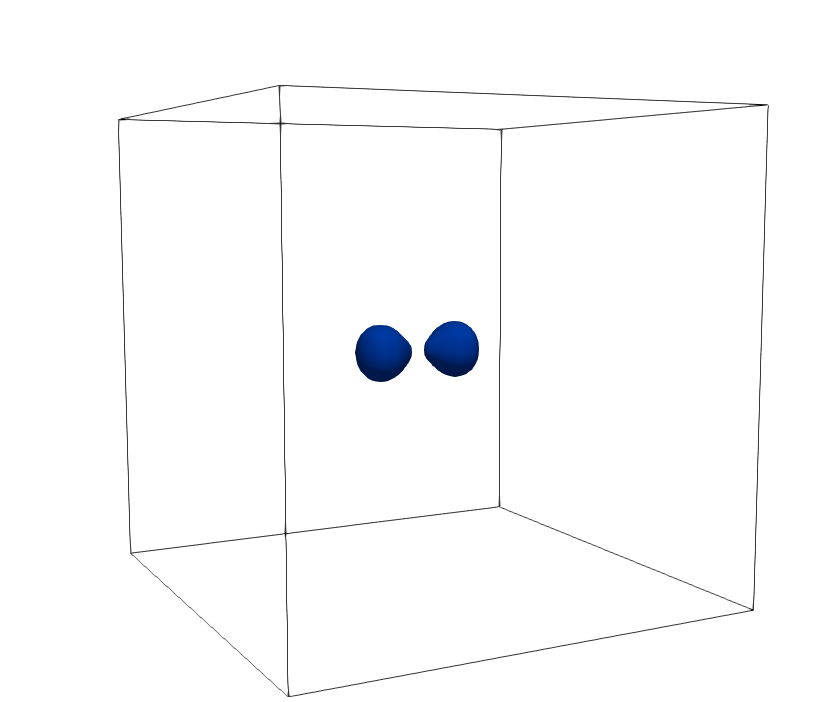}
\includegraphics[angle=-0,width=0.24\textwidth]{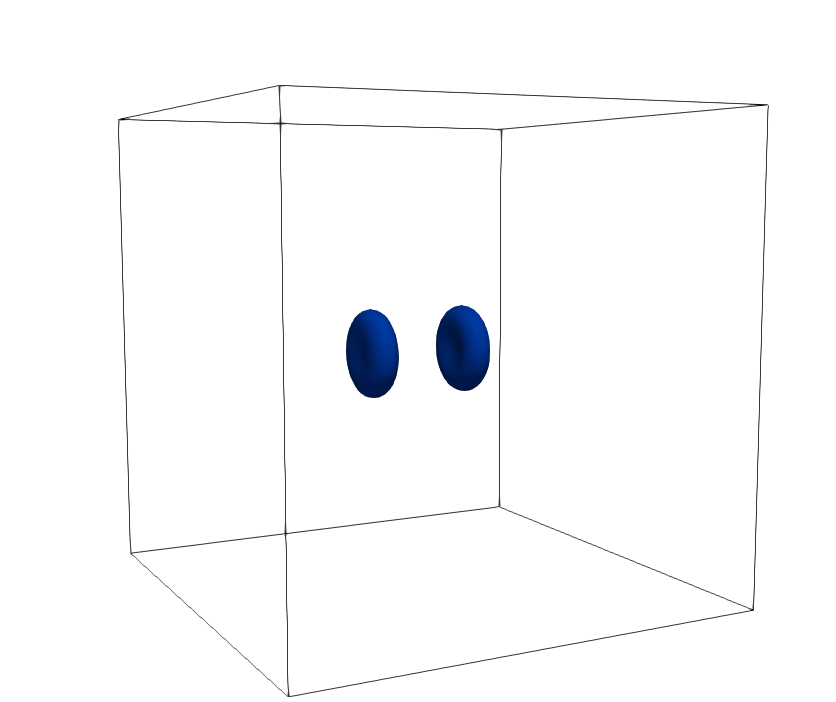}
\caption{(Minus-phase model: $\Omega = (0,8)^3$, $\eps=\frac1{16\pi}$) 
Pinch-off for $\alpha=0.02$ and $S_I=-4$.
We show the evolution at times $t=0$, 1, 1.8, 1.85, 1.9, 2, 3, 4.
The bottom row shows the solution at time $t=0$, 1.85, 1.9, 4 within the domain
$\Omega$.
}
\label{fig:3d16pi_pinchoff}
\end{figure}

Similarly, a simulation for a pinch-off in three dimensions is shown in 
Figure~\ref{fig:3d16pi_pinchoff}. The initial data in the domain
$\Omega=(0,8)^3$ is chosen such that the zero level set 
is a cigar-like shape with total dimension $1.6\times0.8\times0.8$.

\begin{figure}
\center
\includegraphics[angle=-0,width=0.45\textwidth]{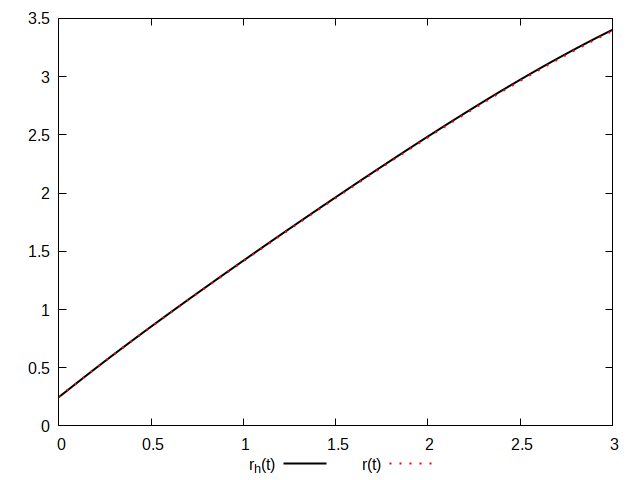}
\includegraphics[angle=-0,width=0.45\textwidth]{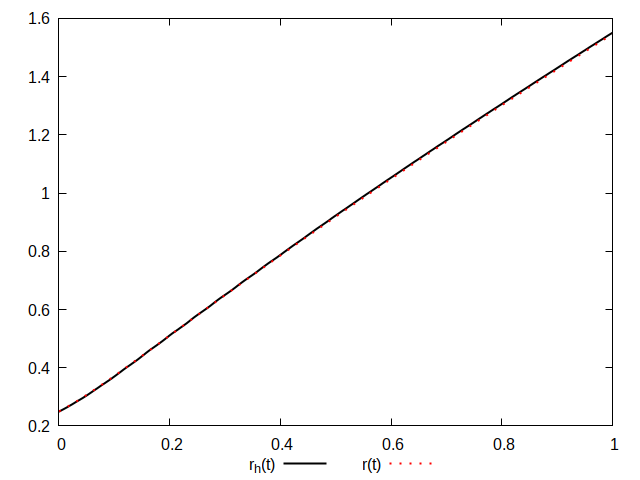}
\caption{(Minus-phase model: $\Omega = \bB^d_4(0)$) 
Comparison of $r(t)$ and $r_h(t)$ for 
$\alpha=0.1$ and $S_I=1$ for $d=2$ (left)
and $d=3$ (right).
}
\label{fig:2f}
\end{figure}

For the sharp interface simulations,
we begin with a single spherical interface inside the domains
$\Omega = \bB^d_4(0)$, $d=\{2,3\}$. For the initial radius we choose $r(0) = 0.25$, 
while the physical parameters are chosen as 
$\alpha=0.1$ and  $S_I=1$.
A comparison between the true radii $r(t_n)$ and the discrete radii $r_h(t_n)$
is shown in Figure~\ref{fig:2f}. We observe that in both two and three dimensions our
numerical algorithm matches the evolution of the true solution obtained from solving the ODE very well.



\begin{figure}
\center
\includegraphics[width=0.3\textwidth]{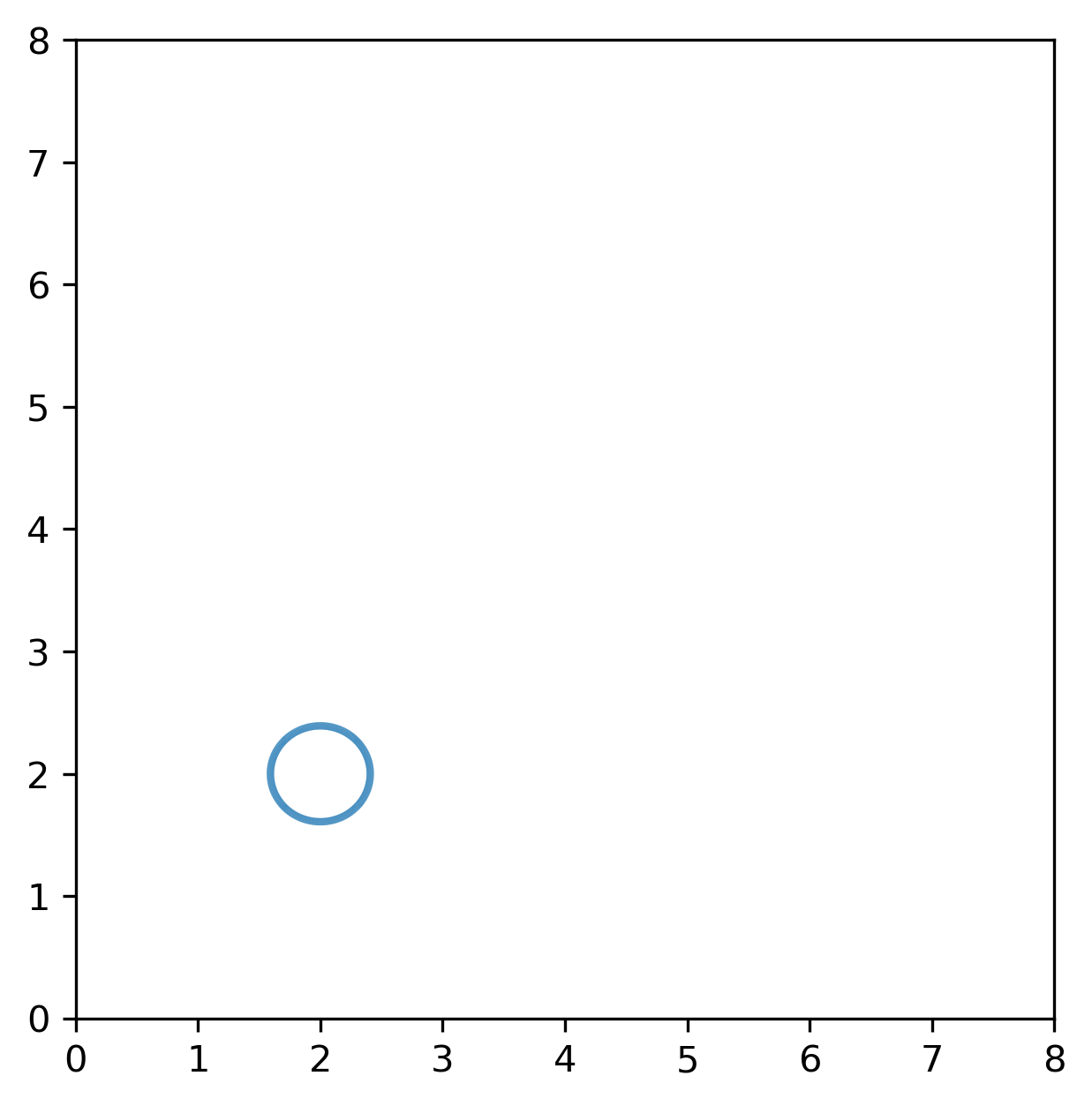}
\includegraphics[width=0.3\textwidth]{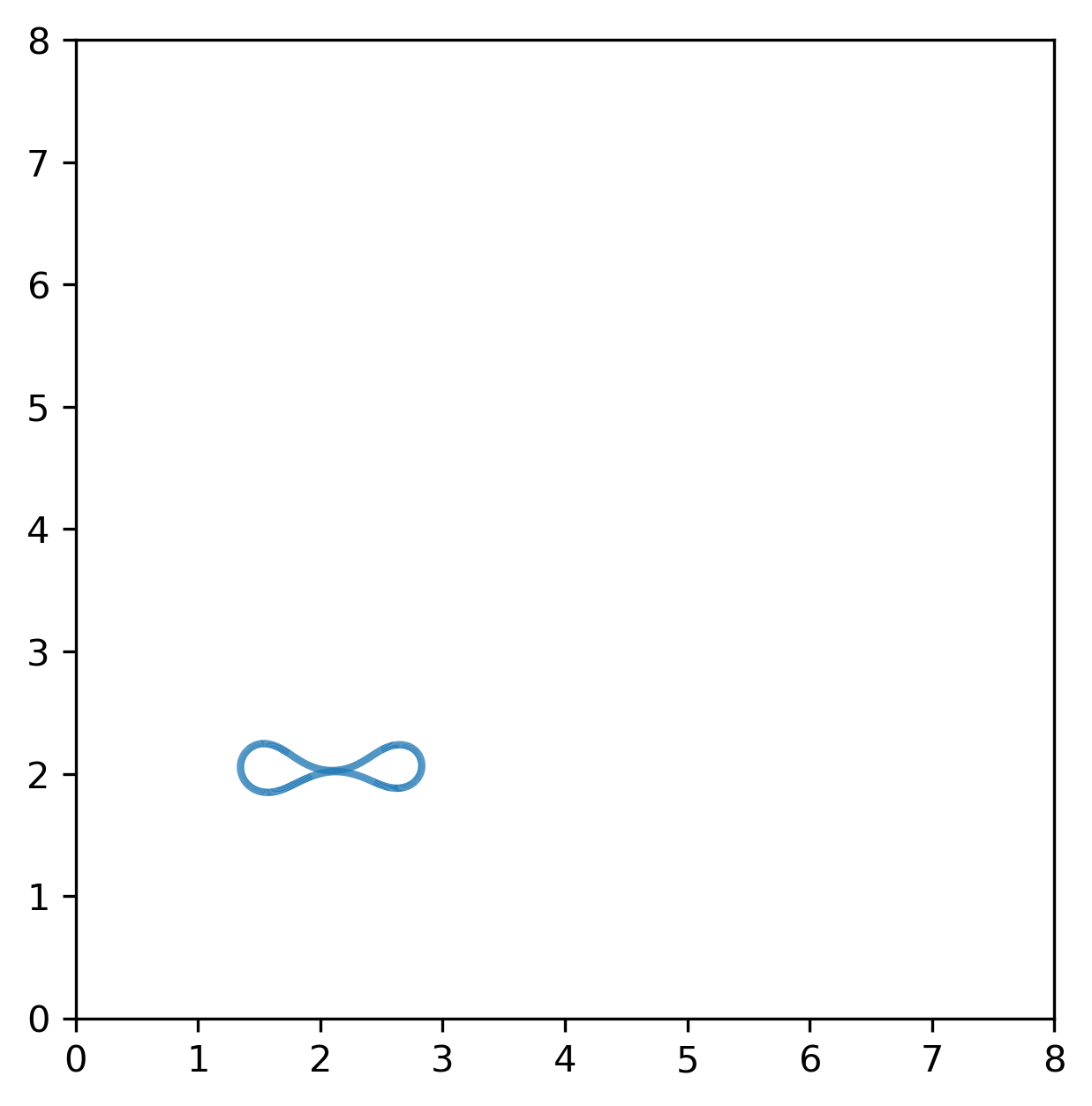}
\includegraphics[width=0.3\textwidth]{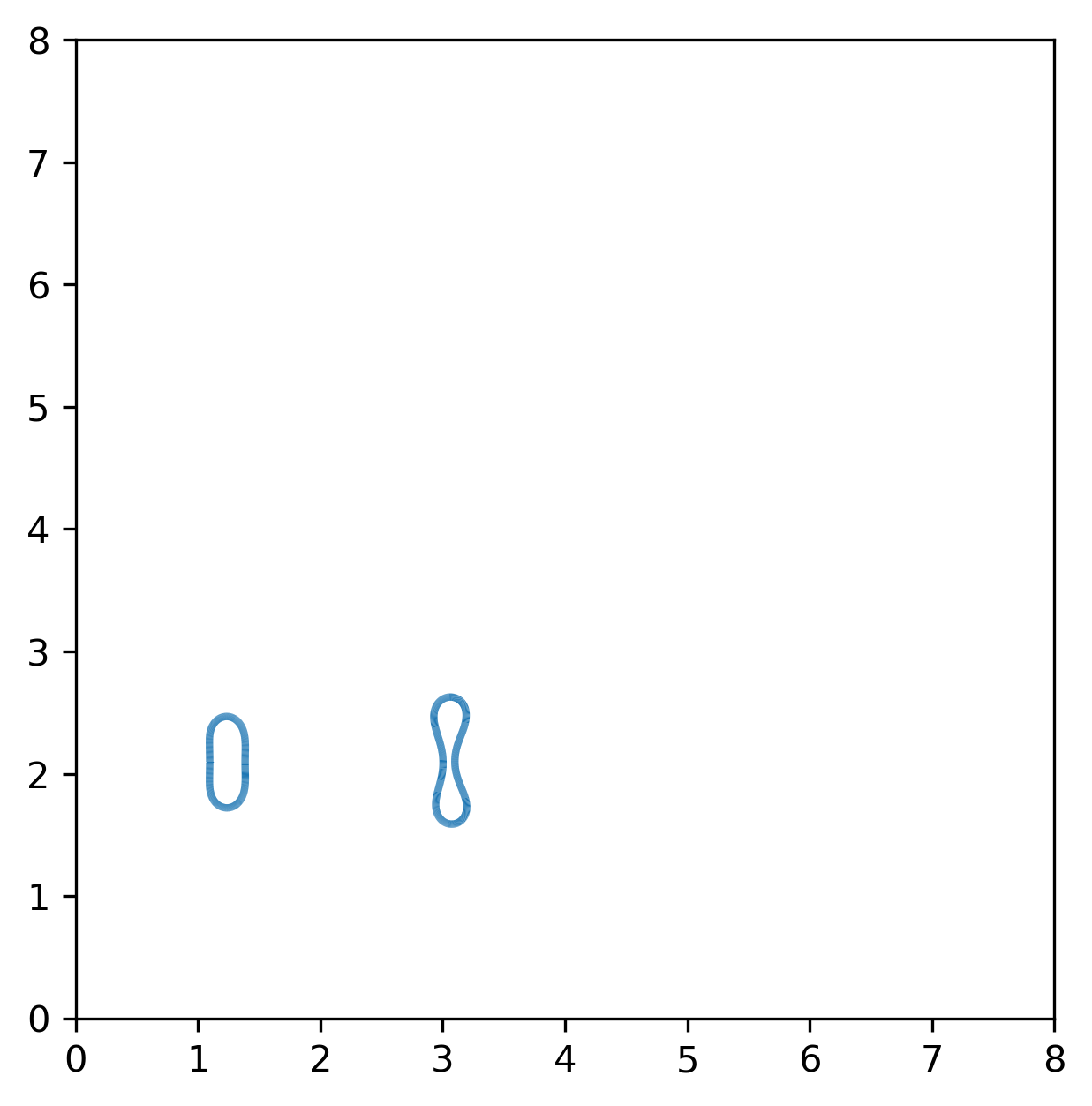}
\includegraphics[width=0.3\textwidth]{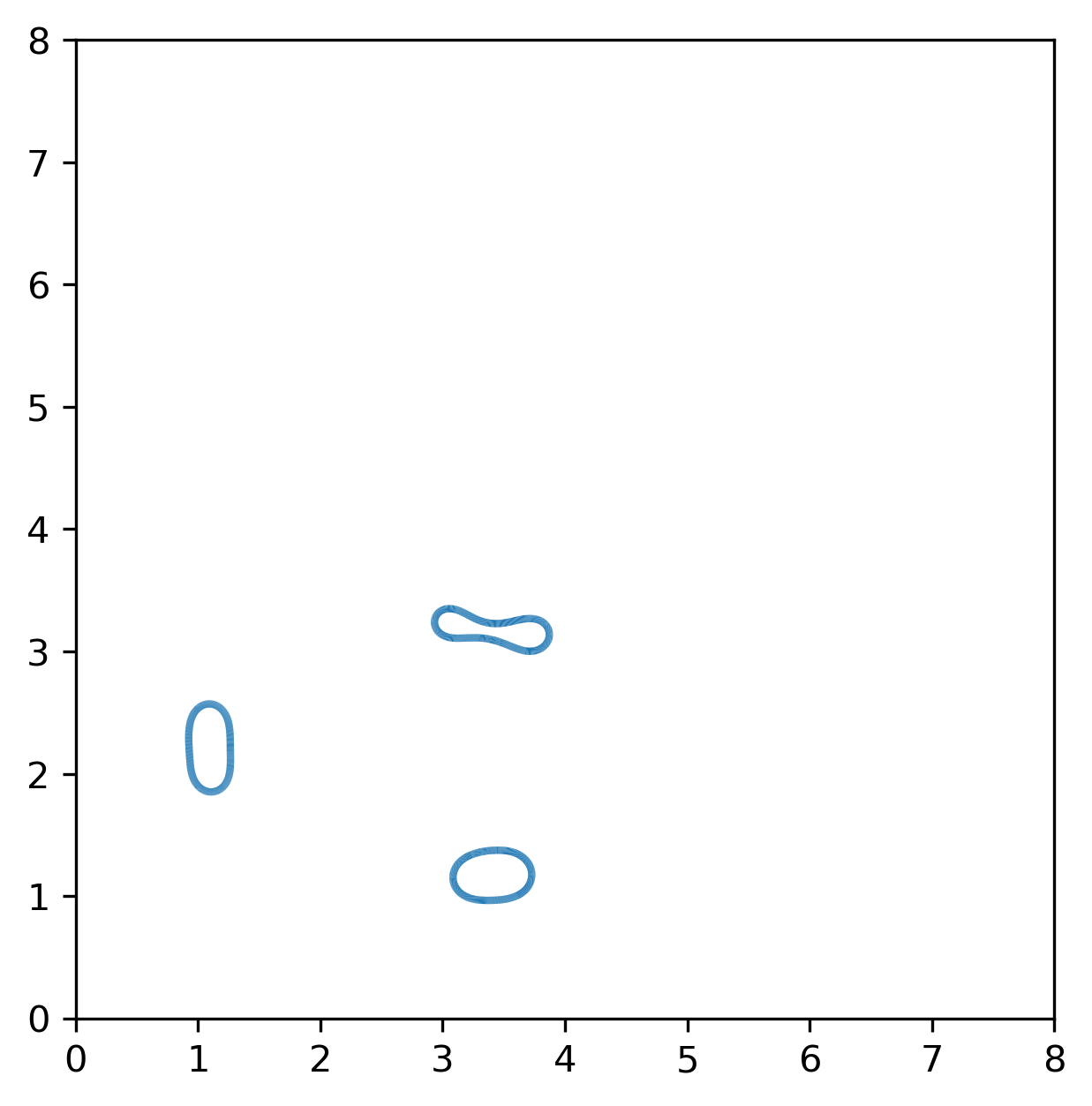}
\includegraphics[width=0.3\textwidth]{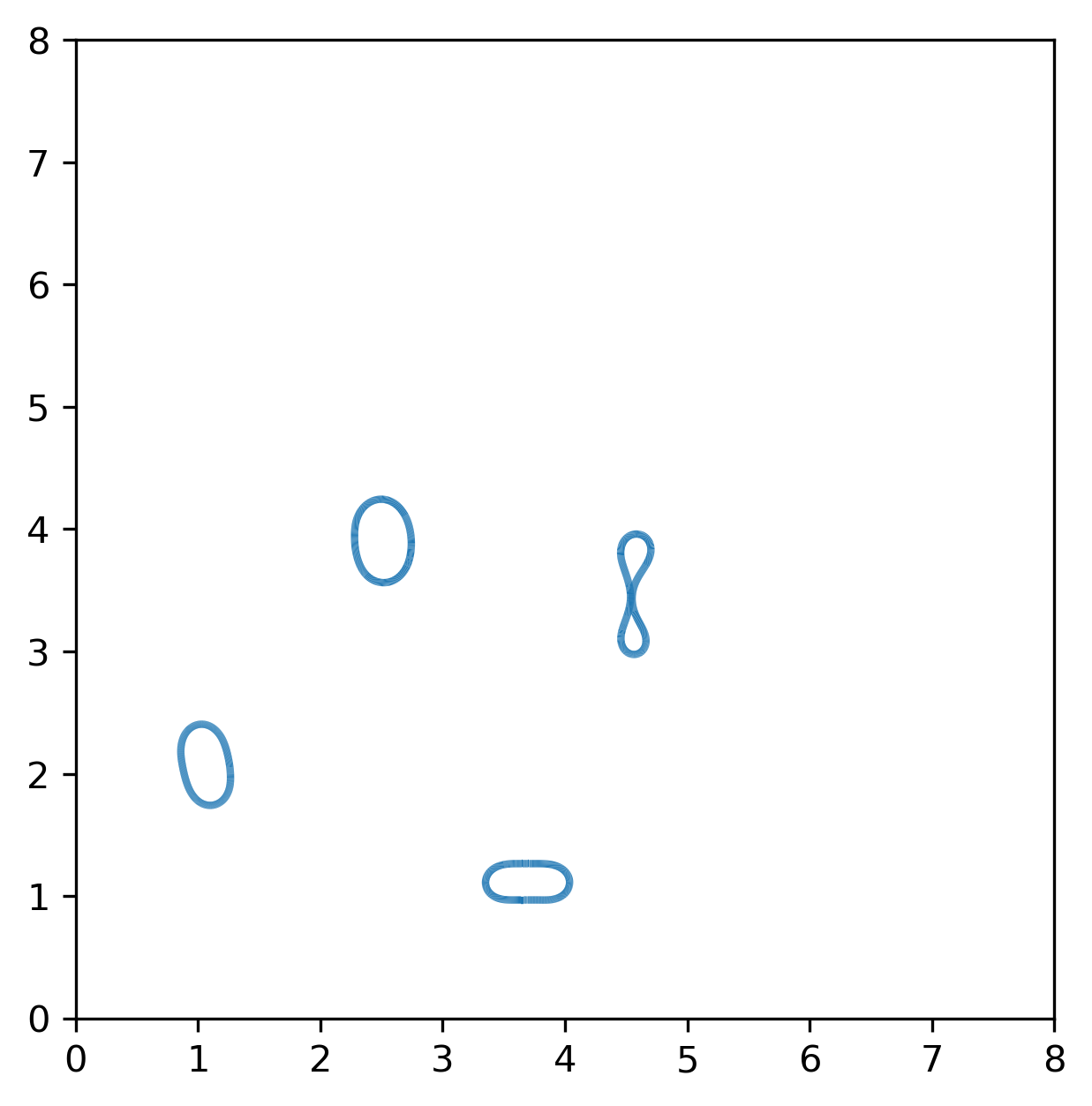}
\includegraphics[width=0.3\textwidth]{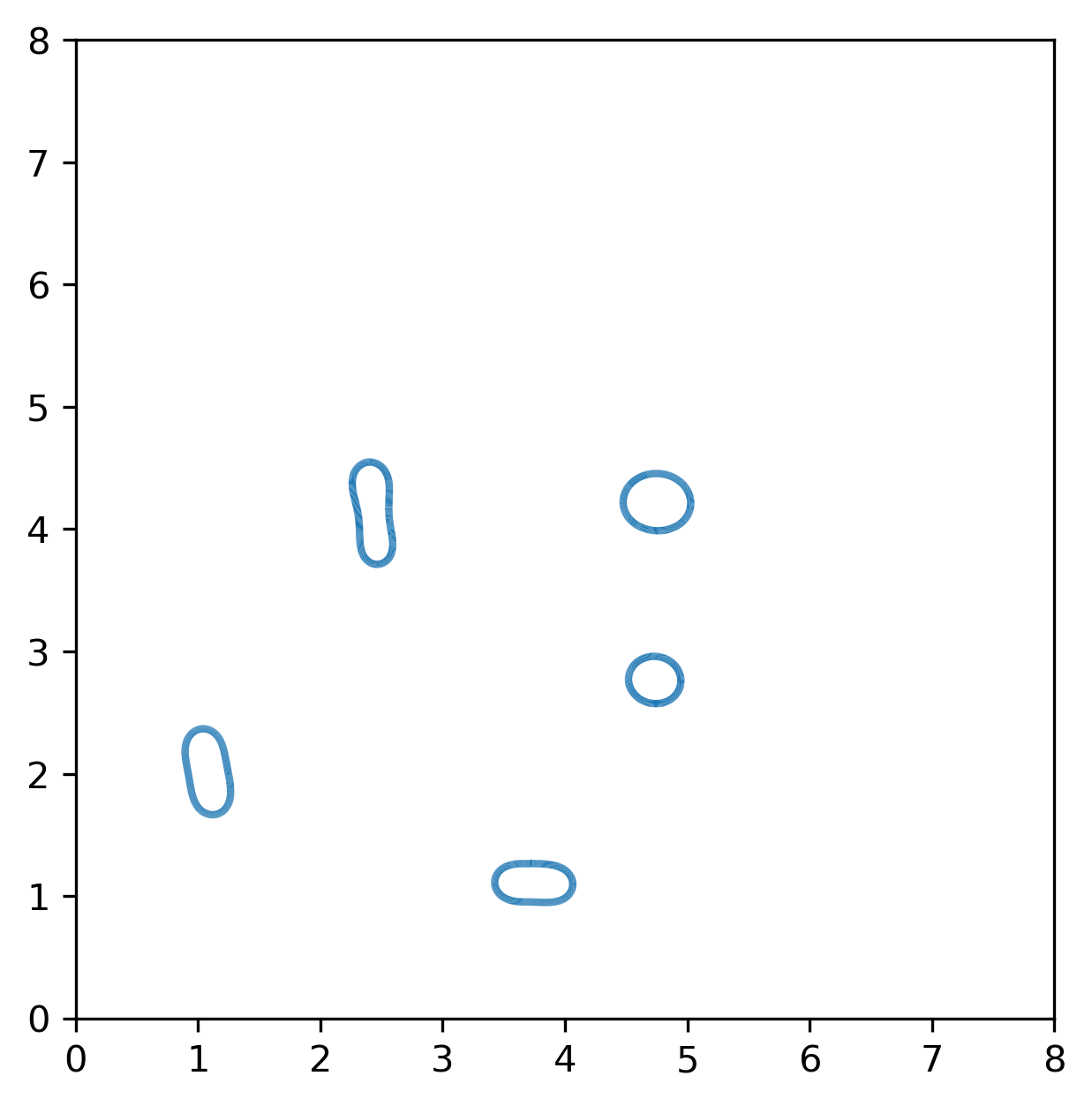}
\includegraphics[width=0.3\textwidth]{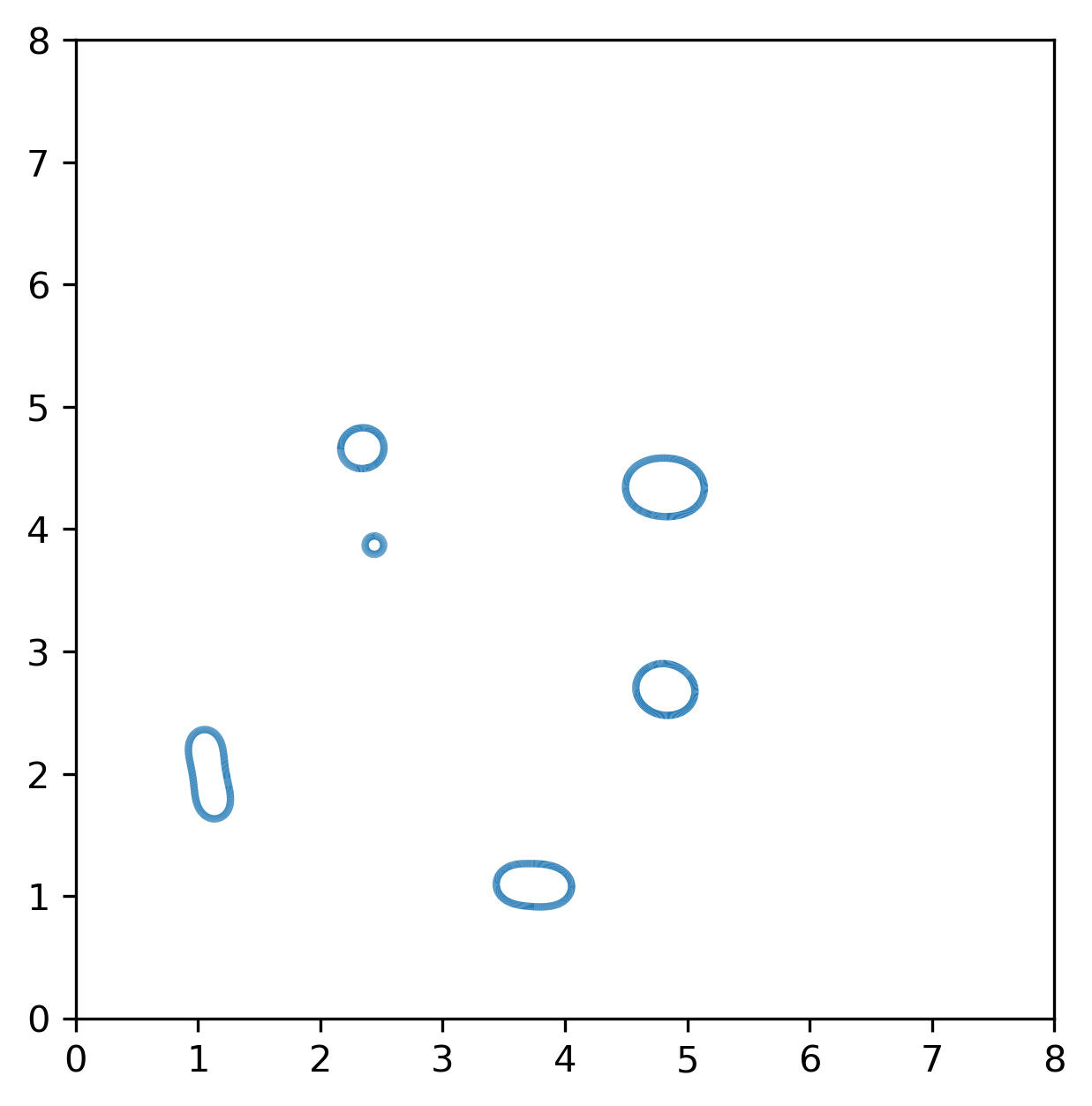}
\includegraphics[width=0.3\textwidth]{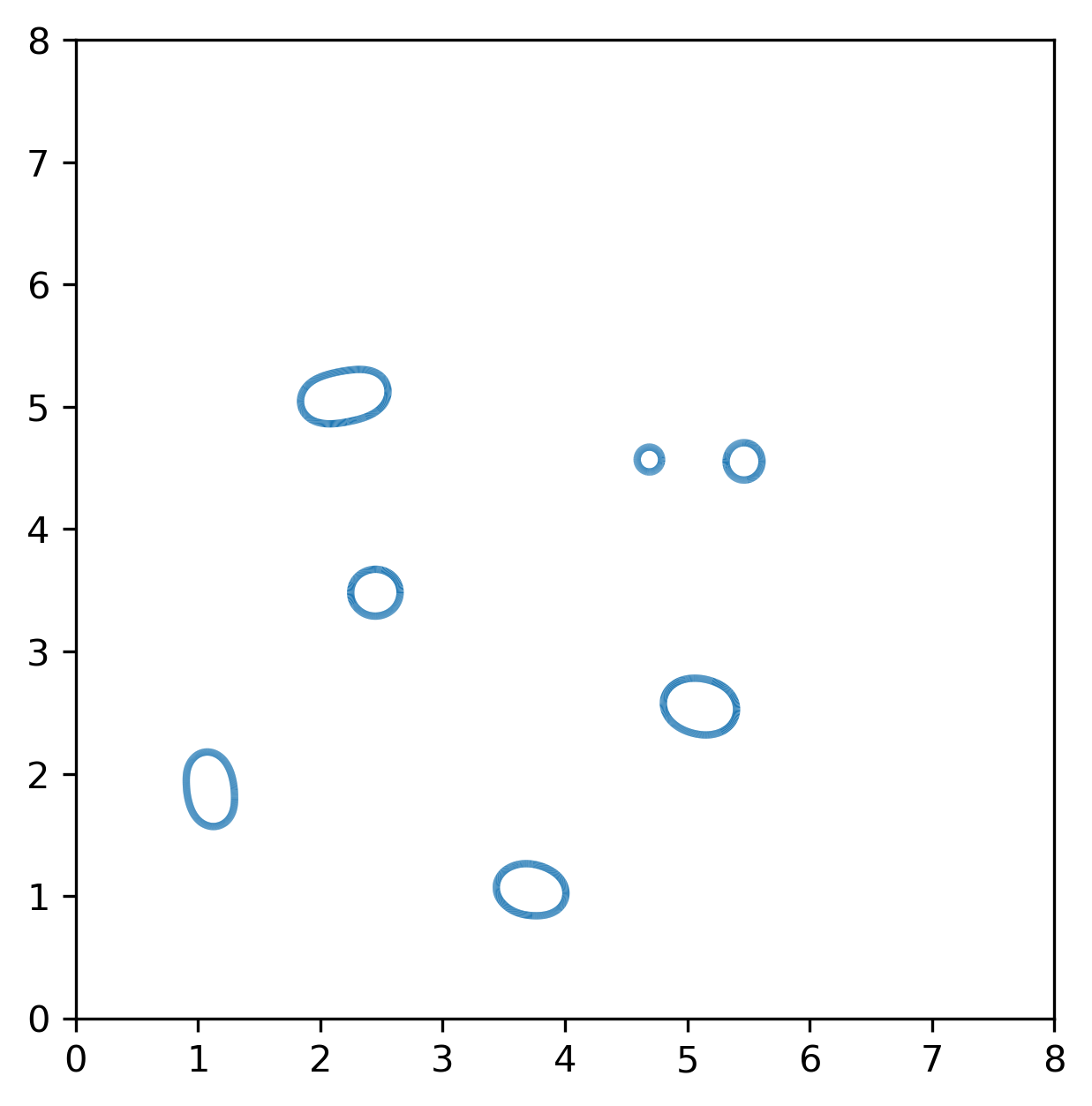}
\includegraphics[width=0.3\textwidth]{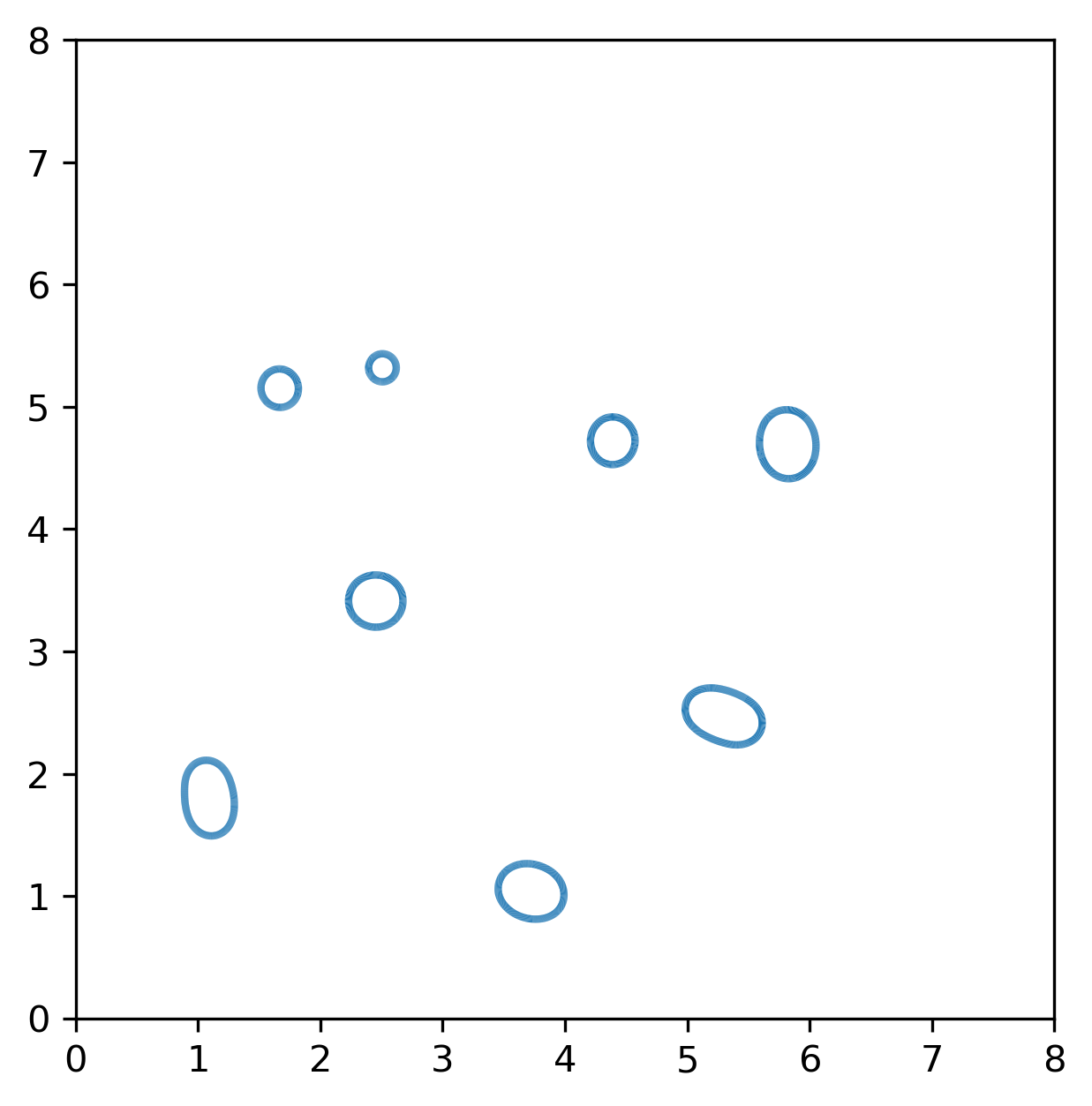}
\includegraphics[width=0.3\textwidth]{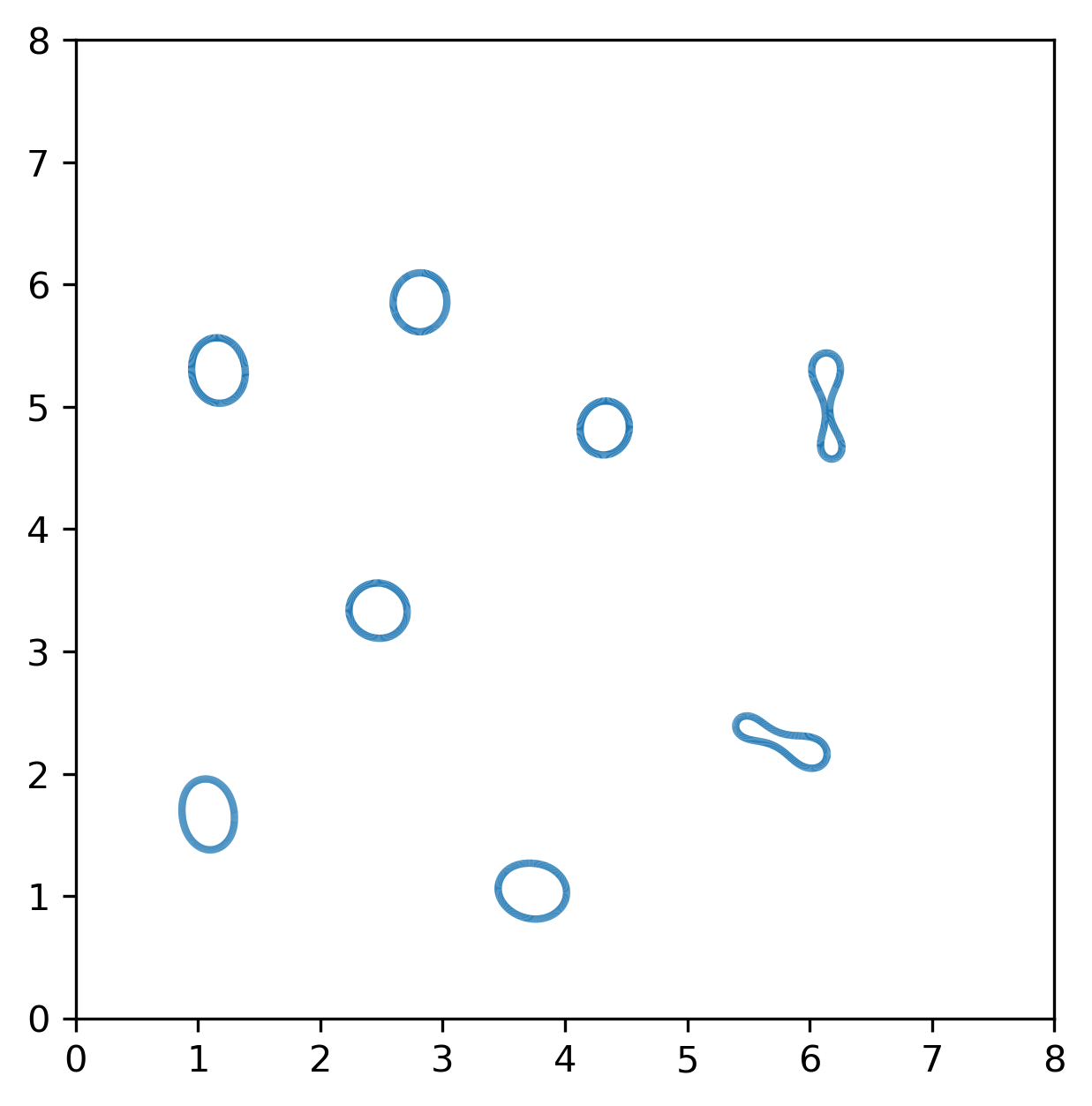}
\includegraphics[width=0.3\textwidth]{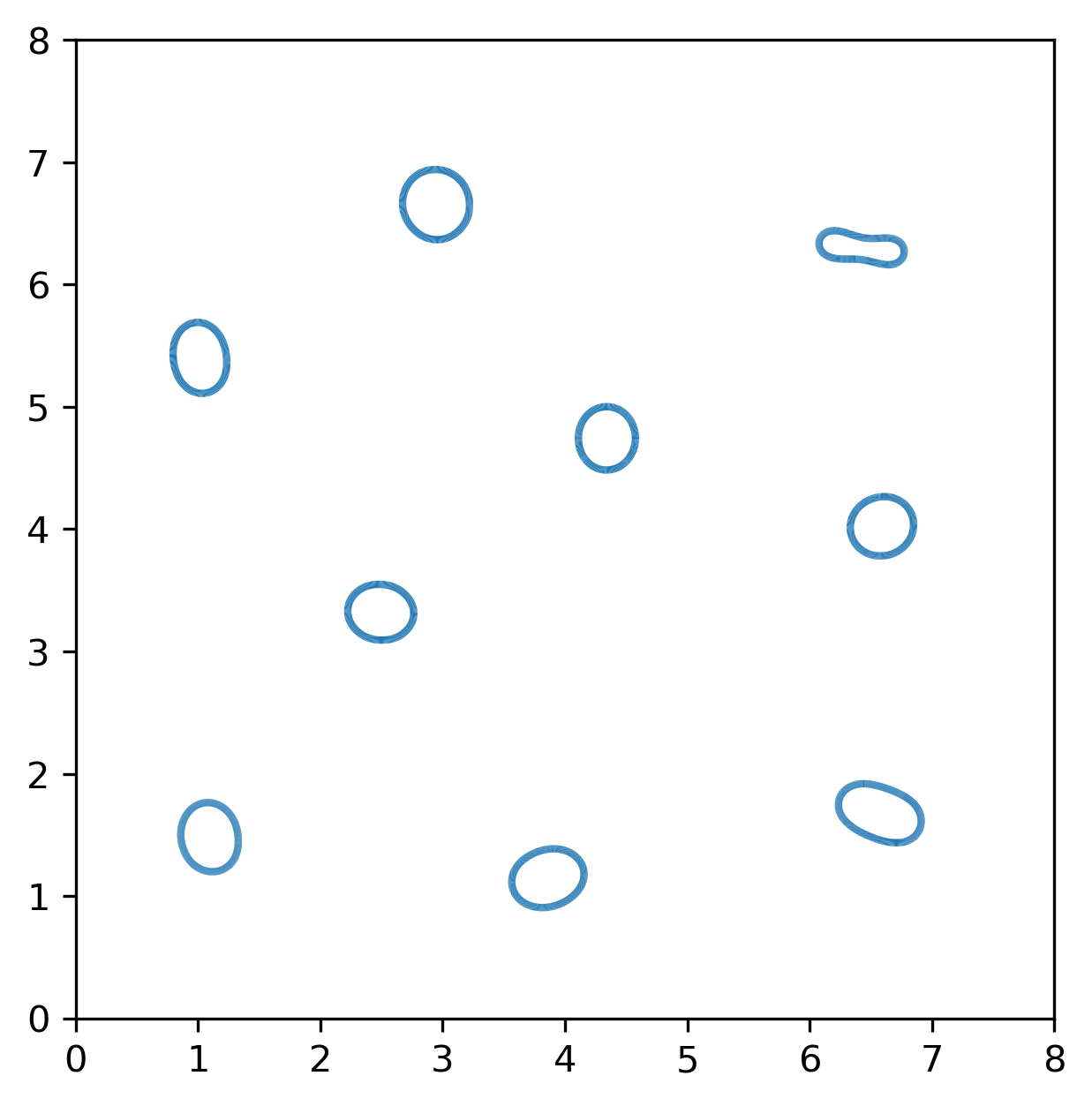}
\includegraphics[width=0.3\textwidth]{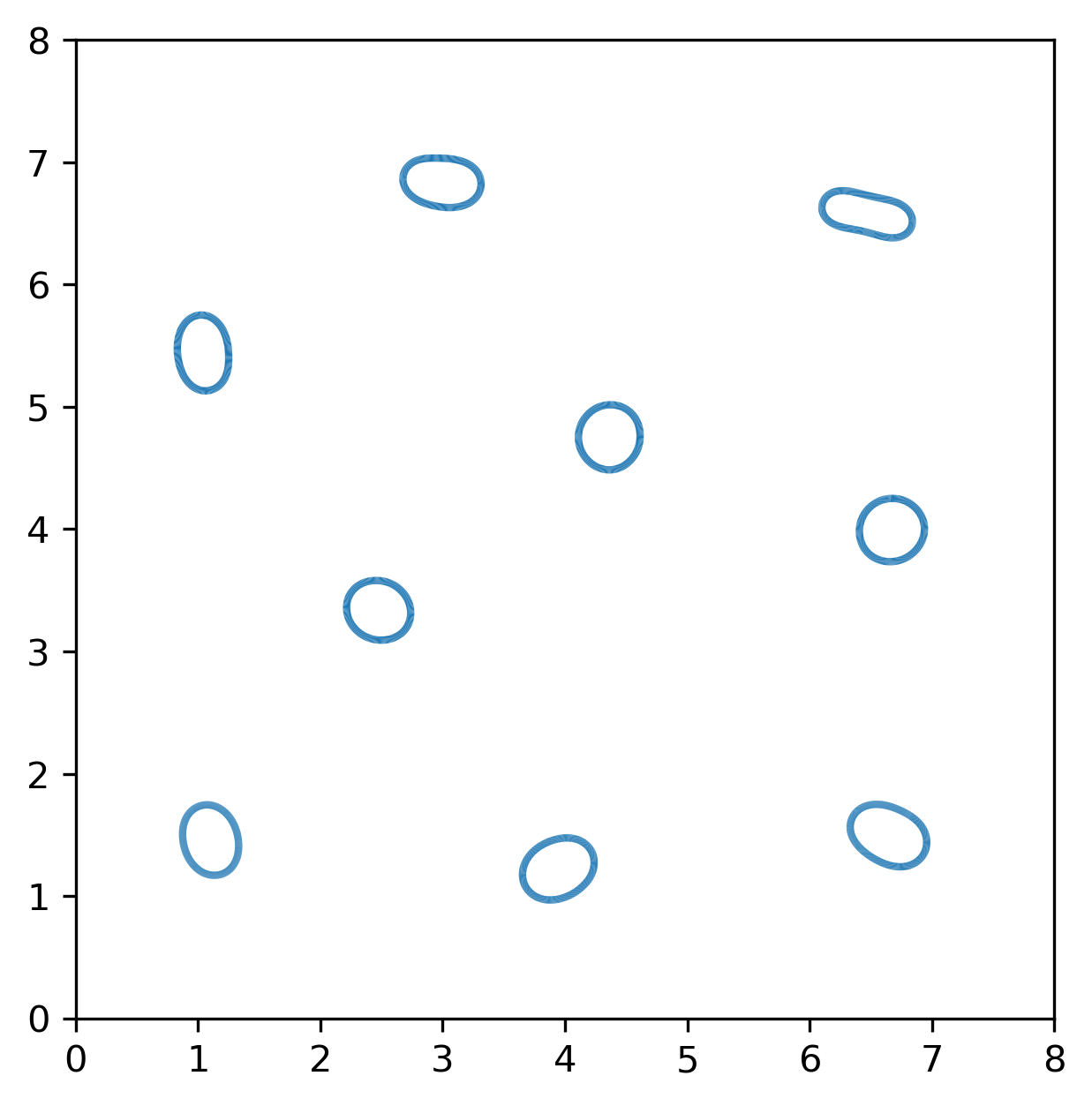}
\caption{(Minus-phase model: $\Omega = (0,8)^2$)
Simulation with pinch-offs for $\alpha= 0.02$ and $S_I=-1.8$.
We show the evolution at times $t=0,4,7,12,18,20,21,24,26,30,40,50$.}
\label{fig:outerradial04pert2}
\end{figure}

For a radial interface with a $\ell=2$ perturbation we show a simulation
in Figure~\ref{fig:outerradial04pert2}. The perturbation is such that 
the total dimensions of the initial interface are $0.816 \times 0.784$,
and we centre the initial interface around the point $(2,2)^{\top}$ within the
domain $\Omega = (0,8)^2$.  
We let $\alpha=0.02$ and $S_I = -1.8$, so that the equilibrium radius is 
$0.3553$.
Throughout the evolution we notice a preference for the growth of the $\ell=2$ 
perturbation, as well as a tendency to keep the size of the cells close to the
equilibrium shape. 
Through several pinch-offs more and more droplets are
created, that visibly interact with each other and affect each other's growth. In the last row of Figure~\ref{fig:outerradial04pert2}, we can observe two examples where some droplets close to the right boundary appear to initiate a pinch-off, but ultimately they remain whole due to the lack of available space or the lack of substrate from the minus-phase. Ultimately, rather than pinching off, they just migrate away from the other droplets.

\begin{figure}
\center
\includegraphics[angle=-0,width=0.24\textwidth]{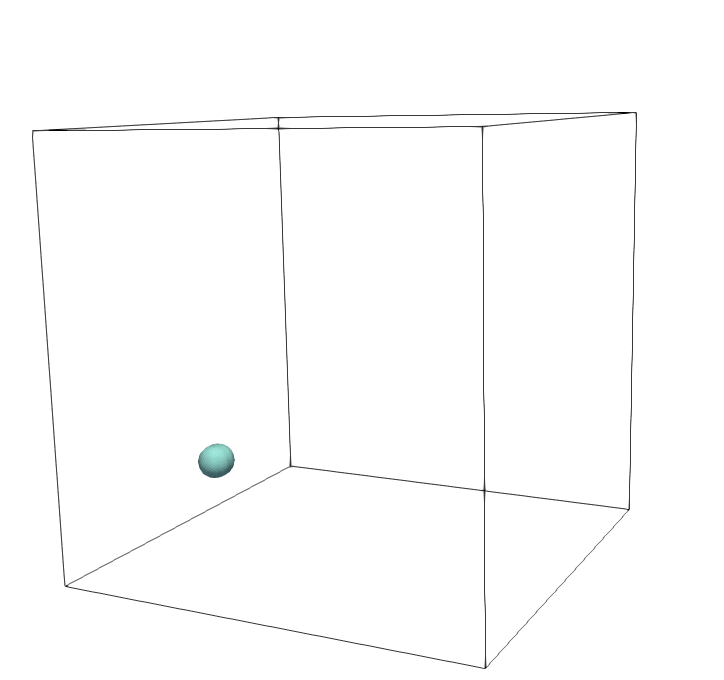}
\includegraphics[angle=-0,width=0.24\textwidth]{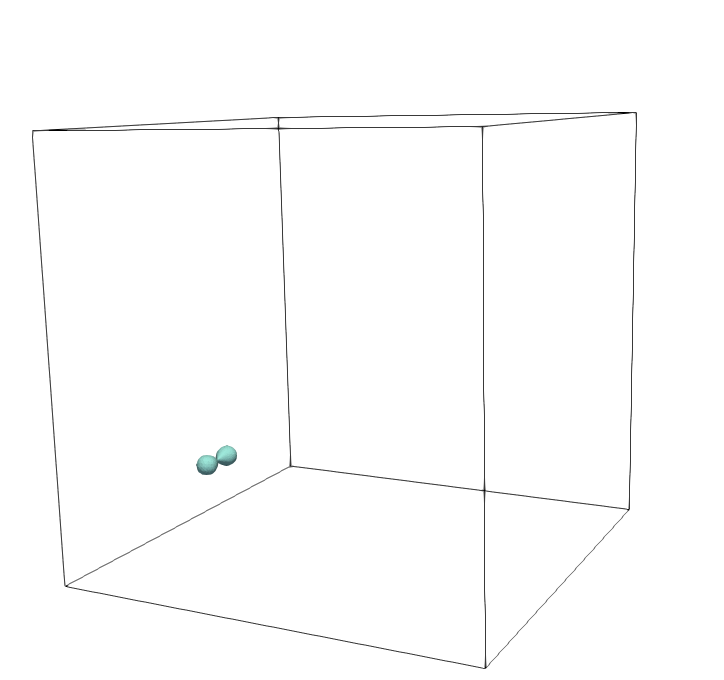}
\includegraphics[angle=-0,width=0.24\textwidth]{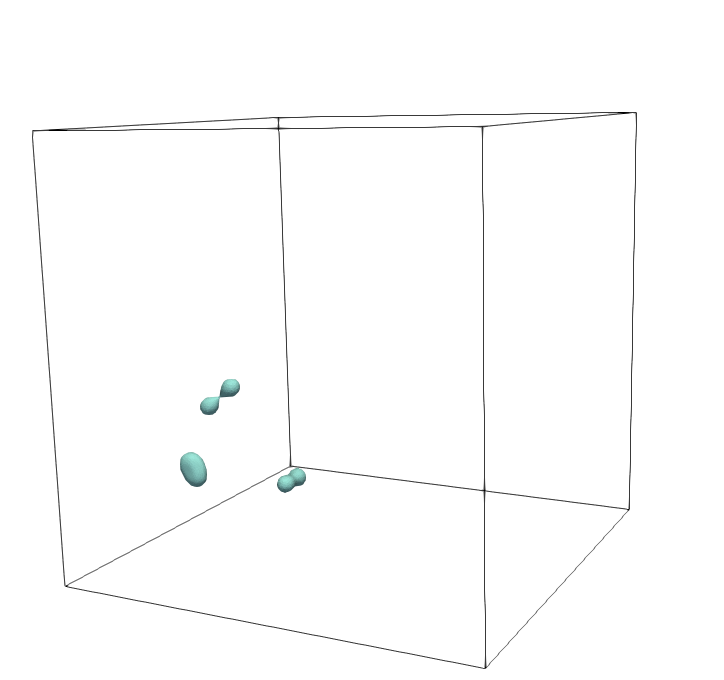}
\includegraphics[angle=-0,width=0.24\textwidth]{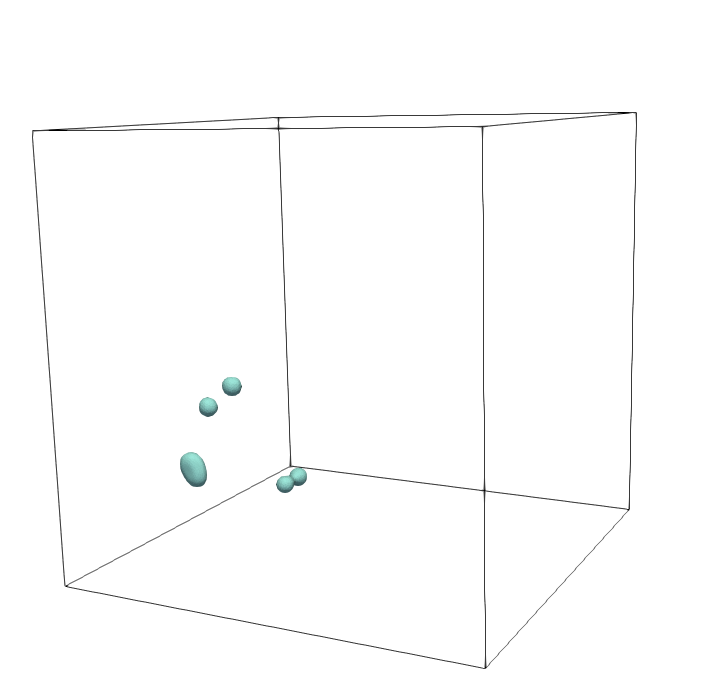}
\includegraphics[angle=-0,width=0.24\textwidth]{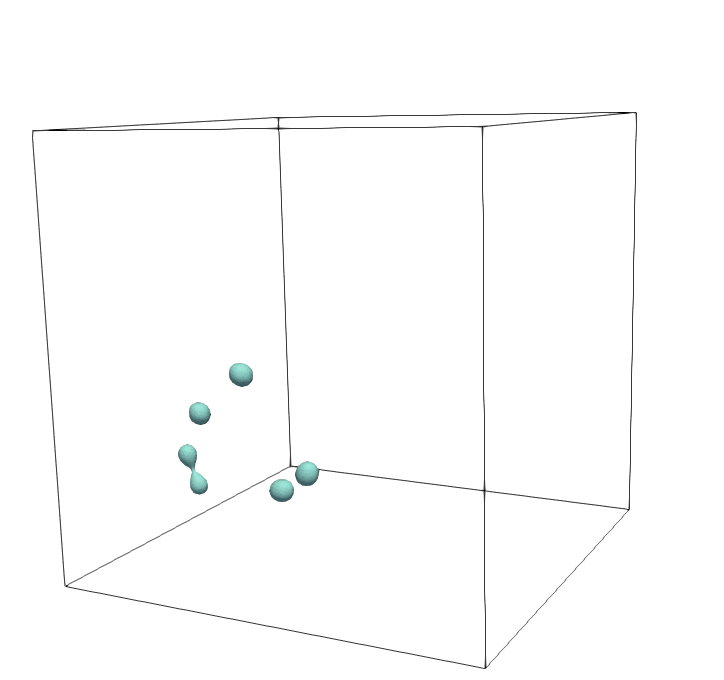}
\includegraphics[angle=-0,width=0.24\textwidth]{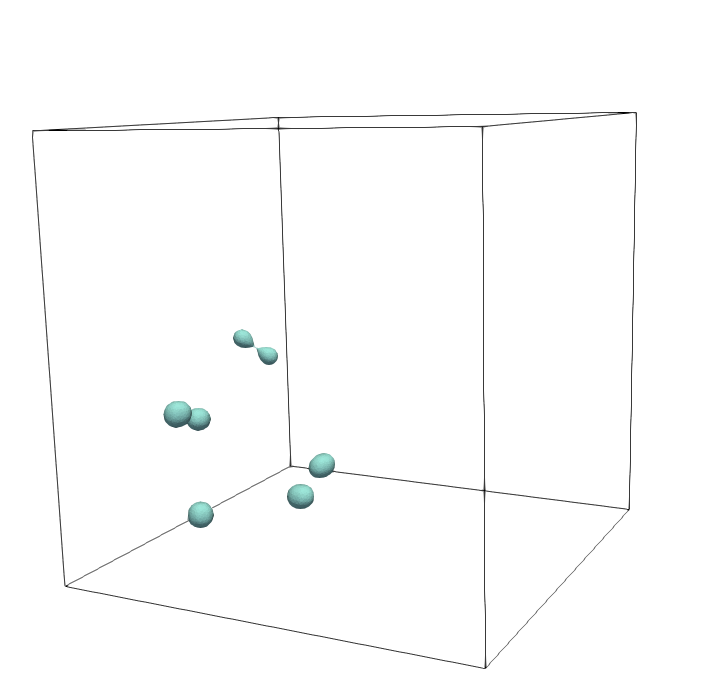}
\includegraphics[angle=-0,width=0.24\textwidth]{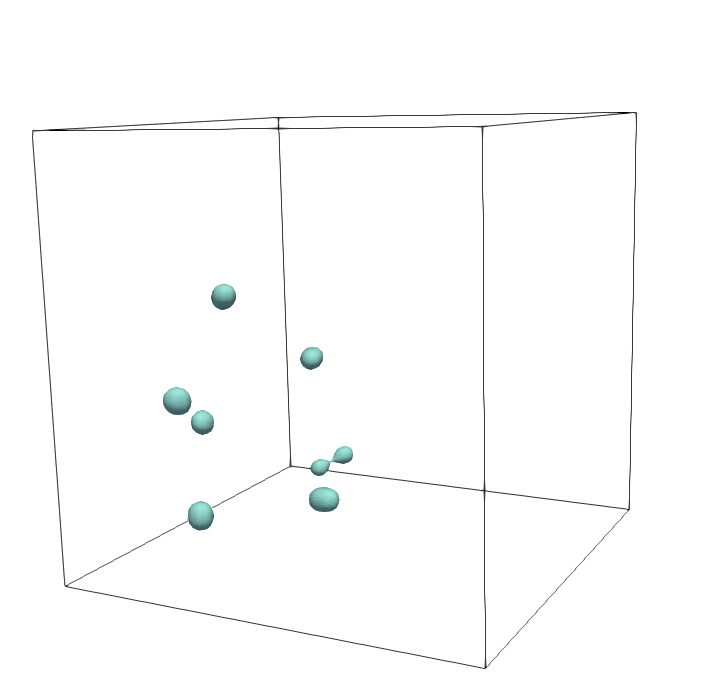}
\includegraphics[angle=-0,width=0.24\textwidth]{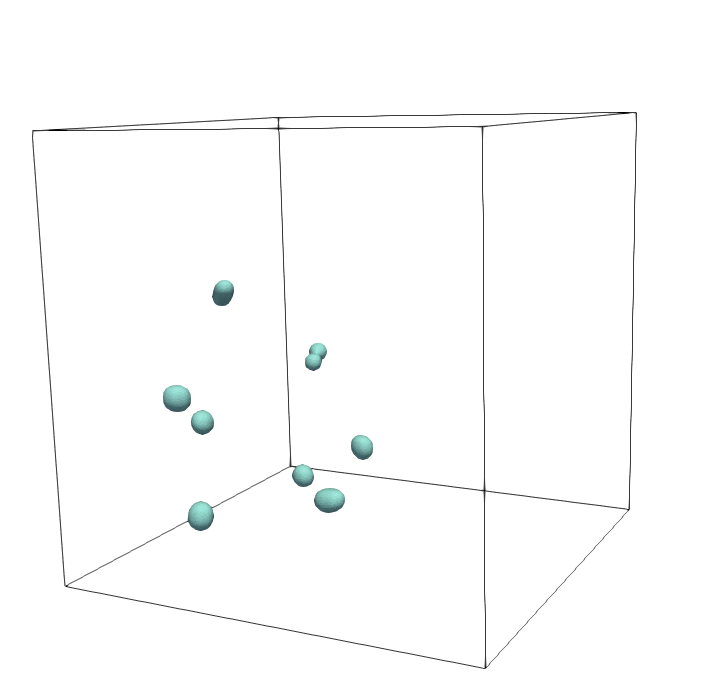}
\includegraphics[angle=-0,width=0.24\textwidth]{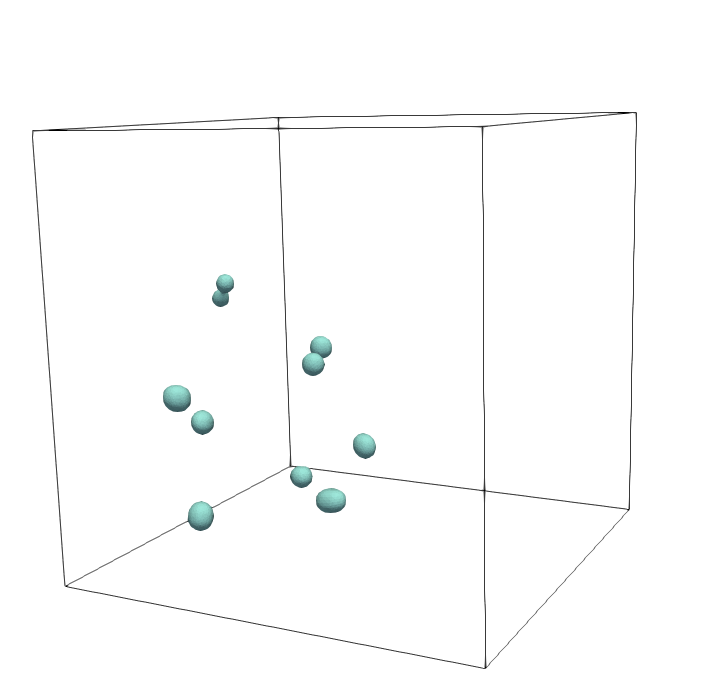}
\includegraphics[angle=-0,width=0.24\textwidth]{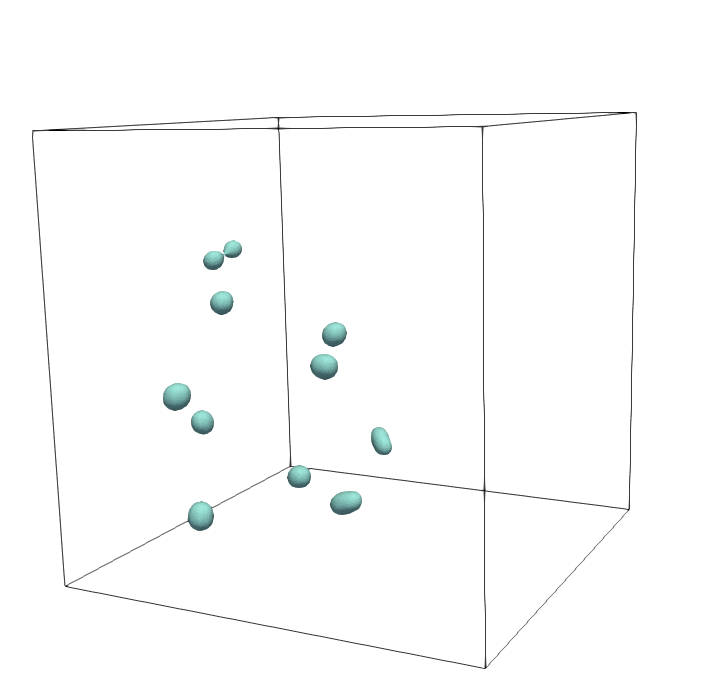}
\includegraphics[angle=-0,width=0.24\textwidth]{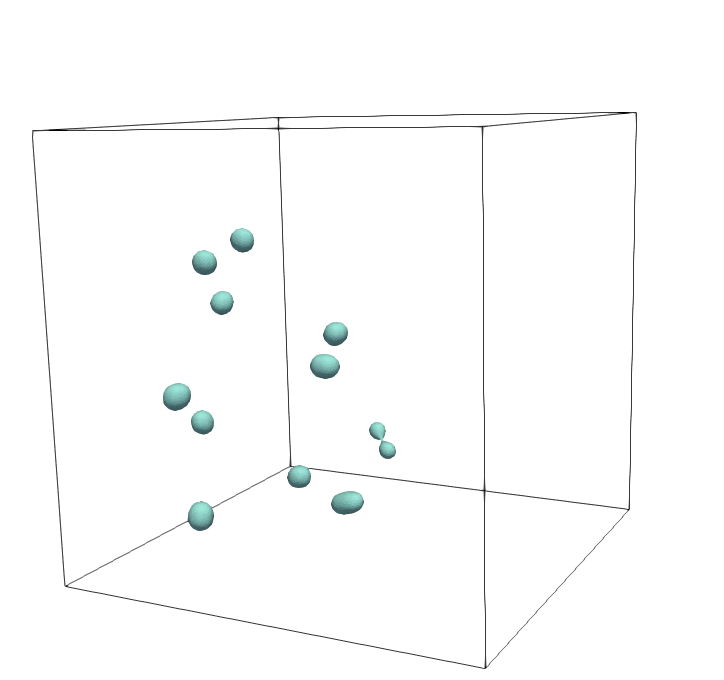}
\includegraphics[angle=-0,width=0.24\textwidth]{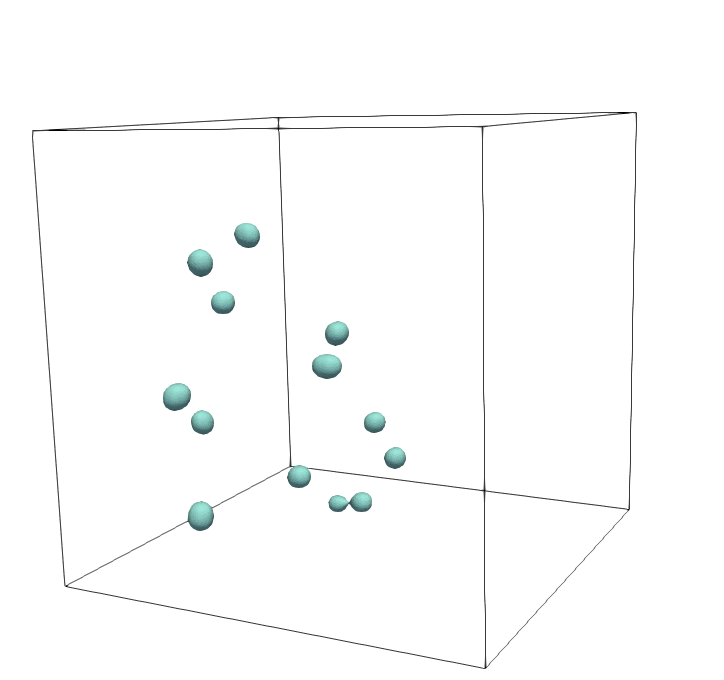}
\includegraphics[angle=-0,width=0.24\textwidth]{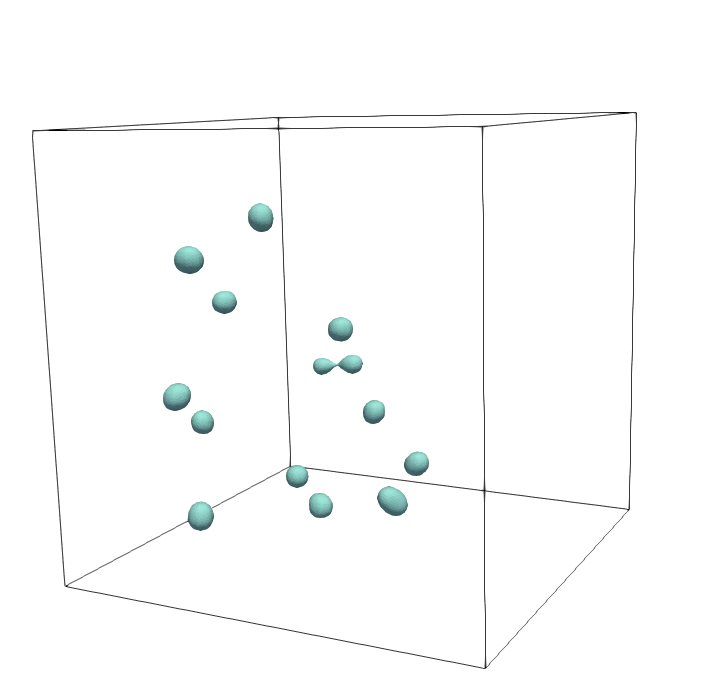}
\includegraphics[angle=-0,width=0.24\textwidth]{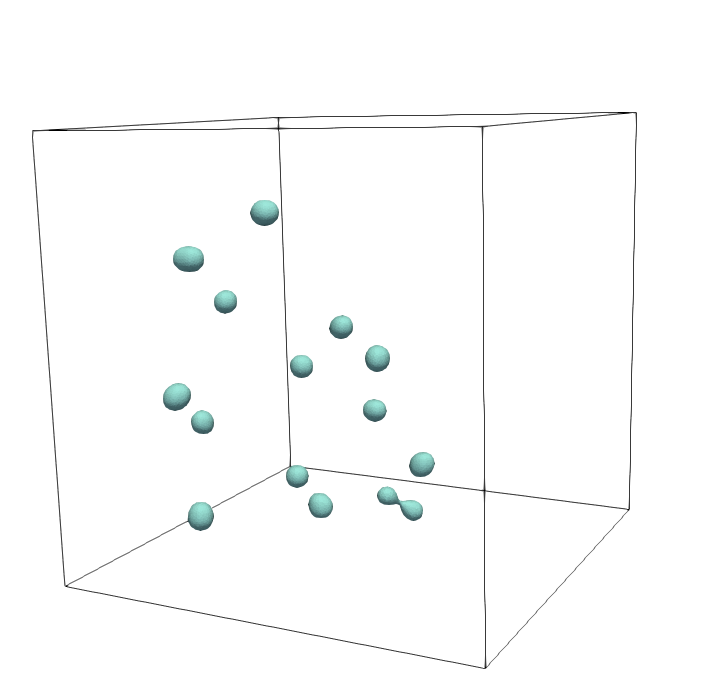}
\includegraphics[angle=-0,width=0.24\textwidth]{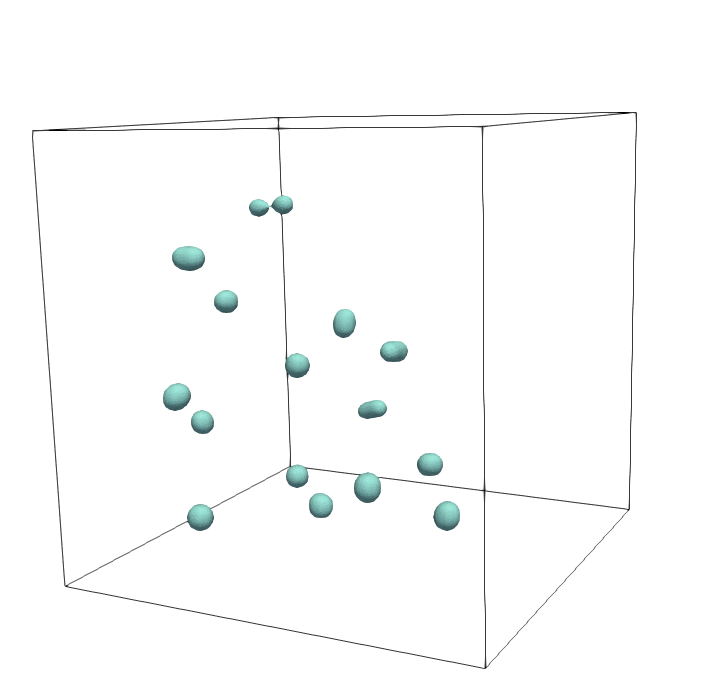}
\includegraphics[angle=-0,width=0.24\textwidth]{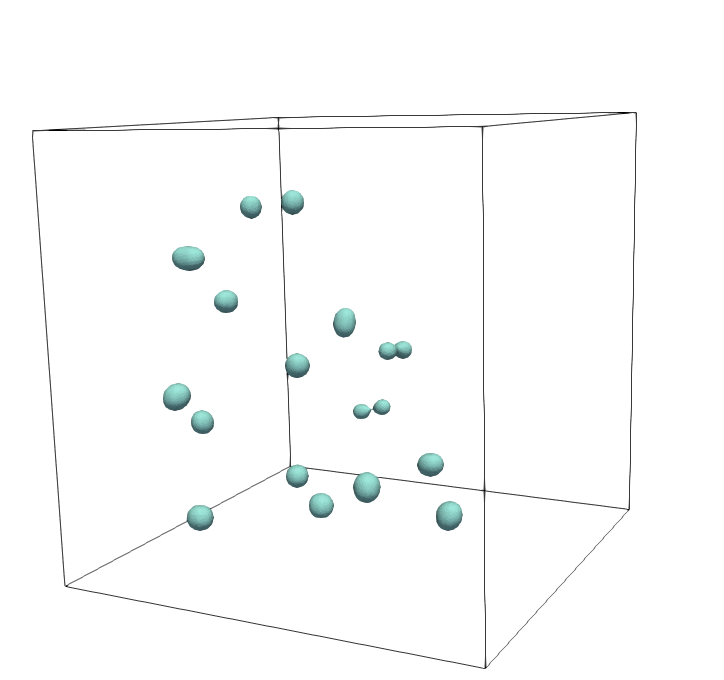}
\includegraphics[angle=-0,width=0.24\textwidth]{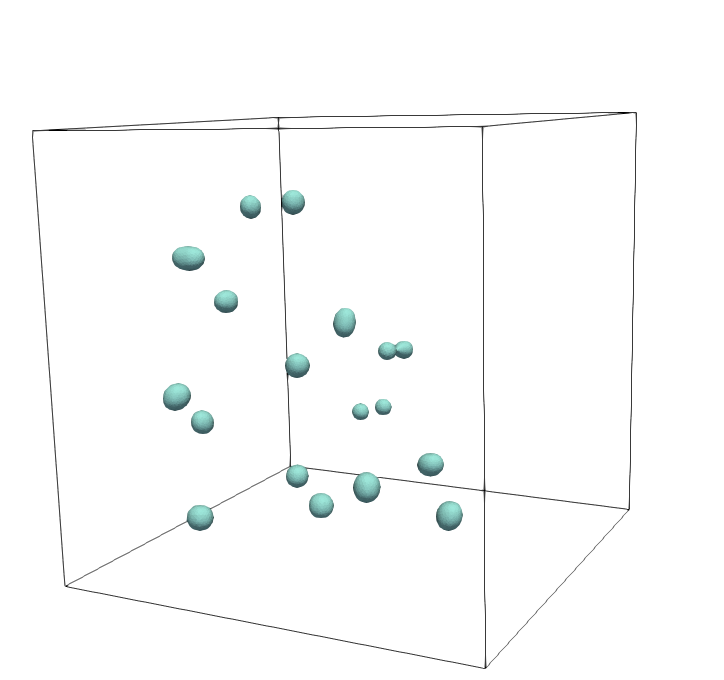}
\includegraphics[angle=-0,width=0.24\textwidth]{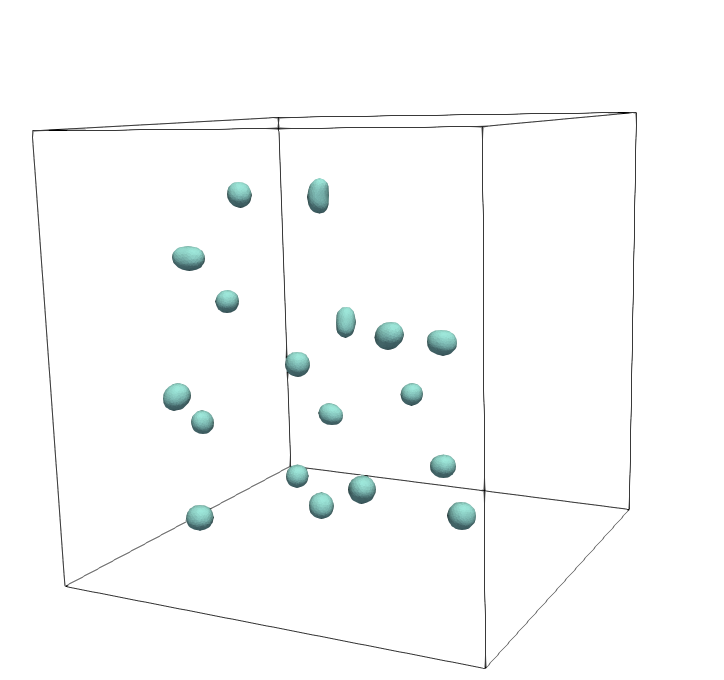}
\caption{(Minus-phase model: $\Omega = (0,8)^3$)
Simulation with pinch-offs for $\alpha=0.02$ and $S_I=-4$.
We show the evolution at times $t=0$, 1.6, 18.13, 18.23, 19.65, 27.51, 38.11,
40.66, 41.37, 49.43, 50.54, 51.6, 58.29, 61.31, 67.32, 68.08, 68.15, 80.
All but the initial and final snap-shot show a time just before a pinch-off 
occurs.
}
\label{fig:3dtopo2offset}
\end{figure}

For a corresponding three-dimensional experiment, we use the parameters
$\alpha=0.02$ and $S_I=-4$, so that the equilibrium radius is $0.2706$.
We start the simulation within $\Omega= (0,8)^3$
from a cigar-like initial shape of total dimension 
$0.7\times 0.6 \times 0.6$, centred at $(2,2,2)^{\top}$. 
Similarly to the two-dimensional simulations
in Figure~\ref{fig:outerradial04pert2}, we observe several pinch-offs in Figure~\ref{fig:3dtopo2offset},
induced by the growth of the mode $\ell=2$ perturbation, and it appears that pinch-off occurs more easily in three dimensions than in two dimensions.

\subsubsection*{Simulations for the phase diagram
Figure~\ref{fig:rad2D_minus}}

\begin{figure}
	\center
	\mbox{
		\includegraphics[width=0.2\textwidth]{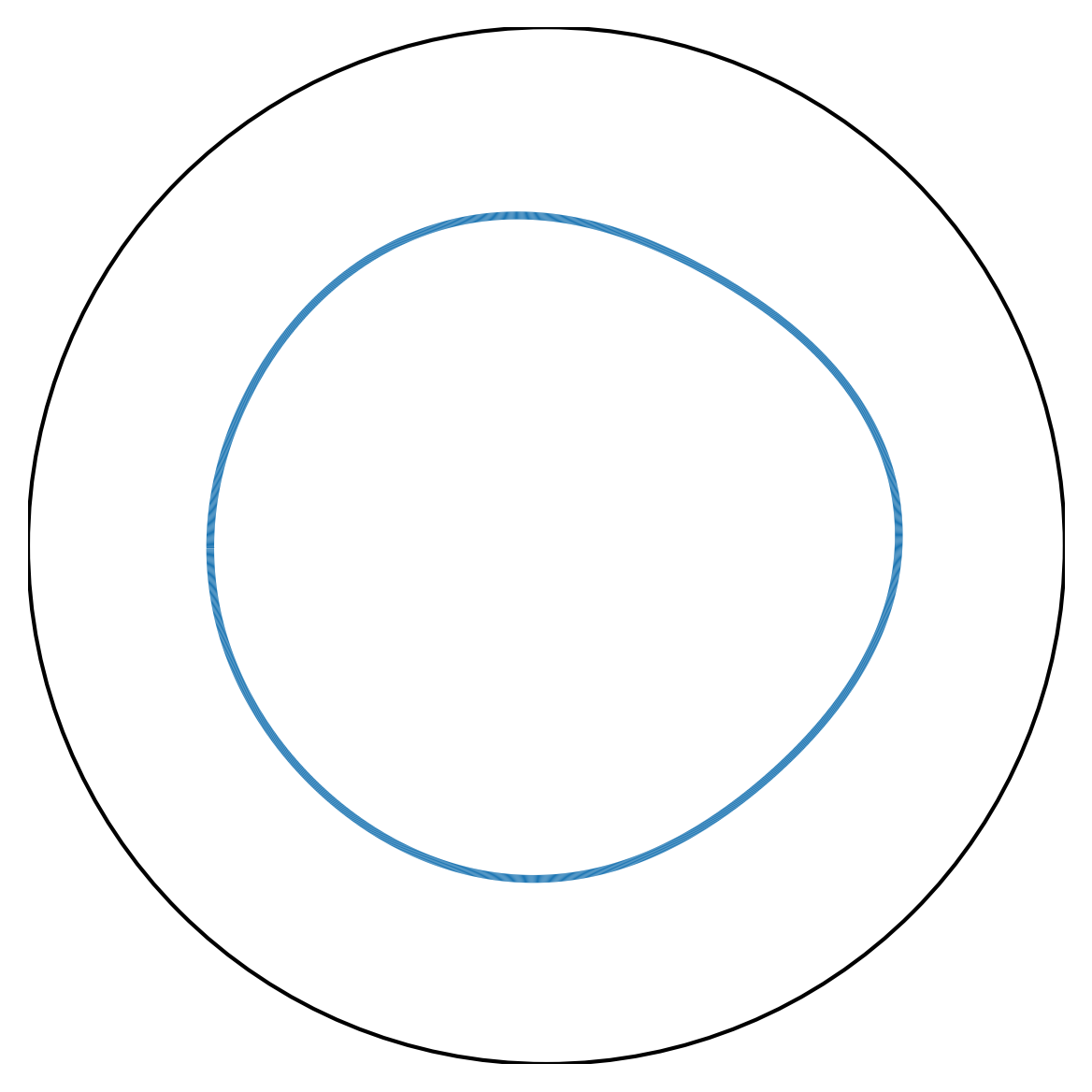}
		\includegraphics[width=0.2\textwidth]{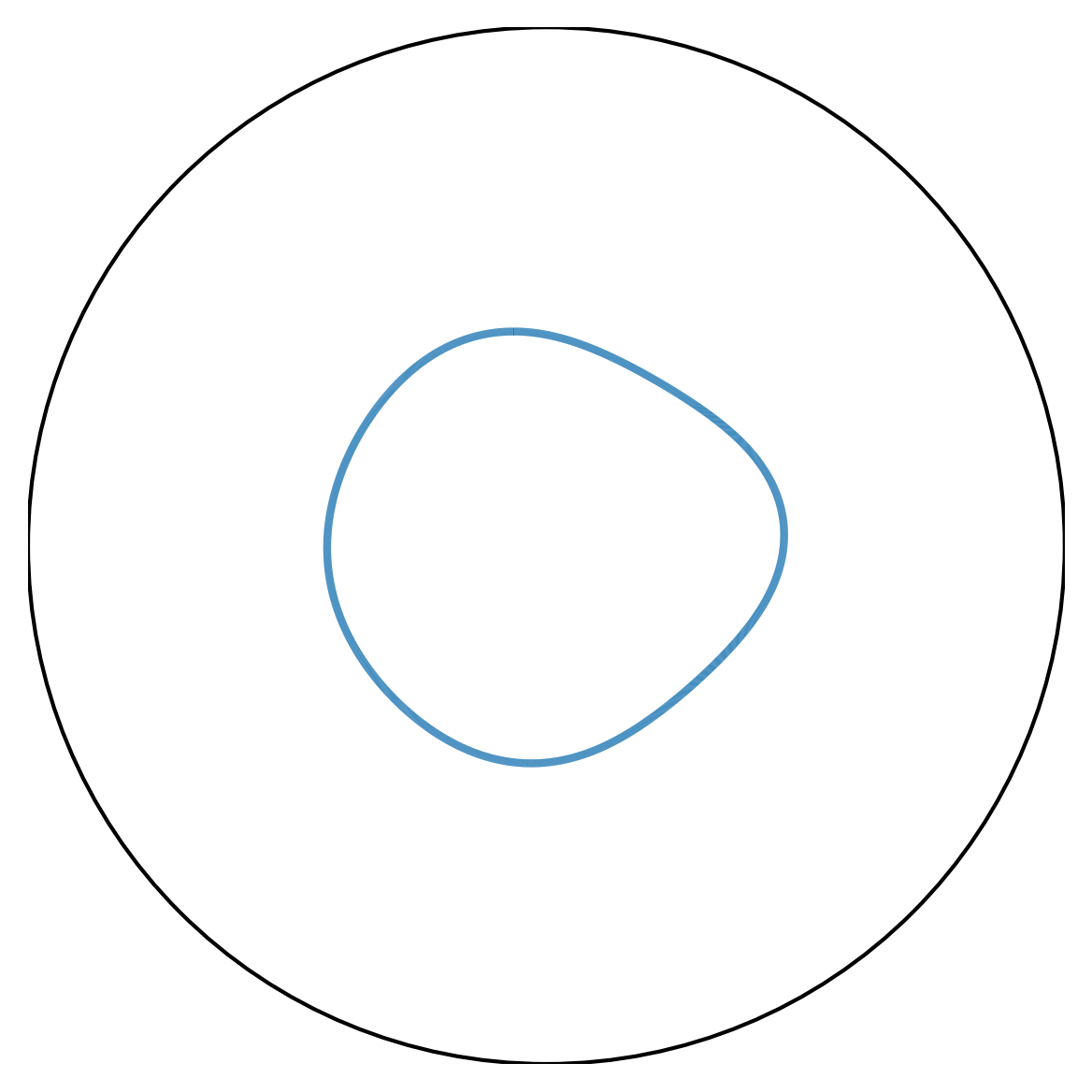}
		\includegraphics[width=0.2\textwidth]{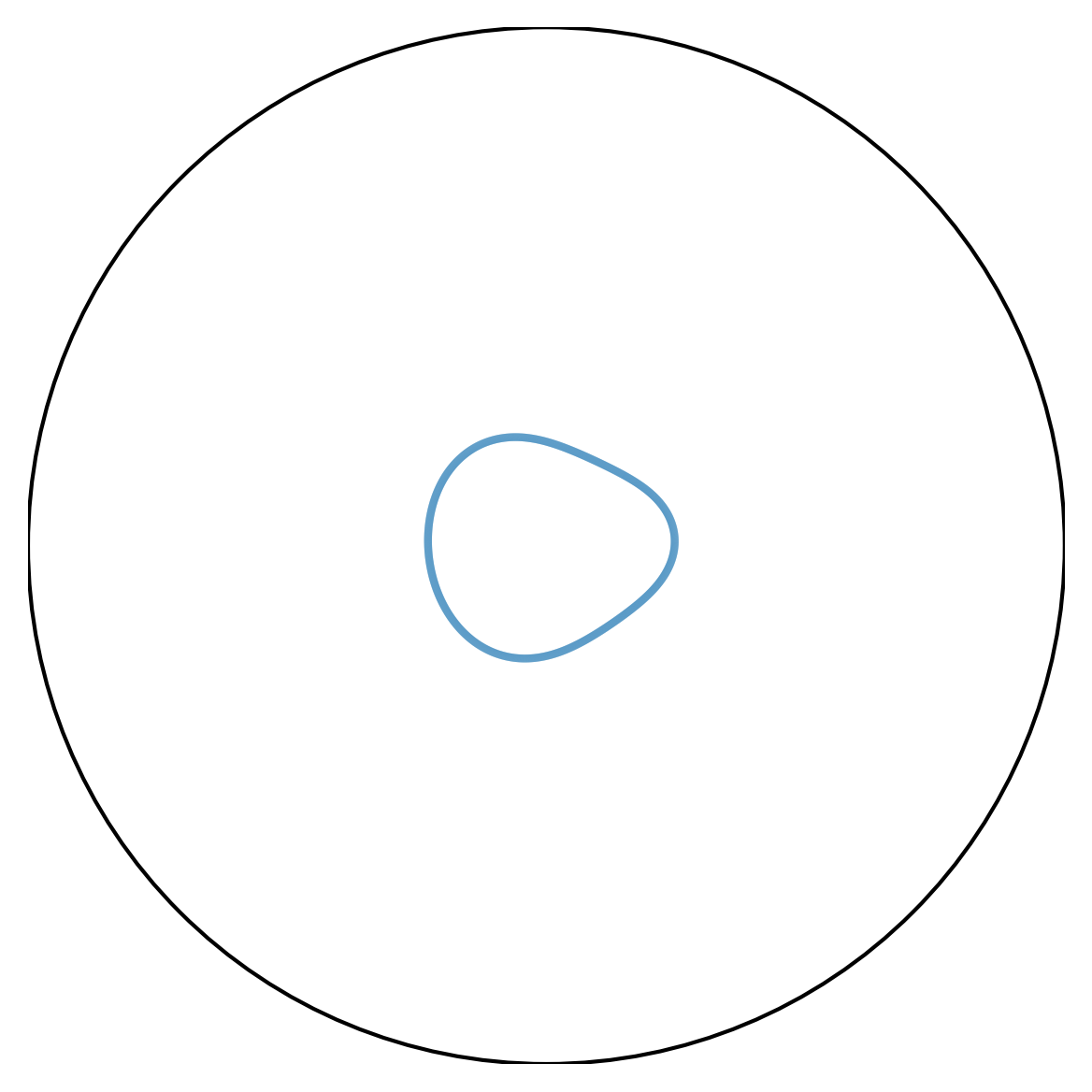}
		\includegraphics[width=0.2\textwidth]{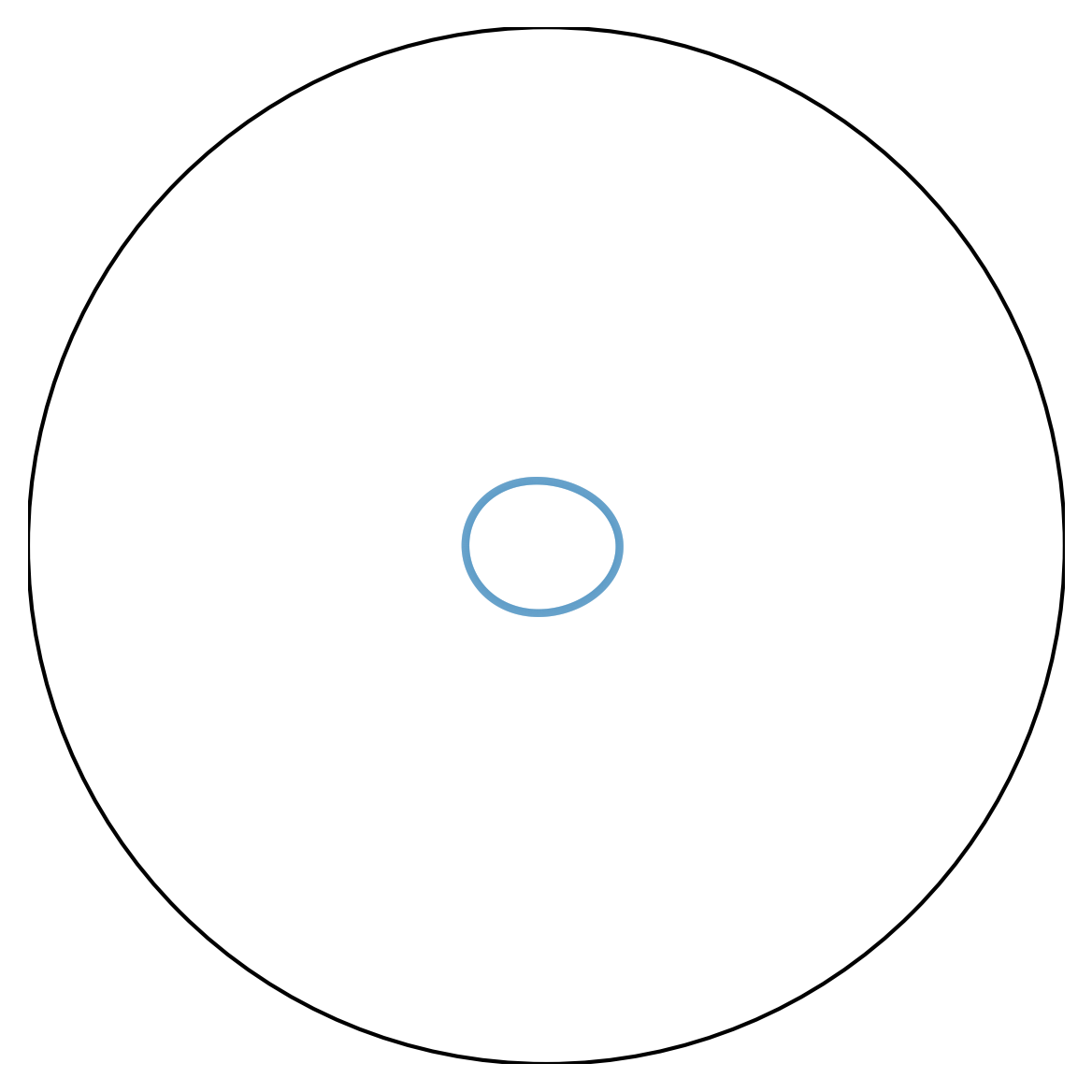}
		\includegraphics[width=0.2\textwidth]{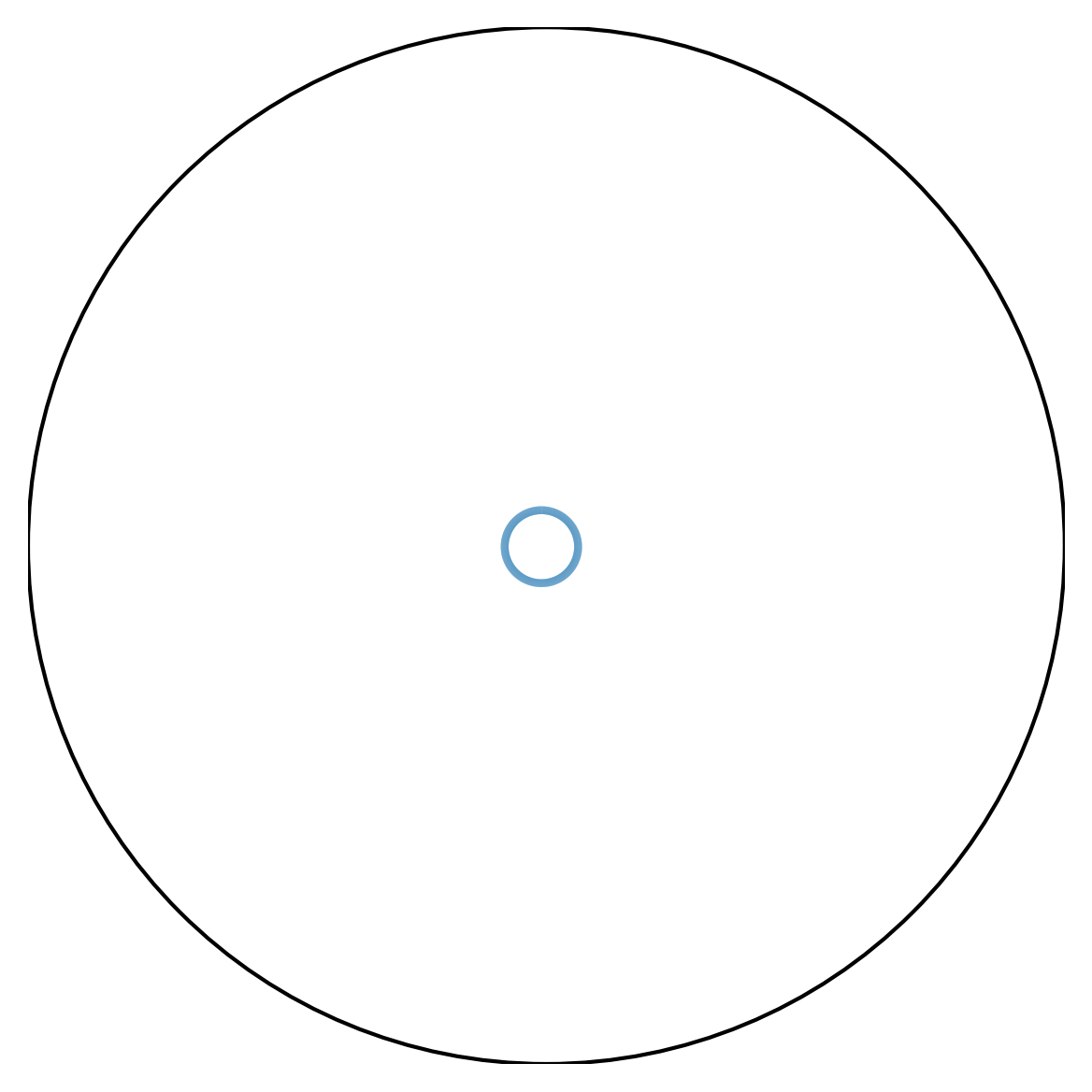}
	}
	\caption{(Minus-phase model: $\Omega = \bB^2_4(0)$)
		Initial data with mode $\ell=2$, 3, 4 perturbations around the radius 2.57,
		for $\alpha= 0.1$ and $S_I=-2$.
		We show the evolution at times $t=0,2,4,5,6$.
		The interface vanishes shortly after $t=6$.}
	\label{fig:pd_a01_SI-2_cp-4}
\end{figure}

\begin{figure}
	\center
	\mbox{
		\includegraphics[width=0.25\textwidth]{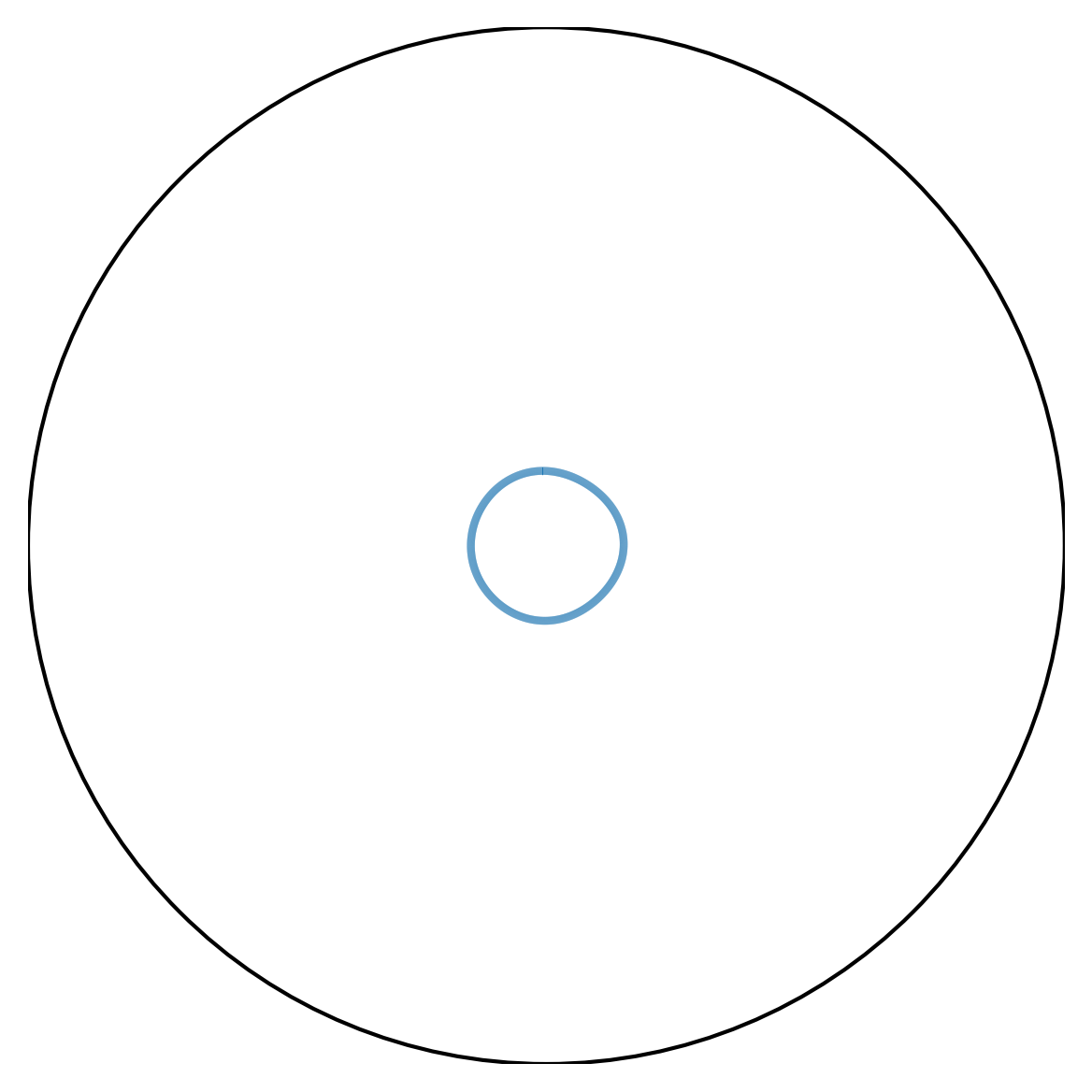}
		\includegraphics[width=0.25\textwidth]{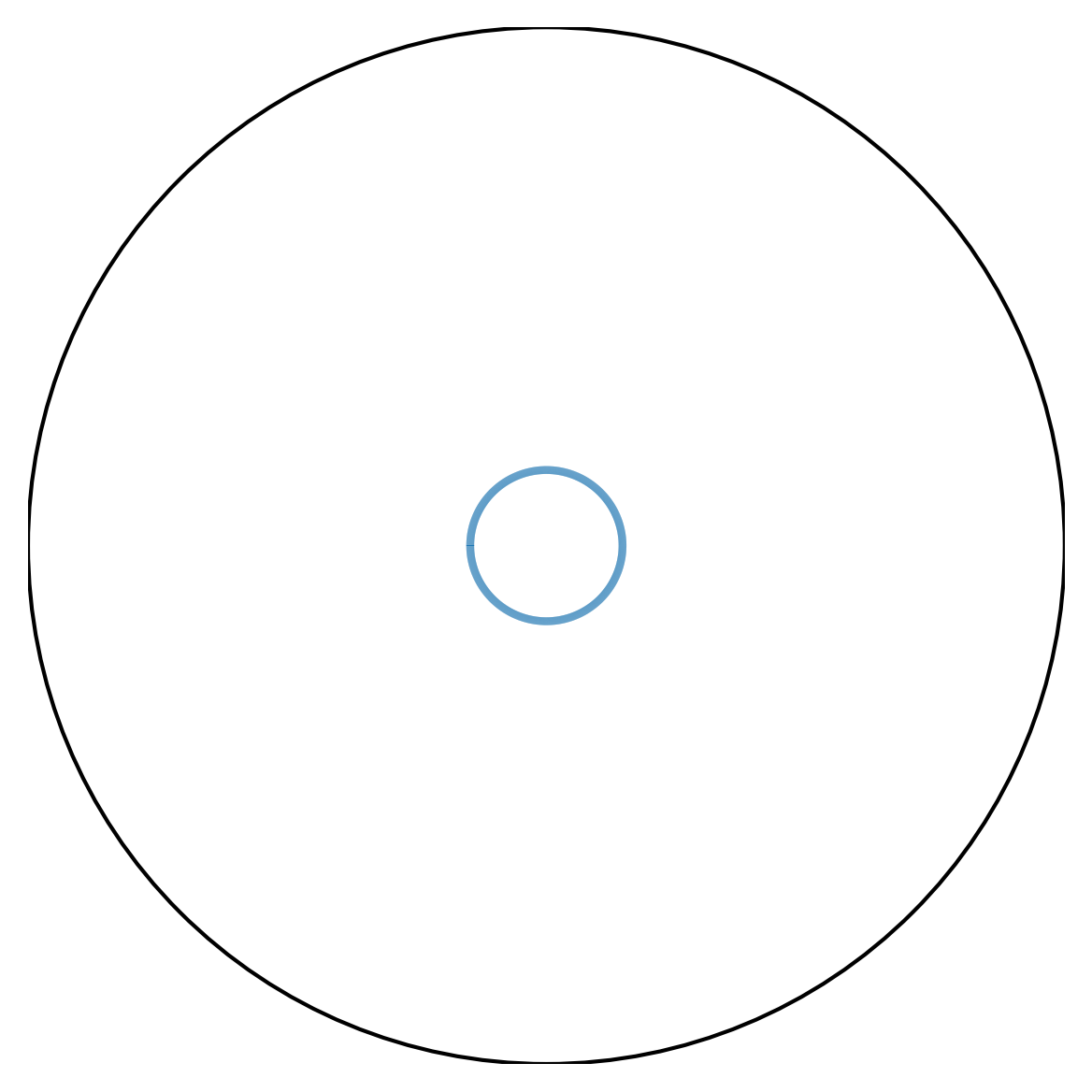}
		\includegraphics[width=0.25\textwidth]{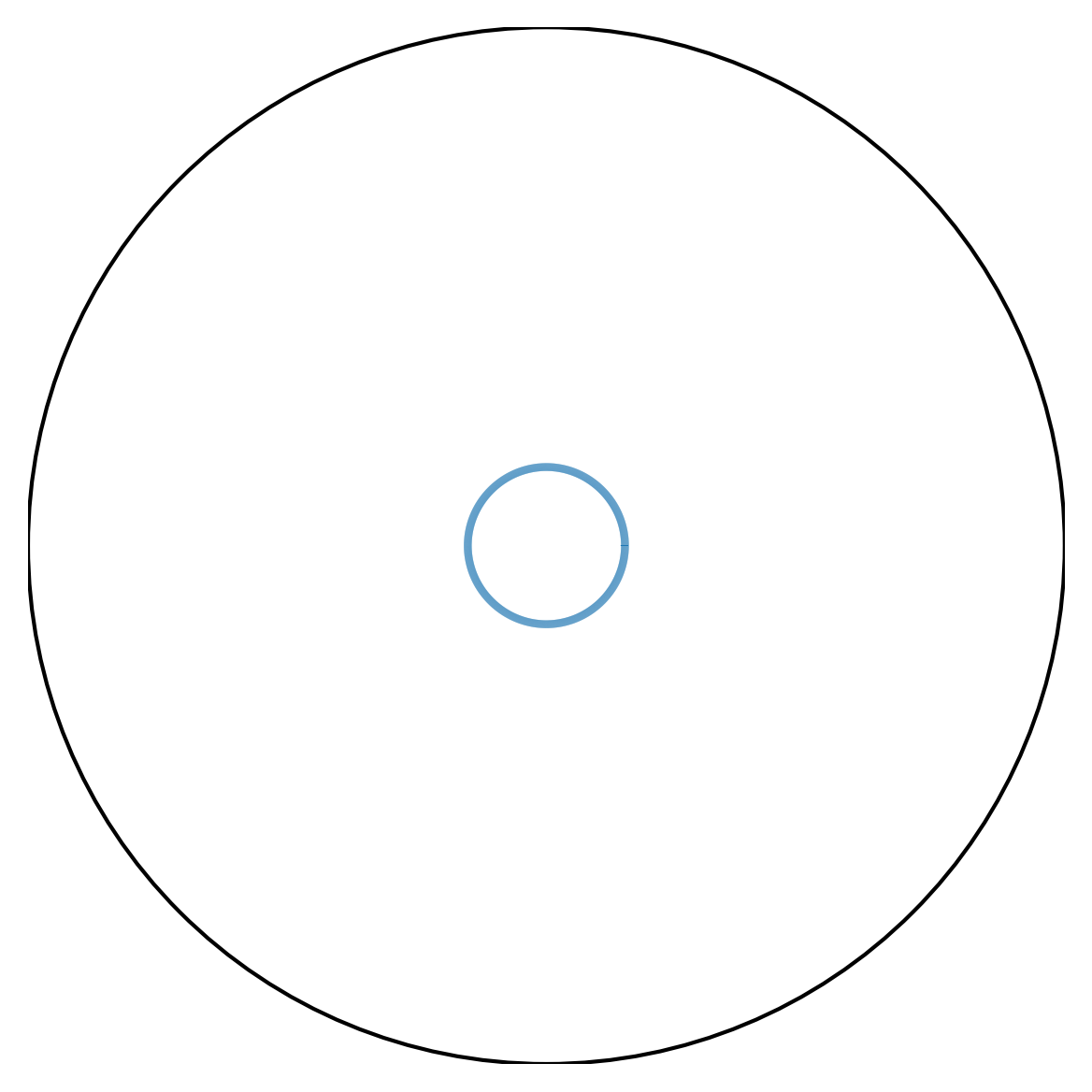}
		\includegraphics[width=0.25\textwidth]{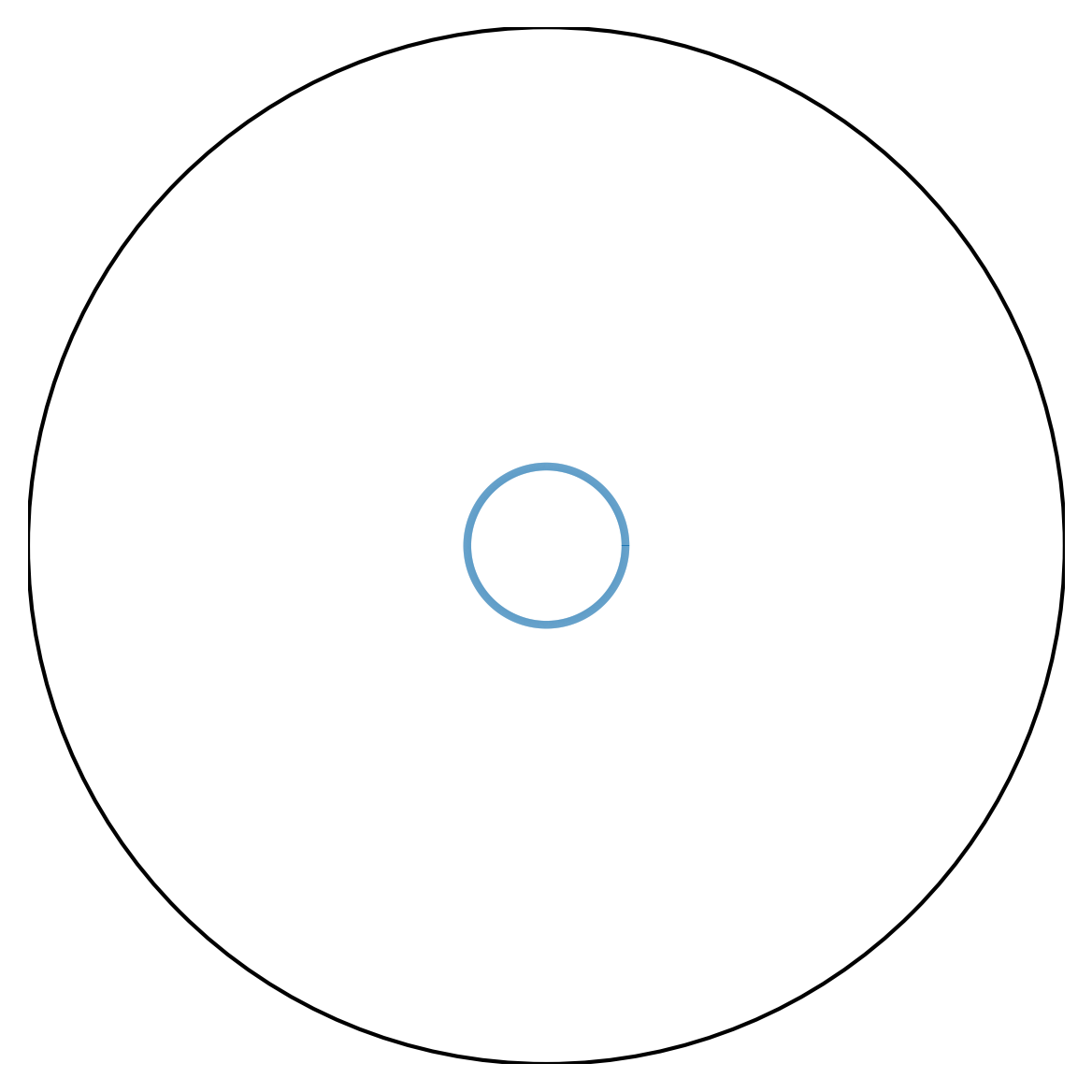}
	}
	\caption{(Minus-phase model: $\Omega = \bB^2_4(0)$)
		Initial data with mode $\ell=2$, 3, 4 perturbations around the radius 0.58,
		for $\alpha= 0.1$ and $S_I=-1.4$.
		We show the evolution at times $t=0,1,10,50$.}
	\label{fig:pd_a01_SI-14_cp-4}
\end{figure}

\begin{figure}
	\center
	\mbox{
		\includegraphics[width=0.2\textwidth]{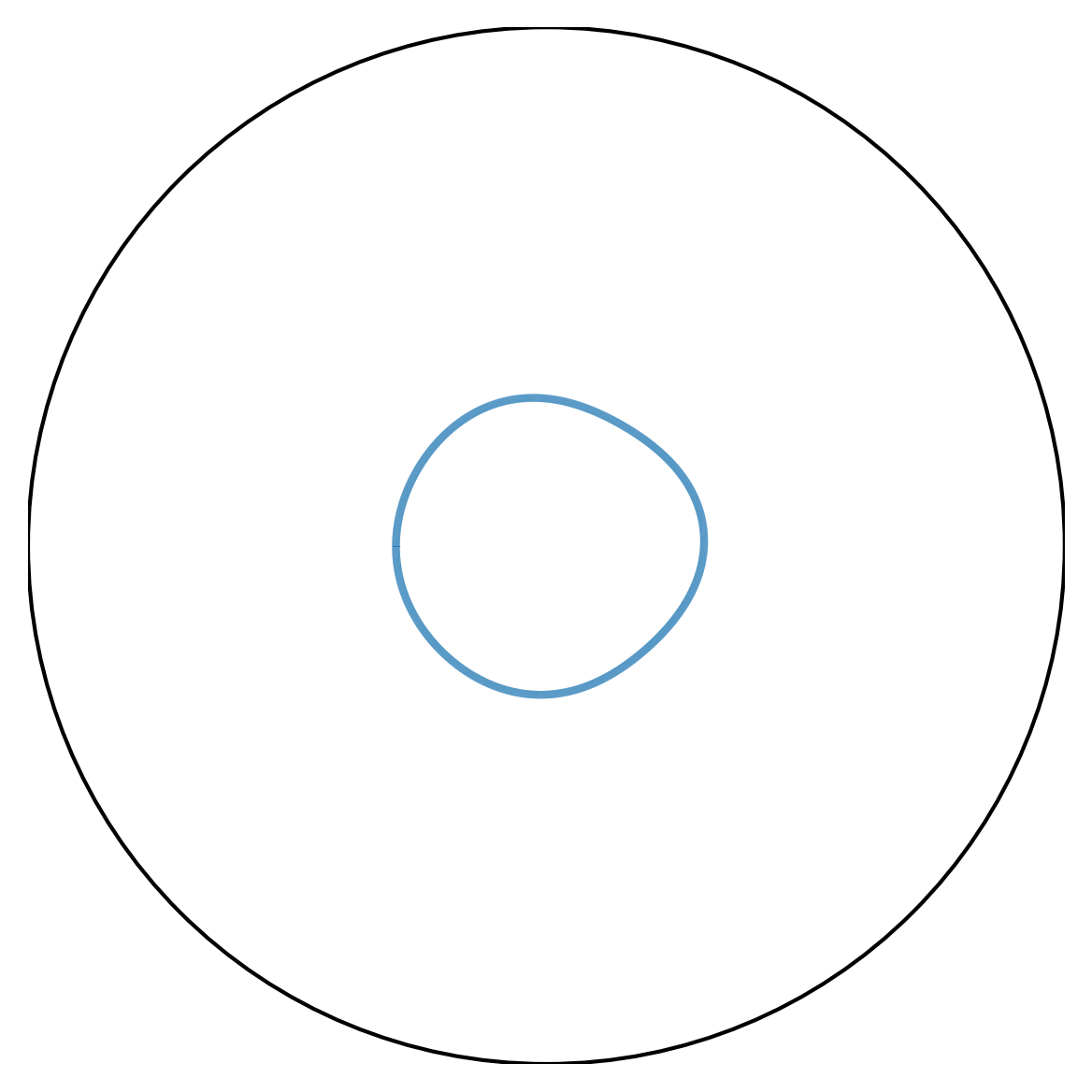}
		\includegraphics[width=0.2\textwidth]{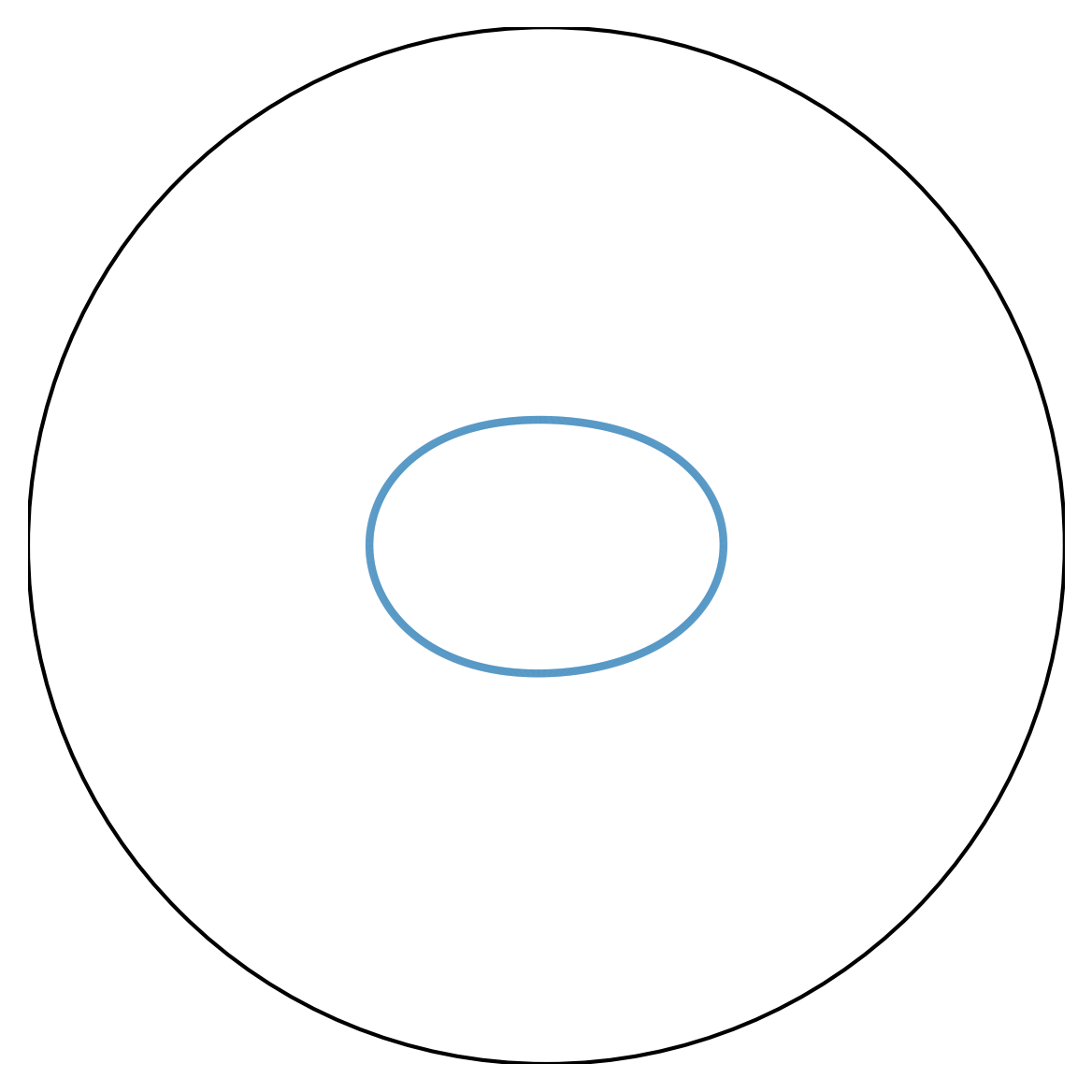}
		\includegraphics[width=0.2\textwidth]{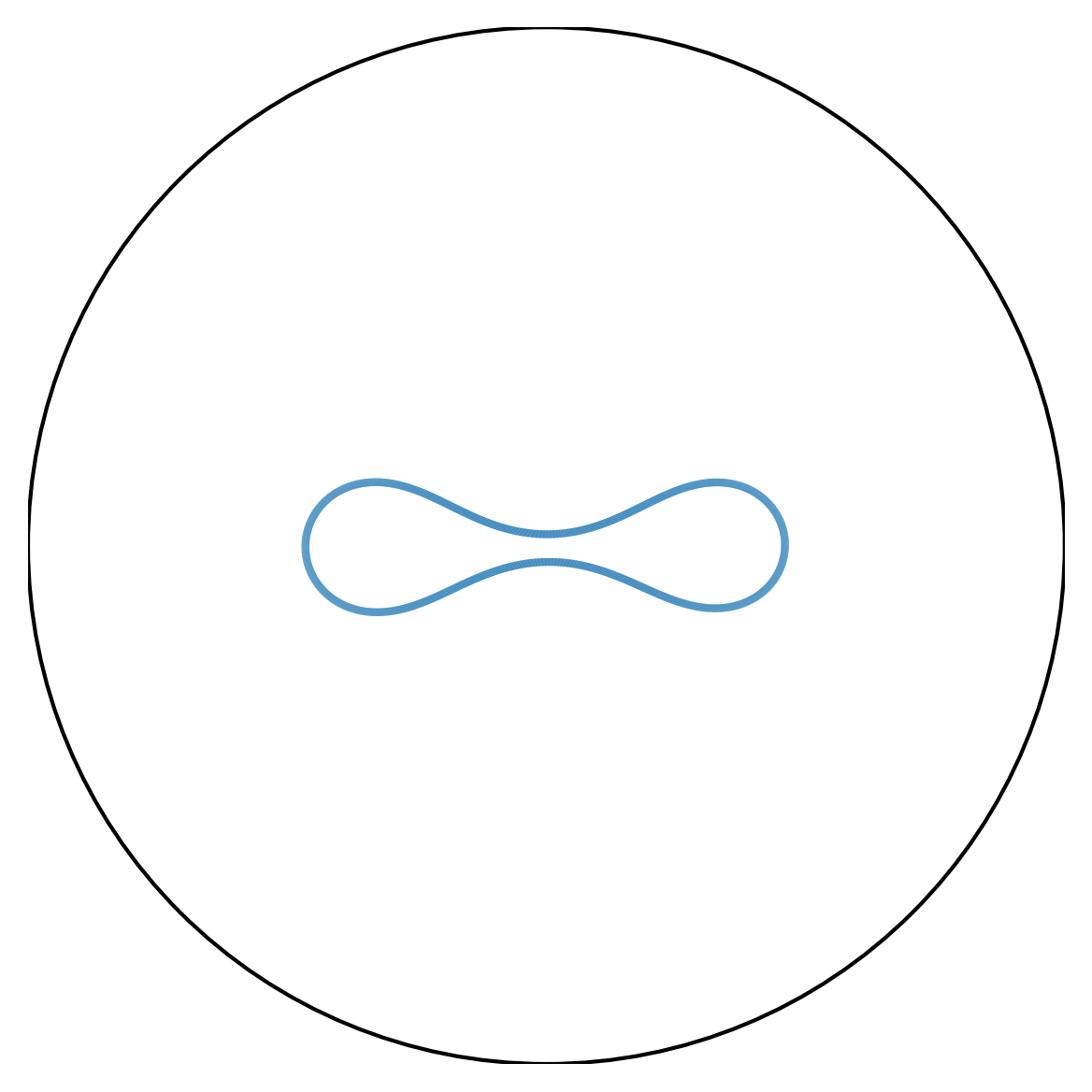}
		\includegraphics[width=0.2\textwidth]{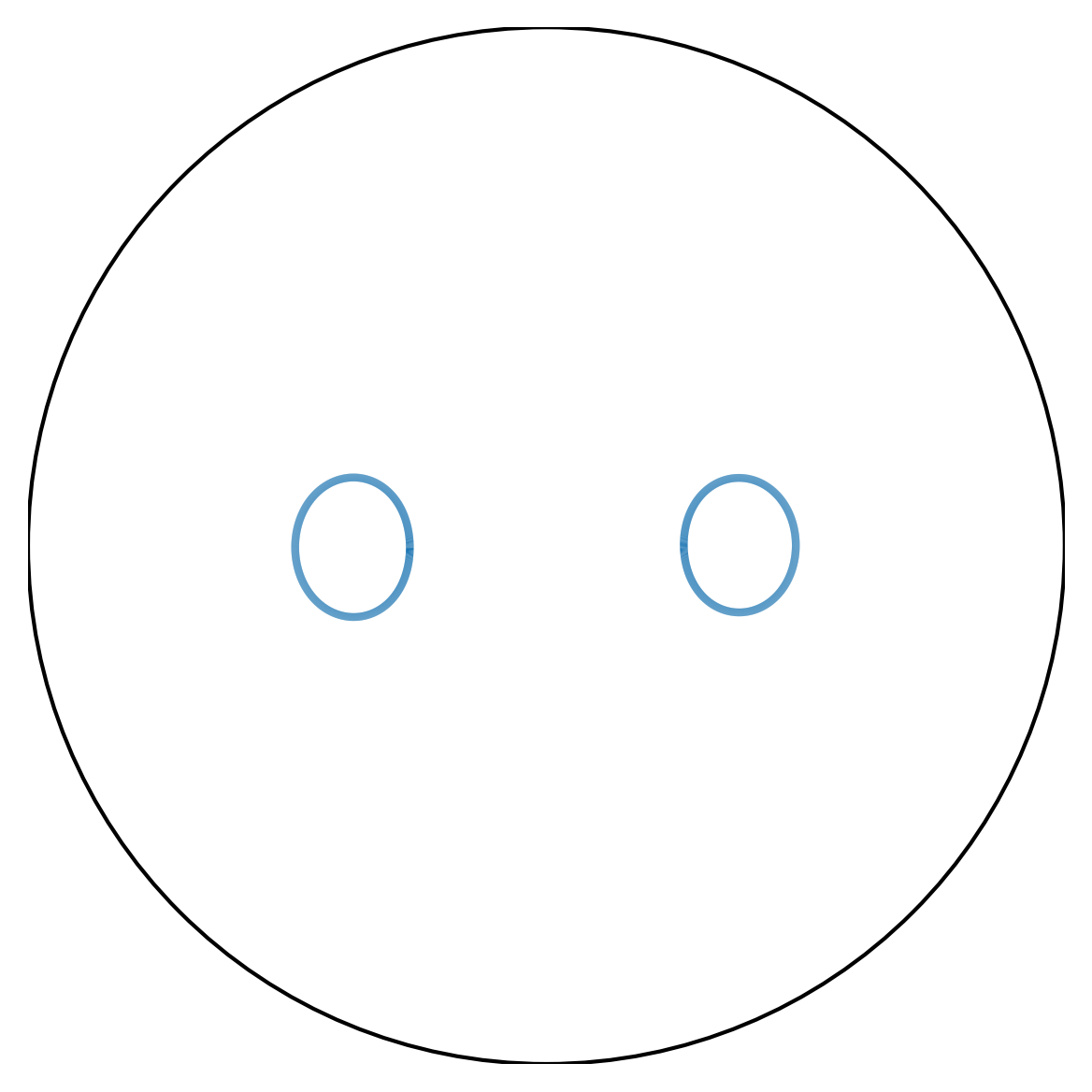}
		\includegraphics[width=0.2\textwidth]{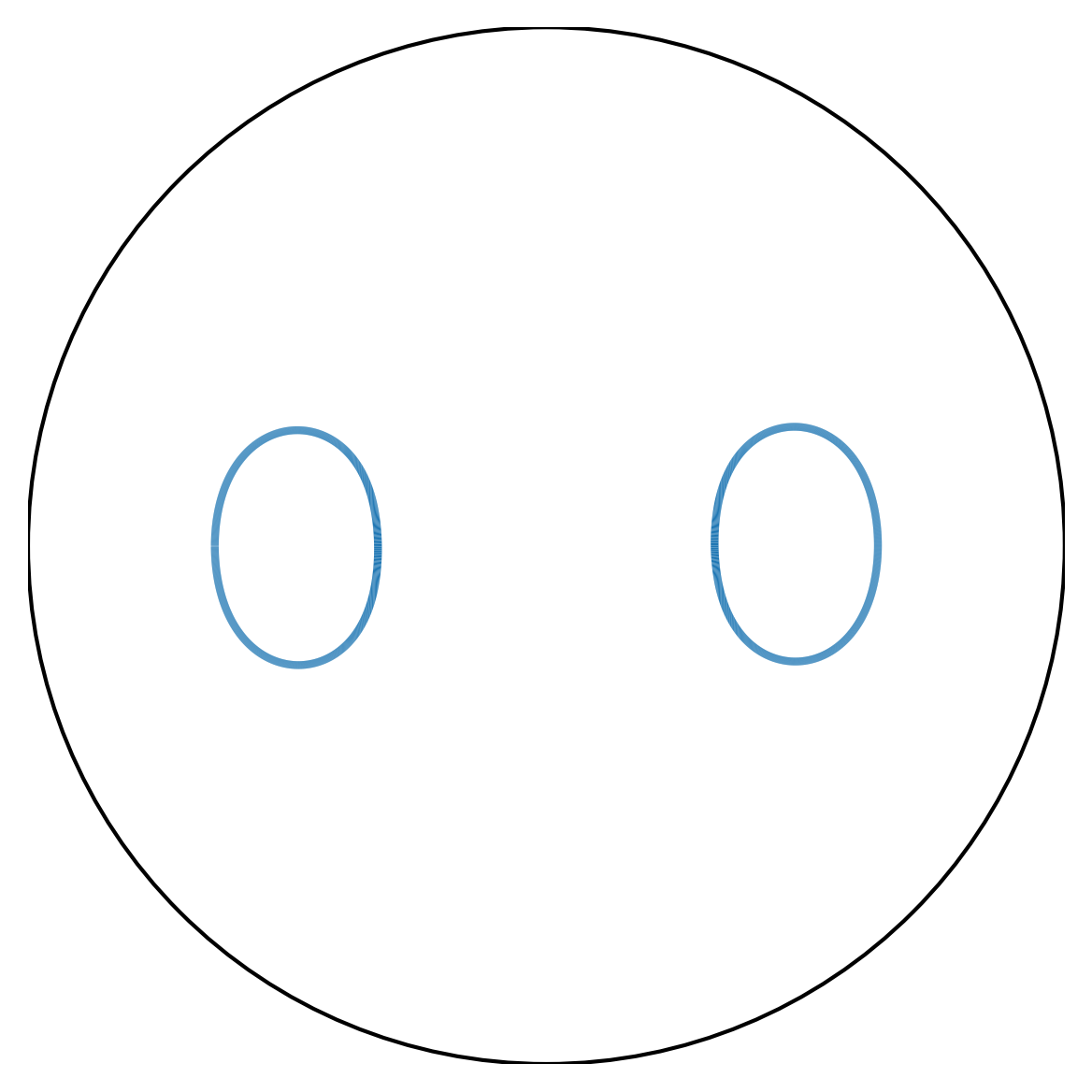}
	}
	\caption{(Minus-phase model: $\Omega = \bB^2_4(0)$)
		Initial data with mode $\ell=2$, 3, 4 perturbations around the radius 1.15,
		for $\alpha= 0.1$ and $S_I=-1.25$.
		We show the evolution at times $t=0,20,33,34,100$.}
	\label{fig:pd_a01_SI-125_cp-4}
\end{figure}

\begin{figure}
\center
\mbox{
\includegraphics[width=0.2\textwidth]{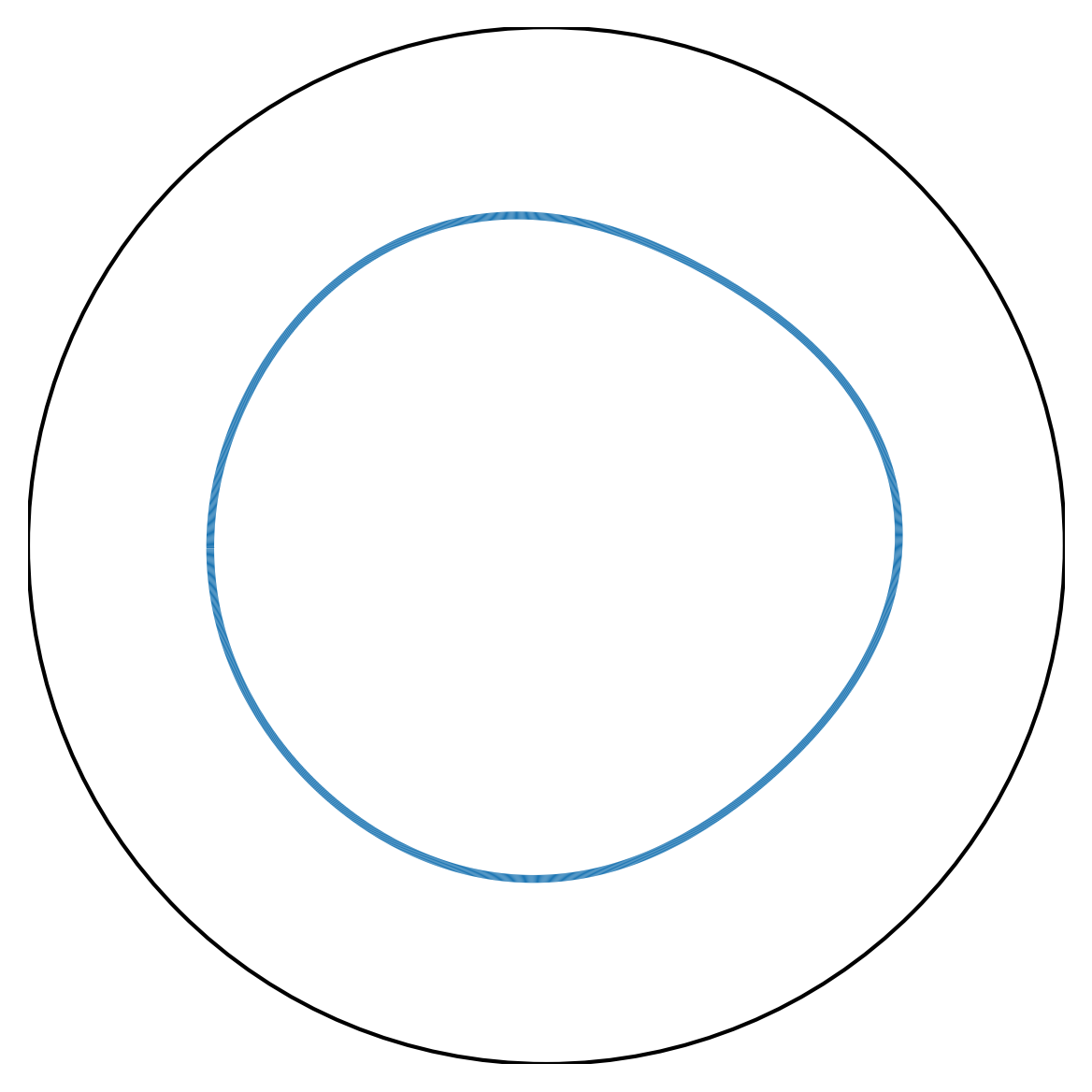}
\includegraphics[width=0.2\textwidth]{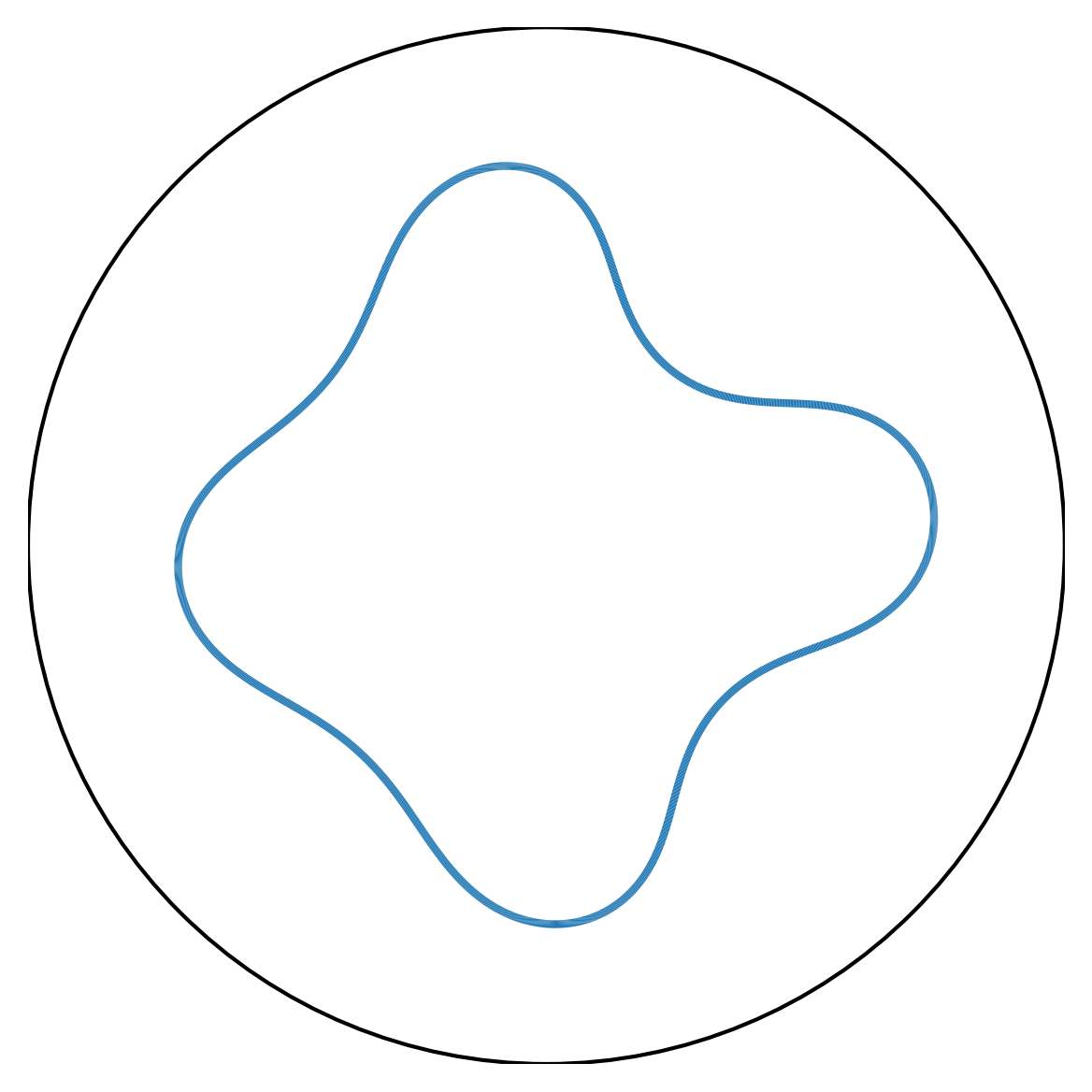}
\includegraphics[width=0.2\textwidth]{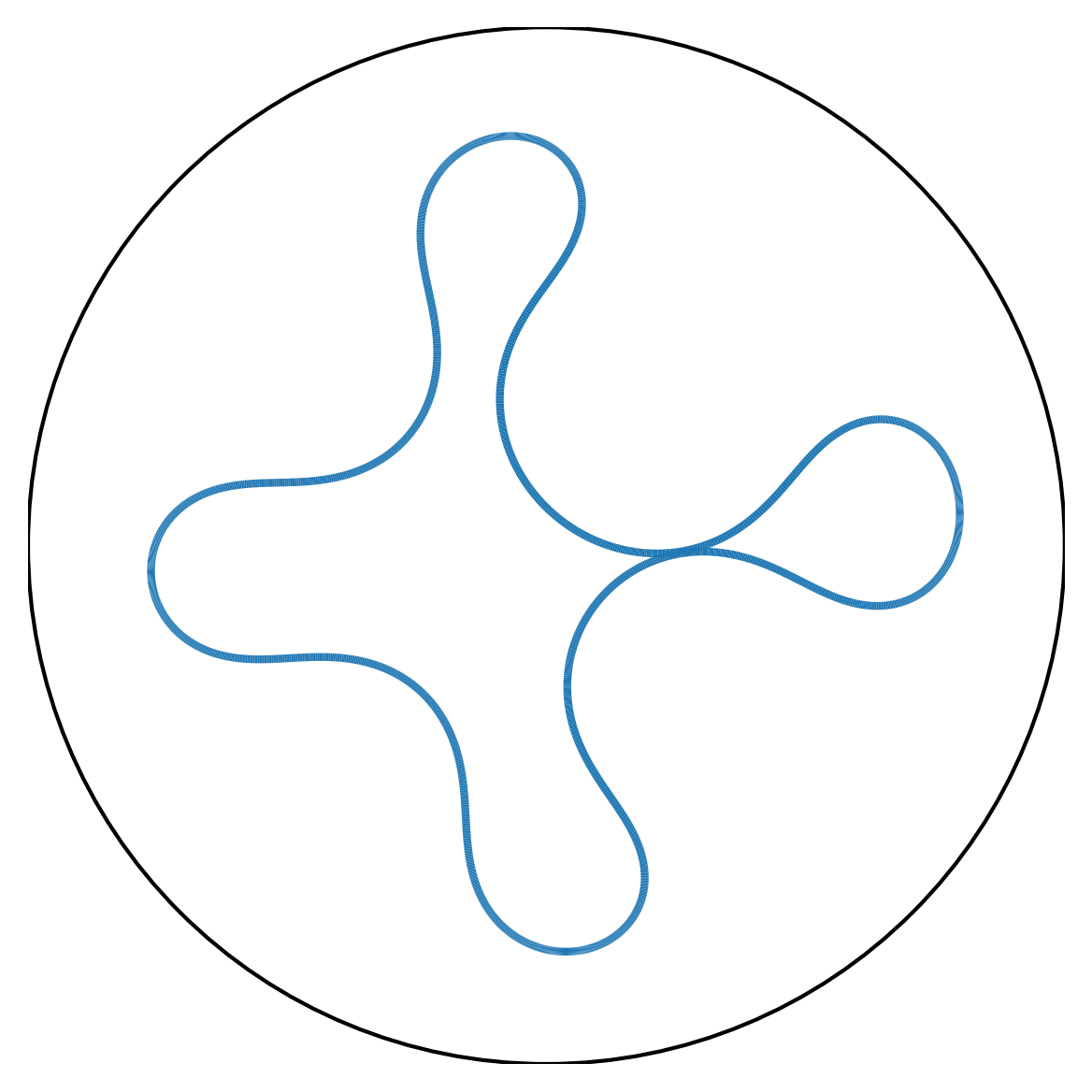}
\includegraphics[width=0.2\textwidth]{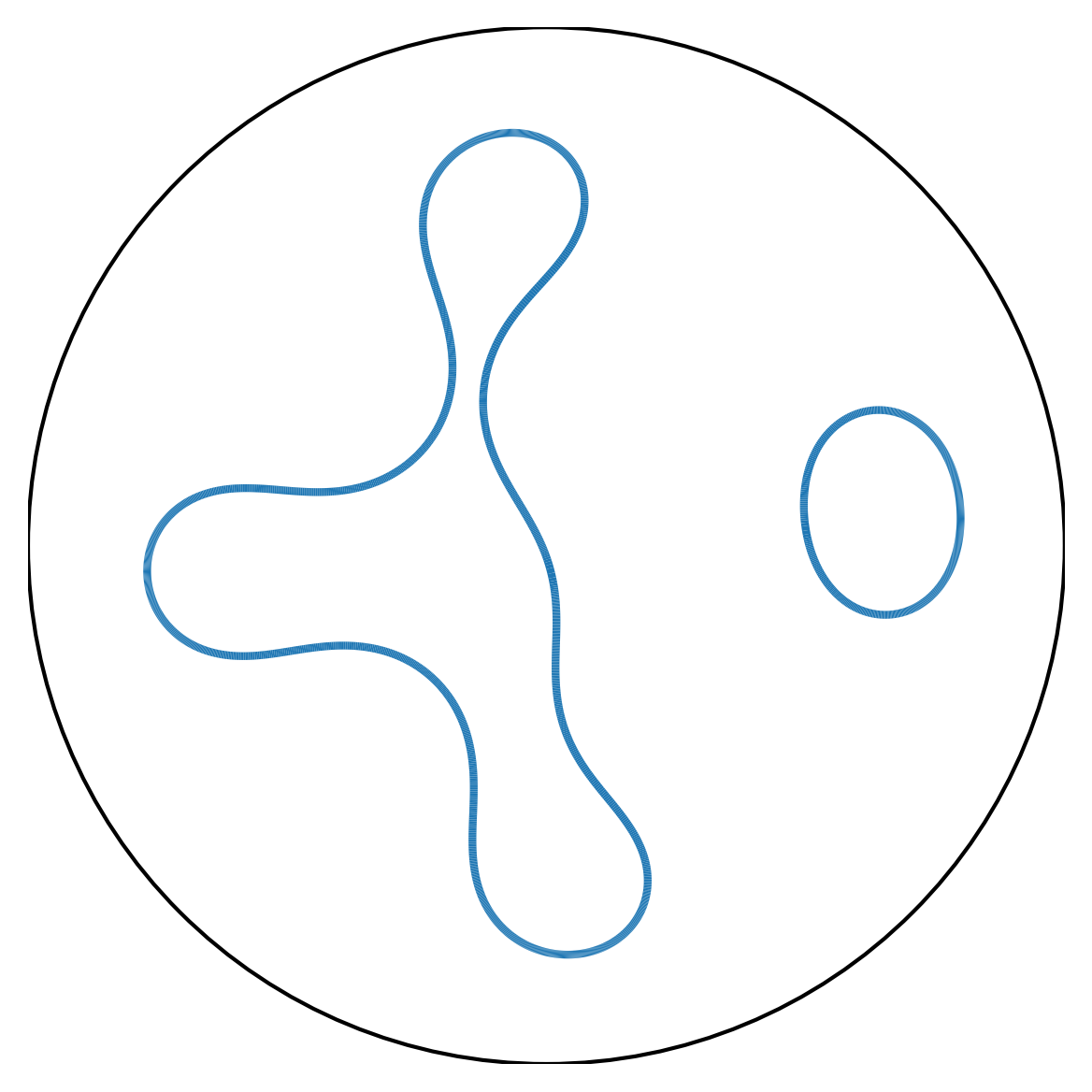}
\includegraphics[width=0.2\textwidth]{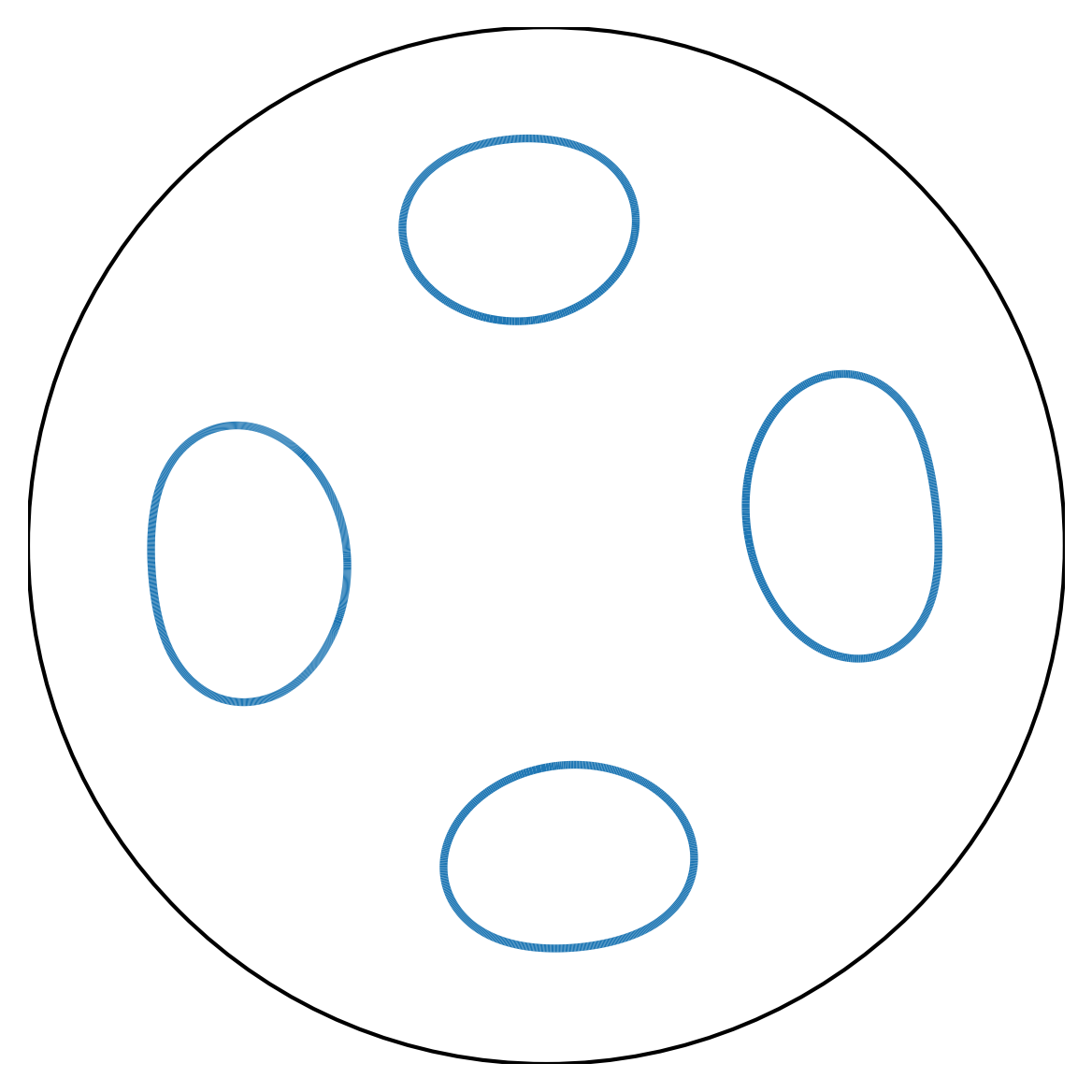}
}
\caption{(Minus-phase model: $\Omega = \bB^2_4(0)$)
Initial data with mode $\ell=2$, 3, 4 perturbations around the radius 2.57,
for $\alpha= 0.1$ and $S_I=-1$.
We show the evolution at times $t=0,40,49,50,80$.}
\label{fig:pd_a01_SI-1_cp-4}
\end{figure}

\begin{figure}
\center
\mbox{
\includegraphics[width=0.2\textwidth]{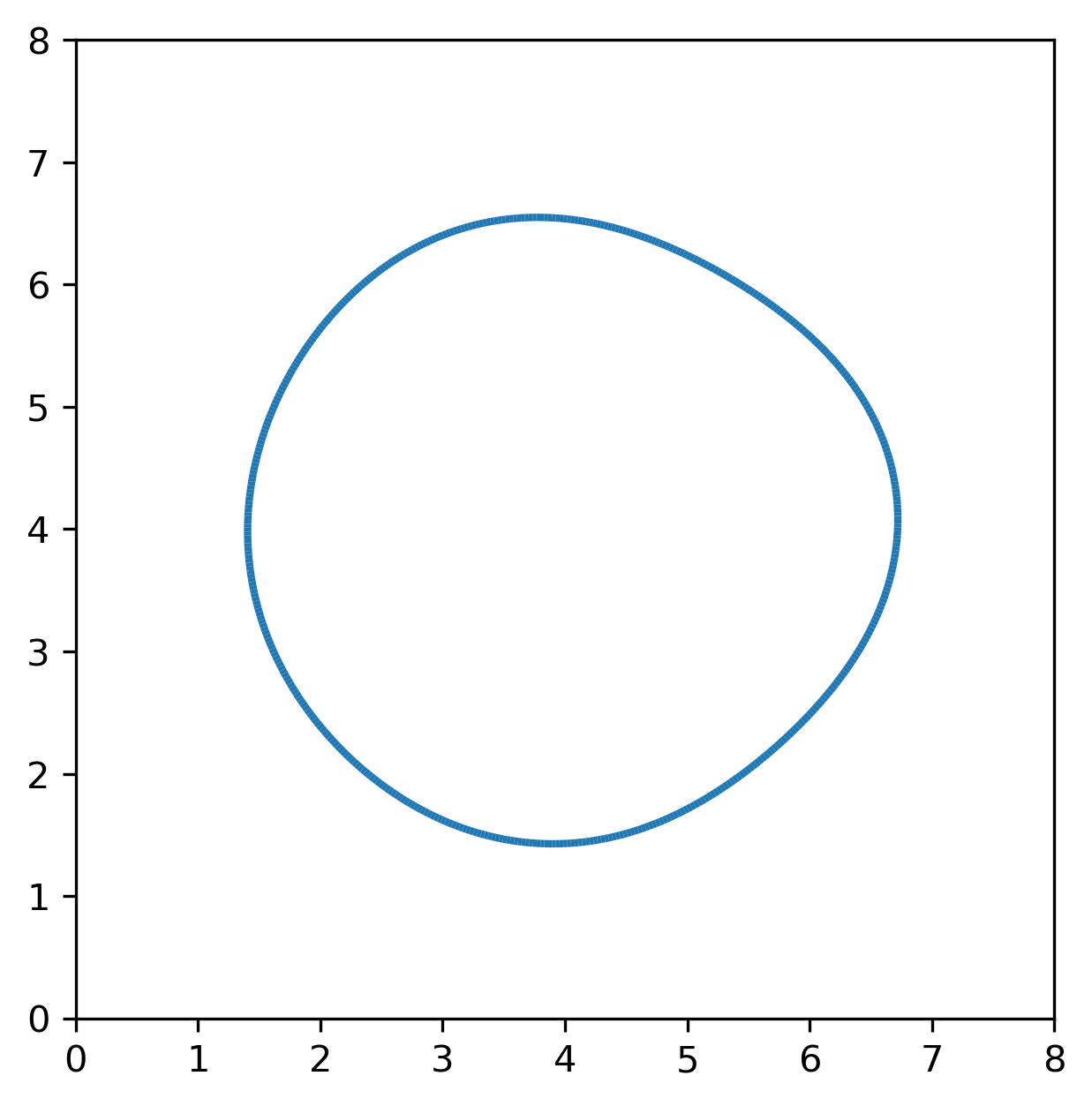}
\includegraphics[width=0.2\textwidth]{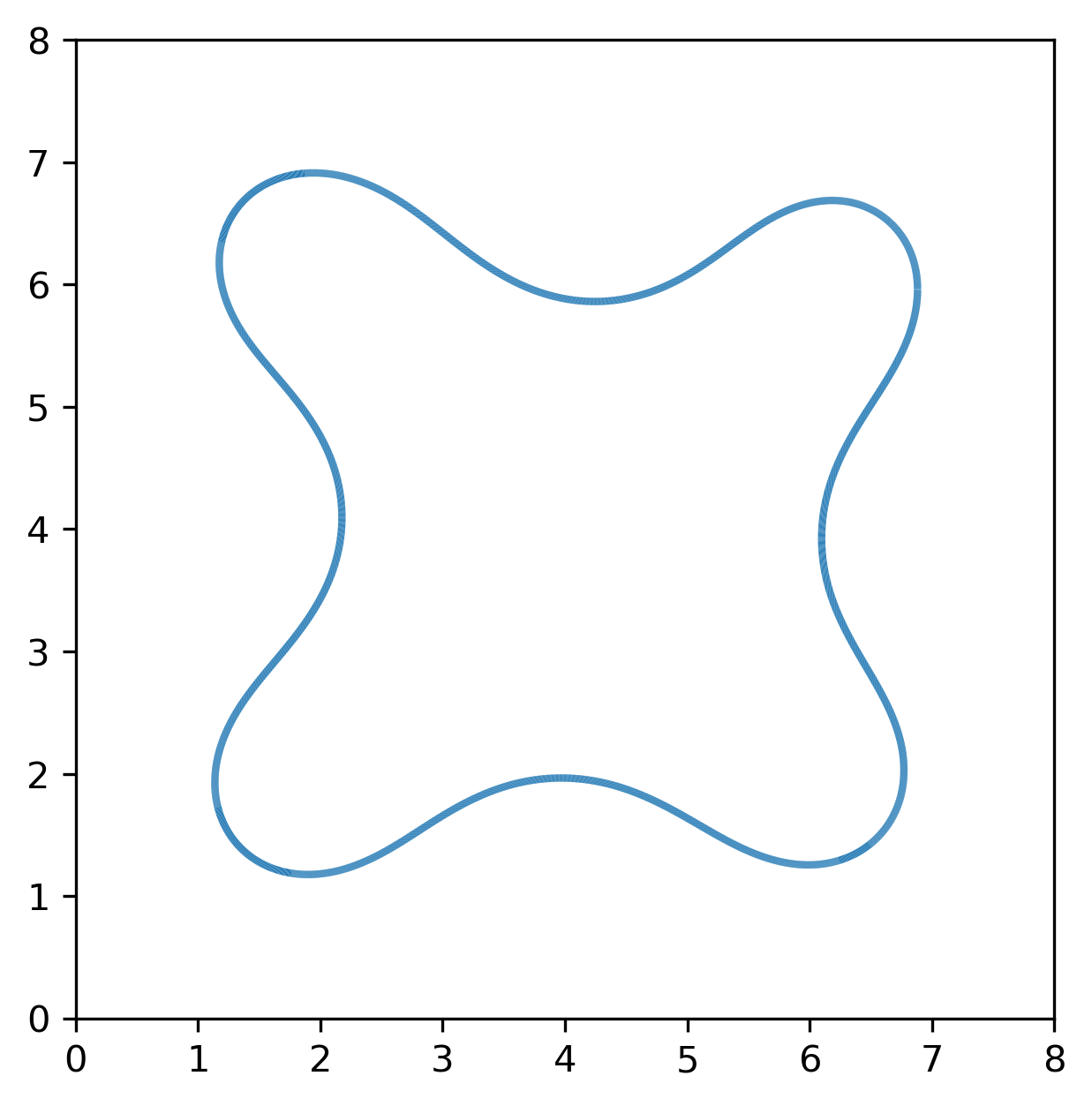}
\includegraphics[width=0.2\textwidth]{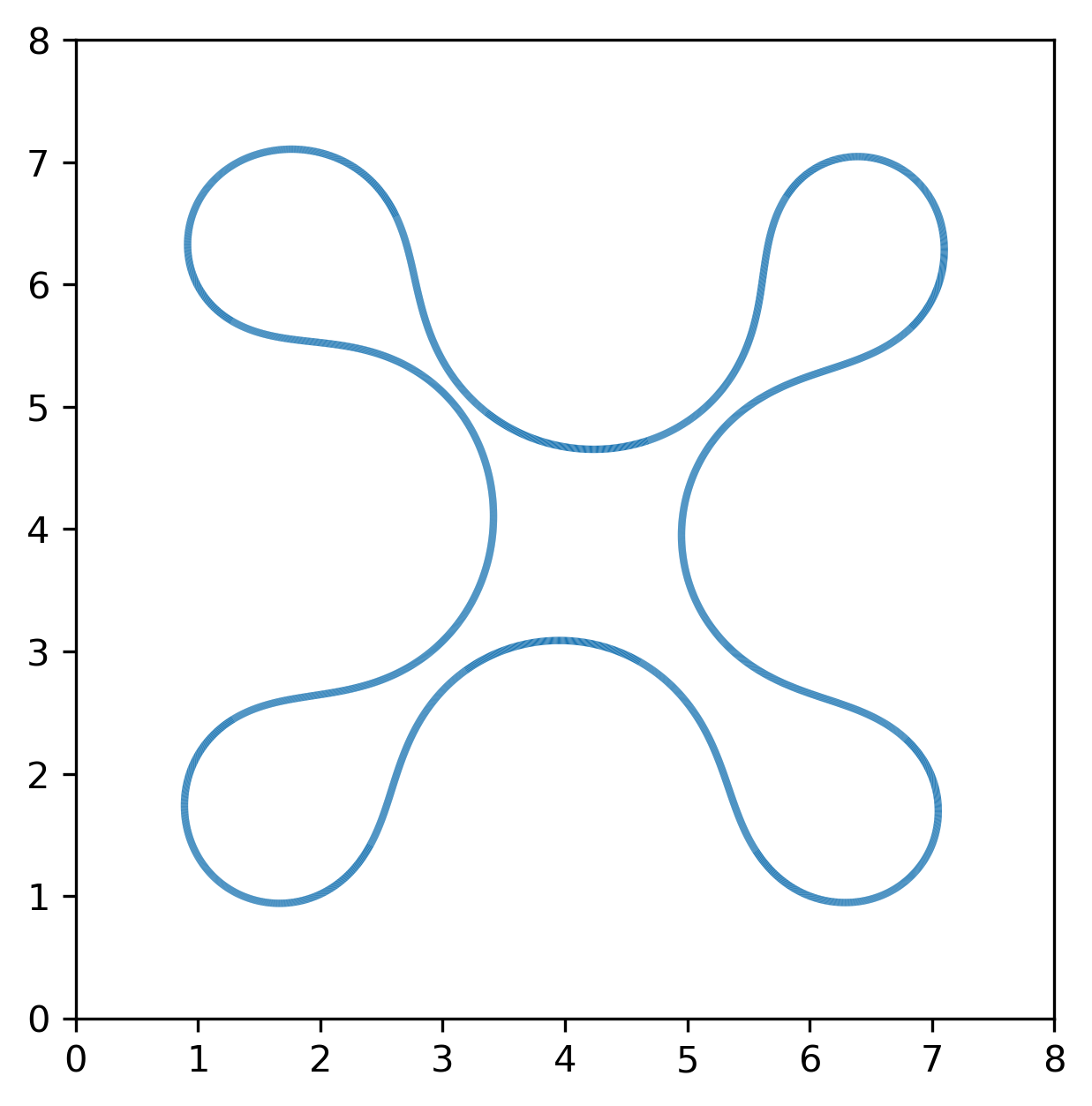}
\includegraphics[width=0.2\textwidth]{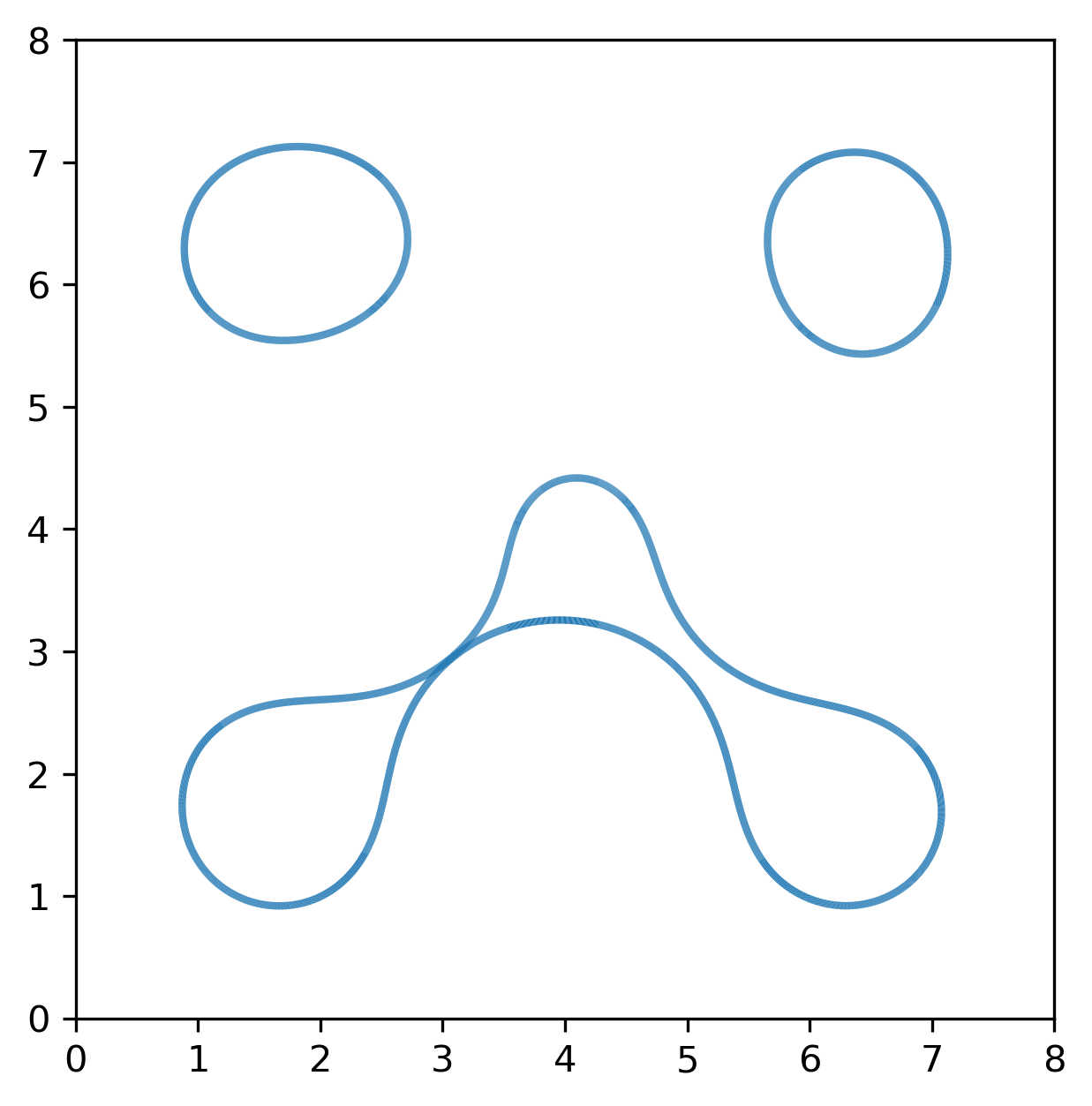}
\includegraphics[width=0.2\textwidth]{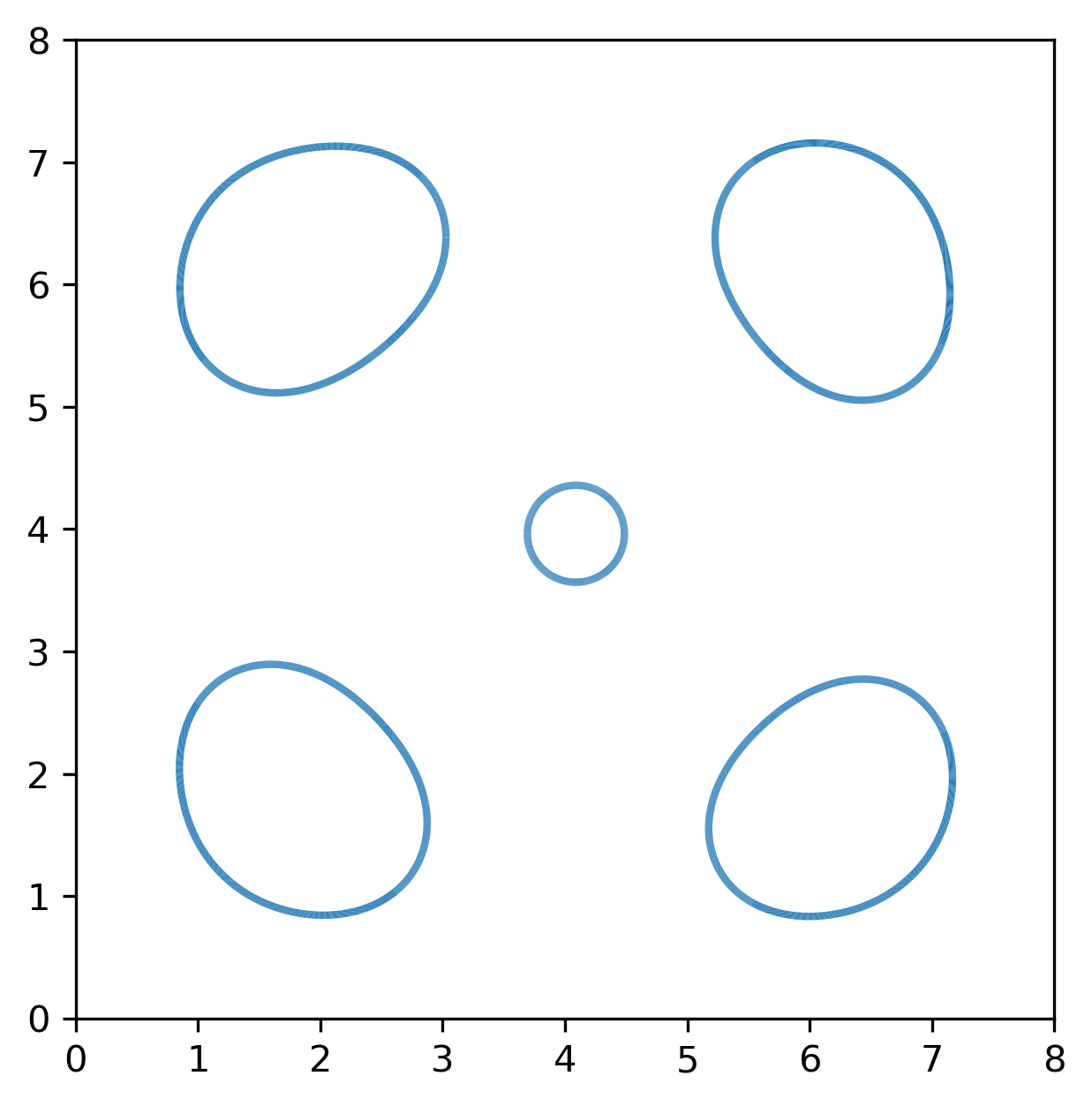}
}
\caption{(Minus-phase model: $\Omega = (0,8)^2$)
Initial data with mode $\ell=2$, 3, 4 perturbations around the radius 2.57,
for $\alpha= 0.1$ and $S_I=-1$.
We show the evolution at times $t=0,20,28,29,40$.}
\label{fig:pd_a01_SI-1_p-4}
\end{figure}

We let $\Omega = \bB^2_4(0)$, 
$\alpha=0.1$ and $S_I=-2, -1.4, -1.25, -1$, starting from a
perturbation of a circle with radius $r_0$ that includes the first 4 modes.
For the final three values of $S_I$ we choose $r_0$ to be the radius
of the stationary solution, so $r_0=2.57, 1.15, 0.58$. For $S_I=-2$ there is 
no stationary radius, and so we let $r_0=2.57$, with the expectation that the
interface will shrink to a point. 
In Figures~\ref{fig:pd_a01_SI-2_cp-4}, \ref{fig:pd_a01_SI-14_cp-4}, \ref{fig:pd_a01_SI-125_cp-4} and \ref{fig:pd_a01_SI-1_cp-4} we show the evolutions for the different setups and observe the expected behaviour: shrinking to extinction, stable evolution to equilibrium shape, growth of mode $\ell=2$, as well as growth of mode $\ell = 4$. Consequently, in all four cases, the observed behaviour agrees with the predictions of the linear stability analysis shown in Figure~\ref{fig:rad2D_minus}. In Figure~\ref{fig:pd_a01_SI-1_p-4} we show the evolution for $S_I = -1$ in $\Omega = (0,8)^2$ to contrast with Figure~\ref{fig:pd_a01_SI-1_cp-4}, namely the influence of the domain boundary on growth of certain perturbations. The behaviour for $S_I = -2, -1.4, -1.25$ when $\Omega = (0,8)^2$ is analogous to Figures~\ref{fig:pd_a01_SI-2_cp-4}, \ref{fig:pd_a01_SI-14_cp-4} and \ref{fig:pd_a01_SI-125_cp-4}.


\footnotesize
\bibliographystyle{abbrv} 
\bibliography{refs} 

@article{Agosti2017,
	author = {Agosti, Abramo and Antonietti, Paola Francesca and Ciarletta, Pasquale and Grasselli, Maurizio and Verani, Marco},
	title = {A {C}ahn--{H}illiard-type equation with application to tumor growth dynamics},
	journal = {Math. Methods Appl. Sci.},
	volume = {40},
	number = {18},
	pages = {7598--7626},
	doi = {https://doi.org/10.1002/mma.4548},
	url = {https://onlinelibrary.wiley.com/doi/abs/10.1002/mma.4548},
	eprint = {https://onlinelibrary.wiley.com/doi/pdf/10.1002/mma.4548},
	year = {2017}
}

@article{DonauSSetal20,
AUTHOR = {Donau, Carsten and  Sp\"ath, Fabian and Sosson, Mrilyne and Kriebisch, Brigitte A.K. and Schnitter, Fabian and  Tena-Solson, Marta and Hyun-Seo Kang and  Salibi, Elia  and  Sattler,Micahel and   Mutschler, Hannes and
  Boekhoven, Job},
TITLE={Active coacervate droplets as a model for membraneless organelles and protocells},
JOURNAL = {Nat. Commun.},
VOLUME ={11}, 
PAGES= {5167},
YEAR= {2020},
}

@article{BergmannBBetal23,
AUTHOR = {Bergmann, A.M. and  Bauermann, J. and Bartolucci, G. and  
 Donau, C. and   Stasi, M. and   Holtmannsp\"otter, A.-L. and  J\"ulicher, F. and Weber, C.A. and Boekhoven, J.},
TITLE={Liquid spherical shells are a non-equilibrium steady state of active droplets},
JOURNAL = {Nat. Commun.},
VOLUME ={14}, 
PAGES= {6552},
YEAR= {2023},
}

@Article{AbelsGG12,
AUTHOR = {Abels, Helmut 
and Garcke, Harald and Gr\"un, G\"unther},
TITLE = {Thermodynamically consistent, frame indifferent diffuse interface
models for incompressible two-phase flows with different densities},
JOURNAL = {Math. Models Methods Appl. Sci.},
FJOURNAL = {Mathematical Models and Methods in Applied Sciences},
VOLUME = {22},
YEAR = {2012},
NUMBER = {3},
PAGES = {1150013},
numpages = {40},
DOI = {10.1142/S0218202511500138},
URL = {https://doi.org/10.1142/S0218202511500138}
}

@Article{AlikakosBC94,
AUTHOR = {Alikakos, Nicholas D. and Bates, Peter W. and Chen, Xinfu},
TITLE = {Convergence of the {C}ahn--{H}illiard equation to the
{H}ele--{S}haw model},
JOURNAL = {Arch. Rational Mech. Anal.},
FJOURNAL = {Archive for Rational Mechanics and Analysis},
VOLUME = {128},
YEAR = {1994},
NUMBER = {2},
PAGES = {165--205},
ISSN = {0003-9527},
CODEN = {AVRMAW},
MRCLASS = {35Q99 (35B25 35C20 35K55 76D10 80A22)},
MRNUMBER = {1308851},
MRREVIEWER = {Peter J. Sternberg},
DOI = {10.1007/BF00375025},
URL = {https://doi.org/10.1007/BF00375025},
}

@preamble{
   "\def\soft#1{\leavevmode\setbox0=\hbox{h}\dimen7=\ht0\advance
    \dimen7 by-1ex\relax\if t#1\relax\rlap{\raise.6\dimen7
    \hbox{\kern.3ex\char'47}}#1\relax\else\if T#1\relax
    \rlap{\raise.5\dimen7\hbox{\kern1.3ex\char'47}}#1\relax
    \else\if d#1\relax\rlap{\raise.5\dimen7\hbox{\kern.9ex
    \char'47}}#1\relax\else\if D#1\relax\rlap{\raise.5\dimen7
    \hbox{\kern1.4ex\char'47}}#1\relax\else\if l#1\relax
    \rlap{\raise.5\dimen7\hbox{\kern.4ex\char'47}}#1\relax
    \else\if L#1\relax\rlap{\raise.5\dimen7\hbox{\kern.7ex
    \char'47}}#1\relax\else\message{accent \string\soft
    \space #1 not defined!}#1\relax\fi\fi\fi\fi\fi\fi} "
}

@Article{BarrettBG99, comment={was Dmob5},
AUTHOR = {Barrett, John W. and Blowey, James F. 
and Garcke, Harald},
TITLE = {Finite element approximation of the {C}ahn--{H}illiard equation
with degenerate mobility},
JOURNAL = {SIAM J. Numer. Anal.},
FJOURNAL = {SIAM Journal on Numerical Analysis},
VOLUME = {37},
YEAR = {1999},
NUMBER = {1},
PAGES = {286--318},
ISSN = {0036-1429},
MRCLASS = {65M60 (35A35)},
MRNUMBER = {1742748},
MRREVIEWER = {Thomas Y. Hou},
DOI = {10.1137/S0036142997331669},
URL = {https://doi.org/10.1137/S0036142997331669},
}

@Article{dendritic,
AUTHOR = {Barrett, John W. and Garcke, Harald and N\"urnberg, Robert},
TITLE = {On stable parametric finite element methods for the {S}tefan problem 
and the {M}ullins--{S}ekerka problem with applications to dendritic growth},
JOURNAL = {J. Comput. Phys.},
FJOURNAL = {Journal of Computational Physics},
volume = {229},
number = 18,
pages = {6270--6299},
year = {2010},
doi = {10.1016/j.jcp.2010.04.039},
URL = {https://doi.org/10.1016/j.jcp.2010.04.039},
ISSN = {0021-9991},
MRCLASS = {80A22 (65M12 65M60 80M10)},
MRNUMBER = {2660305},
MRREVIEWER = {Damian S\l ota},
}

@Article{crystal,
AUTHOR = {Barrett, John W. and Garcke, Harald and N\"urnberg, Robert},
TITLE = {Finite-element approximation of one-sided {S}tefan problems with 
anisotropic, approximately crystalline, {G}ibbs--{T}homson law},
JOURNAL = {Adv. Differ. Equ.},
volume = {18},
number = {3-4},
pages = {383--432},
year = 2013,
url = {http://projecteuclid.org/euclid.ade/1360073021},
ISSN = {1079-9389},
MRCLASS = {80A22 (35R35 65M60 74N05)},
MRNUMBER = {3060200},
}

@incollection{bgnreview, comment = {was: bgnreview19handbook},
AUTHOR = {Barrett, John W. and Garcke, Harald and N\"urnberg, Robert},
TITLE = {Parametric finite element approximations of curvature driven 
interface evolutions},
booktitle = {Handb. Numer. Anal.},
volume = 21,
editor = {Andrea Bonito and Ricardo H. Nochetto},
PAGES = {275--423},
PUBLISHER = {Elsevier},
address = {Amsterdam},
YEAR = {2020},
DOI = {10.1016/bs.hna.2019.05.002},
url = {https://doi.org/10.1016/bs.hna.2019.05.002},
}

@Article{CahnEN-C96, comment = {was: Cahnetal},
AUTHOR = {Cahn, J. W. and Elliott, Charles M. and Novick-Cohen, A.},
TITLE = {The {C}ahn--{H}illiard equation with a concentration dependent
mobility: motion by minus the {L}aplacian of the mean
curvature},
JOURNAL = {European J. Appl. Math.},
FJOURNAL = {European Journal of Applied Mathematics},
VOLUME = {7},
YEAR = {1996},
NUMBER = {3},
PAGES = {287--301},
ISSN = {0956-7925},
MRCLASS = {80A22 (35Q99)},
MRNUMBER = {1401172},
MRREVIEWER = {J. R. Ockendon},
}

@Book{Ciarlet78,
author = {Ciarlet, Philippe G.},
title = {The Finite Element Method for Elliptic Problems},
note = {Studies in Mathematics and its Applications, Vol. 4},
publisher = {North-Holland Publishing Co.},
address = {Amsterdam},
year = 1978,
pages = {xix+530},
isbn = {0-444-85028-7},
mrclass = {65N30},
mrnumber = {0520174},
mrreviewer = {Josef Nedoma}
}

@article{Davis04,
AUTHOR = {Davis, Timothy A.},
TITLE = {Algorithm 832: {UMFPACK} {V}4.3--an unsymmetric-pattern
multifrontal method},
JOURNAL = {ACM Trans. Math. Software},
FJOURNAL = {Association for Computing Machinery. Transactions on
Mathematical Software},
VOLUME = {30},
YEAR = {2004},
NUMBER = {2},
PAGES = {196--199},
ISSN = {0098-3500},
CODEN = {ACMSCU},
MRCLASS = {65F50 (65-04 65F05)},
MRNUMBER = {2075981},
DOI = {10.1145/992200.992206},
URL = {https://doi.org/10.1145/992200.992206},
}

@article{EbenbeckL21,
AUTHOR = {Ebenbeck, Matthias and Lam, Kei Fong},
TITLE = {Weak and stationary solutions to a
{C}ahn--{H}illiard--{B}rinkman model with singular potentials
and source terms},
JOURNAL = {Adv. Nonlinear Anal.},
FJOURNAL = {Advances in Nonlinear Analysis},
VOLUME = {10},
YEAR = {2021},
NUMBER = {1},
PAGES = {24--65},
ISSN = {2191-9496,2191-950X},
MRCLASS = {35Q92 (35D30 76D07 92C50)},
MRNUMBER = {4104274},
DOI = {10.1515/anona-2020-0100},
URL = {https://doi.org/10.1515/anona-2020-0100},
}

@Article{ElliottG96, comment = {was: ElliottAndGarcke96},
AUTHOR = {Elliott, Charles M. and Garcke, Harald},
TITLE = {On the {C}ahn--{H}illiard equation with degenerate mobility},
JOURNAL = {SIAM J. Math. Anal.},
FJOURNAL = {SIAM Journal on Mathematical Analysis},
VOLUME = {27},
YEAR = {1996},
NUMBER = {2},
PAGES = {404--423},
ISSN = {0036-1410},
CODEN = {SJMAAH},
MRCLASS = {35K55 (35Q99 82C26)},
MRNUMBER = {1377481},
MRREVIEWER = {Jing Xue Yin},
DOI = {10.1137/S0036141094267662},
URL = {https://doi.org/10.1137/S0036141094267662},
}

@article{GLSS,
AUTHOR = {Garcke, Harald 
and Lam, Kei Fong and Sitka, Emanuel and Styles, Vanessa},
title = {A {C}ahn--{H}illiard--{D}arcy model for tumour growth with
chemotaxis and active transport},
journal={Math. Models Methods Appl. Sci.},
  volume={26},
  number={06},
  pages={1095--1148},
  year={2016},
  publisher={World Scientific}
}

@Book{GilbargT83, comment = {was: GT},
AUTHOR = {Gilbarg, David and Trudinger, Neil S.},
TITLE = {Elliptic Partial Differential Equations of Second Order},
SERIES = {Grundlehren der Mathematischen Wissenschaften [Fundamental
Principles of Mathematical Sciences]},
VOLUME = {224},
EDITION = {Second},
PUBLISHER = {Springer-Verlag},
ADDRESS = {Berlin},
YEAR = {1983},
PAGES = {xiii+513},
ISBN = {3-540-13025-X},
MRCLASS = {35Jxx (35-01)},
MRNUMBER = {737190},
MRREVIEWER = {O. John},
DOI = {10.1007/978-3-642-61798-0},
URL = {https://doi.org/10.1007/978-3-642-61798-0},
}

@article{Glasner03,
AUTHOR = {Glasner, Karl},
TITLE = {A diffuse interface approach to {H}ele--{S}haw flow},
JOURNAL = {Nonlinearity},
FJOURNAL = {Nonlinearity},
VOLUME = {16},
YEAR = {2003},
NUMBER = {1},
PAGES = {49--66},
ISSN = {0951-7715},
CODEN = {NONLE5},
MRCLASS = {76D27 (35Q35 35R35)},
MRNUMBER = {1950775},
MRREVIEWER = {Darren G. Crowdy},
}

@article{LeeMS16,
AUTHOR = {Lee, Alpha Albert and M{\"u}nch, Andreas and S{\"u}li, Endre},
TITLE = {Sharp-interface limits of the {C}ahn--{H}illiard equation with
degenerate mobility},
JOURNAL = {SIAM J. Appl. Math.},
FJOURNAL = {SIAM Journal on Applied Mathematics},
VOLUME = {76},
YEAR = {2016},
NUMBER = {2},
PAGES = {433--456},
ISSN = {0036-1399},
MRCLASS = {35K35 (35B25 35K59 74N20 80A22 82C26)},
MRNUMBER = {3466205},
DOI = {10.1137/140960189},
URL = {https://doi.org/10.1137/140960189},
}

@Article{MullinsS63,
author = {W. W. Mullins and R. F. Sekerka},
title = {Morphological Stability of a Particle Growing by Diffusion or Heat 
Flow},
journal = {J. Appl. Phys.},
volume = 34,
number = 2,
pages = {323--329},
year = 1963,
doi = {10.1063/1.1702607},
url = {https://doi.org/10.1063/1.1702607},
}

@Book{Alberta,  comment = {also: AlbertaHG},
author = {A. Schmidt and K. G. Siebert},
title = {Design of Adaptive Finite Element Software: The Finite Element Toolbox
{ALBERTA}},
series = {Lect. Notes Comput. Sci. Eng.},
volume = 42,
publisher = {Springer-Verlag},
address = {Berlin},
ISBN = {3-540-22842-X},
MRCLASS = {65-02 (65-04 65N30 65N50 68P05)},
MRNUMBER = {2127659},
year = 2005,
PAGES = {xii+315},
doi = {10.1007/b138692},
url = {https://doi.org/10.1007/b138692},
}

@Article{Simon87, comment = {was: Simon},
author = {J. Simon},
title = {Compact sets in the space {$L^{p}(0,T;B)$}},
journal = {Ann. Mat. Pura Appl.},
volume = 146,
pages = {65--96},
year = 1987
}

@article{bauermann2023formation,
  title={Formation of liquid shells in active droplet systems},
  author={Bauermann, Jonathan and Bartolucci, Giacomo and Boekhoven, Job and Weber, Christoph A and J{\"u}licher, Frank},
  journal={Phys. Rev. Res.},
  volume={5},
  number={4},
  pages={043246},
  year={2023},
  publisher={APS}
}

@article{zwickerostwald, 
title={Suppression of {O}stwald ripening in active emulsions},
author={Zwicker, David and Hyman, Anthony A. and J{\"u}licher, Frank},
journal={Phys. Rev. E},
volume={92},
number={1},
pages={012317},
year={2015}
}

@article{Active_drops,
author = {Zwicker, David and Seyboldt, Rabea and Weber, Christoph A. and
J{\"u}licher, Frank},
title = {Growth and division of active droplets provides a model for 
protocells},
journal = {Nat. Phys.},
fjournal = {Nature Physics},
volume = 13,
number = {4},
pages = {408--413},
year = 2017
}

@book {BDGP,
    AUTHOR = {B\"ansch, Eberhard and Deckelnick, Klaus and Garcke, Harald
              and Pozzi, Paola},
     TITLE = {Interfaces: modeling, analysis, numerics},
    SERIES = {Oberwolfach Seminars},
    VOLUME = {51},
 PUBLISHER = {Birkh\"auser/Springer, Cham},
      YEAR = {[2023] \copyright2023},
     PAGES = {xii+178},
      ISBN = {978-3-031-35549-3; 978-3-031-35550-9},
   MRCLASS = {65M60 (65N30)},
  MRNUMBER = {4676342},
       DOI = {10.1007/978-3-031-35550-9},
       URL = {https://doi.org/10.1007/978-3-031-35550-9},
}

@book{WatsonBessel,
AUTHOR = {G. N. Watson},
TITLE = {A Treatise on the Theory of Bessel Functions},
YEAR = {1995},
PAGES = {814},
PUBLISHER = {Cambridge University Press},
ADDRESS = {Cambridge},
}

@misc{NIST:DLMF,
         key = "{\relax DLMF}",
       title = "{\it NIST Digital Library of Mathematical Functions}",
howpublished = "\url{https://dlmf.nist.gov/}, Release 1.2.1 of 2024-06-15",
         url = "https://dlmf.nist.gov/",
        note = "F.~W.~J. Olver, A.~B. {Olde Daalhuis}, D.~W. Lozier, B.~I. Schneider,
                R.~F. Boisvert, C.~W. Clark, B.~R. Miller, B.~V. Saunders,
                H.~S. Cohl, and M.~A. McClain, eds."}

@article{ESCHER1998267,
title = {A Center Manifold Analysis for the {M}ullins--{S}ekerka Model},
journal = {J. Differential Equations},
volume = {143},
number = {2},
pages = {267--292},
year = {1998},
issn = {0022-0396},
doi = {https://doi.org/10.1006/jdeq.1997.3373},
url = {https://www.sciencedirect.com/science/article/pii/S0022039697933738},
author = {Joachim Escher and Gieri Simonett}
}

@article{OOL,
author = {Garcke, Harald and Lam, Kei Fong and N\"urnberg, Robert and 
Signori, Andrea},
title = {On a {C}ahn--{H}illiard equation for the growth and division of 
chemically active droplets modelling protocells},
journal = {European J. Appl. Math.},
fjournal={European Journal of Applied Mathematics},
volume=37, 
number=5, 
year=2026, 
pages={1106--1136},
url = {https://doi.org/10.1017/S0956792525100211},
doi = {10.1017/S0956792525100211}
}

@article{OOL2,
author = {Garcke, Harald and Lam, Kei Fong and N\"urnberg, Robert and 
Signori, Andrea},
     TITLE = {On a {M}ullins--{S}ekerka model for the growth of active droplets modeling protocells: Stability analysis and numerical computations},
journal = {Math. Models Methods Appl. Sci.},
volume = {Online first},
      YEAR = {2026},
     PAGES = {},
     DOI = {https://doi.org/10.1142/S0218202526500600}
}

@article{ebenbeck2021cahn,
  title={Cahn--{H}illiard--{B}rinkman systems for tumour growth.},
  author={Ebenbeck, Matthias and Garcke, Harald and N{\"u}rnberg, Robert},
  journal={Discrete Contin. Dyn. Syst. Ser.  S},
  volume={14},
  number={11},
  year={2021},
  pages = {3989--4033}
}

@article{garcke2025anisotropic,
  title={The anisotropic {C}ahn--{H}illiard equation with degenerate mobility: existence of weak solutions},
  author={Garcke, Harald and Knopf, Patrik and Signori, Andrea},
  journal={Anal. Appl.},
  year={2025},
  pages={1--23},
  publisher={World Scientific},
}

@article{Carrillo,
	title={Degenerate {C}ahn--{H}illiard systems: {F}rom nonlocal to local},
	author={Carrillo, Jos\'e A. and Elbar, Charles and Skrzeczkowski, Jakub},
	journal={Commun. Contemp. Math.},
	volume = {27},
	pages={2450041},
	year={2025},
	publisher={World Scientific}
}

\end{document}